\documentclass{article} 
\usepackage{amssymb}
\usepackage{amscd}
\usepackage[all]{xy}
\usepackage{amsmath} 
\usepackage[dvipdfmx]{graphicx}
\usepackage{latexsym} 
\usepackage{mathtools}
\usepackage{longtable}
\usepackage{enumerate}
\usepackage{comment}
\usepackage[usenames]{color}  
\allowdisplaybreaks

\usepackage{theorem}
\theoremstyle{plain}
\theorembodyfont{\itshape}
\newtheorem{theorem}{Theorem}[section]
\newtheorem{proposition}[theorem]{Proposition}
\newtheorem{lemma}[theorem]{Lemma}
\newtheorem{corollary}[theorem]{Corollary} 
\newtheorem{result}[theorem]{Result} 
 
\theorembodyfont{\rmfamily}

\newtheorem{example}[theorem]{Example}
\newtheorem{hypothesis}[theorem]{Hypothesis} 
\newtheorem{induction}[theorem]{Induction} 
\newtheorem{notation}[theorem]{Notation}
\newtheorem{problem}[theorem]{Problem}
 
\newtheorem{recipe}[theorem]{Recipe}
\newtheorem{remark}[theorem]{Remark}
\renewcommand{\thefigure}{\arabic{section}.\arabic{figure}}
\renewcommand{\thetable}{\arabic{section}.\arabic{table}}
\makeatletter
\@addtoreset{figure}{section}
\@addtoreset{table}{section}
\makeatother

\def\bC{\mathbb{C}}
\def\bN{\mathbb{N}}
\def\bP{\mathbb{P}}
\def\bQ{\mathbb{Q}}
\def\bR{\mathbb{R}}
\def\bZ{\mathbb{Z}}
\def\cC{\mathcal{C}}
\def\cE{\mathcal{E}}
\def\cF{\mathcal{F}}
\def\cG{\mathcal{G}}
\def\cK{\mathcal{K}} 
\def\cP{\mathcal{P}}
\def\ra{\mathrm{a}}
\def\rb{\mathrm{b}}
\def\rc{\mathrm{c}}
\def\rd{\mathrm{d}}
\def\re{\mathrm{e}}
\def\rf{\mathrm{f}}
\def\rn{\mathrm{n}}
\def\roff{\mathrm{off}}
\def\ron{\mathrm{on}} 
\def\rr{\mathrm{r}}
\def\rt{\mathrm{t}}
\def\rv{\mathrm{v}}
\def\ri{\mathrm{i}}
\def\rii{\mathrm{ii}}
\def\rA{\mathrm{A}}

\def\rC{\mathrm{C}}
\def\rCT{\mathrm{CT}} 
\def\rD{\mathrm{D}}
\def\rE{\mathrm{E}}
\def\rHy{\mathrm{Hy}}
\def\rI{\mathrm{I}} 
\def\rN{\mathrm{N}}
\def\rR{\mathrm{R}}
\def\rS{\mathrm{S}}
\def\rST{\mathrm{ST}}
\def\rII{\mathrm{I\!I}}
\def\diag{\operatorname{diag}}
\def\Ad{\operatorname{Ad}} 
\def\Aut{\operatorname{Aut}} 
\def\Fixi{\operatorname{Fix}^{\mathrm{i}}}
\def\Gal{\operatorname{Gal}}
\def\GL{\operatorname{GL}}
\def\Ker{\operatorname{Ker}}
\def\Mob{\operatorname{M\mbox{\"{o}}b}}
\def\PGL{\operatorname{PGL}}
\def\Pic{\operatorname{Pic}} 
\def\PSL{\operatorname{PSL}}
\def\PSU{\operatorname{PSU}}
\def\Res{\operatorname{Res}}  
\def\Span{\operatorname{Span}}
\def\SL{\operatorname{SL}}
\def\SO{\operatorname{SO}}
\def\SU{\operatorname{SU}}
\def\Tr{\operatorname{Tr}}
\def\bc{\mbox{\boldmath $c$}}
\def\br{\mbox{\boldmath $r$}}
\def\bs{\mbox{\boldmath $s$}}
\def\bu{\mbox{\boldmath $u$}}
\def\bv{\mbox{\boldmath $v$}}
\def\bal{\mbox{\boldmath $\alpha$}}
\def\vD{\varDelta}
\def\ve{\varepsilon}
\def\vG{\varGamma} 
\def\div{ \, | \, }
\def\pd{ \, |\!| \, }
\title{\bf Equivariant Linearization and \\ 
Rotation Domains on K3 Surfaces}   
\author{Katsunori Iwasaki\thanks{Professor Emeritus, Department of Mathematics, 
Hokkaido University, Sapporo 060-0810 Japan, supported by JSPS KAKENHI, JP22K03365. 
{\tt iwasaki@math.sci.hokudai.ac.jp}} }
\date{}   
\begin{document}
\maketitle
\begin{abstract} 
We construct a lot of K3 surface automorphisms of positive entropy having rotation domains 
of ranks $1$ and $2$. 
To carry out this construction, we first lay theoretical foundations concerning equivariant 
linearization of nonlinear maps under resolutions of quotient singularities, 
linear models near exceptional components, Salem numbers and multipliers at periodic points, 
two kinds of fixed point formulas and related indices at exceptional components.   
Then these basic tools are combined with the method of hypergeometric groups 
to enable us to detect various types of rotation domains on K3 surfaces.   
\end{abstract} 
\section{Introduction} \label{sec:intro} 
Given an automorphism $f : X \to X$ of infinite order on a compact complex surface $X$,    
its {\sl Fatou set} $\cF(f)$ is the set of all points $p \in X$ such that 
the family $\{ f^n \}_{n=0}^{\infty}$ is equicontinuous in a neighborhood of $p$, 
where $f^n = f \circ \cdots \circ f$ ($n$ times) stands for the $n$-th iterate of $f$.  
It is the largest open subset of $X$ on which the dynamics of $f$ exhibits regular behaviors.    
Each connected component of $\cF(f)$ is referred to as a {\sl Fatou component}.   
\par
Following Bedford and Kim \cite[\S 1]{BK2}, we rewiew the concept of a rotation domain; 
see also Forn{\ae}ss and Sibony \cite{FS} and Ueda \cite{Ueda}.  
An $f$-invariant Fatou component $U \subset \cF(f)$ is called a {\sl rotation domain of rank} $r$,  
if iterations of $f$ induce an $r$-dimensional real torus action on $U$.    
To be more precise, let $\cG(U)$ be the set of all normal limits 
$g = \displaystyle \lim_{n_j \to \infty} (f|_U)^{n_j} : U \to \overline{U}$ 
of subsequences of $\{ (f|_U)^n \}_{n=0}^{\infty}$.  
By \cite[Theorem 1.2]{BK2}, if the condition  
\begin{equation} \label{eqn:G(U)}
\mbox{for any $g \in \cG(U)$ its image $g(U)$ contans an open set in $\overline{U}$},  
\end{equation}
is satisfied, then $\cG(U)$ is a compact, infinite abelian, Lie group of automorphisms 
on $U$ in the compact-open topology, so the connected component $\cG_0(U)$ of 
the identity in $\cG(U)$ is a real torus of positive dimension. 
We say that $U$ is a {\sl rotation domain of rank} $r$, if $\cG_0(U)$ is a real torus of dimension $r$.  
According to \cite[Theorems 1.3 and 1.5]{BK2} we have $1 \le r \le 4$ and if $X$ is 
K\"{a}hler and $f$ has positive topological entropy, then $r$ is either $1$ or $2$. 
Replacing $f$ by some iterate $f^m$ with $m \ge 2$, 
we can also speak of a {\sl periodic} rotation domain for the map $f$. 
A fixed (or periodic) point lying in a (periodic) rotation domain is called 
a {\sl center} of that rotation domain.  
\par
Bedford and Kim \cite{BK1, BK2} constructed rational surface automorphisms of 
positive entropy having rotation domains of ranks $1$ and $2$, and made a detailed 
study of the constructed mappings.    
In their construction they started with explicit birational self-maps of the 
projective plane $\bP^2$ and performed iterated blowups to get desired automorphisms. 
We also refer to McMullen \cite{McMullen2}, Bedford and Kim \cite{BK1} for rotation domains 
on rational surfaces. 
Now let us turn our attention to rotation domains on K3 surfaces. 
McMullen \cite{McMullen1} constructed K3 surface automorphisms of positive entropy with 
a rotation domain of rank $2$, which he called a {\sl Siegel disk}. 
His examples had Picard number $0$. 
Oguiso \cite{Oguiso} found a similar example of Picard number $8$. 
Iwasaki and Takada \cite{IT2} obtained more automorphisms with Siegel disks that 
covered all possible Picard numbers, that is, even integers between $0$ and $18$. 
Their construction was based on {\sl the method of hypergeometric groups} introduced 
in \cite{IT1}.  
\par
This article develops a detailed study of rotation domains on K3 surfaces, with emphasis 
on rotation domains of rank $1$ or on the coexistence of rank $1$ domains with rank $2$ ones. 
A K3 surface is a simply connected, compact complex surface $X$ admitting a nowhere vanishing 
holomorphic $2$-form $\eta_X$. 
It is K\"{a}hler by Siu \cite{Siu}.       
Let $f : X \to X$ be an automorphism of $X$.  
There exists a complex number $\delta(f)$ of modulus $1$ such that 
\begin{equation} \label{eqn:tn}
f^* \eta_X = \delta(f) \, \eta_X. 
\end{equation}
In this article $\delta(f)$ is referred to as the {\sl twist number} of $f$. 
It was called the {\sl determinant} of $f$ by McMullen \cite{McMullen1} and 
the {\sl special eigenvalue} of $f$ by Iwasaki and Takada \cite{IT1,IT2} respectively. 
Equation \eqref{eqn:tn} implies that $f$ preserves the smooth volume form 
$\eta_X \wedge \overline{\eta}_X$, so that condition \eqref{eqn:G(U)} is satisfied  
by \cite[Proposition 1.1]{BK2}. 
\par
A {\sl Salem number} is a real algebraic integer $\lambda > 1$ which is conjugate to 
$\lambda^{-1}$ and whose remaining conjugates have modulus $1$. 
Its degree is an even integer $d \ge 4$. 
We refer to Smyth \cite{Smyth} for a good survey on Salem numbers.   
The twist number $\delta = \delta(f)$ in \eqref{eqn:tn} is either a root of unity or a 
Galois conjugate to a Salem number $\lambda$ of degree $d \le 22$. 
We always assume that the automorphism $f : X \to X$ satisfies the following two conditions:   
\begin{itemize} 
\setlength{\itemsep}{-1pt}
\item[(A1)] the twist number $\delta = \delta(f)$ is conjugate to a Salem number $\lambda = \lambda(f)$,   
\item[(A2)] the Picard group $\Pic(X)$ is non-degenerate with respect to the intersection form on $H^2(X, \bZ)$.   
\end{itemize}
\par
From (A1) the K3 surface $X$ is non-projective and the automorphism $f$ has topological 
entropy 
$$
h(f) = \log \lambda(f) > 0,  
$$
where $\lambda(f)$ is referred to as the {\sl dynamical degree} of $f$.  
When $X$ is non-projective, $\Pic(X)$ may degenerate, so condition (A2) is 
imposed to rule out this possibility so that $\Pic(X)$ is negative definite. 
Under (A1) and (A2), all irreducible curves in $X$ are $(-2)$-curves 
and the number of them is at most the Picard number $\rho(X)$, the rank of $\Pic(X)$.    
The union $\cE(X)$ of all $(-2)$-curves in $X$ is called the {\sl exceptional set} of $X$.  
Each connected component $E \subset \cE(X)$ is called an {\sl exceptional component}.  
The dual graph of $E$ is a Dynkin diagram $D$ of type $\rA$, $\rD$ or $\rE$. 
The map $f$ permutes the exceptional components and may leave some of them invariant. 
If $E$ is $f$-invariant, then $f$ permutes the $(-2)$-curves in $E$ and induces an automorphism 
of the diagram $D$,  
\begin{equation} \label{eqn:da}
\sigma(f, E) \in \Aut(D). 
\end{equation}
The order of $\sigma(f, E)$ is $1$, $2$, or $3$. 
In the respective cases we say that the pair $(f, E)$, or rather the set $E$ when 
the map $f$ is understood, has {\sl trivial}, {\sl reflectional}, or {\sl tricyclic symmetry}, 
where the tricyclic case occurs only when $D$ is of type $\rD_4$.  
Since the number of exceptional components is finite, all exceptional components 
are left invariant by some iterate $f^m$.  
Now it is natural to pose the following problem.        
\par
\begin{problem} \label{prob:rd} 
Given an $f$-invariant exceptional component $E \subset \cE(X)$, we ask whether 
\begin{enumerate} 
\setlength{\itemsep}{-1pt}
\item[(R1)] $E$ is completely contained in a rotation domain of rank $1$,  
\item[(R2)] $E$ is completely contained in a rotation domain of rank $2$, 
\item[(Hy)] $E$ is hyperbolic to the effect that a stable curve flows into $E$ and an unstable curve emanates from $E$. 
\end{enumerate}
This question can be asked for $f^m$ in place of $f$ and hence makes sense 
for all exceptional components in $X$. 
\end{problem}
\par
One of the main results in this article is the following solution to this problem.    
\begin{theorem}  \label{thm:main1} 
Let $E$ be an $f$-invariant exceptional component. 
If $E$ is of Dynkin type $\rD$ or $\rE$, then $(\rR 1)$ is always the case for $(f, E)$.   
The same is true when $E$ is of Dynkin type $\rA$ with reflectional symmetry.  
On the other hand, when $E$ is of Dynkin type $\rA$ with trivial symmetry, 
all of the three cases $(\rR 1), (\rR 2), (\rHy)$ can occur.      
\end{theorem}
\par
Formally, after a precise formulation, Theorem \ref{thm:main1} is established as Theorems \ref{thm:rd} and \ref{thm:RHAt}. 
\begin{remark} \label{rem:genuine} 
Let $(f, E)$ be a pair of type $\rA$ with trivial symmetry. 
If $f$ is the square of another automorphism $g$ such that $(g, E)$ is of type $\rA$ with reflectional symmetry,  
then Theorem \ref{thm:main1} implies that $(f, E)$ falls in case $(\rR 1)$. 
If $(f, E)$ is not of this form but still lies in case $(\rR 1)$, 
then we say that it is {\sl genuinely} of case $(\rR 1)$.
\end{remark} 
\par
Finding a lot of K3 surface automorphisms illustrating Theorem \ref{thm:main1}
is also a main goal of this article.  
\begin{result} \label{res:main2} 
We give examples of $K3$ surface automorphisms in $\S \ref{sec:mhgg}$ and 
detect rotation domains on them in $\S \ref{sec:rd}$.    
\begin{enumerate}
\setlength{\itemsep}{-1pt}
\item[$(\mathrm{i})$] Table $\ref{tab:Ass}$ gives examples having an    
exceptional component of type $\rA$ with trivial symmetry.  
Proposition $\ref{prop:Ass}$ establishes basic properties of these examples, 
and Theorem $\ref{thm:RHAt}$ shows that they cover all of the three cases 
$(\rR 2), (\rR 1), (\rHy)$.  
We emphasize that Table $\ref{tab:Ass}$ contains an example which is genuinely of case $(\rR 1)$.   
\item[$(\mathrm{ii})$] Tables $\ref{tab:Ars1}$ and $\ref{tab:Ars2}$ give examples having 
exceptional components of type $\rA$ with reflectional symmetry. 
Proposition $\ref{prop:Ars1}$ establishes basic properties of these examples. 
They fall in case $(\rR 1)$ by Theorem $\ref{thm:main1}$.    
\item[$(\mathrm{iii})$] Tables $\ref{tab:DEs1}$ and $\ref{tab:DEs2}$ give examples having 
exceptional components of type $\rD$ or $\rE$. 
Propositions $\ref{prop:DEs}$ and $\ref{prop:DEs2}$ establish basic properties of these examples. 
They fall in case $(\rR 1)$ by Theorem $\ref{thm:main1}$.      
\end{enumerate}
Table $\ref{tab:Ass}$ contains examples having a rotation domain of rank $1$ and one of 
rank $2$ at the same time, the former of which contains an exceptional component of 
type $\rA$ and the latter contains one of type $\rE;$ see Theorem $\ref{thm:RHAt}$.   
\end{result}
\par
While Theorem \ref{thm:main1} and Result \ref{res:main2} are concerned with  
phenomena around an exceptional component, it is more classical to 
find a fixed point $p$ {\sl off} the exceptional set $\cE(X)$ and ask whether $p$ is a 
center of some rotation domain or a hyperboluc fixed point.  
This problem is more interesting when we focus not only on fixed points  
but also on cycles of periodic points, i.e. {\sl periodic cycles}, with various primitive periods.  
We remark that all periodic points off $\cE(X)$ are isolated. 
We now give only a rough and partial sketch of the results in this direction.   
\begin{result} \label{res:main3} 
Among other things we obtain the following kinds of results in $\S \ref{sec:pc}$. 
\begin{enumerate}
\setlength{\itemsep}{-1pt}
\item[$(\mathrm{i})$] An automorphism $f$ in Theorem $\ref{thm:rdpc}$ has periodic cycles 
of primitive periods $1, 7, 12, 14$, which are centers of rotation domains of rank $2$;     
a periodic cycle of primitive period $3$, each point of which is a center of a rotation domains of rank $1$;     
and hyperbolic periodic cycles of primitive periods $2, 5, 9, 10, 11$. 
While those cycles are off $\cE(X)$, the exceptional set $\cE(X)$ itself is an $f$-invariant 
exceptional component of type $\rA_2$ with reflectional symmetry, which therefore is 
contained in a rotation domain of rank $1$ by Theorem $\ref{thm:main1}$.  
\item[$(\mathrm{ii})$] Example $\ref{ex:takada}$ presents two automorphisms which 
have no exceptional set, but still admit a rotation domain of rank $1$. 
These maps are found by Y.~Takada in a conversation with the author.   
\end{enumerate}    
\end{result}     
\par
All examples in Results \ref{res:main2} and \ref{res:main3} arise from the {\sl method of 
hypergeometric groups} developed in \cite{IT1,IT2}. 
A survey of this method is given in \S \ref{sec:mhgg}.   
We remark that any K3 surface automorphism created by this method  
satisfies conditions $(\rA 1)$ and $(\rA 2)$. 
This method is suitable for discussing rotation domains on K3 surfaces.  
\par
In order to detect a rotation domain for a given map $f : X \to X$, it is 
standard to linearize $f$ around a fixed point $p \in X$ into its tangent map 
$(d f)_p : T_p X \to T_p X$ and examine if $p$ is a center of some rotation 
domain by looking at the linear map $(d f)_p$.   
This is typically the Siegel-type linearization of a local biholomorphism      
\begin{equation} \label{eqn:nlm}
F : (\bC^2, 0) \to (\bC^2, 0)    
\end{equation} 
into its Jacobian matrix $J := F'(0) \in \GL_2(\bC)$, where $(\bC^2, 0)$ stands 
for an open neighborhood of the origin $0$ in $\bC^2$, not the whole space $\bC^2$.  
In \eqref{eqn:nlm} we usually assume that the multipliers of $F$ at the origin, i.e. the 
eigenvalues of $J$, have modulus $1$ and satisfy a certain Diophantine 
condition; see \S \ref{sec:equilin} for a brief review. 
\par
In \S \ref{sec:equilin} we generalize this idea to that of {\sl equivariant linearization},  
which enables us to deal with a K3 surface automorphism $f : X \to X$ in 
a neighborhood of an $f$-invariant exceptional component $E$.   
A small neighborhood $V$ of $E$ is a minimal resolution of a complex surface 
$S$ with a Kleinian singularity $p \in S$, and $(S, p)$ is isomorphic to 
the quotient space of $(\bC^2, 0)$ by a finite subgroup $G$ of $\SU(2)$.  
Thus the automorphism $f$ restricted to $V$ lifts to a local biholomorphism $F$  
as in \eqref{eqn:nlm}, which is $G$-{\sl equivariant}. 
In this situation we have the following result.  
\begin{theorem} \label{thm:main4}  
Let $G$ be a nontrivial finite subgroup of $\SU(2)$ and suppose that 
the map $F$ in \eqref{eqn:nlm} is $G$-equivariant. 
Then the Jacobian matrix $J := F'(0) \in \GL_2(\bC)$ normalizes $G$, 
that is, $J g J^{-1} \in G$ for every $g \in G$, and $J$ is diagonalizable 
unless $G$ is a group of order $2$. 
Suppose that $J$ is diagonalizable in this exceptional case and hence in any case.   
Under a Diophantine condition on its multipliers, $F$ is $G$-equivariantly 
linearizable to $J$.    
\end{theorem} 
 \par
The terminology used here is explained in \S \ref{sec:equilin} and  
Theorem \ref{thm:main4} is established formally as Theorem \ref{thm:elinz}. 
Any normalizor $J \in \GL_2(\bC)$ of $G$ acting linearly on $\bC^2$ leads back to 
an automorphism 
\begin{equation} \label{eqn:lm}
f_J : (X, E) \to (X, E)
\end{equation} 
in a vicinity of $E$, which we call 
the {\sl linear model} associated with $J$. 
In \S \ref{sec:nmr} we give a detailed description of how $f_J$ acts on the exceptional component $E$. 
To this end, in \S \ref{sec:rs}, we need to recall explicitly how a minimal resolution of a Kleinian 
singularity takes place, although it is rather technical. 
Now we start with an automorphism $f : (X, E) \to (X, E)$, proceed to a $G$-equivariant map $F$ as  
in \eqref{eqn:nlm}, apply Theorem \ref{thm:main4} to $G$-equivariantly linearize $F$ into its Jacobian matrix $J$, 
and come back to $(X, E)$ with a linear model $f_J$ as in \eqref{eqn:lm}.  
As Proposition \ref{prop:cmid} shows, there is a subtle but not negligible difference 
between $f|_E$ and $f_J|_{E}$, although they are almost the same.  
In any case, the linear model $f_J$ can be used to detect rotation domains for 
the original map $f$; see Lemma \ref{lem:rh-lm}.  
\par
For the linearization of the map $F$ in \eqref{eqn:nlm}, in both classical and $G$-equivariant versions,  
it is necessary to judge whether the multiplier pair $\bal = (\alpha_1, \alpha_2)$ of $F$ at $0$ is multiplicatively 
independent (MI), or multiplicatively dependent (MD) but non-resonant (NR), or resonant; see \S \ref{sec:pp} 
for definitions.  
Writing $\delta := \det J = \alpha_1 \alpha_2$, we express the multipliers $\alpha_1$ and $\alpha_2$ 
as $\delta^{1/2} \alpha^{\pm 1}$ with $\alpha \in \bC^{\times}$.  
At this stage there is a sharp contrast between 
\begin{enumerate}
\setlength{\itemsep}{-1pt}
\item[(a)] the equivariant linearization around an exceptional component of type $\rD$ or $\rE$, 
as well as around an exceptional component of type $\rA$ with reflectional symmetry, 
\item[(b)] the equivariant linearization around an exceptional component of type $\rA$ with trivial symmetry,  
as well as the ordinary linearization around an isolated fixed or periodic point off the exceptional set. 
\end{enumerate}
In case (a) the number $\alpha$ is confined into roots of unity which can be determined explicitly, 
while in case (b) the number $\alpha$ remains unknown and yet to be determined from a {\sl global geometry}  
of the automorphism $f : X \to X$.    
To cope with case (b), we impose an algebraic constraint, or ansatz, on the number $\alpha$; see Hypothesis \ref{hyp:P}. 
The ansatz there is to assume that there exists 
a rational function $\hat{P}(w) \in \bQ(w)$ such that   
\begin{equation*} \label{eqn:P-hat}
\alpha + \alpha^{-1} = \hat{P}(\hat{\tau}), \qquad \hat{\tau} := \delta^{1/2} + \delta^{-1/2},      
\end{equation*}
where $\delta$ is supposed to be a non-real conjugate to a Salem number,  
in anticipation of applications to K3 surface automorphisms.  
In this setting, using some properties of Salem numbers, we establish in Theorem \ref{thm:MDR} 
an exact criterion for the judgement of the trichotomy mentiond at the beginning of this paragraph.  
\par
Now the problem is to determine the function $\hat{P}(w)$ over an 
exceptional component $E$ of type $\rA$ with trivial symmetry, or at a fixed or periodic 
point $p$ off $\cE(X)$. 
In the latter situation, it is also necessary to find such a point $p$ beforehand.   
To this end we use two kinds of Lefschetz-type fixed point formulas (FPFs); one is   
a {\sl topological} FPF \eqref{eqn:saito2}, which comes from S.~Saito's FPF \cite{Saito},  
and the other is a {\sl holomorphic} FPF \eqref{eqn:tt2}, which comes from Toledo and Tong's FPF \cite{TT}.  
While these formulas already appear in our papers \cite{IT1,IT2}, this article develops an 
extensive study of the related indices around exceptional components in \S \ref{sec:fpf}.     
As Theorem \ref{thm:ecA} and its use in \S \ref{sec:rd} suggest, the treatment of exceptional components of 
type $\rA$ with trivial symmetry is quite subtle.   
A rich variety of functions $\hat{P}(w)$ make the phenomena in case (b) more colorful  
than those in case (a).  
\par
All the theoretical foundations in \S\S \ref{sec:equilin}--\ref{sec:fpf} are 
combined with the method of hypergeomteric groups (as reviewed in \S \ref{sec:mhgg}) to 
enable us to detect rotation domains and their periodic versions 
in \S\S \ref{sec:rd}--\ref{sec:pc}.   
We remark that, logically, the contents of \S\S \ref{sec:equilin}--\ref{sec:fpf} 
should not be limitted to K3 surface automorphisms---although applications to them   
always stay in our mind---and the results in those sections should be applicable in wider situations. 
On the other hand, the contents proper to K3 surface automorphisms are presented 
in \S\S \ref{sec:mhgg}--\ref{sec:pc}. 
Appendix \ref{app:exp} gives a detailed study of particular examples, 
which is interesting in itself, but digresses somewhat from the main stream of the article. 
Appendix \ref{app:rfp} provides some data used in \S \ref{sec:pc}, which is too messy 
to be included in the main text.  
This rather long sketch of the article also speaks about its organization,  
so there is no need to repeat it. 
Throughout the paper $\bN := \bZ_{\ge 1}$ stands for the set of all {\sl positive} integers; 
$\ri := \sqrt{-1}$ is the imaginary unit.                           
\section{Equivariant Linearization} \label{sec:equilin} 
Let $\alpha_1$ and $\alpha_2$ be nonzero complex numbers, which are referred 
to as {\sl multipliers}.  
The multiplier pair $\bal = (\alpha_1, \alpha_2)$ is said to be {\sl Diophantine} 
if there exist positive constants $K$ and $c$ such that 
$$
|\alpha_j -\alpha_1^{n_1} \alpha_2^{n_2}| \ge K (n_1 + n_2)^{-c} \qquad \mbox{for any} 
\quad (n_1, n_2) \in N, \quad j = 1, 2,  
$$
where $N := \{ (n_1, n_2) \in \bZ^2 \mid n_1 \ge 0, \, n_2 \ge 0, \, n_1 + n_2 \ge 2 \}$. 
A Siegel-type linearization theorem for nonlinear mappings, e.g. a theorem of 
P\"{o}schel \cite{Poschel} restricted to $2$-dimensions, states the following result.  
\begin{theorem} \label{thm:Poschel}
Let $F : (\bC^2, 0) \to (\bC^2, 0)$ be a local biholomorphism with Jacobian matrix 
$J = F'(0)$ at the origin, and let $\alpha_1$ and $\alpha_2$ be the eigenvalues of $J$.  
If $J$ is diagonalizable and the pair $\bal = (\alpha_1, \alpha_2)$ is 
Diophantine, then there exists a unique local biholomorphism 
$\Phi : (\bC^2, 0) \to (\bC^2, 0)$ such that  
\begin{equation} \label{eqn:poschel}
\Phi \circ F \circ \Phi^{-1} = J, \qquad \Phi'(0) = I. 
\end{equation}
\end{theorem}
\par
If $\bal$ is Diophantine, then it is obviously {\sl non-resonant} (NR for short) 
to the effect that  
$$
\alpha_j -\alpha_1^{n_1} \alpha_2^{n_2} \neq 0 \qquad \mbox{for any} \quad 
(n_1, n_1) \in N, \quad j = 1, 2.  
$$
Transcendence theory \cite{Baker} tells that there is an important 
situation where the converse implication is also true.  
\begin{theorem} \label{thm:tnt} 
Let $\alpha_1$ and $\alpha_2$ be algebraic numbers of modulus $1$. 
If $\bal$ is NR then it is Diophantine.  
\end{theorem}
\par
We shall establish an equivariant version of Theorem \ref{thm:Poschel} for  
finite subgroups of $\SU(2)$.  
To this end, let us recall the classification of nontrivial finite subgroups of $\SU(2)$;  
see e.g. Lamotke \cite{Lamotke} for details.  
\begin{theorem} \label{thm:bclassf} 
Any nontrivial finite subgroup $G$ of $\SU(2)$ is one of the groups 
listed in Table $\ref{tab:bclassf}$, where the first column exhibits the  
Dynkin diagram $D$ of the exceptional set for the minimal resolution of the 
quotient singularity $\bC^2/G$, and binary dihedral groups of order $8$ are distinguished 
from other binary dihedral groups as quaternion groups.   
Any two groups of the same Dynkin type are conjugate to each other in $\SU(2)$.   
\end{theorem}
\begin{table}[h]
\centerline{
\begin{tabular}{llcc}
\hline
$D_{\phantom{n}}$  & groups $G$ & order & $n$ \\
\hline
$\rA_n$ & cyclic groups & $n+1$ & $\ge 1$ \\
\hline
$\rD_4$ & quaternion groups & $8$ & $4$ \\
\hline
$\rD_n$ & binary dihedral groups & $4(n-2)$ & $\ge 5$ \\
\hline
$\rE_6$ & binary tetrahedral groups & $24$ & $6$ \\
\hline
$\rE_7$ & binary octahedral groups & $48$ & $7$ \\
\hline
$\rE_8$ & binary icosahedral groups & $120$ & $8$ \\
\hline
\end{tabular}}
\caption{Classification of nontrivial finite subgroups of $\SU(2)$.} 
\label{tab:bclassf}
\end{table}
\par
Let $F : (\bC^2, 0) \to (\bC^2, 0)$ be a local biholomorphism at the origin and  
let $G$ be a nontrivial finite subgroup of $\SU(2)$. 
We say that $F$ is $G$-{\sl equivariant} 
if for any $g \in G$ and any $z \in \bC^2$ sufficiently close to the origin 
$(F \circ g \circ F^{-1})(z)$ stays in the $G$-orbit through $z$. 
The aim of this section is to establish the following result.    
\begin{theorem} \label{thm:elinz} 
Let $G$ be a nontrivial finite subgroup of $\SU(2)$ and 
$F : (\bC^2, 0) \to (\bC^2, 0)$ be a $G$-equivariant local biholomorphism with 
Jacobian matrix $J := F'(0) \in \GL_2(\bC)$ at the origin. 
Then $J$ normalizes the group $G$, that is, $J g J^{-1} \in G$ for every $g \in G$, 
and $J$ is diagonalizable unless $G$ is a group of order $2$.   
Suppose that 
\begin{enumerate}
\setlength{\itemsep}{-1pt}
\item[$(\mathrm{i})$] $J$ is diagonalizable when $G$ is a group of order $2$, so it is diagonalizable in any case,   
\item[$(\mathrm{ii})$] the pair $\bal = (\alpha_1, \alpha_2)$ is Diophantine, where $\alpha_1$ and $\alpha_2$ 
are the eigenvalues of $J$. 
\end{enumerate} 
Then there exists a unique local biholomorphism 
$\Phi : (\bC^2, 0) \to (\bC^2, 0)$ that linearizes $F$ into $J$ as in equation \eqref{eqn:poschel}. 
The map $\Phi$ commutes with every element of $G$. 
\end{theorem} 
\par
The proof of this theorem is completed at the end of this section.  
We begin with the following lemma, which is helpful in investigating the concept of $G$-equivariance.  
\begin{lemma} \label{lem:free-a} 
For $g \in \SU(2)$ we put $K(g) := |\lambda -1| \ge 0$ for an eigenvalue $\lambda$ of $g$. 
The number $K(g)$ is well defined, i.e. independent of the choice of $\lambda$.  
We have $K(g) = 0$ if and only if $g = 1$, the unit element.   
Moreover,   
\begin{equation} \label{eqn:free-a}
|g z -z| = K(g) \, |z| \qquad \mbox{for any} \quad z \in \bC^2, 
\end{equation}
where $|z|$ is the standard Euclidean norm of $z$. 
In particular the group $\SU(2)$ acts on $\bC^2 \setminus \{ 0 \}$ freely.  
\end{lemma}
{\it Proof}. 
The eigenvalues of $g \in \SU(2)$ are a complex number $\lambda$ of modulus $1$ and 
its reciprocal $\lambda^{-1}$, so we have $|\lambda^{-1} -1| = |\lambda^{-1}(1-\lambda)| = |\lambda -1|$. 
Thus the value of $K(g)$ is independent of the eigenvalue $\lambda$ chosen.   
There exists an orthonormal basis $u$, $v$ of $\bC^2$ such that 
$g u = \lambda u$ and $g v = \lambda^{-1} v$.  
Notice that $g = 1$ if and only if $\lambda = 1$, that is, $K(g) = 0$.   
Writing $z \in \bC^2$ as a linear combination $z = a u + b v$ with coefficients $a$, $b \in \bC$, we have   
$g z -z = (\lambda -1) a u + (\lambda^{-1} -1) b v$ and hence $|g z - z|^2 = |(\lambda-1) a|^2 + |(\lambda^{-1}-1) b|^2  
= |\lambda-1|^2 (|a|^2+|b|^2) = K(g)^2 |z|^2$, which gives formula \eqref{eqn:free-a}.  
If $g \in \SU(2)$ fixes a nonzero vector $z \in \bC^2$, then formula \eqref{eqn:free-a} 
implies $K(g) = 0$, which in turn forces $g = 1$. 
Thus the group $\SU(2)$ acts on $\bC^2 \setminus \{0\}$ freely.  
\hfill $\Box$ 
\begin{lemma} \label{lem:jacobian} 
If $F$ is $G$-equivariant, then its Jacobian matrix $J := F'(0) \in \GL_2(\bC)$ 
at the origin normalizes the group $G$, that is, $J g J^{-1} \in G$ for every $g \in G$, 
and we have   
\begin{equation} \label{eqn:jacobian}
F \circ g \circ F^{-1} = J g J^{-1} \quad 
\mbox{for every} \quad g \in G. 
\end{equation}
\end{lemma}
{\it Proof}. 
Since $F$ is $G$-equivariant, for any $g \in G$ and any $z \in \bC^2 \setminus \{0\}$ 
near the origin there exists an element $\chi(g; z) \in G$ such that 
$(F \circ g \circ F^{-1})(z) = \chi(g; z) z$. 
Note that $\chi(g; z)$ is unique since $G$ acts on $\bC^2 \setminus \{0\}$ freely. 
\par\smallskip
{\bf Claim}. The element $\chi(g; z) \in G$ is independent of $z$ in a connected 
punctured neighborhood of the origin.  
\par\smallskip\noindent
For any $z$, $w \in \bC^2 \setminus \{0\}$ near the origin it follows from 
formula \eqref{eqn:free-a} and $\chi(g; w) \in \SU(2)$ that 
\begin{align*}
K(\chi(g; w)^{-1} \chi(g; z)) \, |z| &= |\chi(g; w)^{-1} \chi(g; z) z - z| =  |\chi(g; z) z - \chi(g; w) z| \\ 
&= |\chi(g; z) z - \chi(g; w) w + \chi(g; w) (w-z)| \\
&\le |\chi(g; z) z - \chi(g; w) w | + | \chi(g; w) (w-z)| \\
&= |(F \circ g \circ F^{-1})(z) -(F \circ g \circ F^{-1})(w)| + |w-z|.  
\end{align*}
Fix $z \neq 0$ and let $w$ tend to $z$. 
Since $(F \circ g \circ F^{-1})(w)$ is continuous in $w$, there exists a positive 
constant $\ve$ such that 
$|(F \circ g \circ F^{-1})(z) -(F \circ g \circ F^{-1})(w)| + |w-z| < K_0 |z|$ 
whenever $|w - z| < \ve$, where $K_0$ is defined by  
$$
K_0 := \min \{ K(h) \mid h \in G \setminus \{1\} \} > 0.  
$$
For some $w$ with $|w-z| < \ve$, if $\chi(g; w) \neq \chi(g; z)$, then 
$\chi(g; w)^{-1} \chi(g; z) \in G \setminus \{1 \}$ and the estimate above yields 
$$
K_0 |z| \le K(\chi(g; w)^{-1} \chi(g; z)) \, |z| \le |(F \circ g \circ F^{-1})(z) -(F \circ g \circ F^{-1})(w)| + |w-z| < K_0 |z|, 
$$
which is a contradiction. 
Thus $\chi(g; w) = \chi(g; z)$ whenever $|w-z| < \ve$. 
Namely $\chi(g; z)$ is a locally constant function of $z$ in a punctured neighborhood of the origin.  
This establishes the Claim.  
\par
By the Claim there exists a unique map $\chi : G \to G$, $g \mapsto \chi(g)$ such that 
$(F \circ g \circ F^{-1})(z) = \chi(g) z$ for every $z \in \bC^2$ near the origin (including $z = 0$).  
Differentiating this equation with respect to $z$ at the origin, we have 
$J g J^{-1} = \chi(g) \in G$ for every $g \in G$.  
This shows that $J \in \GL_2(\bC)$ normalizes $G$ and formula \eqref{eqn:jacobian} holds true. 
\hfill $\Box$ \par\medskip
Lemma \ref{lem:jacobian} tempts us to determine the normalizer $N_{\GL}$ of $G$ in $\GL_2(\bC)$. 
It suffices to know the normalizer $N_{\SL}$ of $G$ in $\SL_2(\bC)$, as  
$N_{\GL} = \bC^{\times} \cdot N_{\SL}$. 
When $G$ is non-abelian, the following lemma is useful.  
\begin{lemma} \label{lem:iwasawa} 
If $G$ is a non-abelian finite subgroup of $\SU(2)$ then $N_{\SL}$ is a finite subgroup of $\SU(2)$.    
\end{lemma}
{\it Proof}. 
We begin by showing the following fact:   
Let $g \in \SL_2(\bC)$ and $u_j \in \SU(2)$ for $j = 1, 2$. 
\par\smallskip  
{\bf Claim 1}. Suppose that $u_1$ and $u_2$ do not commute.   
If  $g u_j g^{-1} \in \SU(2)$ for $j = 1, 2$, then $g \in \SU(2)$. 
\par\smallskip\noindent
Gram-Schmidt orthogonalization allows us to write $g = t h$ for an  
$h \in \SU(2)$ and an upper triangular matrix 
$$
t = 
\begin{pmatrix} r & s \\[1mm]  
0 & r^{-1} 
\end{pmatrix} \qquad \mbox{with} \quad  
r > 0, \quad s \in \bC. 
$$
Put $v_j := h u_j h^{-1} \in \SU(2)$ for $j = 1, 2$. 
Then $t v_j t^{-1} = g u_j g^{-1} \in \SU(2)$ for $j = 1, 2 $. 
Since $u_1$ and $u_2$ are non-commutative, so are $v_1$ and $v_2$, 
hence  at least one of them, say $v_j$, is not diagonal.  
Here we have  
$$
t v_j t^{-1} 
= 
\begin{pmatrix}
\alpha + r^{-1} s \beta & -s \{ r(\alpha-\bar{\alpha}) + s \beta \} 
- r^2 \bar{\beta} \\[1mm]
r^{-2} \beta & \bar{\alpha} 
\end{pmatrix} 
\quad 
\mbox{for} 
\quad 
v_j = 
\begin{pmatrix*}[r] 
\alpha & - \bar{\beta} \\[1mm] 
\beta & \bar{\alpha} 
\end{pmatrix*},  
$$
where $|\alpha|^2 + |\beta|^2 = 1$ and $\beta \neq 0$.  
In view of $t v_j t^{-1} \in \SU(2)$, the $(1, 1)$-entry of this 
matrix must be the complex conjugate to the $(2, 2)$-entry. 
This forces $\alpha + r^{-1} s \beta = \alpha$, hence $s = 0$ and 
$$
t v_j t^{-1} 
= 
\begin{pmatrix}
\alpha & - r^2 \bar{\beta} \\[1mm]
r^{-2} \beta & \bar{\alpha} 
\end{pmatrix} 
\in \SU(2).   
$$
So the $(1, 2)$-entry of this matrix must be the negative of the complex 
conjugate to the $(2, 1)$-entry. 
This forces $r^2 = r^{-2}$ and so $r = 1$, as $r$ is real and positive.  
Thus $t = I$ and $g = h \in \SU(2)$, which proves the claim. 
\par
Claim $1$ implies that $N := N_{\SL}$ is a subgroup of $\SU(2)$, since $G$ is non-abelian.  
Consider the homomorphism 
$$
\Ad_G : N \to \Aut(G), \,\, n \mapsto \Ad_G(n), \quad 
\mbox{defined by} \quad \Ad_G(n) g := n g n^{-1} \quad \mbox{for} \quad g \in G.    
$$
\par
{\bf Claim 2}. We have $\Ker(\Ad_G) = \{ \pm I\}$. 
\par\smallskip\noindent
Indeed, the inclusion $\{ \pm I\} \subset \Ker(\Ad_G)$ is obvious. 
Conversely, let $h$ be any element of $\Ker(\Ad_G)$. 
Then $h$ commutes with every element in $G$. 
However, any non-abelian finite subgroup of $\SU(2)$, in particular 
the group $G$ under consideration, acts on $\bC^2$ irreducibly. 
Thus Schur's lemma implies that $h$ must be a scaler matrix. 
Since $\det(h) = 1$, we have $h = \pm I$ and hence  
$\Ker(\Ad_G) \subset \{ \pm I\}$.  
This proves Claim $2$.  
\par
Let $\Ad_G(N) \subset \Aut(G)$ be the image of   
$\Ad_G : N \to \Aut(G)$. 
Claim $2$ yields the short exact sequence 
$$
\begin{CD} 
1 @>>> \{ \pm I \} @> \mathrm{inclusion} >> N @> \Ad_G >> \Ad_G(N) @>>> 1. 
\end{CD}
$$
Since $G$ is a finite group, its automorphism group $\Aut(G)$ is also finite, 
and so is its subgroup $\Ad_G(N)$. 
Now the exact sequence shows that $N$ is also a finite group 
with order $|N| = 2 |\Ad_G(N)| \le 2 |\Aut(G)| < \infty$. \hfill $\Box$
\begin{remark} \label{rem:bsupgroup} 
By the {\sl trivial subgroups} of a binary group we mean $\{ I\}$, $\{\pm I \}$ 
and the group itself.  
\begin{enumerate}
\setlength{\itemsep}{-1pt}
\item Any binary octahedral group contains exactly two nontrivial normal 
subgroups, which are a quaternion group and a binary tetrahedral group.  
\item No binary octahedral group is a proper subgroup of a finite subgroup of $\SU(2)$. 
Similarly, no binary icosahedral group is a proper subgroup of a finite subgroup of $\SU(2)$.  
\end{enumerate}
\end{remark}
\begin{theorem} \label{thm:normalizer} 
If $G$ be a nontrivial finite subgroup of $\SU(2)$, then its normalizer $N_{\GL}$ in $\GL_2(\bC)$ 
is the nonzero scalar multiples of its normalizer $N_{\SL}$ in $\SL_2(\bC)$, that is,  
$N_{\GL} = \bC^{\times} \cdot N_{\SL}$, where $N_{\SL}$ is given as follows.  
\begin{enumerate}
\setlength{\itemsep}{-1pt} 
\item Suppose that $G$ is a cyclic group of order $n+1$. 
If $n = 1$ then $N_{\SL} = \SL_2(\bC)$.  
If $n \ge 2$ and $G$ is generated by $A(\ve)$ with $\ve = \exp(2\pi \ri/(n+1))$, then 
$N_{\SL} = \{ A(a) \mid a \in \bC^{\times} \} \cup \{ B(b) \mid b \in \bC^{\times}\}$, 
where  
\begin{equation*}
A(a) := 
\begin{pmatrix} 
a & 0 \\[1mm] 
0 & a^{-1} 
\end{pmatrix}, 
\qquad  
B(b) := 
\begin{pmatrix} 
0 & b \\[1mm] -b^{-1} & 0 
\end{pmatrix}.   
\end{equation*}
\item If $G$ is a binary dihedral group of order $4(n-2)$ with $n \ge 5$, 
then $N_{\SL}$ is the unique binary dihedral group of order $8(n-2)$ that contains 
$G$ as a subgroup of index $2$.  
\item If $G$ is a quaternion group or a binary tetrahedral group, then 
$N_{\SL}$ is the unique binary octahedral group that contains $G$ as a normal 
subgroup. 
\item If $G$ is a binary octahedral group or a binary icosahedral group, 
then $N_{\SL}$ is the group $G$ itself.  
\end{enumerate}
Every matrix in $N_{\GL}$ is diagonalizable, unless $G$ is a group of order $2$.    
\end{theorem}
{\it Proof}. 
Assertion (1).  
The case $n = 1$ is obvious. 
Suppose that $n \ge 2$. 
A direct calculation shows that 
$$
g \, A(\ve) \, g^{-1} = 
\begin{pmatrix} 
a d \ve - b c \ve^{-1} & -a b (\ve - \ve^{-1}) \\[1mm]
c d (\ve - \ve^{-1}) & a d \ve^{-1} - b c \ve 
\end{pmatrix} 
\quad \mbox{for} \quad 
g = 
\begin{pmatrix} 
a & b \\[1mm]
c & d
\end{pmatrix}, 
\quad a d - b c = 1.   
$$
Notice that $\ve - \ve^{-1}$ is nonzero, since $\ve = \exp( 2 \pi \ri/(n+1))$ 
with $n \ge 2$ is different from $\pm 1$.  
If $g \in N_{\SL}$ then $g \, A(\ve) \, g^{-1}$ must be diagonal, so we have 
$a b = 0 = c d$. 
When $b=0$, we have $a d = 1$ and $c = 0$, hence 
$a \in \bC^{\times}$, $g = A(a)$ and $g \, A(\ve) \, g^{-1} = A(\ve) \in G$. 
When $a =0$, we have $b c = -1$ and $d = 0$, hence 
$b \in \bC^{\times}$, $g = B(b)$ and $g \, A(\ve) \, g^{-1} = A(\ve)^{-1} \in G$. 
Thus $N_{\SL} = \{ A(a) \mid a \in \bC^{\times} \} \cup \{ B(b) \mid b \in \bC^{\times}\}$.    
\par
Assertion (2).  
We may assume that $G$ is generated by $A(\ve)$ and $B(1)$, where 
$\ve = \exp(\pi \ri/(n-2))$. 
In this case $G$ has a total of $4(n-2)$ elements  
$A(\ve^j)$ and $B(\ve^j)$, $j = 1, \dots, 2(n-2)$, with all  
$B(\ve^j)$ being of order $2$.  
If $g \in N_{\SL}$ then $g \, A(\ve) \, g^{-1} \in G$ must be diagonal, 
because it is an element of order $2(n-2) \ge 6 > 2$.  
Then the same reasoning as in the proof of assertion (1) shows that $g$ 
must be either $A(a)$ or $B(a)$ for some $a \in \bC^{\times}$. 
When $g = A(a)$, we have $g \, A(\ve^j) \, g^{-1} = A(\ve^j)$ and 
$g \, B(\ve^j) \, g^{-1} = B(a^2 \ve^j)$.  
 When $g = B(a)$, we have $g \, A(\ve^j) \, g^{-1} = A(\ve^{-j})$ and 
$g \, B(\ve^j) \, g^{-1} = B(a^2 \ve^{-j})$. 
In either case we have $g \, B(\ve^j) \, g^{-1} \in G$ for $j = 1, \dots, 2(n-2)$ 
if and only if $a^2$ is a $2(n-1)$-th root of unity, that is, $a$ is a $4(n-2)$-th root of unity.  
Thus we have 
$$
N_{\SL} = \{ A(a) \mid a^{4(n-2)} = 1 \} \cup \{ B(a) \mid a^{4(n-2)} = 1 \}. 
$$
This means that $N_{\SL}$ is a binary dihedral group of order $8(n-2)$ that contains $G$. 
Let $G'$ be any binary dihedral group of order $8(n-2)$ that contains $G$. 
Then $G$ is of index $2$ in $G'$ and hence a normal subgroup of $G' \subset \SL_2(\bC)$.  
This fact forces $G' \subset N_{\SL}$ and so $G' = N_{\SL}$, because $G'$ and $N_{\SL}$ 
are of the same order.  
Therefore $N_{\SL}$ is the unique binary dihedral group of order $8(n-2)$ that contains $G$. 
\par
Assertion (3).  
We begin by the case where $G$ is a quaternion group. 
Take a binary octahedral group $O_0$. 
It contains a quaternion group $G_0$ as a normal subgroup by 
Remark \ref{rem:bsupgroup}.(1).  
Since any two quaternion groups are conjugate to each other in $\SU(2)$, 
there exists a conjugation $\varphi : \SU(2) \to \SU(2)$ that converts $G_0$ to $G$.  
Let $O \subset \SU(2)$ be the image of $O_0$ under the map $\varphi$. 
Then $O$ is a binary octahedral group that contains $G$ as a normal subgroup. 
Hence $O$ is a subgroup of the normalizer $N_{\SL}$ of $G$ in $\SL_2(\bC)$.  
Lemma \ref{lem:iwasawa} implies $O \subset N_{\SL} \subset \SU(2)$, 
since $G$ is non-abelian. 
This inclusion forces $O = N_{\SL}$ by Remark \ref{rem:bsupgroup}.(2).  
We also have the uniqueness of the binary octahedral group that contains 
$G$ as a normal subgroup, because any such group coincides with the 
normalizer $N_{\SL}$. 
This proves assertion (3) for quaternion groups. 
Assertion (3) for binary tetrahedral groups can be established just in the same manner.  
\par
Assertion (4). 
The proof of this assertion is analogous to that of assertion (3) or even simpler than that.     
Indeed, if $G$ is a binary octahedral group or a binary icosahedral group, 
then Lemma \ref{lem:iwasawa} implies $G \subset N_{\SL} \subset \SU(2)$. 
This inclusion forces $G = N_{\SL}$ by Remark \ref{rem:bsupgroup}.(2). 
\par
The final statement in the theorem follows from assertion (1) when $G$ is a cyclic group 
of order $\ge 3$, since $A(a)$ is diagonal and $B(b)$ is diagonalizable.  
It follows from Lemma \ref{lem:iwasawa} when $G$ is a non-cyclic group and 
hence a non-abelian group by Theorem \ref{thm:bclassf}, 
since every matrix in $\SU(2)$ is diagonalizable. \hfill $\Box$    
\begin{lemma} \label{lem:el} 
Suppose that the $G$-equivariant map $F$ is linearized into its Jacobian matrix 
$J$ by a local biholomorphism $\Phi$ as in formula \eqref{eqn:poschel}.  
Let $\Psi : (\bC^2, 0) \to (\bC^2, 0)$ be the local holomorphism  
defined by 
\begin{equation} \label{eqn:el2}
\Psi := \frac{1}{|G|} \sum_{g \in G} \Phi^g,    
\end{equation}
where $\Phi^g := g^{-1} \circ \Phi \circ g$ for $g \in G$. 
If $\Psi^g$ is defined in a similar manner, then $\Psi$ satisfies the equations    
\begin{equation} \label{eqn:el3} 
\Psi \circ F \circ \Psi^{-1} = J, \qquad \Psi'(0) = I, \qquad 
\Psi^g = \Psi \quad \mbox{for every} \quad g \in G. 
\end{equation}
\end{lemma}
{\it Proof}. 
Putting $\chi(g) := J g J^{-1}$ for $g \in G$, we begin by showing that  
\begin{equation} \label{eqn:el4} 
\Psi = \frac{1}{|G|} \sum_{g \in G} \Phi^{\chi(g)}, \qquad     
\Phi^{\chi(g)} \circ F = J \circ \Phi^g \quad \mbox{for} \quad g \in G.  
\end{equation}
By Lemma \ref{lem:jacobian} we have $\chi \in \Aut(G)$ and hence the first equality in \eqref{eqn:el4}.    
The second equality follows from     
\begin{align*}
\Phi^{\chi(g)} \circ F 
&= \chi(g)^{-1} \circ \Phi \circ \chi(g) \circ F 
= \chi(g^{-1}) \circ \Phi \circ (F \circ g \circ F^{-1}) \circ F   
= \chi(g^{-1}) \circ (\Phi \circ F) \circ g \\
&= (J g^{-1} J^{-1}) \circ (J \circ \Phi) \circ g = J \circ (g^{-1} \circ \Phi \circ g) 
= J \circ \Phi^g,  
\end{align*}   
where formulas \eqref{eqn:jacobian} and \eqref{eqn:poschel} are used in the second and 
fourth equalities. 
Formulas \eqref{eqn:el4} then lead to    
\begin{align*} 
\Psi \circ F = \frac{1}{|G|} \sum_{g \in G} \Phi^{\chi(g)} \circ F 
= \frac{1}{|G|} \sum_{g \in G} J \circ \Phi^g  
= J \circ \left( \frac{1}{|G|} \sum_{g \in G} \Phi^g \right) = J \circ \Psi.  
\end{align*}
It is evident from $\Phi(0) = I$ in \eqref{eqn:poschel} and definition \eqref{eqn:el2} that $\Psi(0) = I$. 
This in particular implies that $\Psi$ is a local biholomorphism around the origin, so the 
equation $\Psi \circ F = J \circ \Psi$ becomes $\Psi \circ F \circ \Psi^{-1} = J$.  
Finally we have 
\begin{align*}
\Psi^g &= g^{-1} \circ \left(  \frac{1}{|G|} \sum_{h \in G} \Phi^h \right) \circ g 
=  \frac{1}{|G|} \sum_{h \in G} g^{-1} \circ \Phi^h \circ g 
= \frac{1}{|G|} \sum_{h \in G} \Phi^{h g} = \Psi    
\end{align*}
for every $g \in G$. 
Therefore all equations in \eqref{eqn:el3} are verified. 
\hfill $\Box$ \par\medskip
{\it Proof of Theorem $\ref{thm:elinz}$}. 
The first half of Theorem \ref{thm:elinz} follows from Lemma \ref{lem:jacobian} 
and the last part of Theorem \ref{thm:normalizer}.  
Now $J$ is always diagonalizable by assumption (i) in Theorem \ref{thm:elinz}.  
By Theorem \ref{thm:Poschel} together with assumption (ii), there exists a unique 
local biholomorphism $\Phi : (\bC^2, 0) \to (\bC^2, 0)$ satisfying condition \eqref{eqn:poschel}.  
The first two equations in \eqref{eqn:el3} show that the map $\Psi$ in definition \eqref{eqn:el2} 
is also such a local biholomorphism.  
So the uniqueness assertion in Theorem \ref{thm:Poschel} forces $\Phi = \Psi$. 
This together with the third equation in \eqref{eqn:el3} implies that $\Phi^g = \Phi$, 
that is, $\Phi$ commutes with every element $g \in G$.  \hfill $\Box$ 
\section{Resolutions of Singularities} \label{sec:rs}
For a nontrivial finite subgroup $G$ of $\SU(2)$ the orbit space $S := \bC^2/G$ 
has a Kleinian singularity at the origin.    
Constructing a minimal resolution $X$ of $S$ is a classical topic in singularity theory and 
a detailed explanation of it can be found e.g. in Lamotke \cite{Lamotke}. 
We review it in a manner suitable for our purpose.  
Explicit description of an atlas on $X$, although somewhat technical, is necessary 
for our calculations in later sections.        
\subsection{Cyclic Cases} \label{ss:cc1}
Let us consider the case where $G$ is a cyclic group of order $n+1$ with $n \in \bN$. 
We may assume that $G$ is generated by the diagonal matrix     
$g = \diag(\zeta^{-1}, \, \zeta)$, where $\zeta = \exp(2 \pi \ri/(n+1))$.  
The group $G$ acts on the polynomial ring $\bC[w_1, w_2]$ in a natural manner 
and the ring $\bC[w_1, w_2]^G$ of $G$-invariant polynomials is generated by 
\begin{equation} \label{eqn:xyz}
x := w_1 w_2,  \qquad y := w_1^{n+1}, \qquad z := w_2^{n+1},  
\end{equation}
which satisfy the relation $x^{n+1} = y z$. 
So the orbit space $S := \bC^2/G$ is realized as an affine algebraic surface 
$$
S = \{ (x, y, z) \in \bC^3 \mid x^{n+1} - y z = 0 \},  
$$
which has a Kleinian singularity of Dynkin type $\rA_n$ at the origin $0 = (0, 0, 0) \in S$.  
\par
Let $X$ be the set of all points $p = ((x, y, z), [a_1: b_1], \dots, [a_n : b_n]) \in \bC^3 \times (\bP^1)^n$ 
which satisfy the equations 
$$ 
a_i b_{i+1} x = a_{i+1} b_i, \qquad 0 \le i \le n, 
$$
with $a_0 := 1$, $a_{n+1} := y$, $b_0 := z$, $b_{n+1} := 1$ by convention.  
Then a minimal resolution of $S$ is given by
\begin{equation*} \label{eqn:proj}
X \to S, \quad p = ((x, y, z), [a_1: b_1], \dots, [a_{n} : b_{n}]) \mapsto (x, y, z).    
\end{equation*} 
Its exceptional set $E$, that is, the inverse image of the origin under this map,  
is the union of $(-2)$-curves  
\begin{equation*} \label{eqn:Ei}
E_i := \{ ((0,0,0), \overbrace{[1:0], \dots, [1:0]}^{i-1}, [a_i : b_i], \overbrace{[0:1], \dots, [0:1]}^{n-i}) 
\mid [a_i:b_i] \in \bP^1 \}, \qquad 1 \le i \le n.    
\end{equation*}
The dual graph of $E$ is a Dynkin diagram of type $\rA_n$, that is, 
\begin{equation} \label{eqn:dA}
\xymatrix@C=15pt{
E_1 \ar@{-}[r] & E_2 \ar@{-}[r] & \cdots \ar@{-}[r] & E_i \ar@{-}[r] & \cdots \ar@{-}[r] & E_{n-1} \ar@{-}[r] & E_n}  
\end{equation} 
where the consecutive $(-2)$-curves $E_i$ and $E_{i+1}$ meet transversally in a single point 
\begin{equation} \label{eqn:pi}
p_i := ((0,0,0), \overbrace{[1:0], \dots, [1:0]}^{i}, \overbrace{[0:1], \dots, [0:1]}^{n-i}), \qquad 1 \le i \le n-1.
\end{equation} 

\par
The smooth surface $X$ is covered by the complex coordinate system $\{ (W_i, \varphi_i) \}_{i=1}^{n+1}$ 
consisting of      
\begin{gather*} \label{eqn:W} 
W_i := \{ p = ((x, y, z), [a_1: b_1], \dots, [a_{n} : b_{n}]) \in X 
\mid a_1 \cdots a_{i-1} b_i \cdots b_{n} \neq 0 \}, 
\\[1mm] 
\varphi_i : W_i \to \bC^2, \quad p \mapsto (u_i, v_{i-1}) := (a_i/b_i, b_{i-1}/a_{i-1}), \qquad 1 \le i \le n+1.  
\end{gather*}
A little calculation shows that the inverse mapping of $\varphi_i$ is given by  
\begin{equation*} \label{eqn:resol2-2}
\varphi_i^{-1} : \bC^2 \to W_i, \quad  
(u_i, v_{i-1}) \mapsto ((x, y, z), [1: v_1], \dots, [1 : v_{i-1}], [u_i : 1], \dots, [u_{n} : 1]),   
\end{equation*} 
where $v_j := u_i^{i-j-1} v_{i-1}^{i-j}$ for $1 \le j \le i-1$,  
$u_j := u_i^{j-i+1} v_{i-1}^{j-i}$ for $i \le j \le n$, and   
\begin{equation} \label{eqn:resol2-1a} 
x := u_i v_{i-1}, \qquad y := u_i^{n-i+2} v_{i-1}^{n-i+1}, \qquad z := u_i^{i-1} v_{i-1}^i, 
\qquad 1 \le i \le n+1.   
\end{equation}
Notice also that the inversion of formula \eqref{eqn:resol2-1a} is given by 
\begin{equation} \label{eqn:resol2-3}
u_i = x^i z^{-1} = x^{i-n-1} y, \qquad v_{i-1} = x^{1-i} z = x^{n+2-i} y^{-1}, \qquad 
1 \le i \le n+1. 
\end{equation} 
\par
For $1 \le i < j \le n+1$ the intersection of $W_i$ and $W_j$ is described as follows:      
\begin{itemize} 
\item If $j = i+1$ then $\varphi_i(W_i \cap W_{i+1}) = \{ u_i \neq 0 \}$ and $\varphi_{i+1}(W_i \cap W_{i+1}) = \{ v_i \neq 0 \}$. 
\item If $j \ge i+2$ then $\varphi_i(W_i \cap W_j) = \{ u_i v_{i-1} \neq 0 \}$ and $\varphi_j(W_i \cap W_j) = \{ u_j v_{j-1} \neq 0 \}$.      
\end{itemize}
For $1 \le i \neq j \le n+1$ the coordinate transformation from $(W_i, \varphi_i)$ to $(W_j, \varphi_j)$ is given by   
$$
\varphi_j \circ \varphi_i^{-1} :  \varphi_i(W_i \cap W_j) \to \varphi_j(W_i \cap W_j), \qquad   
u_j = u_i^{j-i+1} v_{i-1}^{j-i}, \quad  v_{j-1} = u_i^{i-j} v_{i-1}^{i-j+1}.   
$$
\par
We remark that definition \eqref{eqn:pi} makes sense for $i = 0$ and for $i = n$, allowing to define  
two additional points $p_0$ and $p_n$.   
If we consider the affine curves $E_0 := \{ u_1 = 0 \}$ in $W_1$ and $E_{n+1} := \{ v_n = 0\}$ in $W_{n+1}$,   
then 
\begin{equation} \label{eqn:p0pn}
E_0 \cap E_1 = \{ p_0 \} \leftrightarrow (u_1, v_0) = (0, 0), \qquad 
E_n \cap E_{n+1} = \{ p_n \} \leftrightarrow (u_{n+1}, v_n) = (0, 0).  
\end{equation} 
The points $p_0$ and $p_n$ as well as the affine curves $E_0$ and $E_{n+1}$ also play important    
roles in \S \ref{ss:ncc1} and \S \ref{sec:nmr}. 
\par
Another remark is that the construction so far carries over, when we restrict the entire $(w_1, w_2)$-space 
$\bC^2$ to a subspace $\bC \times U$, where $U := \{ w_2 \in \bC \mid |w_2| < r \}$ with $r > 0$, redefine 
the cyclic quotient sigularity as $S := (\bC \times U)/G$ and modify the minimal resolution $X \to S$ accordingly.  
In \S \ref{ss:ncc1} we need a description of this construction in terms of a pushforward lift on a cotangent bundle.     
For a local biholomorphism of complex manifolds $f : M \to N$, its {\sl pushforward lift} $f_* : T^* M \to T^* N$ 
is defined by $(f_*)_p := (g^*)_{f(p)} : T_p^* M \to T_{f(p)}^* N$ at $p \in M$, where $g : (N, f(p)) \to (M, p)$ is 
the local inverse to the local biholomorphism $f : (M, p) \to (N, f(p))$ around $p \in M$.   
The pushforward lift $f_*$ is again a local biholomorphism; if $f$ is a global biholomorphism then so is $f_*$. 
We will always think of a space as a subspace of its cotangent bundle, upon identifying it with 
the zero section of the bundle.  
In particular a point on the space is a point in its cotangent bundle.    
\par
Back in the current situation, under the identification $\bC \times U = T^* U$, $(w_1, w_2) \leftrightarrow w_1 \, d w_2$, 
the action of $g = \diag(\zeta^{-1}, \zeta) \in G$ on $\bC \times U$ is just the pushforward lift $\gamma_* : T^* U \to T^* U$ 
of the rotation $\gamma : U \to U$, $w_2 \mapsto \zeta w_2$. 
The quotient map $U \to V := U/\langle \gamma \rangle$ is realized by    
$V := \{ z \in \bC \mid |z| < r^{n+1} \}$ and $w_2 \mapsto z := w_2^{n+1}$.  
It restricts to an unramiled $\langle \gamma \rangle$-covering $U \setminus \{ 0 \} \to V \setminus \{ 0 \}$, 
so the pushforward lift $T^*(U \setminus \{ 0 \}) \to T^*(V \setminus \{ 0 \})$, 
$w_1 \, d w_2 \mapsto \xi \, d z$, induces an isomorphism 
$(T^*(U \setminus \{0 \}))/G \to T^*(V \setminus \{ 0 \})$.  
This gives an identification 
\begin{equation} \label{eqn:cotangent} 
T^*(V \setminus \{ 0 \}) = (T^*(U \setminus \{0 \}))/G = S \setminus \{ z = 0 \},   
\end{equation}
which avoids the singularity $0 \in S$. 
Thus $T^*(V \setminus \{ 0 \})$ can be thought of as an open subset of the space $X$.   
\begin{lemma} \label{lem:zs} 
The zero section of $T^*(V \setminus \{ 0 \})$ coincides with the punctured curve $E_0 \setminus \{ p_0 \}$ in $X$.  
\end{lemma} 
{\it Proof}. 
Formula \eqref{eqn:resol2-1a} with $i = 1$ reads $x = u_1 v_0$ and $z = v_0$,  
so the definition of pushforward lift together with formula \eqref{eqn:xyz} tells us that   
the fiber coordinate $\xi$ of the cotangent bundle $T^*(V \setminus \{ 0 \})$ is calculated as  
 $$
 \xi = \dfrac{w_1}{(d z/dw_2)} = \frac{w_1}{(n+1) w_2^n} = \frac{w_1 w_2}{(n+1) w_2^{n+1}} 
 = \frac{x}{(n+1) z} = \frac{u_1 v_0}{(n+1) v_0} = \frac{u_1}{n+1}. 
 $$
Note also that $z = v_0 \neq 0$ on $T^*(V \setminus \{ 0 \}) = S \setminus \{ z = 0 \}$. 
Thus the zero section $\{ \xi = 0 \}$ of $T^*(V \setminus \{0 \})$ coincides with 
the curve $E_0 = \{ u_1 = 0 \}$ with the point $p_0 \leftrightarrow (u_1, v_0) = (0, 0)$ removed.   
\hfill $\Box$
\subsection{Non-Cyclic Cases} \label{ss:ncc1} 
Any finite non-cyclic subgroup of $\SU(2)$ is a binary polyhedral group by Theorem \ref{thm:bclassf}.   
To deal with this case, we begin with a general construction for the standard $\GL_2(\bC)$-action on $\bC^2$. 
First we consider the action of the subgroup $\{ \pm 1\} \subset \GL_2(\bC)$.      
The quotient map $\pi : \bC^2 \to Y := \bC^2/\{ \pm 1 \}$ is realized by   
\begin{equation*} \label{eqn:S2}
\pi : \bC^2 \to Y :=  \{ (y_1, y_2, y_3) \in \bC^3 \mid y_1 y_2 = y_3^2 \},  
\quad (w_1, w_2) \mapsto (y_1, y_2, y_3) := (w_1^2, w_2^2, w_1 w_2).   
\end{equation*}  
The group $\GL_2(\bC)$ acts on $Y\subset \bC^3$ through the group homomorphism   
$$
\rho : \GL_2(\bC) \to \GL_3(\bC), \quad 
\begin{pmatrix} a & b \\ c & d \end{pmatrix} 
\mapsto 
\begin{pmatrix} 
a^2 & b^2 & 2 a b \\ c^2 & d^2 & 2 c d \\ a c & b d & ad + b c 
\end{pmatrix},      
$$
where $\rho(g)$ preserves $Y$ for every $g \in \GL_2(\bC)$. 
This action descends to an action of $\GL_2(\bC)/ \{\pm 1\}$, since the kernel of $\rho$ is $\{ \pm 1\}$. 
Now the map $\pi : \bC^2 \to Y$ is $\GL_2(\bC)$-equivariant. 
\par
Consider the blowup $\varSigma_1 \to \bC^2$ at the origin in $\bC^2$.  
Explicitly, $\varSigma_1$ is realized by 
$$
\varSigma_1 := \{ (w, [z] ) = ((w_1, w_2), [z_1 : z_2]) \in \bC^2 \times \bP^1 \mid z_1 w_2 = z_2 w_1 \},  
$$
with the blowup morphism being the first projection $(w, [z]) \mapsto w$.   
The group $\GL_2(\bC)$ acts on $\varSigma_1$ diagonally by 
$$
g (w, [z]) := (g w, [g z]) \qquad \mbox{for} \quad g \in \GL_2(\bC), \,\, (w, [z]) \in \varSigma_1.  
$$ 
Moreover the quotient space $\varSigma_2 := \varSigma_1/\{ \pm 1\}$ and the quotient map 
$\phi : \varSigma_1 \to \varSigma_2$ are realized by 
\begin{align*}
\varSigma_2 &:= \{ ((y_1, y_2), [z_1 : z_2]) \in \bC^2 \times \bP^1 \mid z_1^2 y_2 = z_2^2 y_1 \}, \\[1mm]
\phi &: \varSigma_1 \to \varSigma_2, \quad  
((w_1, w_2), [z_1 : z_2]) \mapsto ((y_1, y_2), [z_1 : z_2]) \qquad \mbox{with} \quad (y_1, y_2) := (w_1^2, w_2^2).    
\end{align*} 
Then the map $\phi : \varSigma_1 \to \varSigma_2$ is $\GL_2(\bC)$-equivariant,  
where $\GL_2(\bC)$ acts on $\varSigma_2$ by  
\begin{equation} \label{eqn:actSig2}
g \left( 
\begin{pmatrix} y_1 \\ y_2 \end{pmatrix}, 
\begin{bmatrix} z_1 \\ z_2 \end{bmatrix} \right) 
:= 
\left( 
\begin{pmatrix} a^2 y_1 + b^2 y_2 + 2 a b \, y_3 \\ c^2 y_1 + d^2 y_2 + 2 c d \, y_3  \end{pmatrix}, 
\begin{bmatrix} a z_1 + b z_2 \\ c z_1 + d z_2 \end{bmatrix} \right) 
\quad \mbox{for} \quad 
g = \begin{pmatrix} a & b \\ c & d \end{pmatrix} \in \GL_2(\bC). 
\end{equation}
This action naturally descends to an action of $\GL_2(\bC)/\{ \pm 1\}$. 
Thanks to the relation $z_1^2 y_2 = z_2^2 y_1$, a well-defined map 
$\psi : \varSigma_2 \to Y$, $((y_1, y_2), [z_1 : z_2]) \mapsto (y_1, y_2, y_3)$ 
is defined by $y_3 := (z_2/z_1) \, y_1$ if $z_1 \neq 0$ and $y_3 := (z_1/z_2) \, y_2$ if 
$z_2 \neq 0$, so that there exists a $\GL_2(\bC)$-equivariant commutative diagram 
\begin{equation} \label{cd:SigS}
\begin{CD}
\varSigma_1 @> {\scriptstyle \mathrm{blowup}} >> \bC^2 \\
@V {\phi} VV @VV {\pi} V \\
\varSigma_2 @> {\psi} >> Y.  
\end{CD}
\end{equation}
\par
The surface $\varSigma_2$ together with the second projection $\varSigma_2 
\to \bP^1$, $(y, [z]) \mapsto [z]$ is a line bundle over $\bP^1$ isomorphic to the 
holomorphic cotangent bundle $T^* \bP^1$.  
Let $g \in \GL_2(\bC)$ act on $\bP^1$ as a projective linear transformation 
in the usual manner.  
Then under the identification $T^* \bP^1 = \varSigma_2$ the pushforward 
lift $g_* : T^* \bP^1 \to T^* \bP^1$ is given by  
\begin{equation} \label{eqn:actcot}
g_* \left( 
\begin{pmatrix} y_1 \\ y_2 \end{pmatrix}, 
\begin{bmatrix} z_1 \\ z_2 \end{bmatrix} \right) 
=  
\left( \frac{1}{\det g}
\begin{pmatrix} a^2 y_1+ b^2 y_2 + 2 a b \, y_3 \\ c^2 y_1 + d^2 y_2 + 2 c d \, y_3  \end{pmatrix}, 
\begin{bmatrix} a z_1 + b z_2 \\ c z_1 + d z_2 \end{bmatrix} \right) 
\quad \mbox{for} \quad 
g = \begin{pmatrix} a & b \\ c & d \end{pmatrix} \in \GL_2(\bC).  
\end{equation}
In particular the actions \eqref{eqn:actSig2} and \eqref{eqn:actcot} coincide 
when they are restricted to the subgroup $\SL_2(\bC)$. 
\par
Consider a diagonalizable matrix $g \in \GL_2(\bC)$ with eigenvalues $\beta \alpha^{\pm 1}$. 
Let $l \subset \bC^2$ be an eigenline of $g$ corresponding to the eigenvalue $\beta \alpha$. 
If $\alpha \neq \pm 1$ then $l$ is unique; otherwise, $g$ is a scalar matrix 
and $l$ is an arbitrary line through the origin. 
The line $l$ designates a unique point $q \in \bP^1$, which is fixed by $[g] \in \PGL_2(\bC)$.  
\begin{lemma} \label{lem:actfix}
In the setting mentioned above the matrix $g \in \GL_2(\bC)$ acts on $\varSigma_2 = T^* \bP^1$ by 
\begin{equation} \label{eqn:actfix} 
g : (w_1, \, w_2) \mapsto (\beta (\beta \alpha^2) w_1, \, \beta (\beta \alpha^2)^{-1} w_2) 
\quad \mbox{around the fixed point $q \in \bP^1 \subset T^* \bP^1$},  
\end{equation}
in terms of an inhomogeneous coordinate $w_2$ on $\bP^1$ around $q$ and the associated 
fiber coordinate $w_1$ on $T^* \bP^1$.   
\end{lemma}
{\it Proof}. 
Taking a conjugate to $g$ in $\GL_2(\bC)$ if necessary, we may assume from the beginning that 
$g$ is a diagonal matrix $\diag( \beta \alpha, \, \beta \alpha^{-1})$ and $q$ is the point $[1 : 0] \in \bP^1$.  
Then formula \eqref{eqn:actSig2} reads 
$$
g((y_1, y_2), \, [z_1 : z_2]) 
= (((\beta \alpha)^2 y_1, \, (\beta \alpha^{-1})^2 y_2), \, [\beta \alpha z_1 : \beta \alpha^{-1} z_2]). 
$$
We can now take $w_2 := z_2/z_1$ as an inhomogeneous coordinate of $\bP^1$ around $q = [1 : 0]$ 
and $w_1 := y_1$ as the associated fiber coordinate on the cotangent bundle $T^* \bP^1$ near the point $q$. 
In terms of the local chart $(w_1, w_2)$ we have 
$g (w_1, \, w_2) = (\beta^2 \alpha^2 w_1, \, \alpha^{-2} w_2)$ and hence formula \eqref{eqn:actfix}. 
\hfill $\Box$ \par\medskip 
\begin{table}[h]
\centerline{
\begin{tabular}{ccccccccccccc}
\hline
       & &          &          &          & &            &            &             &  &        & &            \\[-3mm]  
type & & $n^{\re}$ & $n^{\rf}$ & $n^{\rv}$ & & $|P^{\re}|$ & $|P^{\rf}|$ & $|P^{\rv}|$  & & $|\vG|$ & & $n$ \\[1mm]
\hline
       & &          &          &          & &            &            &             &  &        & &            \\[-3mm]  
$\rD_4$ & & $1$ & $1$ & $1$ & & $2$ & $2$ & $2$  & &  $4$ &                 &  $4$    \\[1mm]
$\rD_n$ & & $1$ & $n-3$ & $1$ & & $n-2$ & $2$ & $n-2$  & &  $2(n-2)$ & & $\ge 5$ \\[1mm]
$\rE_6$ & & $1$ & $2$ & $2$ & & $6$ & $4$ & $4$ & & $12$ & & $6$ \\[1mm]
$\rE_7$ & & $1$ & $2$ & $3$ & & $12$ & $8$ & $6$ & & $24$ & & $7$ \\[1mm]
$\rE_8$ & & $1$ & $2$ & $4$ & & $30$ & $20$ & $12$ && $60$ & & $8$ \\[1mm]
\hline
\end{tabular}}
\caption{Numerical data concerning exceptional orbits.} 
\label{tab:excep}
\end{table}
Let $G$ be a binary polyhedral group in $\SU(2)$.  
Then $G$ has center $\{ \pm 1\}$ and the quotient group $\vG := G/\{ \pm 1\}$ 
is a regular polyhedral group in $\PSU(2) \cong \SO(3)$. 
We have an identification $S := \bC^2/G = Y/\vG$ via the map $\pi$ in diagram \eqref{cd:SigS}. 
Together with the identification $\varSigma_2/\vG = (T^* \bP^1)/ \vG$ the map 
$\psi$ in \eqref{cd:SigS} induces a morphism    
$$
\psi : 
(T^* \bP^1)/\vG = \varSigma_2/\vG \to Y/ \vG = S.  
$$  
This map is composed with a minimal resolution of $(T^* \bP^1)/\vG$ to get a 
minimal resolution of $S$, where $\vG$ acts on $T^* \bP^1$ as the pushforward 
lift of the usual action of $\vG$ on $\bP^1$ as projective linear transformations.     
A point $q \in \bP^1$ is said to be {\sl exceptional} if the isotropy subgroup 
$\vG_q := \{ \gamma \in \vG \mid \gamma(q) = q \}$ is nontrivial. 
All exceptional points in $\bP^1$ are partitioned into three $\vG$-orbits $P^{\re}$, $P^{\rf}$, $P^{\rv}$, which 
arise respectively from the mid-edge points, face-centers, vertices of a regular polyhedron having $\vG$ as its rotational symmetry group.        
For $\nu = \re, \rf, \rv$, if $q \in P^{\nu}$ then $\vG_{q}$ is a cyclic group of order $n^{\nu} + 1$ and $P^{\nu}$ 
is a set of cardinality $|\vG|/(n^{\nu} +1)$, where the number $n^{\nu}$ and related quantities are given 
in Table \ref{tab:excep}.  
The quotient space $\mathbf{P}^1 := \bP^1/\vG$ is again a projective line and the quotient map 
$\bP^1 \to \mathbf{P}^1$ restricts to an unramified $\vG$-covering 
$\bP^1 \setminus (P^{\re} \cup P^{\rf} \cup P^{\rv}) \to \mathbf{P}^1 \setminus \{ p^{\re}, p^{\rf}, p^{\rv} \}$, 
where the point $p^{\nu} \in \mathbf{P}^1$ is the image of $P^{\nu}$. 
So its pushforward lift is again an unramified $\vG$-covering,  
\begin{equation} \label{eqn:TP}
T^*(\bP^1 \setminus (P^{\re} \cup P^{\rf} \cup P^{\rv})) \to T^*(\mathbf{P}^1 \setminus \{ p^{\re}, p^{\rf}, p^{\rv} \}),  
\end{equation}
which induces an isomorphism $T^*(\bP^1 \setminus (P^{\re} \cup P^{\rf} \cup P^{\rv}))/\vG \to 
T^*(\mathbf{P}^1 \setminus \{p^{\re}, p^{\rf}, p^{\rv} \})$.       
For each $\nu = \re, \rf, \rv$, take a representative point $q^{\nu} \in P^{\nu}$ and write $\vG^{\nu} := \vG_{q^{\nu}}$. 
Let $U^{\nu}$ be a sufficiently small $\vG^{\nu}$-invariant disk in $\bP^1$ centered at $q^{\nu}$. 
The pushforward lift $T^* U^{\nu} \to T^* \bP^1$ of the inclusion $U^{\nu} \hookrightarrow \bP^1$ induces 
an isomorphism of $(T^* U^{\nu})/\vG^{\nu}$ onto an open subset of $(T^* \bP^1)/ \vG$. 
From these observations we have open embeddings 
\begin{equation} \label{eqn:oemb}
T^*(\mathbf{P}^1 \setminus \{ p^{\re}, p^{\rf}, p^{\rv} \}) \hookrightarrow (T^*\bP^1)/\vG 
\hookleftarrow (T^* U^{\nu})/ \vG^{\nu}, \qquad \nu = \re, \rf, \rv. 
\end{equation}
\par
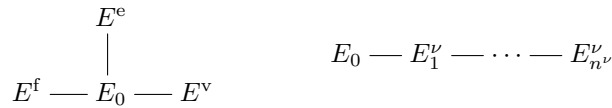
\begin{figure}[h]
\centerline{
\begin{minipage}{0.3\linewidth}
\[
\xymatrix@C=14pt@R=14pt{
            & E^{\re} \ar@<-0.3ex>@{-}[d]  &  \\
 E^{\rf} \ar@{-}[r]  & E_0 \ar@{-}[r] & E^{\rv}}
\]
\end{minipage} \hspace{-8mm}
\begin{minipage}{0.33\linewidth}
\[
\xymatrix@C=12pt@R=12pt{ 
E_0 \ar@{-}[r] & E_1^{\nu} \ar@{-}[r] & \cdots \ar@{-}[r] & E_{n^{\nu}}^{\nu}}
\]
\end{minipage}}
\caption{Dynkin diagram of type $\rD$ or $\rE$ and its arms.} 
\label{fig:dDE}
\end{figure}

Since $S^{\nu} := (T^* U^{\nu})/\vG^{\nu}$ is a cyclic quotient singularity of order 
$n^{\nu} + 1$, the recipe in \S \ref{ss:cc1} gives a minimal resolution $X^{\nu} \to S^{\nu}$.  
Various objects associated with this resolution are denoted by the same symbols as those 
in \S \ref{ss:cc1} followed by a superscript $\nu = \re, \rf, \rv$.  
For example the exceptional set $E^{\nu} \subset X^{\nu}$ is the union of $(-2)$-curves 
$E^{\nu}_i$, $1 \le i \le n^{\nu}$.      
A minimal resolution $X$ of $(T^* \bP^1)/\vG$ and hence of $S := \bC^2/G$ is obtained by gluing three pieces 
$X^{\nu}$, $\nu = \re, \rf, \rv$, to $T^*(\mathbf{P}^1 \setminus \{ p^{\re}, p^{\rf}, p^{\rv} \})$ 
via the mappings \eqref{eqn:oemb} along the open sets 
$$
S^{\nu} \setminus \{ z^{\nu} = 0 \} = (T^*(U^{\nu} \setminus \{ q^{\nu} \}))/\vG^{\nu} = 
T^*(V^{\nu} \setminus \{ p^{\nu}\}), 
$$ 
where $V^{\nu}$ is the image of $U^{\nu}$ under the projection $\bP^1 \to \mathbf{P}^1$ 
and the identification \eqref{eqn:cotangent} is employed.  
Let $E_0$ be the closure in $X$ of the zero section of $T^*(\mathbf{P}^1 \setminus \{ p^{\re}, p^{\rf}, p^{\rv} \})$.  
Here the present $E_0$ is different from $E_0$ in \S \ref{ss:cc1}, the latter being $E_0^{\nu}$ 
in the current situation.       
Then $E_0$ is also a $(-2)$-curve in $X$ and the union $E := E_0 \cup E^{\re} \cup E^{\rf} \cup E^{\rv}$ 
constitutes the exceptional set of the minimal resolution $X \to S$.   
Lemma \ref{lem:zs} implies 
\begin{equation} \label{eqn:E} 
E_0^{\nu} \subset E_0, \qquad E^{\nu}_1 \cap E_0 = \{ p^{\nu}_0 \}, \qquad \nu = \re, \rf, \rv,  
\end{equation}
and the dual graph of $E$ is given by the left diagram in Figure \ref{fig:dDE}, 
while the right diagram shows its arm $E^{\nu}$. 
\section{Normalizers on Minimal Resolutions} \label{sec:nmr} 
Let $G$ be a nontrivial finite subgroup of $\SU(2)$ and $X \to S$ be a minimal 
resolution of the quotient singularity $S := \bC^2/G$ as constructed in \S \ref{sec:rs}.    
The exceptional set $E \subset X$ is a union of $(-2)$-curves, the dual graph of which 
is a Dynkin diagram $D$ of type $\rA$, $\rD$ or $\rE$, where the correspondence $G \mapsto D$ is 
given in Table \ref{tab:bclassf}. 
Let $N_{\GL}$ be the normalizer of $G$ in $\GL_2(\bC)$. 
For each $J \in N_{\GL}$ the linear automorphism $J : \bC^2 \to \bC^2$ descends to a map $S \to S$, 
which in turn lifts to a holomorphic automorphism $f_J : X \to X$. 
It permutes the $(-2)$-curves in $E$ and induces a diagram automorphism $\sigma_J : D \to D$. 
According to their Dynkin types, the automorphism groups $\Aut(D)$ are classified in Table \ref{tab:auto}.        
The aim of this section is to describe how the map $f_J \in \Aut(X)$ acts on the exceptional set $E$ 
in connection with the action of $\sigma_J \in \Aut(D)$ on the diagram $D$.    
\begin{table}[h] 
\centerline{
\begin{tabular}{cccccccl}
\hline 
$D$ & $\rA_1$ & $\rA_n$ & $\rD_4$ & $\rD_n$ & $\rE_6$ & $\rE_7$ & $\rE_8$ \\
\hline
$\Aut(D)$ & $1$ & $\bZ_2$ & $S_3$ & $\bZ_2$ & $\bZ_2$ & $1$ & $1$ \\
\hline
$n$ & $1$ & $\ge 2$ & $4$ & $\ge 5$ & $6$ & $7$ & $8$ \\ 
\hline
\end{tabular}} 
\caption{Dynkin diagrams $D$ and their automorphism groups $\Aut(D)$.} 
\label{tab:auto}
\end{table}
\begin{remark} \label{rem:fJ} 
For $K \in N_{\SL}$ let $\langle K \rangle$ denote the class of $K$ in the quotient group $N_{\SL}/G$. 
If $J \in N_{\GL}$ is expressed as $J = \beta K$ with $\beta \in \bC^{\times}$ and 
$K \in N_{\SL}$, then $f_J$ depends only on $\beta$ and $\langle K \rangle$, since 
so does the map $S \to S$. 
Theorem \ref{thm:normalizer} and Table \ref{tab:auto} show that when $G$ is of type $\rD$ or $E$ 
we have an isomorphism $N_{\SL}/G \cong \Aut(D)$. 
\end{remark}   
\subsection{Cyclic Cases} \label{ss:cc2}
Let $G$ be the cyclic group of order $n+1$ as in \S \ref{ss:cc1}, where $n \in \bN$.  
The notation there is retained in this subsection.     
\begin{lemma} \label{lem:ida1} 
If $J \in N_{\GL}$ is a diagonal matrix of the form 
\begin{equation} \label{eqn:fa}
J = \beta K \qquad \mbox{with} \quad 
K = \diag( \alpha, \, \alpha^{-1}) \in N_{\SL},   
\end{equation}
where $\alpha$, $\beta \in \bC^{\times}$, then the induced automorphism 
$f_J : X \to X$ maps each $W_i$ onto itself, sending  
\begin{equation} \label{eqn:f3}
(u_i, v_{i-1}) \mapsto (u_i', v_{i-1}') := (\beta^{2 i-n-1} \alpha^{n+1} u_i, \, \beta^{n+3-2 i} \alpha^{-n-1} v_{i-1}),   
\qquad 1 \le i \le n+1. 
\end{equation}
It preserves each $E_i$, $0 \le i \le n+1$, and fixes each $p_i$, $0 \le i \le n$.  
In particular $\sigma_J \in \Aut(\rA_n)$ is trivial.  
The multipliers of $f_J$ at $p_i$ are $\beta^{n+1-2 i} \alpha^{-n-1}$ along $E_i$ and 
$\beta^{2 i-n+1} \alpha^{n+1}$ along $E_{i+1}$, so $\det (d f_J)_{p_i} = \beta^2 = \det J$.  
\end{lemma}
{\it Proof}. 
In view of formula \eqref{eqn:xyz} the matrix $J$ induces an automorphism 
\begin{equation} \label{eqn:f2a}
S \to S, \quad  (x, y, z) \mapsto (x', y', z') := (\beta^2 x, \, (\beta \alpha)^{n+1} y, \, (\beta \alpha^{-1})^{n+1} z). 
\end{equation}
The map $f_J$ is then calculated as the composite of the following three maps  
$$
\begin{CD} 
(u_i, v_{i-1}) @> \eqref{eqn:resol2-1a} >> (x, z) @> \eqref{eqn:f2a} >> (x', z') 
@> (\ref{eqn:resol2-3}') >> (u_i', v_{i-1}'),   
\end{CD}
$$
where $(\ref{eqn:resol2-3}')$ is the map \eqref{eqn:resol2-3} with $(u_i, v_{i-1})$ and $(x, y)$ 
replaced by $(u_i', v_{i-1}')$ and $(x', y')$.  
Explicitly we have   
\begin{align*}
u_i' &= (x')^i (z')^{-1} = (\beta^2 x)^i \{ (\beta \alpha^{-1})^{n+1}  z \}^{-1} 
= \beta^{2 i-n-1} \alpha^{n+1} x^i z^{-1} = \beta^{2 i-n-1} \alpha^{n+1} u_i,  
\\[1mm] 
v_{i-1}' &= (x')^{1-i} z' = (\beta^2 x)^{1-i} (\beta \alpha^{-1})^{n+1} z 
= \beta^{n+3-2 i} \alpha^{-n-1} x^{1-i} z = \beta^{n+3-2 i} \alpha^{-n-1} v_{i-1}, 
\end{align*}
which proves formula \eqref{eqn:f3}. 
Since $E_i \cap W_{i+1} = \{ u_{i+1} = 0\}$, $E_{i+1} \cap W_{i+1} = \{ v_i = 0 \}$ and 
$p_i$ has coordinates $(u_{i+1}, v_i) = (0, 0)$ for $0 \le i \le n$, the remaining assertions 
follow from formula \eqref{eqn:f3}. \hfill $\Box$ 
\begin{lemma} \label{lem:ida2} 
Let $n \ge 2$. 
If $J \in N_{\GL}$ is an anti-diagonal matrix of the form 
\begin{equation} \label{eqn:fb} 
J = \beta K \qquad \mbox{with} \quad 
K = \begin{pmatrix} 0 & \ri \gamma \\ \ri \gamma^{-1} & 0 \end{pmatrix} \in N_{\SL},   
\end{equation}
where $\beta$, $\gamma \in \bC^{\times}$, then the induced automorphism 
$f_J : X \to X$ maps each $W_i$ onto $W_{n+2-i}$, sending  
\begin{equation} \label{eqn:f4}
(u_i, v_{i-1}) \mapsto (u_{n+2-i}', v_{n+1-i}') := ((\ri \beta)^{n+3-2 i} \gamma^{n+1} v_{i-1}, \, (\ri \beta)^{2 i-n-1} \gamma^{-n-1} u_i),    
\quad 1 \le i \le n+1.  
\end{equation}
It swaps $E_i \leftrightarrow E_{n+1-i}$, $1 \le i \le n$, and 
$p_i \leftrightarrow p_{n-i}$, $1 \le i \le n-1$, so $\sigma_J \in \Aut(\rA_n) \cong \bZ_2$ is an involution.      
\end{lemma}
{\it Proof}. 
In view of formula \eqref{eqn:xyz} the matrix $J$ induces an automorphism 
\begin{equation} \label{eqn:f2b}
S \to S, \quad (x, y, z) \mapsto (x', y', z') := ((\ri \beta)^2 x, \, (\ri \beta \gamma)^{n+1} z, \, (\ri \beta \gamma^{-1})^{n+1} y).  
\end{equation}
The map $f_J$ is then calculated as the composite of the following three maps  
$$
\begin{CD} 
(u_i, v_{i-1}) @> \eqref{eqn:resol2-1a} >> (x, z) @> \eqref{eqn:f2b} >> (x', y') 
@> (\ref{eqn:resol2-3}'') >> (u_{n+2-i}', v_{n+1-i}'),   
\end{CD}
$$
where $(\ref{eqn:resol2-3}'')$ is the map \eqref{eqn:resol2-3} with $(u_i, v_{i-1})$ and $(x, y)$ 
replaced by $(u_i', v_{i-1}')$ and $(x', y')$ and then followed by the change of index $i \mapsto n+2 - i$.   
Explicitly we have   
\begin{align*}
u_{n+2-i}' &= (x')^{(n+2-i)-n-1} y' = \{ (\ri \beta)^2 x \}^{1-i} (\ri \beta \gamma)^{n+1} z 
= (\ri \beta)^{n+3-2 i} \gamma^{n+1} x^{1-i} z =  (\ri \beta)^{n+3-2 i} \gamma^{n+1} v_{i-1}, \\[1mm] 
v_{n+1-i}' 
&= (x')^{n+2-(n+2-i)} (y')^{-1} = \{ (\ri \beta)^2 x \}^i \{ (\ri \beta \gamma)^{n+1} z\}^{-1} 
= (\ri \beta)^{2 i-n-1} \gamma^{-n-1} x^i z^{-1} = (\ri \beta)^{2 i-n-1} \gamma^{-n-1} u_i,    
\end{align*}
which proves formulas \eqref{eqn:f4}. 
The remaining assertions follow from formula \eqref{eqn:f4}. \hfill $\Box$ 
\begin{figure}[h]
\begin{minipage}{.45\linewidth}
\centerline{\includegraphics*[width=51mm,clip]{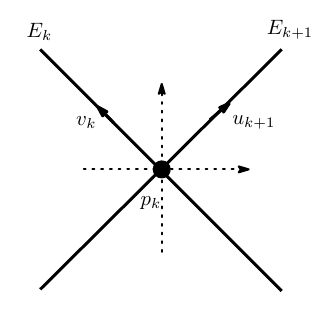}} 
\caption{The case of $n = 2 k$ with $k \ge 1$.} 
\label{fig:excep1} 
\end{minipage}
\hspace{3mm}
\begin{minipage}{.45\linewidth}
\centerline{\includegraphics*[width=75mm,clip]{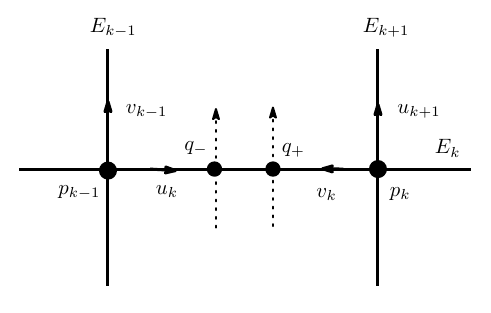}}  
\caption{The case of $n = 2 k-1$ with $k \ge 2$.}
\label{fig:excep2} 
\end{minipage}
\end{figure}
\begin{corollary} \label{cor:ida2}
In the situation of Lemma $\ref{lem:ida2}$ the following assertions hold. 
\begin{enumerate}
\item If $n = 2k$ with $k \ge 1$ then $f_J$ swaps $E_k \leftrightarrow E_{k+1}$ 
while fixing the point $p_k$. 
This is the unique fixed point of $f_J$ on $E$. 
The tangent map $(d f_J)_{p_k} : T_{p_k} X \to T_{p_k} X$ has eigenvalues $\pm \ri \beta$, 
whose eigen-directions are shown by the dotted arrows in Figure $\ref{fig:excep1}$.  
In particular we have $\det (d f_J)_{p_k} = \beta^2 = \det J$. 
\item If $n = 2k-1$ with $k \ge 2$ then $f_J$ swaps $(E_{k-1}, p_{k-1}) \leftrightarrow 
(E_{k+1}, p_k)$, preserves $E_k$, and has two fixed points $q_{\pm}$ on $E_k$. 
They are the only fiexd points of $f_J$ on $E$. 
At each of them the tangent map $(d f)_{q_{\pm}} : T_{q_{\pm}} X \to T_{q_{\pm}} X$ has 
eigenvalue $-1$ along $E_k$ and eigenvalue $- \beta^2$ in a normal direction to $E_k$ 
as shown by the dotted arrows in Figure $\ref{fig:excep2}$.  
In particular we have $\det (d f_J)_{q_{\pm}} = \beta^2 = \det J$.    
\end{enumerate}
\end{corollary}
{\it Proof}. 
Assertion (1). 
Formula \eqref{eqn:f4} with $i = k+1$ shows that $f_J$ maps $W_{k+1}$ onto itself,  
sending $(u_{k+1}, v_k) \mapsto (u_{k+1}', v_k')$ via   
$u_{k+1}' = (\ri \beta) \gamma^{2 k+1} v_k$ and $v_k' = (\ri \beta) \gamma^{-2 k-1} u_{k+1}$.  
This expression readily leads to the assertion, since 
$E_k = \{ u_{k+1} = 0\}$, $E_{k+1} = \{ v_k = 0\}$ and $p_k$ corresponds to $(u_{k+1}, v_k) = (0, 0)$. 
\par
Assertion (2). 
Formula \eqref{eqn:f4} with $i = k$ shows that $f_J$ maps $W_k$ onto $W_{k+1}$, sending 
$(u_k, v_{k-1}) \mapsto (u_{k+1}', v_k')$ via $u_{k+1}' = -\beta^2 \gamma^{2 k} v_{k-1}$ and 
$v_k' = \gamma^{-2 k} u_k$. 
Thus $f_J$ swaps $(E_{k-1}, p_{k-1}) \leftrightarrow (E_{k+1}, p_k)$,  
since $E_{k-1} = \{ u_k = 0\}$ and $E_{k+1} = \{ v_k = 0\}$ while $p_{k-1}$ and $p_k$ 
correspond to $(u_k, v_{k-1}) = (0, 0)$ and 
$(u_{k+1}, v_k) = (0, 0)$, respectively. 
In terms of the local coordinates $(u_k, v_{k-1})$ the map $f_J$ is represented as 
$$
u_k' = \gamma^{2 k} u_k^{-1}, \qquad v_{k-1}' = - \beta^2 \gamma^{- 2k} \cdot u_k^2 v_{k-1}. 
$$  
Hence $f_J$ preserves $E_k = \{ v_{k-1} = 0\}$ and has two fixed points $q_{\pm} \in E_k$ 
corresponding to $(u_k, v_{k-1}) = (\pm \gamma^k, 0)$.  
In the current coordinates the tangent map $(d f_J)_{q_{\pm}} : T_{q_{\pm}} X \to T_{q_{\pm}} X$ is 
represented by the Jacobian matrix $\diag(-1, \, -\beta^2)$.  
From this observation the assertion about its eigenvalues follows immediately. 
\hfill $\Box$ \par\medskip
When $n = 1$, we have $N_{\GL} = \GL_2(\bC)$ by Theorem \ref{thm:normalizer}.(1), 
so $J \in N_{\GL}$ is not always diagonalizable. 
Whether $J$ is diagonalizable or not can be judged in terms of a property of the  
M\"{o}bius transformation $f_J|_E$ on $E = E_1 \cong \bP^1$. 
A M\"{o}bius transformation is said to be {\sl parabolic} if it has only one fixed point.    
\begin{lemma} \label{lem:ndiag} 
When $n = 1$, a matrix $J \in N_{\GL} = \GL_2(\bC)$ is not diagonalizable if and only if $f_J|_E$ is parabolic.    
\end{lemma}
{\it Proof}. 
Let $\beta \alpha^{\pm 1} \in \bC^{\times}$ be the eigenvalues of $J$.  
Up to conjugacy in $\GL_2(\bC)$ we may assume 
$$
J = \beta \begin{pmatrix} \alpha & \alpha b_1 \\[1mm] 0 & \alpha^{-1} \end{pmatrix}, \qquad \mbox{so that} \qquad 
\begin{pmatrix} w_1' \\[1mm] w_2' \end{pmatrix} = J \begin{pmatrix} w_1 \\[1mm] w_2 \end{pmatrix} 
= \beta \begin{pmatrix} \alpha(w_1 + b_1 w_2) \\[1mm] \alpha^{-1} w_2 \end{pmatrix},    
$$
where $b_1 \in \bC$. 
Equations \eqref{eqn:xyz} with $n = 1$ become $x = w_1 w_2$, $y = w_1^2$, $z = w_2^2$. 
Thus we have 
$$
x' = w_1' w_2' = \beta^2 (w_1 w_2 + b_1 w_2^2) = \beta^2 (x + b_1 z), \qquad 
z' = (w_2')^2 = \beta^2 \alpha^{-2} \, w_2^2 = \beta^2 \alpha^{-2} z. 
$$
Formula \eqref{eqn:resol2-3} with $i = 1$ reads $u_1 = x z^{-1}$ and $v_0 = z$, 
so the map $f_J$ is calculated as
\begin{equation} \label{eqn:fJn1} 
u_1' = x' (z')^{-1} = \alpha^2(x z^{-1} + b_1) = \alpha^2(u_1 + b_1), \qquad 
v_0' = z' = \beta^2 \alpha^{-2} \, z = \beta^2 \alpha^{-2} \, v_0, 
\end{equation}
in terms of the local coordinates $(u_1, v_0)$, where $E$ is the closure of $\{v_0 = 0\}$.  
Hence $J$ is not diagonalizable, if and only if $\alpha^2 = 1$ and $b_1 \neq 0$, if and only if 
$f_J|_E$ is a parabolic transformation $u_1 \mapsto u_1 + b_1$. \hfill $\Box$ 
\subsection{Non-Cyclic Cases} \label{ss:ncc2} 
We proceed to the case where $G$ is a binary polyhedral group as in \S\ref{ss:ncc1}. 
The notation there is retained in this subsection. 
Let $J \in N_{\GL}$ and write $J = \beta K$ with $\beta \in \bC^{\times}$ and $K \in N_{\SL}$ as in Remark \ref{rem:fJ}. 
The projective linear transformation $[J] = [K]$ on $\bP^1$ induced by $J$ or $K$ permutes the exceptional orbits 
$P^{\re}$, $P^{\rf}$, $P^{\rv}$, while every element of $\vG = G/\{ \pm 1\}$ leaves them invariant.   
Thus the permutation of them depends only on the class $\langle K \rangle \in N_{\SL}/G$. 
Upon passing to the resolution, the map $f_J : X \to X$ restricts to a projective linear transformation on 
$E_0$ permuting the three points $p_0^{\re}$, $p_0^{\rf}$, $p_0^{\rv} \in E_0$; see also \eqref{eqn:E}. 
The class $\langle K \rangle$, the permutation of $P^{\re}$, $P^{\rf}$, $P^{\rv}$ and that of 
$p_0^{\re}$, $p_0^{\rf}$, $p_0^{\rv}$ are all identified with the Dynkin automorphism $\sigma_J \in \Aut(D)$. 
The map $f_J$ depends only on $\beta \in \bC^{\times}$ and $\sigma_J = \langle K \rangle \in \Aut(D) = N_{\SL}/G$, 
so we can take a various representative $K \in N_{\SL}$ of $\sigma_J$ without changing the map $f_J$.   
The order of $\sigma_J$ is $1$, $2$, or $3$, and each case is discussed separately.  
\begin{lemma} \label{lem:order1} 
Suppose that $\sigma_J$ is the unit element in $\Aut(D)$. 
Then the map $f_J$ preserves each of $E^{\re}$, $E^{\rf}$, $E^{\rv}$ and fixes $E_0$ pointwise.  
For each $\nu = \re, \rf, \rv$ and $0 \le i \le n^{\nu}$, the point $p_i^{\nu}$ is a fixed point of $f_J$ 
having multipliers $(\det J)^{i+1}$ along $E_{i+1}^{\nu}$ and $(\det J)^{-i}$ along $E_i^{\nu}$.  
\end{lemma}
{\it Proof}.  
Since $\sigma_J = \langle K \rangle$ is the unit element in $\Aut(D) = N_{\SL}/G$, we have  
$J = \beta K$ with $K \in G$ and hence we can take $K = I$ as a representative of the class $\sigma_J$, so that $J = \beta I$. 
For each $\nu = \re, \rf, \rv$, take a point $q^{\nu} \in P^{\nu}$. 
Since $q^{\nu}$ is fixed by $[\beta I] \in \PGL_2(\bC)$, Lemma \ref{lem:actfix} is available with 
$q = q^{\nu}$ and $\alpha = 1$ in formula \eqref{eqn:actfix}. 
This enables us to apply Lemma \ref{lem:ida1} with $\alpha = \beta$ to the resolution $X^{\nu} \to (T^* U^{\nu})/\vG^{\nu}$ 
of a cyclic quotient singularity of order $n^{\nu}+1$. 
The result is that $f_J$ preserves $E^{\nu}$ and its multipliers at the fixed points $p_i^{\nu} \in E^{\nu}$, 
$0 \le i \le n^{\nu}$, are $\beta^{2(i+1)} = (\det J)^{i+1}$ along $E_{i+1}^{\nu}$ and 
$\beta^{-2 i} = (\det J)^{-i}$ along $E_i^{\nu}$. 
The map $f_J$ fixes $p_0^{\re}$, $p_0^{\rf}$, $p_0^{\rv}$ and hence fixes the entire    
$E_0$ pointwise, because a projective linear transformation fixing three points is the identity. 
\hfill $\Box$
\begin{lemma} \label{lem:order2}
Suppose that $\sigma_J$ is of order $2$ in $\Aut(D)$. 
Then $D$ must be of type $\rD$ or $\rE_6$ and the map $f_J$ preserves just 
one of $E^{\re}$, $E^{\rf}$, $E^{\rv}$, say $E^{\nu}$, and swaps the remaining two. 
For each $0 \le i \le n^{\nu}$ the point $p_i^{\nu}$ is a fixed pont of $f_J$ having  
multipliers $-(\det J)^{i+1}$ along $E_{i+1}^{\nu}$ and $-(\det J)^{-i}$ along $E_i^{\nu}$. 
Moreover $f_J$ fixes a unique point $p_{\infty}$ on $E_0 \setminus \{ p_0^{\re}, p_0^{\rf}, p_0^{\rv} \}$, 
at which $f_J$ has multipliers $-1$ along $E_0$ and $- (\det J)$ along $T_{p_{\infty}}^*E_0$. 
\end{lemma} 
{\it Proof}.  
(1) Case of type $\rD_n$ with $n \ge 5$.  
In this case $G$ is a binary dihedral group of order $4(n-2)$ and $N_{\SL}$ is the unique 
binary dihedral group $H$ that contains $G$ as an index $2$ subgroup. 
Notice that $P^{\rf}$ is special to the effect that $n^{\rf} \neq n^{\re} = n^{\rv}$ (see Table \ref{tab:excep}).    
Let $Q^{\re}$, $Q^{\rf}$, $Q^{\rv}$ denote the exceptional orbits for the dihedral group $\vD := H/\{ \pm I \}$. 
Then $Q^{\rf} = P^{\rf}$ and $Q^{\rv} = P^{\re} \sqcup P^{\rv}$, while  
$Q^{\re}$ is a non-exceptional orbit of $\vG = G/\{ \pm I \}$. 
\par
(1a) Take a point $q^{\rf} \in Q^{\rf} = P^{\rf}$.  
Its isotropy subgroup in $\vD$ is a cyclic group of order $2(n^{\rf}+1)$, 
whose binary lift is a cyclic subgroup of order $4(n^{\rf}+1)$ in $H$.  
Let $K^{\rf}$ be a generator of the latter group.  
Then $[K^{\rf}] \in \PSL_2(\bC)$ swaps $P^{\re}$ and $P^{\rv}$, so $K^{\rf}$ lies in $H \setminus G$ 
and represents the nontrivial class $\sigma_J \in H/G \cong \bZ_2$. 
We may take $K^{\rf}$ as a representative of $\sigma_J$, so that $J = \beta K^{\rf}$. 
The eigenvalues of $K^{\rf}$ are a primitive $4(n^{\rf}+1)$-st root of unity, say $\zeta$, and its 
reciprocal $\zeta^{-1}$. 
Since $q^{\rf}$ is fixed by $[\beta K^{\rf}] \in \PGL_2(\bC)$, Lemma \ref{lem:actfix} is available with 
$q = q^{\rf}$ and $\alpha = \zeta$ in formula \eqref{eqn:actfix}.  
This enables us to apply Lemma \ref{lem:ida1} with $\alpha = \beta \zeta^2$ to the resolution 
$X^{\rf} \to (T^* U^{\rf})/\vG^{\rf}$ of a cyclic quotient singularity of order $n^{\rf}+1$. 
The result is that $f_J$ preserves $E^{\rf}$ and has multipliers $\beta^{2(i+1)} \zeta^{2(n^{\rf}+1)} 
= - \beta^{2(i+1)} = - (\det J)^{i+1}$ along $E_{i+1}^{\rf}$ and $\beta^{-2 i} \zeta^{-n^{\rf}-1} 
= - \beta^{-2 i} = (\det J)^{-i}$ along $E_i^{\rf}$ at the fixed points $p_i^{\rf} \in E^{\rf}$, $0 \le i \le n^{\rf}$.   
Moreover $f_J$ swaps $E^{\re} \leftrightarrow E^{\rv}$ as $[\beta K^{\rf}]$ swaps $P^{\re} \leftrightarrow P^{\rv}$. 
\par
(1b) Take a point $q^{\re} \in Q^{\re}$. 
Its isotropy subgroup in $\vD$ is a group of order $2$,  
whose binary lift is a cyclic subgroup of order $4$ in $H$. 
Let $K^{\re}$ be a generator of the latter group.  
Then $[K^{\re}] \in \PSL_2(\bC)$ swaps $P^{\re} \leftrightarrow P^{\rv}$, so $K^{\re}$ lies in $H \setminus G$ 
and represents the nontrivial class $\sigma_J \in H/G \cong \bZ_2$. 
We may take $K^{\re}$ as an alternative representative of $\sigma_J$, so that $J = \beta K^{\re}$.  
The eigenvalues of $K^{\re}$ are $\pm \ri$. 
Since $q^{\re}$ is fixed by $[\beta K^{\re}] \in \PGL_2(\bC)$, Lemma \ref{lem:actfix} is available with 
$q = q^{\re}$ and $\alpha = \ri$ in formula \eqref{eqn:actfix}. 
It implies that the multipliers of $\beta K^{\re}$ at $q^{\re}$ are $-1$ along the base space $\bP^1$ and 
$-\beta^2 = - (\det J)$ along the fiber $T^*_{q^{\re}} \bP^1$. 
Let $p_{\infty}$ be the image of $Q^{\re}$ under the quotient map $\bP^1 \setminus 
(P^{\re} \cup P^{\rf} \cup P^{\rv}) \to \mathbf{P}^1 \setminus \{ p^{\re}, p^{\rf}, p^{\rv} \}$. 
Since the map \eqref{eqn:TP} is a local biholomorphism, we observe that $p_{\infty}$ is a fixed point 
of $f_J$ in $E_0 \setminus \{ p_0^{\re}, p_0^{\rf}, p_0^{\rv} \}$ having multipliers 
$-1$ along $E_0$ and $-(\det J)$ along $T^*_{p_{\infty}} E_0$. 
\par
(2) Case of type $\rE_6$. 
In this case $G$ is a binary tetrahedral group and $N_{\SL}$ is the unique binary octahedral group $H$ 
that contains $G$ as an index $2$ subgroup. 
Notice that $P^{\re}$ is special to the effect that $n^{\re} \neq n^{\rf} = n^{\rv}$ (see Table \ref{tab:excep}). 
Let $Q^{\re}$, $Q^{\rf}$, $Q^{\rv}$ denote the exceptional orbits for the octahedral group 
$\vD := H/\{ \pm I \}$. 
Then $Q^{\rv} = P^{\re}$ and $Q^{\rf} = P^{\rf} \sqcup P^{\rv}$, while  
$Q^{\re}$ is a non-exceptional orbit of $\vG = G/\{ \pm I \}$. 
Take a point $q^{\rv} \in Q^{\rv} = P^{\re}$.  
Its isotropy subgroup in $\vD$ is a cyclic group of order $4$, whose binary lift is a cyclic subgroup of order $8$ in $H$.  
Let $K^{\rv}$ be a generator of the latter group. 
We can take $K^{\rv}$ as a representative of the nontrivial class $\sigma_J \in H/G \cong \bZ_2$, so $J = \beta K^{\re}$.  
As in (1a) we can show that $f_J$ preserves $E^{\re}$ and has multipliers $- (\det J)^{i+1}$ along 
$E_{i+1}^{\re}$ and $-(\det J)^{-i}$ along $E_i^{\re}$ at the fixed points $p_i^{\re} \in E^{\re}$, $0 \le i \le n^{\re}=1$.  
Moreover $f_J$ swaps $E^{\rf}$ and $E^{\rv}$.  
Next, take a point $q^{\re} \in Q^{\re}$. 
Its isotropy subgroup in $\vD$ is a group of order $2$, whose binary lift is a cyclic subgroup of order $4$ in $H$. 
Let $K^{\re}$ be a generator of the latter group. 
We can take $K^{\re}$ as an alternative representative of $\sigma_J$, so $J = \beta K^{\re}$. 
Let $p_{\infty}$ be the image of $Q^{\re}$ under the quotient map $\bP^1 \setminus 
(P^{\re} \cup P^{\rf} \cup P^{\rv}) \to \mathbf{P}^1 \setminus \{ p^{\re}, p^{\rf}, p^{\rv} \}$. 
As in (1b) we can show that $p_{\infty}$ is a fixed point of $f_J$ in $E_0 \setminus \{ p_0^{\re}, p_0^{\rf}, p_0^{\rv} \}$ 
having multipliers $-1$ along the base $E_0$ and $-(\det J)$ along the fiber $T^*_{p_{\infty}} E_0$.  
\par
(3) Case of type $\rD_4$.   
In this case $G$ is a quaternion group and $N_{\SL}$ is the unique binary octahedral group $H$ 
that contains $G$ as a normal subgroup.  
Notice that $P^{\re}$, $P^{\rf}$, $P^{\rv}$ are ``even" as $n^{\re} = n^{\rf} = n^{\rv}$ (see Table \ref{tab:excep}).   
Let $Q^{\re}$, $Q^{\rf}$, $Q^{\rv}$ be the exceptional orbits for the octahedral group 
$\vD := H/\{ \pm I \}$. 
Then $Q^{\rv} = P^{\re} \sqcup P^{\rf} \sqcup P^{\rv}$, while $Q^{\re}$ is the disjoint union of three 
non-exceptional  orbits $Q^{\re}_{\re}$, $Q^{\re}_{\rf}$, $Q^{\re}_{\rv}$ of the Klein $4$-group 
$\vG = G/ \{ \pm I\}$.  
Corresponding to this decomposition there are three dihedral groups $\vD^{\re}$, $\vD^{\rf}$, $\vD^{\rv}$ in 
$\vD$ that contain $\vG$ as an index $2$ subgroup, where $\vD^{\re}$ preserves $P^{\re}$ and swaps 
$P^{\rf} \leftrightarrow P^{\rv}$, and so on.  
For $\nu = \re, \rf, \rv$, the binary lift $H^{\nu}$ of $\vD^{\nu}$ is a binary dihedral group in 
$H$ that contains $G$ as an index $2$ subgroup. 
Given $J = \beta K$ with $\sigma_J = \langle K \rangle$ of order $2$, we have $K \in H^{\nu} \setminus G$ 
for a unique $\nu = \re, \rf, \rv$.  
Our claim can then be proved by the same argument as that in case (1) with $H^{\nu}$ in place of $H$.    
Meanwhile, $Q^{\rf}$ is the disjoint union of two non-exceptional $\vG$-orbits $Q^{\rf}_{\pm}$, which are 
antipodal to each other on the octahedron.    
Discussion about them is left to Lemma \ref{lem:order3}. \hfill $\Box$
\begin{lemma} \label{lem:order3} 
Suppose that $\sigma_J$ is of order $3$ in $\Aut(D)$. 
Then $D$ must be of type $\rD_4$ and the map $f_J$ permutes $E^{\re}$, $E^{\rf}$, $E^{\rv}$ cyclically. 
Moreover $f_J$ fixes exactly two points $p_{\pm}$ on $E$, which lie in $E_0 \setminus 
\{ p_0^{\re}, p_0^{\rf}, p_0^{\rv} \}$ and at which $f_J$ has multipliers 
$\omega^{\pm 1}$ along $E_0$ and $ \omega^{\mp 1} (\det J)$ along 
$T_{p_{\pm}}^*E_0$, where $\omega$ is a primitive cube root of unity.    
\end{lemma} 
{\it Proof}. 
It is clear from Table \ref{tab:auto} that $D$ must be of type $\rD_4$ and $G$ is a quaternion group. 
We are in the situation of case (3) in the proof of Lemma \ref{lem:order2}. 
Take a pair of antipodal points $q^{\rf}_{\pm} \in Q^{\rf}_{\pm}$. 
Their (common) isotropy subgroup in $\vD$ is a cyclic group of order $3$, whose binary lift is a cyclic 
group of order $6$. 
A generator $K^{\rf}$ of the latter group is a representative of the given class $\sigma_J \in H/G \cong S_3$  
of order $3$, so that $J = \beta K^{\rf}$.  
Since $K^{\rf}$ has eigenvalues $\omega^{\mp 1/2}$ and $[K^{\rf} ] \in \PGL_2(\bC)$ fixes $q^{\rf}_{\pm}$,   
Lemma \ref{lem:actfix} is available with $q = q^{\rf}_{\pm}$ and $\alpha = \omega^{\mp 1/2}$ in formula \eqref{eqn:actfix}. 
It implies that the multipliers of $\beta K^{\re}$ at $q^{\rf}_{\pm}$ are $\omega^{\pm 1}$ along  
$\bP^1$ and $\omega^{\mp 1} \beta^2 = \omega^{\mp 1} (\det J)$ along $T^*_{q^{\rf}_{\pm}} \bP^1$. 
Let $p_{\pm}$ be the image of $Q^{\rf}_{\pm}$ under the quotient map $\bP^1 \setminus 
(P^{\re} \cup P^{\rf} \cup P^{\rv}) \to \mathbf{P}^1 \setminus \{ p^{\re}, p^{\rf}, p^{\rv} \}$. 
Since the map \eqref{eqn:TP} is a local biholomorphism, $p_{\pm}$ are fixed points  
of $f_J$ in $E_0 \setminus \{ p_0^{\re}, p_0^{\rf}, p_0^{\rv} \}$ having multipliers 
$\omega^{\pm 1}$ along $E_0$ and $\omega^{\mp 1} (\det J)$ along the fiber $T^*_{p_{\pm}} E_0$. 
Since $[K^{\rf} ] \in \PGL_2(\bC)$ permutes $P^{\re}$, $P^{\rf}$, $P^{\rv}$ cyclically, $f_J$ permutes 
$E^{\re}$, $E^{\rf}$, $E^{\rv}$ cyclically. 
\hfill $\Box$ 
\section{Linear Models near Exceptional Components} \label{sec:lmnec}
The purpose of this section is to introduce linear model maps around the 
exceptional set of a minimal resolution of a Kleinian singularity. 
They are used in \S \ref{sec:rd} to detect rotation domains around 
exceptional components.  
Let $G$ be a nontrivial finite subgroup of $\SU(2)$.   
The orbit space $S := \bC^2/G$ is realized by an algebraic surface 
$$
S = \{ x = (x_1, x_2, x_3) \in \bC^3 \mid u(x) = 0 \} \qquad 
\mbox{with Poincar\'{e} residue} \quad 
\eta_S := \left. (d x/ d u) \right|_S, 
$$  
where $u(x)$ is a certain polynomial of $x = (x_1, x_2, x_3)$ and $d x := d x_1 \wedge d x_2 \wedge d x_3$. 
The standard $2$-form $\eta_{\bC^2} := d w_1 \wedge d w_2$ on $\bC^2$ descends to $\eta_S$ 
(up to a nonzero constant multiple) via the quotient map $\bC^2 \to S = \bC^2/G$ off the singular point $p \in S$. 
Let $(X, E) \to (S, p)$ be a minimal resolution of $(S, p)$ with exceptional set $E$.  
Then $\eta_S$ lifts to a nonwhere vanishing holomorphic $2$-form $\eta_X$ on $(X, E)$, 
which we call the {\sl satandard $2$-form} on it.  
\par
Let $J \in N_{\GL}$ be a normalizer of $G$ in $\GL_2(\bC)$.      
As indicated in the diagram 
\begin{equation} \label{cd:lm}
\begin{CD} 
J \curvearrowright (\bC^2, 0) + \eta_{\bC^2} @.         f_J \curvearrowright (X, E) + \eta_X \\
@V \mbox{quotient by $G$} VV    @VV \mbox{minimal resolution} V \\
(\bC^2, 0)/G @= (S, p) + \eta_S,  
\end{CD}
\end{equation}
the linear action of $J$ on $(\bC^2, 0)$ descends to $(S, p)$ and then lifts to 
an automorphism $f_J$ on $(X, E)$, whose structure was investigated in \S\S \ref{ss:cc2}--\ref{ss:ncc2}. 
The twist relation $J^* \eta_{\bC^2} = (\det J) \, \eta_{\bC^2}$ transfers to 
$f_J^* \eta_X = (\det J) \, \eta_X$.     
We refer to $f_J$ as the {\sl linear model} associated with the normalizer $J \in N_{\GL}$. 
\par
The other way around, we start with a nonlinear biholomorphism $f : (X, E) \to (X, E)$ 
with twist relation $f^* \eta_X = \delta \, \eta_X$ for a complex constant $\delta = \delta(f)$ of modulus $1$.   
This map descends to $(S, p)$ and then lift to a $G$-equivariant nonlinear biholomorphism  
$F : (\bC^2, 0) \to (\bC^2, 0)$ through the diagram
\begin{equation} \label{cd:rtlm}
\begin{CD} 
f \curvearrowright (X, E) + \eta_X  @. F \curvearrowright (\bC^2, 0) + \eta_{\bC^2}  \\
@V \mbox{collapsing $E$}  VV    @VV \mbox{lift by $G$} V \\
(S, p) + \eta_S   @= (\bC^2, 0)/G. 
\end{CD}
\end{equation}
Let $J := F'(0) \in N_{\GL}$ be the Jacobian matrix of $F$ at the origin.  
Then the map $f_J$ is referred to as the {\sl linear model} of $f$ around $E$. 
We then have $\det J = \delta$, $F^* \eta_{\bC^2} = \delta \, \eta_{\bC^2}$  
and $f_J^* \eta_X = \delta \, \eta_X$, so $f$ and $f_J$ have the same twist number.   
The map $f$ permutes the $(-2)$-curves in $E$ and induces an automorphism 
$\sigma(f) \in \Aut(D)$ of the Dynkin diagram $D$ of $E$, just as $f_J$ induces 
$\sigma_J \in \Aut(D)$. 
\par
We are interested in how $f$ and $f_J$ are different from each other on $E$. 
It is measured by the mappings 
\begin{equation} \label{eqn:hH}
h := f_J^{-1} \circ f : (X, E) \to (X, E), \qquad 
H := J^{-1} \circ F : (\bC^2, 0) \to (\bC^2, 0).  
\end{equation} 
which pass to each other through diagram \eqref{cd:rtlm}.  
It follows from $J := F'(0)$ and formula \eqref{eqn:jacobian} that 
\begin{equation} \label{eqn:cmid}
H'(0) = I, \qquad H \circ g = g \circ H  \qquad \mbox{for any} \quad g \in G. 
\end{equation}
It is possible to determine the map $h|_E$ from a bit of information about $H$.   
First we consider the case where $G$ is a cyclic group of order $n+1$ as  
in \S\ref{ss:cc1}, using the notation there. 
Recall that $u_i$ (resp. $v_i$) is an inhomogeneous coordinate on 
$E_i \cong \bP^1$ around the point $p_{i-1}$ (resp. $p_i$) for $i = 1, \dots, n$. 
We write $H = (H_1(w), H_2(w))$ in terms of the standard coordinate system $w = (w_1, w_2)$ on $(\bC^2, 0)$.  
The following proposition shows that there is a subtle but not negligible  
difference between $f|_E$ and $f_J|_E$, although they are almost the same.  
\begin{proposition} \label{prop:cmid}
Let $G$ be the cyclic group of order $n+1$ generated by the matrix 
$g = \diag(\zeta^{-1}, \zeta)$, where $\zeta := \exp( 2 \pi \ri /(n+1))$.   
For $n = 1$ tha map $h|_E$ is the identity on $E = E_1$.   
For $n \ge 2$ it is described as follows:     
\begin{enumerate}
\setlength{\itemsep}{-1pt} 
\item[$(1)$] $h|_{E_1}$ maps $E_1$ onto itself by $u_1 \mapsto u_1 + b_1$ or equivalently by  
$v_1 \mapsto v_1/(1 + b_1 v_1)$, where $b_1$ is the coefficient of $w_2^n$ in the Taylor series   
expansion of $H_1(w)$ at the origin $w = 0$.    
\item[$(i)$] If $n \ge 3$, then $h|_{E_i}$ is the identity mapping on $E_i$ for each $i = 2, \dots, n-1$.  
\item[$(n)$] $h|_{E_n}$ maps $E_n$ onto itself by $v_n \mapsto v_n + b_2$ or equivalently by  
$u_n \mapsto u_n/(1 + b_2 u_n)$, where $b_2$ is the coefficient of $w_1^n$ in the Taylor series  
expansion of $H_2(w)$ at the origin $w = 0$.    
\end{enumerate}    
In particular $f$ and $f_J$ induce the same Dynkin automorphism, that is, $\sigma(f) = \sigma_J 
\in \Aut(\rA_n)$. 
\end{proposition} 
{\it Proof}. 
First we discuss the main case of $n \ge 2$. 
Forced by the constraints in \eqref{eqn:cmid}, the Taylor series expansions of  
$w_1' =H_1(w)$ and $w_2' = H_2(w)$ at $w = 0$ must be of the forms:    
\begin{align*}
w_1' &= w_1 \left\{  \sum_{j=0}^{\infty} \sum_{k=0}^{\infty} a_{j k}^{(1 1)} \, (w_1 w_2)^j  w_1^{k (n+1)}  
+ \sum_{j=0}^{\infty} \sum_{k=1}^{\infty} a_{j k}^{(1 2)} \, (w_1 w_2)^j w_2^{k (n+1)} 
+ \frac{w_2^n}{w_1} \sum_{j=0}^{\infty} b_j^{(1)} \, w_2^{j(n+1)} \right\} \\
&= w_1 \left\{ \sum_{j=0}^{\infty} \sum_{k=0}^{\infty} a_{j k}^{(1 1)} \, x^j y^k 
+ \sum_{j=0}^{\infty} \sum_{k=1}^{\infty} a_{j k}^{(1 2)} \, x^j  z^k  
+ \frac{z}{x} \sum_{j=0}^{\infty} b_j^{(1)} \, z^j \right\}, \qquad 
a_{0 0}^{(1 1)} = 1, \quad b_0^{(1)} = b_1; 
\\
w_2' &=  w_2 \left\{ \sum_{j=0}^{\infty} \sum_{k=1}^{\infty} a_{j k}^{(2 1)} \, (w_1 w_2)^j w_1^{k (n+1)} 
+ \sum_{j=0}^{\infty} \sum_{k=0}^{\infty} a_{j k}^{(2 2)} \, (w_1 w_2)^j  w_2^{k (n+1)}  
+ \frac{w_1^n}{w_2} \sum_{j=0}^{\infty} b_j^{(2)} \, w_1^{j(n+1)} \right\} \\
&= w_2 \left\{ \sum_{j=0}^{\infty} \sum_{k=1}^{\infty} a_{j k}^{(2 1)} \, x^j y^k 
+ \sum_{j=0}^{\infty} \sum_{k=0}^{\infty} a_{j k}^{(2 2)} \, x^j z^k  
+ \frac{y}{x} \sum_{j=0}^{\infty} b_j^{(2)} \, y^j \right\}, \qquad  
a_{0 0}^{(2 2)} = 1, \quad b_0^{(2)} = b_2, 
\end{align*}
where formula \eqref{eqn:xyz} is used in the second and fourth lines.   
In other words, we have  
\begin{equation} \label{eqn:wd}
w_1' = w_1 \left\{ A_1(x, y, z) + (z/x) \, B_1(z) \right\},  \qquad 
w_2' = w_2 \left\{ A_2(x, y, z) + (y/x) \, B_2(y) \right\},  
\end{equation}
where $A_1 = A_1(x, y, z)$, $A_2 = A_2(x, y, z)$, $B_1 = B_1(z)$ and $B_2 = B_2(y)$ are 
convergent power series in the indicated variables such that 
$A_1(0, 0, 0) = A_2(0, 0, 0) = 1$, $B_1(0) = b_1$ and $B_2(0) = b_2$. 
Regarding $x$, $y$, $z$ as holomorphic functions on $(X, E)$, we see that 
$A_1$, $A_2$, $B_1$ and $B_2$ are holomorphic functions on $(X, E)$ such that 
\begin{equation} \label{eqn:ABE}
A_1 \equiv 1, \qquad A_2 \equiv 1, \qquad B_1 \equiv b_1, \qquad B_2 \equiv b_2 
\qquad \mbox{on} \quad E.  
\end{equation}  
\par
Let $S \to S$, $(x, y, z) \mapsto (x', y', z')$ be the map induced from $H : w \mapsto w'$. 
It follows from \eqref{eqn:xyz} and \eqref{eqn:wd} that 
\begin{align*}
x' &= x \left\{ A_1 A_2 + (y/x) \, A_1 B_2 + (z/x) \, A_2 B_1 + x^{n-1} \, B_1 B_2 \right\}, 
\\[1mm] 
y' &= y \left\{ A_1 + (z/x) \, B_1 \right\}^{n+1}, \qquad 
z' = z \left\{ A_2 + (y/x) \, B_2 \right\}^{n+1}. 
\end{align*}
Substituting these equations into the prime version of formula \eqref{eqn:resol2-3} and 
using formula \eqref{eqn:resol2-1a}, we have  
\begin{equation} \label{eqn:udvd}
u_i' = (x')^i (z')^{-1} = u_i \, C^i \, D^{-n-1}, \qquad 
v_{i-1}' = (x')^{1-i} z' =v_{i-1} \, C^{1-i} \, D^{n+1}, \qquad 
1 \le i \le n+1,  
\end{equation}
where the functions $C = C(u_i, v_{i-1})$ and $D = D(u_i, v_{i-1})$ are defined by 
\begin{align*}
C &:= A_1 A_2 + u_i^{n-i+1} v_{i-1}^{n-i} \, A_1 B_2 + u_i^{i-2} v_{i-1}^{i-1} \, A_2 B_1 
+ u_i^{n-1} v_{i-1}^{n-1} \, B_1 B_2, 
\\[1mm] 
D &:= A_2 + u_i^{n-i+1} v_{i-1}^{n-i} B_2. 
\end{align*}  
\par
Putting $v_{i-1} = 0$ in the first formula of \eqref{eqn:udvd} and taking \eqref{eqn:ABE} into account, 
we observe the following: 
For $n \ge 2$ and $i = 1$ we have $u_1' = u_1 + b_1$ on $E_1$; 
for $n \ge 3$ and $2 \le i \le n-1$ we have $u_i' = u_i$ on $E_i$; 
for $i = n \ge 2$ we have $u_n' = u_n/(1 + b_2 u_n)$ on $E_n$. 
In a similar manner, putting $u_i = 0$ in the second formula of \eqref{eqn:udvd} and taking \eqref{eqn:ABE} 
into account, we observe the following: 
For $n \ge 2$ and $i = 2$ we have $v_1' = v_1/(1 + b_1 v_1)$ on $E_1$; 
for $n \ge 3$ and $3 \le i \le n$ we have $v_{i-1}' = v_{i-1}$ on $E_{i-1}$; 
for $n \ge 2$ and $i = n+1$ we have $v_n' = v_n + b_2$ on $E_n$. 
\par
Next we consider the case $n = 1$. 
Formula \eqref{eqn:xyz} now reads $x = w_1 w_2$, $y = w_1^2$, $z = w_2^2$  
and \eqref{eqn:wd} turns into  
\begin{alignat*}{2}
w_1' &= w_1 \{ A_1(x, y, z) + (w_2^3/w_1) \, B_1(w_2^2) \} = w_1 \{ A_1(x, y, z) + (z^2/x) \, B_1(z) \}, 
\quad & A_1(0,0,0) &= 1, \\[1mm]
w_2' &= w_1 \{ A_2(x, y, z) + (w_1^3/w_2) B_1(w_1^2) \} = w_1 \{ A_2(x, y, z) + (y^2/x) \, B_2(y) \},  
\quad & A_2(0,0,0) &= 1. 
\end{alignat*}
where $A_1$, $A_2$, $B_1$, $B_2$ are convergent power series in the indicated variables. 
A chief difference from formula \eqref{eqn:wd} lies in the terms $z^2/x$ and $y^2/x$ in front 
of $B_1(z)$ and $B_2(y)$; the values of $B_1(0)$ and $B_2(0)$ are not relevant here.  
We have $u_1 = x z^{-1}$ and $v_0 = z$ from formula \eqref{eqn:resol2-3}.  
Proceeding as in the case of $n \ge 2$, we obtain $x' = x \, C$ and $z' = z \, D^2$, hence 
$u_1' = x' (z')^{-1} = x z^{-1} \, C D^{-2} = u_1 \, C D^{-1}$ and $v_0' = z' = z D^2 = v_0 D^2$, 
where    
$$
C := A_1 A_2 + u_1^3 v_0 \, A_1 B_2 + u_1^{-1} v_0 \, A_2 B_1 + u_1^2 v_0^2 \, B_1 B_2, 
\qquad D := A_2 + u_1^3 v_0 \, B_2.  
$$
Putting $v_0 = 0$, we have $u_1' = u_1$ on $E$, that is, $h|_E$ is the identity transformation. 
\hfill $\Box$ 
\begin{remark} \label{rem:cmid} 
When $G$ is a non-cyclic finite subgroup of $\SU(2)$, the construction of a minimal resolution in 
\S \ref{ss:ncc1} and Proposition \ref{prop:cmid} applied to each arm $E^{\nu}$ in Figure \ref{fig:dDE} 
show that the map $h|_E$ fixes each $(-2)$-curve in $E$ pointwise except for the terminal 
curves $E_{n^{\nu}}^{\nu}$, $\nu = \re, \rf, \rv$, on which $h|_E$ behaves as   
$v_{n^{\nu}}^{\nu} \mapsto v_{n^{\nu}}^{\nu} + b^{\nu}$ in terms of the inhomogeneous 
coordinate $v_{n^{\nu}}^{\nu}$ around the point $p_{n^{\nu}}^{\nu} \in E_{n^{\nu}}^{\nu}$, 
where $b^{\nu} \in \bC$ is a suitable constant.  
The formal proof of this fact is omitted, because it is not used in what follows.      
\end{remark}
\par
We say that $f$ is {\sl equivariantly linearizable to its linear model} $f_J$, if the map $F$ 
in diagram \eqref{cd:rtlm} is $G$-equivariantly linearizable to $J$ around $0$.    
In this case $f$ is biholomorphically conjugate to $f_J$ in a neighborhood of $E$,  
so $f$ and $f_J$ exhibit the same dynamics there. 
By definition, the {\sl multipliers of} $f$ {\sl at} $E$ are those of $F$ at $0$,  
that is, the eigenvalues of $J$. 
When $f$ is understood, they are also referred to as the {\sl multipliers of} $E$. 
It is natural to ask when $f$ is equivariantly linearizable to its linear model. 
We give a partial answer to this question.   
\begin{proposition} \label{prop:eltlm} 
Let $\bal = (\alpha_1, \alpha_2) = (\beta \alpha, \beta \alpha^{-1})$ 
be the multipliers of $E$.  
If $\alpha$ is a root of unity while $\beta$ is an algebraic number of modulus $1$, 
but not a root of unity, then $f$ is equivariantly linearizable to its linear model.      
\end{proposition}
{\it Proof}. 
First we show that $\bal$ is NR. 
Suppose the contrary that $\alpha_j - \alpha_1^{n_1} \alpha_2^{n_2} = 0$ 
holds for some $(n_1, n_2) \in N$ and $j = (3 \mp 1)/2$, 
where the set $N$ is defined at the beginning of \S \ref{sec:equilin}.   
Then we have $\beta^{n_1 + n_2-1} = \alpha^{n_2-n_1 \pm 1}$, 
hence $\beta^{k(n_1 + n_2-1)} = \alpha^{k(n_2-n_1 \pm 1)} = 1$ for some $k \in \bN$. 
This forces $n_1+n_2-1 = 0$, which contradicts $(n_1, n_2) \in N$. 
Since $\alpha_1$ and $\alpha_2$ are algebraic numbers, 
$\bal$ is Diophantine by Theorem \ref{thm:tnt}. 
Theorem \ref{thm:elinz} then concludes that $F$ is $G$-equivariantly 
linearizable to $J = F'(0)$, that is, $f$ is equivariantly linearizable to its linear model $f_J$.  
\hfill $\Box$ \par\medskip
There are important cases to which Proposition \ref{prop:eltlm} is applicable.         
\begin{theorem} \label{thm:lm}
Let $(X, E)$ be a minimal resolution of a Kleinian singularity with exceptional set $E$  
and let $f : (X, E) \to (X, E)$ be an automorphism with twist number $\delta$ of modulus $1$.  
Suppose that $\delta$ is an algebraic number which is not a root of unity. 
Then $f$ is equvariantly linearizable to its linear model $f_J$, if either     
\begin{enumerate}
\setlength{\itemsep}{-1pt}
\item[$(\ri)$] $E$ is of Dynkin type $\rA_n$, $n \ge 2$, and $\sigma(f) \in \Aut(\rA_n) \cong \bZ_2$ is nontrivial; or   
\item[$(\rii)$] $E$ is of Dynkin type $\rD$ or $\rE$, in which case $\sigma(f)$ may be arbitrary,     
\end{enumerate}
\end{theorem}
{\it Proof}. 
Let $\beta \alpha^{\pm 1}$ be the multipliers of $E$, that is, the eigenvalues of $J$, where $\beta = \delta^{1/2}$. 
In view of Proposition \ref{prop:eltlm}, it suffices to show that $\alpha$ is a root of unity. 
In case (i), by Proposition \ref{prop:cmid} we have $\sigma(f) = \sigma_J$, so a comparison 
of Lemmas \ref{lem:ida1} and \ref{lem:ida2} shows that $J$ must be of the form \eqref{eqn:fb}, 
{\sl not} of the form \eqref{eqn:fa}. 
This observation tells us that $\alpha = \ri$. 
In case $(\rii)$ we have $J = \delta^{1/2} K$ with $K \in N_{\SL}$, where $N_{\SL}$ is a finite 
subgroup of $\SU(2)$ by Theorem \ref{thm:normalizer}. 
So the eigenvalues $\alpha^{\pm 1}$ of $K$ are roots of unity and those of $J$ are $\beta \alpha^{\pm 1}$ 
with $\beta = \delta^{1/2}$.      
\hfill $\Box$ \par\medskip
Theorem \ref{thm:lm} relies only on a {\sl semi-local} property of the map $f$ 
in a neighborhood of the exceptional set $E$.    
In order to get similar results when $E$ is of Dynkin type $\rA$ and $\sigma(f)$ is trivial, however, 
we need more elaborate machinery involving {\sl  global geometry} of the map $f : X \to X$.  
This is one of the main issues in this article, to which a large part of 
\S\S \ref{sec:pp}--\ref{sec:rd} is devoted.  
\section{Salem Numbers and Multipliers at Periodic Points} \label{sec:pp}
Let $f : X \to X$ be an automorphism of a complex surface $X$. 
Given a positive integer $m$, let $p \in X$ be a periodic point of $f$ with period 
$m$, that is, a fixed point of the $m$-th iterate $f^m$.       
The eigenvalues $\alpha_1$ and $\alpha_2$ of the holomorphic tangent map 
$(d f^m)_p : T_p X \to T_p X$ are referred to as the {\sl multipliers} of $f^m$ at $p$.  
The multiplier pair $\bal = (\alpha_1, \alpha_2)$ is said to be 
{\sl multiplicatively independent} (MI for short) if any relation  
\begin{equation} \label{eqn:MI}
\mbox{ $\alpha_1^{n_1} = \alpha_2^{n_2}$ \, with \, $(n_1, n_2) \in \bZ^2$ \, 
implies \, $n_1 = n_2 = 0$}.  
\end{equation}
Otherwise, $\bal$ is {\sl multiplicatively dependent} (MD for short). 
Clearly, MI implies NR, or equivalently resonance implies MD. 
If the power $\bal^k := (\alpha_1^k, \alpha_2^k)$ is MD (resp. NR) for some $k \in \bZ_{\ge 2}$     
then $\bal$ itself is MD (resp. NR).  
Conversely, if $\bal$ is MD then $\bal^k$ remains MD for every $k \in \bZ_{\ge 2}$.   
Even if $\bal$ is NR, however, $\bal^k$ may be resonant for some $k \in \bZ_{\ge 2}$.   
In this sense, MD and MI are stable under taking powers, but resonance and NR are not.  
\par
Now let $f : X \to X$ be a K3 surface automorphism with twist number $\delta = \delta(f)$.     
Suppose that $f$ satisfies condition $(\rA1)$ in \S \ref{sec:intro}, which says that  
$\delta$ is a non-real conjugate to a Salem number $\lambda$.     
Let $p \in X$ be a periodic point of $f$ with period $m$ and let $\bal = (\alpha_1, \alpha_2)$ 
be the multiplier pair of $f^m$ at $p$.  
The tangent map $(d f^m)_p : T_p X \to T_p X$ has determinant $\delta^m$, so  
the multipliers satisfy $\alpha_1 \alpha_2 = \delta^m$. 
This allows us to write 
\begin{equation} \label{eqn:alpha}
\alpha_1 = \delta^{m/2} \alpha, \qquad 
\alpha_2 = \delta^{m/2} \alpha^{-1}, \qquad 
\alpha \in \bC^{\times},       
\end{equation}
where $\alpha^{\pm 1}$ may be swapped and the branch of $\delta^{1/2}$ is chosen so that it becomes conjugate to the 
{\sl positive} square root $\lambda^{1/2}$ of $\lambda$. 
If $\lambda^{1/2}$ is again a Salem number, then the branch of $\delta^{1/2}$ required is unique; 
otherwise, either branch of it is alright.  
The goal of this section is to establish a {\sl criterion} for deciding whether 
$\bal$ is MI, or MD but NR, or resonant (Theorem \ref{thm:MDR}).      
We shall work with the following hypothesis, or ansatz,            
\begin{hypothesis} \label{hyp:P}  
There exists a rational function $\hat{P}(w) \in \bQ(w)$ with rational coefficients such that 
\begin{equation}  \label{eqn:P1}
\alpha + \alpha^{-1} = \hat{P}( \hat{\tau}), \qquad \hat{\tau} := \delta^{1/2} + \delta^{-1/2}, 
\end{equation}
where $\hat{P}(w)$ is an {\sl even} function for even $m$ and an {\sl odd} function for odd $m$ respectively.    
\end{hypothesis}
\par
Imposing Hypothesis \ref{hyp:P} is not unreasonable, because applying a holomorphic 
Lefschetz-type fixed point formula to a K3 surface automorphism often gives rise 
to such a function $\hat{P}(w)$ as in \eqref{eqn:P1}; 
see \S\S \ref{sec:rd}--\ref{sec:pc} for details.  
Our main goal in this section is Theorem \ref{thm:MDR} below.  
It is motivated by applications to dynamics on K3 surfaces, but the theorem itself is 
logically independent of K3 surface automorphisms.   
\par
To formulate the result we introduce some terminology and notation. 
A polynomial $u(z)$ of even degree, say $2 n$, is said to be {\sl palindromic} 
if $u(z) = z^{2 n} \, u(z^{-1})$. 
If this is the case, then there exists a unique polynomial $U(w)$ of degree $n$, 
called the {\sl trace polynomial} of $u(z)$, such that $u(z) = z^n \, U(z + z^{-1})$. 
The most basic examples are the trace polynomials $H_n(w)$ of $z^{2n} + 1$, 
which are the unique polynomials $H_n(w) \in \bZ[w]$ such that  
\begin{equation} \label{eqn:H}
H_n(z + z^{-1}) = z^n + z^{-n}, \qquad n \in \bZ_{\ge 0}.   
\end{equation}
Although of the same notation, $H_1(w) = w$ and $H_2(w) = w^2-2$ here have nothing 
to do with those in \S \ref{sec:lmnec}.   
In terms of the $n$-th Chebyshev polynomial $T_n(w)$ of the first kind we have 
$H_n(w) = 2 T_n(w/2)$; see e.g. Rivlin \cite{Rivlin}.  
Notice that $H_{m n}(w) = H_m(H_n(w)) = H_n(H_m(w))$ holds for all $m, n \in \bZ_{\ge 0}$.  
\begin{remark} \label{rem:Q}  
In Hypothesis \ref{hyp:P}, thanks to its parity assumption, the function $\hat{P}(w)$ can be expressed as 
\begin{equation} \label{eqn:P2}
\hat{P}(w) = P(H_2(w)) \quad \mbox{when $m$ is even}; \qquad  
\hat{P}(w) = w^{-1} P(H_2(w)) \quad \mbox{when $m$ is odd},  
\end{equation}
for some rational function $P(w) \in \bQ(w)$, where $H_2(w) = w^2-2$. 
If we define a function $Q(w) \in \bQ(w)$ by  
\begin{equation} \label{eqn:Q2}
Q(w) := P^2(w) -2 \quad \mbox{when $m$ is even}; \qquad 
Q(w) := (w+2)^{-1} P^2(w) - 2 \quad \mbox{when $m$ is odd}, 
\end{equation}
then applying $H_2(\,\cdot\,)$ to equation \eqref{eqn:P1} gives the following equation, 
free from the square root of $\delta$,    
\begin{equation} \label{eqn:Q1}
\alpha^2 + \alpha^{-2} = Q(\tau), \qquad \tau := H_2(\hat{\tau}) = \delta + \delta^{-1}. 
\end{equation}
From definition \eqref{eqn:Q2} we have $Q(w) \ge -2$ on the interval $w > -2$, 
although $Q(w)$ may have some poles there. 
\end{remark} 
\par
For each integer $n \ge 3$, the $n$-th cyclotomic polynomial $\rC_n(z)$ is palindromic 
of even degree $\phi(n) \ge 2$, where $\phi(n)$ is Euler's totient, so we can 
speak of its trace polynomial $\rCT_n(w)$, that is, the $n$-th {\sl cyclotomic 
trace polynomial}. 
The special case of $n = 1, 2$ needs some care.  
If we employ the nonstandard convention 
\begin{equation} \label{eqn:conv}
\rC_1(z) = (z-1)^2, \quad \rC_2(z) := (z+1)^2, \quad \rCT_1(w) := w-2, 
\quad \rCT_2(w) := w+2, \quad \phi(1) = \phi(2) = 2,   
\end{equation}
then all cases $n \in \bN$ can be treated in a unified manner. 
Convention \eqref{eqn:conv} is valid only within this section. 
\par
A {\sl reciprocal number} is an algebraic integer $\xi \neq \pm 1$ 
that is conjugate to its reciprocal $\xi^{-1}$ over $\bQ$.  
In this case the sum $\eta := \xi + \xi^{-1}$ is called the {\sl trace} of $\xi$. 
The minimal polynomial of $\xi$ is a palindromic polynomial of even degree,   
and its trace polynomial is the minimal polynomial of the trace $\eta$.     
Typical examples of reciprocal numbers are roots of unity (other than $\pm 1$) 
and Salem numbers, whose traces are {\sl cyclotomic traces} and {\sl Salem traces}, respectively. 
In passing from reciprocal numbers to their traces, the polynomials $H_n(w)$ arise naturally. 
Taking powers $\xi \mapsto \xi^n$ of a reciprocal number $\xi$ induces the 
operations $\eta \mapsto H_n(\eta)$ on its trace $\eta := \xi + \xi^{-1}$.       
\begin{theorem} \label{thm:MDR} 
Let $\bal = (\alpha_1, \alpha_2)$ be a multiplier pair as in formula \eqref{eqn:alpha}, 
where $\delta$ is a non-real conjugate to a Salem number $\lambda$. 
Suppose that $\alpha$ satisfies equation \eqref{eqn:P1} in Hypothesis $\ref{hyp:P}$.  
Let $\rho := \lambda + \lambda^{-1}$ be the trace of $\lambda$.  
\begin{enumerate}
\setlength{\itemsep}{-1pt} 
\item[$(\rI)$] 
If $Q(\rho) \in [-2, \, 2]$, then $\bal$ is MD exactly when there exists an integer  
$k \in \bN$ such that 
\begin{equation} \label{eqn:MDI} 
\rCT_k(Q(\rho)) = 0, \qquad \phi(k) \, | \, (\deg \lambda),  
\end{equation}
where $\phi(k)$ is Euler's totient and $\deg \lambda$ is the degree 
of $\lambda$ as an algebraic number.  
In this case $\bal$ is NR. 
\item[$(\rII)$] If $Q(\rho) > 2$, then $\bal$ is MD exactly when there 
exist a coprime $(j, k) \in \bN^2$ and a Salem trace $\sigma$ such that 
\begin{equation} \label{eqn:MDII}
H_j(\sigma) = Q(\rho), \qquad \rho = H_k(\sigma).   
\end{equation}
In this case, $(j, k)$ and $\sigma$ are unique, and $\bal$ 
is NR if and only if $(j, k)$ satisfies the conditions:     
\begin{enumerate}
\setlength{\itemsep}{-1pt} 
\item[$(\rN 1)$] $\vD := k m - j \ge 2$, 
\item[$(\rN 2)$] $\vD$ does not divide $m$, 
\item[$(\rN 3)$] if $\vD$ is even, $\vD$ divides $2 m$ and $(2 m)/\vD$ is odd, 
then $\hat{P}(\hat{\rho}) < 0$,   
\end{enumerate}
where $\hat{\rho} :=\lambda^{1/2} + \lambda^{-1/2}$ with 
$\lambda^{1/2}$ being the positive square root of $\lambda$. 
\end{enumerate}
\end{theorem}
\par
When $\vD = 2$, $(\rN 2)$ implies that $m$ must be odd, 
so $(\rN 3)$ simply requires $\hat{P}(\hat{\rho}) < 0$.   
In this theorem we must always keep in mind that $j$ and $k$ should be coprime, 
even when this constraint is not mentioned explicitly.         
\begin{remark} \label{rem:MDR} 
All criteria in Theorem $\ref{thm:MDR}$ are stated only in terms of  the 
Salem number $\lambda$, period $m$, rational function $\hat{P}(w)$ and 
a couple of quantities uniquely determined by them.      
\begin{enumerate}
\setlength{\itemsep}{-1pt} 
\item In case $(\rI)$ the division relation $\phi(k) \, | \, (\deg \lambda)$ in 
condition \eqref{eqn:MDI} confines the integer $k$ into finite possibilities.  
For example, if the Salem number $\lambda$ is the dynamical degree 
of a K3 surface automorphism, then $\deg \lambda$ is an even integer between $4$ and $22$.  
There are exactly $43$ numbers $k \in \bN$ such that $\phi(k) \le 22$, 
which are given in Table \ref{tab:EulerT}. 
Hence any number $k$ satisfying condition \eqref{eqn:MDI}, if it exists, lies in this table. 
\begin{table}[h]
\centerline{
\begin{tabular}{clccl}
\hline
              &        & &              &       \\[-3mm]
$\phi(k)$ & $k$  & $\phantom{a}$ & $\phi(k)$ & $k$ \\[1mm]
\hline
              &        & &             &       \\[-3mm]
$2$ & $1$, $2$ \, by convention \eqref{eqn:conv} & & $12$ & $13$, $21$, $26$, $28$, $36$, $42$ \\[1mm]
$2$ & $3$, $4$, $6$ & & $14$ & none \\[1mm]
$4$ & $5$, $8$, $10$, $12$ & & $16$ & $17$, $32$, $34$, $40$, $48$, $60$ \\[1mm]
$6$ & $7$, $9$, $14$, $18$ & & $18$ & $19$, $27$, $38$, $54$ \\[1mm]
$8$ & $15$, $16$, $20$, $24$, $30$ & & $20$ & $25$, $33$, $44$, $50$, $66$ \\[1mm]
$10$ & $11$, $22$ & & $22$ &  $23$, $46$ \\[1mm]  
\hline
\end{tabular}} 
\caption{Those numbers $k \in \bN$ which satisfy $\phi(k) \le 22$.} 
\label{tab:EulerT} 
\end{table} 
\item In case $(\rII)$, once $\lambda$ and $\hat{P}(w)$ are given explicitly, it is only a task of 
finite procedures to search for those triples $(j, k; \sigma)$ which satisfy condition 
\eqref{eqn:MDII}.  
In fact, there exists a bound $(j_0, k_0) \in \bN^2$ depending only on  
$\rho$ and $Q(w)$ such that the pair $(j, k) \in \bN^2$ in condition \eqref{eqn:MDII} 
must satisfy $j \le j_0$ and $k \div k_0$. 
Moreover, the Salem trace $\sigma$ can be determined from $(\rho, k)$ uniquely and concretely.  
We refer to Remark \ref{rem:MDR2} for the details of how $(j_0, k_0; \sigma)$ is determined.   
If $\bal$ is NR then $j$ must also satisfy $ j \le k m -2$ by criterion (N1).     
\end{enumerate} 
\end{remark}
\par
An easy consequence of Theorem \ref{thm:MDR} is the following corollary, which 
already appears in \cite[Lemma 7.1]{IT2}. 
\begin{corollary} \label{cor:MDR} 
In the situation of Theorem $\ref{thm:MDR}$, $(\ri)$ if $Q(\rho)$ is not an algebraic integer, 
or $(\rii)$ if $Q(\rho') > 2$ for some Galois conjugate $\rho'$ to $\rho$ but different from $\rho$ 
itself,  then the pair $\bal = (\alpha_1, \alpha_2)$ is MI.  
\end{corollary}
\par
Theorem \ref{thm:MDR} is established at the end of \S \ref{ss:MD} after providing  
three lemmas (Lemmas \ref{lem:I}, \ref{lem:mMD} and \ref{lem:mNR}).     
\subsection{Non-Resonance} \label{ss:NR}  
In a general setting we consider when an MD pair $\bal = (\alpha_1, \alpha_2)$ is NR. 
The result of this subsection, Proposition \ref{prop:NR}, is used in \S \ref{ss:MD} 
where Salem numbers come to play.     
We say that a relation $\alpha_1^{N_1} = \alpha_2^{N_2}$ for $\bal$ is  
{\sl primitive}, if for any relation $\alpha_1^{n_1} = \alpha_2^{n_2}$ there exists 
an integer $k$ such that $n_1 = k N_1$ and $n_2 = k N_2$.   
If at least one of $\alpha_1$ and $\alpha_2$ is a root of unity, then 
$\bal$ is obviously resonant.   
So we assume that  
\begin{equation} \label{eqn:NRU} 
\mbox{neither of $\alpha_1$ and $\alpha_2$ is a root of unity}. 
\end{equation}
\begin{lemma} \label{lem:MD-NR}
Suppose that $\bal$ is MD and satisfies condition \eqref{eqn:NRU}. 
Then $\bal$ admits a unique primitive relation $\alpha_1^{N_1} = \alpha_2^{N_2}$ 
such that $N_1 \ge 1$ and $N_2 \neq 0$.  
Moreover the pair $\bal$ is NR if and only if either  
\begin{equation} \label{eqn:MD-NR} 
N_1 \ge 2, \quad N_2 \ge 2; \qquad \mbox{or} \qquad N_1 = N_2 = 1. 
\end{equation}
\end{lemma}    
{\it Proof}. 
Let $\alpha_1^{n_1} = \alpha_2^{n_2}$ be any nontrivial relation for $\bal$.  
Then neither of $n_1$ and $n_2$ is zero by condition \eqref{eqn:NRU}.  
We may take $n_1$ to be positive upon replacing $(n_1, n_2)$ by    
$(-n_1, -n_2)$ when $n_1$ is negative. 
Let $\alpha_1^{N_1} = \alpha_2^{N_2}$ be the nontrivial relation with the 
least positive value of $N_1 \ge 1$. 
Notice that $N_2$ is nonzero.     
Given any relation $\alpha_1^{n_1} = \alpha_2^{n_2}$, Euclidean algorithm 
allows us to write $n_1 = k N_1 + r$ with $0 \le r < N_1$ to have   
$$
\alpha_2^{n_2} = \alpha_1^{n_1} = \alpha_1^{k N_1+ r} 
= \alpha_1^{r} (\alpha_1^{N_1})^k = \alpha_1^r (\alpha_2^{N_2})^k 
 = \alpha_1^r \alpha_2^{k N_2}, \quad 
\mbox{that is}, \quad \alpha_1^r = \alpha_2^{n_2-k N_2}. 
$$
The minimality of $N_1$ forces $r = 0$ and $\alpha_2^{n_2-k N_2} = 1$. 
The former equation $r = 0$ yields $n_1 = k N_1$, while the latter equation 
implies $n_2 = k N_2$, because $\alpha_2$ is not a root of unity.  
Thus $\alpha_1^{N_1} = \alpha_2^{N_2}$ is a primitive relation such that 
$N_1 \ge 1$ and $N_2 \neq 0$. 
The uniqueness of such a primitive relation is clear. 
\par
Next we show that $\bal$ is resonant exactly when one of the following 
three conditions is satisfied:  
$$
(\mathrm{a}) \quad N_2 \le -1; \qquad 
(\mathrm{b}) \quad N_1 = 1, \quad N_2 \ge 2; \qquad 
(\mathrm{c}) \quad N_1 \ge 2, \quad N_2 =1.  
$$
In case (a) the pair $\bal$ is resonant, because the relation $\alpha_1^{N_1} = \alpha_2^{N_2}$ 
yields $\alpha_1 = \alpha_1^{N_1+1} \alpha_2^{|N_2|}$ with $(N_1+1)+|N_2| \ge 3$. 
Suppose that $N_2 \ge 1$ and $\bal$ is resonant. 
Then there exist integers $n_1 \ge 0$ and $n_2 \ge 0$ with $n_1 + n_2 \ge 2$
such that either $(\mathrm{i})$ $\alpha_1 = \alpha_1^{n_1} \alpha_2^{n_2}$  
or $(\mathrm{ii})$ $\alpha_2 = \alpha_1^{n_1} \alpha_2^{n_2}$.  
In case $(\mathrm{i})$ we have $\alpha_1^{1-n_1} = \alpha_2^{n_2}$, 
so there exists an integer $k$ such that $1-n_1 = k N_1$, i.e. $n_1 = 1- k N_1$ 
and $n_2 = k N_2$.   
Then we have $2 \le n_1 + n_2 = 1 + k (N_2-N_1)$, that is, $1 \le k(N_2-N_1)$, 
which implies $k \neq 0$. 
Since $0 \le n_2 = k N_2$ and $1 \le N_2$, we have $1 \le k$, so 
$0 \le n_1 = 1 - k N_1 \le 1-1 \cdot 1 = 0$, thus $k = N_1 = 1$. 
The inequality $1 \le k(N_2-N_1)$ now becomes $2 \le N_2$. 
Hence case (i) implies case (b). 
Similarly case $(\mathrm{ii})$ implies case (c).  
The converse implications are obvious.  
The logical negation of $(\mathrm{a}) \vee (\mathrm{b}) \vee (\mathrm{c})$ 
is condition \eqref{eqn:MD-NR}, which is a necessary and sufficient 
condition for $\bal$ to be NR.  \hfill $\Box$
\begin{proposition} \label{prop:NR} 
Suppose that $\bal = (\alpha_1, \alpha_2)$ is MD and  
$\alpha_1 \alpha_2$ is not a root of unity. 
Then the two conditions  
\begin{enumerate}
\setlength{\itemsep}{-1pt}
\item[$(\rD 1)$] there exists a relation $\alpha_1^{n_1} = \alpha_2^{n_2}$  
such that $n_1 n_2 \ge 1$,  
\item[$(\hat{\rD} 1)$]  any nontrivial relation $\alpha_1^{m_1} = \alpha_2^{m_2}$ 
satisfies $m_1 m_2 \ge 1$,    
\end{enumerate} 
are equivalent. 
Moreover, $\bal$ is NR exactly when $(\rD 1)$, or $(\hat{\rD} 1)$, and   
the following condition are satisfied:  
\begin{enumerate}
\setlength{\itemsep}{-1pt}
\item[$(\rD 2)$] $\alpha_1 \neq \alpha_2^{n}$ and $\alpha_2 \neq \alpha_1^{n}$ 
for any integer $n \ge 2$.      
\end{enumerate} 
\end{proposition}   
{\it Proof}. 
The implication $(\hat{\rD}1) \Rightarrow (\rD 1)$ is obvious. 
We show the converse $(\rD 1) \Rightarrow (\hat{\rD} 1)$. 
Notice that (D1) implies condition \eqref{eqn:NRU}. 
Indeed, if one of $\alpha_1$ and $\alpha_2$ is a root of unity, then so is the other 
by (D1) and hence $\alpha_1 \alpha_2$ is also a root of unity, contrary to the assumption. 
By Lemma \ref{lem:MD-NR} the pair $\bal $ admits a primitive relation 
$\alpha_1^{N_1} = \alpha_2^{N_2}$ with $N_1 \ge 1$ and $N_2 \neq 0$.  
Then in (D1) we have $n_1 = k N_1$, $n_2 = k N_2$ and  
$1 \le n_1 n_2 = k^2 N_1 N_2$ for some $k \in \bZ$, hence 
$N_1 \ge 1$ and $N_2 \ge 1$.  
So, for any nontrivial relation $\alpha_1^{m_1} = \alpha_2^{m_2}$ we have 
$m_1 = j N_1$, $m_2 = j N_2$ and $m_1 m_2 = j^2 N_1 N_2 \ge 1$ for some 
$j \in \bZ \setminus \{ 0\}$, which verifies $(\hat{\rD} 1)$. 
Now assume (D2) in addition to (D1).   
By (D2) the primitive relation $\alpha_1^{N_1} = \alpha_2^{N_2}$ just mentioned must satisfy condition \eqref{eqn:MD-NR}, 
so Lemma \ref{lem:MD-NR} implies that $\bal$ is NR .  
Conversely, suppose that $\bal$ is NR. 
Then $\bal$ satisfies condition \eqref{eqn:NRU}, so Lemma \ref{lem:MD-NR} tells us 
that $\bal$ admits a primitive relation $\alpha_1^{N_1} = \alpha_2^{N_2}$ 
satisfying condition \eqref{eqn:MD-NR}.  
In particular it satisfies $N_1 N_2 \ge 1$ and serves as a relation satisfying (D1). 
Condition (D2) is also a consequence of the property 
\eqref{eqn:MD-NR} for $(N_1, N_2)$. \hfill $\Box$
\subsection{Some Properties of Salem Numbers} \label{ss:recipr} 
\par
As a preliminary to the next section we provide some basic properties of Salem numbers.   
\begin{lemma} \label{lem:salem}  
For Salem numbers as defined in $\S \ref{sec:intro}$, the following properties hold.   
\begin{itemize} 
\setlength{\itemsep}{-1pt}
\item[$(\rS 1)$] If $\lambda > 1$ is an algebraic integer all of whose 
conjugates other than $\lambda$ lie in the closed unit disk with at least 
one conjugate on the unit circle, then $\lambda$ is a Salem number.  
\item[$(\rS 2)$] If $\lambda$ is a Salem number of degree $d$, then so is 
the power $\lambda^n$ and $\bQ(\lambda) = \bQ(\lambda^n)$ for any $n \in \bN$.    
\item[$(\rS 3)$] If $\lambda$ is a Salem number, then there exists a unique 
Salem number $\lambda_0$, called the $\mathrm{primitive}$ of $\lambda$, such that 
any Salem number in $\bQ(\lambda)$ is a positive power of $\lambda_0$.  
\end{itemize}
\end{lemma} 
\par
For these statements and their proofs we refer to  
Smyth \cite[Lemmas 1, 2 and Proposition 3]{Smyth}.  
In accordance with (S3), if $\lambda_0$ is the primitive of a Salem umber $\lambda$,  
then the Salem trace $\rho_0 := \lambda_0 + \lambda_0^{-1}$ is referred to as the 
{\sl primitive} of the Salem trace $\rho := \lambda + \lambda^{-1}$. 
We have $\rho = H_{k_0}(\rho_0)$ if $\lambda$ is the $k_0$-th power of $\lambda_0$.  
\begin{lemma} \label{lem:primitive} 
Let $\rho$ and $\sigma$ be Salem traces and suppose that $\rho = H_k(\sigma)$ with $k \in \bN$. 
If $\rho_0$ is the primitive of $\rho$ and $\rho = H_{k_0}(\rho_0)$ with $k_0 \in \bN$,  
then $k_0$ is divisible by $k$, the Salem trace $\rho_0$ is the primitive of $\sigma$, 
and $\sigma = H_{k_0/k}(\rho_0)$. 
\end{lemma}
{\it Proof}. 
Let $\lambda$, $\mu$ and $\lambda_0$ be the Salem numbers whose traces are 
$\rho$, $\sigma$ and $\rho_0$ respectively. 
The relation $\rho = H_k(\sigma)$ then means $\lambda + \lambda^{-1} = \mu^k + \mu^{-k}$, 
which implies $\lambda = \mu^k$ since $\lambda > 1$ and $\mu^k > 1$. 
Similarly we have $\lambda = \lambda_0^{k_0}$. 
From $(\rS 2)$ in Lemma \ref{lem:salem} we have $\bQ(\mu) = \bQ(\mu^k) = \bQ(\lambda)$, 
so $(\rS3)$ shows that $\lambda_0$ is also the primitive of $\mu$, hence $\mu = \lambda_0^m$ 
for some $m \in \bN$.  
Thus we have $\lambda_0^{k_0} = \lambda = \mu^k = \lambda_0^{m k}$ and hence $k_0 = m k$. 
Therefore $k_0$ is divisible by $k$ with quotient $k_0/k = m$, and we have 
$\sigma = \mu + \mu^{-1} = \lambda_0^m + \lambda_0^{-m} = H_m(\rho_0) = H_{k_0/k}(\rho_0)$.  
\hfill $\Box$
\begin{remark} \label{rem:MDR2} 
This is a continuation of Remark \ref{rem:MDR}.(2). 
The number $k_0 \in \bN$ mentioned there is given as the integer   
such that $\rho = H_{k_0}(\rho_0)$, where $\rho_0$ is the primitive 
of the Salem trace $\rho$. 
By Lemma \ref{lem:primitive} the integer $k \in \bN$ in condition \eqref{eqn:MDII} divides $k_0$ and 
determines the Salem trace $\sigma$ as $\sigma = H_{k_0/k}(\rho_0)$. 
From it we have $\sigma \ge \rho_0$. 
Note that $H_j(w)$ is strictly increasing in $j \in \bN$ and in $w > 2$.      
Take $j_0 \in \bN$ to be the least integer such that $H_{j_0}(\rho_0) \ge Q(\rho)$, 
which certainly exists since $H_j(\rho_0) \to \infty$ as $j \to \infty$ strictly increasingly.   
For any $j > j_0$ we have $H_j(\sigma) > H_{j_0}(\sigma) \ge H_{j_0}(\rho_0) \ge Q(\rho)$. 
So the integer $j \in \bN$ in condition \eqref{eqn:MDII} must satisfy $j \le j_0$ and $k \div k_0$.  
\end{remark}
\begin{lemma} \label{lem:power} 
If two Salem numbers $\lambda$ and $\mu$ are related by 
$\lambda^n = \mu^m$ for some $m$, $n \in \bN$, then there exist a unique Salem 
number $\nu$ and a unique coprime pair $(j, k) \in \bN^2$ such that 
$\lambda = \nu^j$ and $\mu = \nu^k$. 
\end{lemma}
{\it Proof}. 
From $(\rS 2)$ in Lemma \ref{lem:salem} we have 
$\bQ(\mu) = \bQ(\mu^m) = \bQ(\lambda^n) = \bQ(\lambda)$. 
Then by $(\rS 3)$ the Salem numbers $\lambda$ and 
$\mu$ have the same primitive $\lambda_0$ and admit power representations  
$\lambda = \lambda_0^s$ and $\mu = \lambda_0^t$ for some $s$, $t \in \bN$.   
Letting $\nu := \lambda_1^g$ with $g := \gcd\{ s, t\}$, we have 
$\lambda = \nu^j$ and $\mu = \nu^k$, where $j := s/g$ and $k := t/g$ are coprime. 
The uniqueness of $\nu$, $j$ and $k$ is obvious from the structure of positive powers 
of a Salem number.  \hfill $\Box$
\subsection{Multiplicative Dependence and Non-Resonance} \label{ss:MD}
We give three lemmas, the first of which is prepared for case $(\rI)$ of Theorem \ref{thm:MDR}, 
while the second and third ones are for case $(\rII)$.    
They immediately enable us to establish Theorem \ref{thm:MDR} at the end of this subsection.  
\par
Let $\pi \in \Gal(\overline{\bQ}/\bQ)$ be a Galois automorphism such that $\pi(\delta) = \lambda$.   
It converts equation \eqref{eqn:Q1} into 
\begin{equation} \label{eqn:zeQ}
\zeta + \zeta^{-1} = Q(\rho), \qquad \zeta := \pi(\alpha^2).
\end{equation} 
We note that $\bal$ is MD if and only if $\bal^2 := (\alpha_1^2, \alpha_2^2)$ is MD.   
A nontrivial relation for $\bal^2$ is given by  
\begin{equation} \label{eqn:alal}
\alpha_1^{2 n_1} = \alpha_2^{2 n_2} 
\end{equation} 
where $n_1, n_2 \in \bZ$ not both zero.  
Under representation \eqref{eqn:alpha} this relation is equivalent to 
\begin{equation} \label{eqn:aldel}
\alpha^{2(n_1+n_2)} = \delta^{m(n_1-n_2)},  \tag{$\ref{eqn:alal}'$}
\end{equation}
where $n_1 + n_2$ must be nonzero, for otherwise $n_1-n_2$ is nonzero and equation 
\eqref{eqn:aldel} becomes $\delta^{m(n_1-n_2)} = 1$, which contradicts the fact 
that $\delta$ is not a root of unity. 
We may and shall assume that 
\begin{equation} \label{eqn:n1n2}
 (n_1 + n_2)(n_1-n_2) \ge 0 \qquad \mbox{with} \quad n_1 + n_2 \neq 0,  
\end{equation}
after swapping $n_1$ and $n_2$ (or swapping $\alpha^{\pm 1}$) if necessary.  
Since $\delta$ is a non-real conjugate to the Salem number $\lambda$, equation 
\eqref{eqn:aldel} implies that $\alpha$ is an algebraic integer of modulus $1$.  
Moreover $\pi$ converts \eqref{eqn:aldel} into equation  
\begin{equation} \label{eqn:zela}
\zeta^{n_1+n_2} = \lambda^{m(n_1-n_2)}. 
\end{equation}
\begin{lemma} \label{lem:I}
Suppose that $Q(\rho) \in [-2, \, 2]$. 
\begin{enumerate}
\setlength{\itemsep}{-1pt}  
\item If $\bal$ is MD, then $\rCT_k(Q(\rho)) = 0$ and $\phi(k) \, | \,(\deg \lambda)$ for some $k \in \bN$.  
\item If $\rCT_k(Q(\rho)) = 0$ for some $k \in \bN$, then $\bal$ is MD and NR.  
\end{enumerate}
\end{lemma}
{\it Proof}. 
Since $Q(\rho) \in [-2, \, 2]$, equation \eqref{eqn:zeQ} forces $|\zeta| = 1$. 
With convention \eqref{eqn:conv}, equation \eqref{eqn:zeQ} yields   
\begin{equation} \label{eqn:CTkQ}
\rC_k(\zeta) = \zeta^{\phi(k)/2} \, \rCT_k(\zeta + \zeta^{-1}) = \zeta^{\phi(k)/2} \, \rCT_k(Q(\rho))    
\qquad \mbox{for any} \quad k \in \bN.   
\end{equation}
\par
(1) Suppose that $\bal$ is MD and let \eqref{eqn:alal} be a nontrivial relation for $\bal^2$.  
Equation \eqref{eqn:zela} yields $\lambda^{m(n_1-n_2)} = |\zeta|^{n_1 + n_2} = 1$. 
Since $\lambda > 1$, we must have $n_1- n_2 = 0$, so equation \eqref{eqn:aldel} becomes $\alpha^{4 n_1} = 1$. 
Thus $\alpha$ is a root of unity and so is $\zeta = \pi(\alpha^2)$.  
Let $k \in \bN$ be the integer such that $\zeta$ is a primitive $k$-th root of unity. 
Since $\rC_k(\zeta) = 0$, equation \eqref{eqn:CTkQ} yields $\rCT_k(Q(\rho)) = 0$. 
In the field extensions $\bQ \subset \bQ(Q(\rho)) \subset \bQ(\rho)$ we have
$$
[\bQ(\rho) : \bQ] = \deg \rho = (\deg \lambda)/2, \qquad  
[\bQ(Q(\rho)) : \bQ] = \deg (\zeta + \zeta^{-1}) = \phi(k)/2.  
$$
Thus $[\bQ(\rho) : \bQ(Q(\rho))] = (\deg \lambda)/\phi(k)$ is an integer, that is, $\phi(k)$ divides $\deg \lambda$. 
\par
(2) Suppose that $\rCT_k(Q(\rho)) = 0$. 
Then equation \eqref{eqn:CTkQ} implies $\rC_k(\zeta) = 0$, so $\zeta^k = 1$  
and hence $\alpha^{2 k} = 1$ by $\alpha^2 = \pi^{-1}(\zeta)$. 
Now a nontrivial relation \eqref{eqn:aldel} holds with $n_1 = n_2 = k$, thus $\bal$ is MD.     
In \eqref{eqn:alpha} the number $\alpha$ is a root of unity, 
while $\delta^{m/2}$ is not, so an argument in the proof of Proposition \ref{prop:eltlm} shows that $\bal$ is NR.    
\hfill $\Box$
\begin{lemma} \label{lem:mMD} 
Suppose that $Q(\rho) > 2$. 
\begin{enumerate}
\setlength{\itemsep}{-1pt}  
\item Suppose that $\bal$ is MD. 
Then  there exist a Salem number $\mu$ and a coprime pair $(j, k) \in \bN^2$ 
such that $\lambda = \mu^k$ and for any Galois automorphism $\pi$ sending 
$\delta \mapsto \lambda$ we have $\mu^j = \pi(\alpha^2)$ up to swapping of $\alpha^{\pm 1}$.    
Moreover equations \eqref{eqn:MDII} hold with Salem trace $\sigma := \mu + \mu^{-1}$.  
\item If equations \eqref{eqn:MDII} hold for some Salem trace $\sigma$ and some 
coprime pair $(j, k) \in \bN^2$, then $\bal$ is MD.   
\end{enumerate} 
\end{lemma}
{\it Proof}. Since $Q(\rho) > 2$, equation \eqref{eqn:zeQ} implies that 
$\zeta$ is real, positive, and $\zeta \neq 1$. 
\par
(1) Suppose that $\bal$ is MD and let \eqref{eqn:alal} be a nontrivial relation for $\bal^2$. 
Using $0 < \zeta \neq 1$, $1 < \lambda$ and \eqref{eqn:n1n2} in equation 
\eqref{eqn:zela}, we have $(n_1 + n_2)(n_1 - n_2) \ge 1$ and $\zeta > 1$. 
Thus equation \eqref{eqn:zeQ} determines the number $\zeta = \pi(\alpha^2)$ uniquely and independently of $\pi$. 
Let $\varpi$ be any Galois automorphism such that $|\varpi(\alpha^2)| > 1$. 
Since $\varpi$ takes \eqref{eqn:aldel} into equation 
$(\varpi(\alpha^2))^{n_1+n_2} = (\varpi(\delta))^{m(n_1-n_2)}$, we have $|\varpi(\delta)| > 1$. 
Thus $\varpi(\delta)$ is conjugate to the Salem number $\lambda$ and has modulus $|\varpi(\delta)| > 1$.  
Such a number must be the Salem number $\lambda$ itself, so $\varpi(\delta) = \lambda$, 
hence $\varpi$ is a Galois automorphism such that $\varpi(\delta) = \lambda$. 
This forces $\varpi(\alpha^2) = \zeta$. 
Therefore $\zeta$ is the unique conjugate to $\alpha^2$ with modulus $|\zeta| > 1$, 
while $\alpha^2$ itself is an algebraic integer with modulus $1$.  
To put it in another way, $\zeta > 1$ is an algebraic integer all of whose conjugates 
other than $\zeta$ lie in the closed unit disk, with its conjugate $\alpha^2$ on the unit circle.    
It follows from $(\rS 1)$ in Lemma \ref{lem:salem} that $\zeta$ is a Salem number.     
Lemma \ref{lem:power} can now be applied to equation \eqref{eqn:zela} to conclude 
that there exist a unique Salem number $\mu$ and a unique coprime pair $(j, k) \in \bN^2$ 
such that $\pi(\alpha^2) = \zeta = \mu^j$ and $\lambda = \mu^k$. 
We have $\rho := \lambda + \lambda^{-1} = \mu^k + \mu^{-k} = H_k(\sigma)$, 
the second equation in \eqref{eqn:MDII}, and equation \eqref{eqn:zeQ} becomes 
$Q(\rho) = \zeta + \zeta^{-1} =\mu^j + \mu^{-j} = H_j(\mu + \mu^{-1}) = H_j(\sigma)$, 
the first equation in \eqref{eqn:MDII}. 
\par
(2)  Suppose that equations \eqref{eqn:MDII} hold for a Salem trace $\sigma$ and a  
coprime pair $(j, k) \in \bN^2$. 
Let $\mu$ be the Salem number whose trace is $\sigma$, i.e. $\mu + \mu^{-1} = \sigma$. 
From the second equation in \eqref{eqn:MDII} we have $\lambda + \lambda^{-1} = \rho = H_k(\sigma) = \mu^k + \mu^{-k}$.  
Since $\lambda > 1$ and $\mu^k > 1$, we have $\lambda = \mu^k$. 
The first equation in \eqref{eqn:MDII} and equation \eqref{eqn:zeQ} yield   
$\mu^j + \mu^{-j} = H_j(\sigma) = Q(\rho) = \zeta + \zeta^{-1}$. 
So we have $\mu^j = \zeta = \pi(\alpha^2)$ up to possible inversion of $\alpha$.      
Put $n_1 := k m + j$ and $n_2 := k m - j$. 
Then $n_1 + n_2 = 2 k m \ge 1$, $n_1-n_2 = 2 j \ge 1$ and 
$(\pi(\alpha^2))^{n_1 + n_2} = (\mu^j)^{2 k m} = (\mu^k)^{2 j m} =\lambda^{m(n_1-n_2)} = (\pi(\delta))^{m(n_1-n_2)}$. 
The last equation is brought back to \eqref{eqn:aldel} or equivalently to \eqref{eqn:alal}   
by the inverse Galois automorphism $\pi^{-1}$. 
Therefore $\bal$ is MD.  
\hfill $\Box$
\begin{lemma} \label{lem:mNR} 
Suppose that $Q(\rho) > 2$ and $\bal$ is MD. 
Then $\bal$ is NR if and only if the pair $(j, k) \in \bN^2$ in Lemma $\ref{lem:mMD}.(1)$  
satisfies all of the conditions $(\rN 1), (\rN 2), (\rN 3)$ in Theorem $\ref{thm:MDR}$.     
\end{lemma}
{\it Proof}. 
Since $\bal$ is MD, as in \eqref{eqn:alal} or \eqref{eqn:aldel} for $\bal^2$, we can think of any nontrivial 
relation for $\bal$ in the form  
\begin{equation} \label{eqn:alaldel}
\alpha_1^{n_1} = \alpha_2^{n_2}, \qquad \mbox{or equivalently} \qquad 
\alpha^{n_1+ n_2} = \delta^{m(n_1-n_2)/2}.   
\end{equation}
Here we may and shall assume \eqref{eqn:n1n2}.  
From the way the branch of $\delta^{1/2}$ is chosen, there exists a Galois automorphism 
$\pi$ such that $\pi(\delta^{1/2}) = \lambda^{1/2}$; see the sentence containing 
formula \eqref{eqn:alpha}. 
Then $\pi(\delta) = \lambda = \mu^k$, so we have $\pi(\alpha^2) = \mu^j$ 
by a statement in Lemma \ref{lem:mMD}.(1).  
Thus we have $\pi(\alpha) = \ve \, \mu^{j/2}$ for some sign $\ve = \pm 1$, where   
$\mu^{1/2}$ is the positive square root of $\mu > 1$.  
Then $\pi$ converts equations \eqref{eqn:P1} and \eqref{eqn:alaldel} into 
\begin{equation} \label{eqn:mujk}
\ve \cdot H_j (\hat{\sigma}) = \hat{P}(\hat{\rho} ), 
\qquad 
\ve^{n_1 + n_2} \, \mu^{j (n_1 + n_2)/2} = \mu^{m k (n_1-n_2)/2}, 
\end{equation} 
where $\hat{\sigma} := \mu^{1/2} + \mu^{-1/2} > 2$ so that $H_j(\hat{\sigma}) > 2$.  
Applying $H_2(\,\cdot\,)$ to the first equation in \eqref{eqn:mujk} leads to 
the first equation in \eqref{eqn:MDII}, which was already established in Lemma \ref{lem:mMD}.   
So conditions \eqref{eqn:mujk} are equivalent to       
\begin{equation} \label{eqn:jk}
\ve \cdot \hat{P}( \hat{\rho} ) > 0, \qquad \ve^{n_1+n_2} = 1, \qquad 
j (n_1+ n_2) = k m (n_1-n_2).  
\end{equation}
\par
Suppose that $\bal$ is NR. 
It then follows from $(\hat{\rD} 1)$ in Proposition \ref{prop:NR} that $n_1 n_2 \ge 1$. 
This together with the last equation in \eqref{eqn:jk} yields $1 \le |n_1 - n_2| < |n_1 + n_2|$ 
and $j < k m$, that is, $\vD := k m - j \ge 1$.  
Conversely, assuming that $\vD \ge 1$, we put $n_1 := k m + j \ge 2$ and 
$n_2 := \vD = k m -j \ge 1$ to obtain $n_1 + n_2 = 2 k m$ and $n_1 - n_2 = 2 j$. 
Then the last two equations in \eqref{eqn:jk} are fulfilled. 
The first inequality in \eqref{eqn:jk} is also true, if $\ve = \pm 1$ is the sign of 
$\hat{P}(\hat{\rho}) \in \bR^{\times}$. 
Thus we can get a relation for $\bal$ as in \eqref{eqn:alaldel} such that $n_1 n_2 \ge 1$, 
so condition $(\rD 1)$ in Proposition \ref{prop:NR} is satisfied. 
Now that $\bal$ is MD and $\alpha_1 \alpha_2 = \delta^m$ is not a root of unity, 
Proposition \ref{prop:NR} is available in the current situation.       
Therefore $\bal$ is NR if and only if $\vD \ge 1$ and $\bal$ satisfies condition 
$(\rD 2)$ in Proposition \ref{prop:NR}.  
The latter condition is equivalent to say (up to swapping of $n_1$ and $n_2$) that 
{\sl no} pair $(n_1, \ve)$ of integer $n_1 \ge 2$ and sign $\ve = \pm 1$ can satisfy condition \eqref{eqn:jk} with $n_2 = 1$, 
that is,   
\begin{equation} \label{eqn:jk2} 
\ve \cdot \hat{P}( \hat{\rho} ) > 0, \qquad \ve^{n_1+1} = 1, \qquad 
j (n_1 + 1) = k m (n_1-1).  \tag{$\ref{eqn:jk}'$}
\end{equation}
The other way around, we ask when \eqref{eqn:jk2} admits a solution $(n_1, \ve)$. 
Hereafter we assume $\vD \ge 1$.   
\par\smallskip
{\bf Claim 1}. 
Condition \eqref{eqn:jk2} is satisfied by an {\sl odd} integer $n_1 \ge 2$ and a sign 
$\ve = \pm 1$ if and only if $\vD \div m$. 
\par\smallskip
Suppose that \eqref{eqn:jk2} is satisfied by an odd integer $n_1 \ge 2$ and a sign 
$\ve = \pm 1$.  
Then the last equation in \eqref{eqn:jk2} reads $j (n_1+1)/2 = k m (n_1-1)/2$, 
where the positive integers $(n_1+1)/2$ and $(n_1-1)/2$ are coprime  
as they differ only by $1$.    
So there exists an integer $l \in \bN$ such that $j = l (n_1-1)/2$ 
and $k m = l (n_1+1)/2$. 
Subtracting the former equation from the latter yields $\vD := k m-j = l$,  
showing that $\vD$ divides both $j$ and $k m$. 
As $j$ and $k$ are coprime, $\vD$ must divide $m$.  
Conversely, suppose that $\vD$ divides $m$ and put 
$n_1 = 2 k m' -1$ with $m' := m/\vD \in \bN$. 
Clearly $n_1$ is odd and $n_1 \ge 1$.  
More strictly we have $n_1 \ge 2$, for otherwise $1 = n_1 = 2 k m'-1$ would lead 
to $k = m' = 1$, $m = \vD$ and hence $j = k m - \vD = m - \vD = 0$, 
contradicting $j \ge 1$. 
We have also   
$$
\vD (n_1-1)/2 = \vD(k m' -1) = k m -\vD = j, 
\qquad 
\vD (n_1+1)/2 = \vD k m' = k m,   
$$
which readily lead back to the last equation in \eqref{eqn:jk2}.   
We can determine the sign $\ve = \pm 1$ so that the inequality 
$\ve \cdot P( \hat{\rho} ) > 0$ holds true. 
Equation $\ve^{n_1+1} = 1$ is automatically fulfilled as $n_1 + 1$ is even. 
This proves Claim $1$. 
\par\smallskip
{\bf Claim 2}.    
Condition \eqref{eqn:jk2} is satisfied by an {\sl even} integer $n_1 \ge 2$ and 
a sign $\ve = \pm 1$ if and only if 
\begin{equation} \label{eqn:S3}
\mbox{$\vD$ is even, $\vD$ divides $2 m$, 
$(2 m)/\vD$ is odd, and $\hat{P}(\hat{\rho}) > 0$}.  
\end{equation}  
\par
Suppose that \eqref{eqn:jk2} is satisfied by an even integer $n_1 \ge 2$ 
and a sign $\ve = \pm 1$.   
The equation $\ve^{n_1 + 1} = 1$ forces $\ve = 1$, since $n_1 + 1$ is odd. 
Now the first inequality in \eqref{eqn:jk2} reads $\hat{P}(\hat{\rho}) > 0$.    
In the last equation in \eqref{eqn:jk2} we notice that $n_1 + 1$ and $n_1 -1$ 
are coprime, as $n_1 \pm 1$ are odd integers which differ by $2$.  
So there exists an integer $l \in \bN$ such that $j = l (n_1-1)$  
and $k m = l (n_1+1)$.  
Subtracting the former equation from the latter yields $\vD = 2 l$, 
which shows that $\vD$ is even and $\vD/2$ divides both of $j$ and $k m$.   
As $j$ and $k$ are coprime, $\vD/2$ must divide $m$, that is, 
$\vD$ must divide $2 m$.  
We have also $(2m)/\vD = m/l = (n_1 + 1)/k \in \bN$ with $n_1 + 1$ odd. 
This shows that $(2 m)/\vD$ must be an odd integer.  
Therefore we have \eqref{eqn:S3}. 
Conversely, suppose that \eqref{eqn:S3} holds. 
Notice that $k$ must be odd, for otherwise $j = k m - \vD$ would be even, 
so $j$ and $k$ could not be coprime. 
Take the odd integer $k' := (2m)/\vD \in \bN$ and put $n_1 := k k'-1$.   
Obviously, $n_1$ is even and $n_1 \ge 0$.   
If $n_1 = 0$, then $k = k' = 1$, so $2 m = \vD = k m -j = m-j$ and hence 
$j = -m \le -1$, contradicting $j \ge 1$.  
Thus $n_1$ is even and $n_1 \ge 2$. 
Since $\vD$ is even, we can write $\vD = 2 l$ for some $l \in \bN$. 
The definition of $n_1$ then leads to $k m = l (n_1+1)$ and 
$j = k m - \vD = l(n_1 + 1) - 2 l = l(n_1-1)$, which give 
the last equation in \eqref{eqn:jk2}.  
The first inequality and the second equation in \eqref{eqn:jk2} are 
fulfilled by putting $\ve = 1$. 
This proves Claim $2$.  
\par
Logical negation of $\vD \div m$ in Claim $1$ gives condition $(\rN 2)$,  
which together with $\vD \ge 1$ yields $(\rN 1)$. 
Finally, the logical negation of condition \eqref{eqn:S3} in Claim $2$ leads to condition $(\rN 3)$.   
So the lemma is established.  \hfill $\Box$ \par\medskip
{\sl Proof of Theorem $\ref{thm:MDR}$}. 
Case $(\rI)$ follows from Lemma \ref{lem:I}, while case $(\rII)$ from 
Lemmas \ref{lem:mMD} and \ref{lem:mNR}. \hfill $\Box$
\section{Fixed Point Formulas} \label{sec:fpf} 
Two kinds of fixed point formulas (FPFs for short) are important tools in our analysis. 
They were already used in our previous articles \cite{IT1,IT2}, but the main theme of this section is  
a thorough study of topological and holomorphic local indices around exceptional sets,  
especially around exceptional components of type $\rA$ with trivial symmetry 
(\S \ref{ss:linec}, Theorem \ref{thm:ecA}), as well as components of other types 
or other symmetries (\S \ref{ss:oec}). 
\subsection{Reformulations of Fixed Point Formulas} \label{ss:reform}
First we reformulate the FPFs used in \cite{IT1,IT2}.  
Let $f : X \to X$ be a K3 surface automorphism satisfying conditions (A1) and (A2) in \S \ref{sec:intro}. 
Let $\Fixi(f)$ denote the set of all isolated fixed points of $f$, and $N(f)$ the number of 
$(-2)$-curves in $X$ fixed pointwise by $f$. 
From \cite[Proposition 4.2]{IT2} we have a {\sl topological} Lefschetz-type FPF,      
\begin{equation} \label{eqn:saito} 
L(f) := 2 + \Tr( f^*|H^2(X, \bC)) = \sum_{p \in \Fixi(f)} \mu(f, p) + 2 N(f),    
\end{equation}
where $L(f)$ is the Lefschetz number of $f$ and the local index $\mu(f, p)$ of $f$ at $p \in \Fixi(f)$ is defined by 
$$
\mu(f, p) := \dim_{\bC} (\bC\{ z \}/\mathfrak{a}), \qquad \mbox{where} \quad 
\mathfrak{a} := (z_1-f_1(z), \, z_2-f_2(z)), 
$$
with $(f_1, f_2)$ being the local representation of $f$ in terms of a local chart 
$z = (z_1, z_2)$ around $p$ : $z =(0, 0)$. 
The right-hand side of equation \eqref{eqn:saito}  is the sum of local contributions of isolated fixed points 
and fixed $(-2)$-curves. 
Given an $f$-invariant exceptional component $E \subset \cE(X)$, we put 
\begin{equation} \label{eqn:mufE}
\mu(f, E) := \sum_{p \in \Fixi(f, E)} \mu(f, p) + 2 N(f, E),  
\end{equation}
where $\Fixi(f, E)$ is the set of all {\sl isolated} fixed points on $E$ and $N(f, E)$ is the number 
of fixed $(-2)$-curves in $E$.  
In is convenient to divide the right-hand side of formula \eqref{eqn:saito} into {\sl on} and {\sl off} parts:   
\begin{equation} \label{eqn:muof}
\mu_{\ron}(f) := \sum_{E\, : \, \mathrm{inv.}} \mu(f, E), \qquad 
\mu_{\roff}(f) := \sum_{p \, : \, \roff} \mu(f, p),   
\end{equation}
where the fomer sum extends over all $f$-invariant exceptional components $E \subset \cE(X)$, 
while the latter sum extends over all fixed points $p$ off $\cE(X)$; note that they are isolated. 
FPF \eqref{eqn:saito} is then represented as  
\begin{equation} \label{eqn:saito2}
L(f) = \mu_{\ron}(f) + \mu_{\roff}(f).  
\end{equation}
\par
From \cite[Proposition 4.3]{IT2} we also have a {\sl holomorphic} Lefechetz-type FPF,  
\begin{equation} \label{eqn:tt} 
\mathcal{L}(f) := 1+ \delta^{-1}  = \sum_{p \in \Fixi(f)} \nu(f, p) + N(f) \, \dfrac{\delta + 1}{(\delta-1)^2},  
\end{equation} 
where $\mathcal{L}(f)$ is the holomorphic Lefschetz number of $f$, $\delta = \delta(f)$ is the twist number of $f$ 
and the holomorphic local index $\nu(f, p)$ is given by the Grothendieck residue   
$$
\nu(f, p) := \Res_p \, \dfrac{d z_1 \wedge d z_2}{(z_1-f_1(z))(z_2-f_2(z))}. 
$$
The right-hand side of \eqref{eqn:tt}  is the sum of local contributions 
of isolated fixed points and fixed $(-2)$-curves. 
\par
It is convenient to rewrite FPF \eqref{eqn:tt} in terms of 
$\hat{\tau} := \delta^{1/2} + \delta^{-1/2}$. 
We define $\hat{\nu}(f, p) := \delta^{1/2} \, \nu(f, p)$ and  
\begin{subequations} \label{eqn:nuhat}
\begin{align}
\hat{\nu}(f, E) &:= \sum_{p \in \Fixi(f, E)} \hat{\nu}(f, p) + N(f, E) \, 
\delta^{1/2} \dfrac{\delta+1}{(\delta-1)^2},  \label{eqn:nuhat1}
\\[1mm]   
\hat{\nu}_{\ron}(f) &:= \sum_{E \, : \, \mathrm{inv.}} \hat{\nu}(f, E), \qquad 
\hat{\nu}_{\roff}(f) := \sum_{p \, : \, \roff} \hat{\nu}(f, p),   \label{eqn:nuhat2}
\end{align}
\end{subequations}
where the summation rules in \eqref{eqn:nuhat2} are the same as those in \eqref{eqn:muof}. 
At a {\sl simple} isolated fixed point $p$ we have 
$$
\hat{\nu}(f, p) = \dfrac{1}{\hat{\tau}-(\alpha + \alpha^{-1})}, \qquad 
\hat{\tau} := \delta^{1/2} + \delta^{-1/2},  
$$
when $f$ has multipliers $\delta^{1/2} \alpha^{\pm 1}$ at this point.  
FPF \eqref{eqn:tt} is now rewritten as   
\begin{equation} \label{eqn:tt2} 
\hat{\tau} = \hat{\nu}_{\ron}(f) + \hat{\nu}_{\roff}(f).  
\end{equation} 
\par
Let $E \subset \cE(X)$ be an $f$-invariant exceptional component of Dynkin type $D$. 
The map $f$ then permutes the $(-2)$-curves in $E$ and induces an automorphism of the Dynkin diagram $D$,  
\begin{equation} \label{eqn:sigma} 
\sigma(f, E) \in \Aut(D).
\end{equation}
We say that the pair $(f, E)$, or the set $E$ when $f$ is understood, has 
{\sl trivial}, {\sl reflectional}, or {\sl tricyclic symmetry}, if $\sigma(f, E)$ has order $1$, $2$, or $3$, respectively.  
Note that $\sigma(f, E)$ is written $\sigma(f)$ in \S \ref{sec:lmnec} with $E$ being understood.     
\begin{notation} \label{not:iec} 
For a Dynkin diagram $D$ of type $\rA$, $\rD$ or $\rE$, when there is no fear of confusion, 
we often denote an $f$-invariant exceptional component of type $D$ by the same symbol $D$. 
Moreover, if $D$ has trivial, reflectional, or tricyclic symmetry, then we write it as 
$D^{\rt}$, $D^{\rr}$, or $D^{\rc}$, respectively, where $D^{\rc}$ occurs only when $D = \rD_4$.  
\end{notation} 
\par
Although FPFs \eqref{eqn:saito2} and \eqref{eqn:tt2} are formulated for K3 surface 
automorphisms $f$ satisfying conditions $(\rA 1)$ and $(\rA 2)$, the quantities 
$\mu(f, E)$, $\hat{\nu}(f, E)$, $\sigma(f, E)$ etc. are {\sl semi-local} objects which make sense 
for the map $f$ restricted to a small neighborhood of $E$, and hence for a map on 
a minimal resolution of a Kleinian singularity.
\subsection{Exceptional Components of Type A with Trivial Symmetry} \label{ss:linec}
Back to the situation in \S\ref{sec:lmnec}, let $(X, E)$ be a minimal resolution 
of a Kleinian singularity with exceptional set $E$, and $f : (X, E) \to (X, E)$ an 
automorphism with twist number $\delta$.  
We shall evaluate the numbers $\mu(f, E)$ and $\hat{\nu}(f, E)$ defined by formulas \eqref{eqn:mufE} 
and \eqref{eqn:nuhat1} respectively, especially when $E$ is of type $\rA$ with trivial symmetry. 
\par
Suppose that $E$ is of Dynkin type $\rA_n$ with trivial symmetry.   
Let $\bal = (\alpha_1, \, \alpha_2) = (\delta^{1/2} \alpha, \, \delta^{1/2} \alpha^{-1})$ be 
the multipliers of $E$, that is, the eigenvalues of the Jacobian matrix $J := F'(0) \in N_{\GL}$, 
where $F : (\bC^2, 0) \to (\bC^2, 0)$ is the automorphism induced from $f$ via diagram \eqref{cd:rtlm}.  
Let us consider the geometric progression 
\begin{equation} \label{eqn:gamma}
\gamma_j := \alpha_1^j \alpha_2^{j-n-1} = \delta^{(2 j-n-1)/2} \alpha^{n+1}, \qquad 0 \le j \le n+1,    
\end{equation}
with common ratio $\alpha_1 \alpha_2 = \delta$. 
Suppose that $\delta$ is not a root of unity. 
Then the map $j \mapsto \gamma_j$ is one-to-one, so there is no or only one index 
$i$ such that $\gamma_i = 1$.  
Therefore we have the following trichotomy.    
\begin{enumerate}
\setlength{\itemsep}{-1pt} 
\item[(C1)] $\gamma_j \neq 1$ for all indices $j = 0, 1, \dots, n+1$; 
\item[(C2)] $\gamma_i = 1$ for a unique index $i \in\{1, \dots, n \}$;   
\item[(C3)] $\gamma_0 = 1$ or $\gamma_{n+1} = 1$. 
\end{enumerate}
\par
In order to give an explicit formula for $\hat{\nu}(f, E)$, in addition to the polynomials $H_n(w)$ in definition \eqref{eqn:H}, 
we need to use the polynomials $K_n(w) \in \bZ[w]$ uniquely determined by the equations  
$$
z^{n+1} - z^{-(n+1)} = (z - z^{-1}) \, K_n(z + z^{-1}),  \qquad n \in \bZ_{\ge 0}. 
$$
In terms of the $n$-th Chebyshev polynomial $U_n(w)$ of the second kind we have $K_n(w) = U_n(w/2)$;  
see \cite{Rivlin}.   
\begin{remark} \label{rem:HK} 
We note the following ``even-odd" formulas which are used somewhere or elsewhere;    
\begin{alignat*}{2} 
H_{2 n}(w) &= H_n(H_2(w)), \qquad & H_{2 n +1}(w) &= w \{ K_n(H_2(w)) - K_{n-1}(H_2(w)) \}, \\[1mm] 
K_{2 n}(w) &= K_n(H_2(w)) + K_{n-1}(H_2(w)), \qquad & K_{2 n+1}(w) &= w K_n(H_2(w)). 
\end{alignat*}
\end{remark} 
\par
The main result of this section is stated as follows, where we often denote $E$ by  
$\rA_n^{\rt}$ following Notation \ref{not:iec}. 
\begin{theorem} \label{thm:ecA} 
Let $f : (X, E) \to (X, E)$ be an automorphism on a minimal resolution of a Kleinian 
singularity of Dynkin type $\rA_n$, $n \in \bN$.  
Suppose that the twist number $\delta$ of $f$ is not a root of unity and that the exceptional set 
$E$ has trivial symmetry and multipliers 
$\bal = (\alpha_1, \, \alpha_2) = (\delta^{1/2} \alpha, \, \delta^{1/2} \alpha^{-1})$. 
Then we have
\begin{alignat}{3}
\mu(f, \rA_n^{\rt}) &= n + 1, \qquad & 
\hat{\nu}(f, \rA_n^{\rt}) &= 
\dfrac{K_n(\hat{\tau})}{H_{n+1}(\hat{\tau}) -  \{ \alpha^{(n+1)} +\alpha^{-(n+1)} \} } 
\qquad && \mbox{in case $(\rC 1)$},    
\label{eqn:Ai} \\[2mm]
\mu(f, \rA_n^{\rt}) &= n + 1,  \qquad &  
\hat{\nu}(f, \rA_n^{\rt}) &= 
\dfrac{K_n(\hat{\tau})}{(H_2(\hat{\tau})-2) \, K_{i-1}(\hat{\tau}) \, K_{n-i}(\hat{\tau})} 
\qquad && \mbox{in case $(\rC 2)$},   
\label{eqn:Aii}    
\end{alignat} 
where $\hat{\tau} := \delta^{1/2} + \delta^{-1/2}$.   
In case $(\rC 3)$ we have $\mu(f, \rA_n^{\rt}) \ge n+2$.    
\end{theorem}
\begin{remark} \label{rem:ecA0} 
For the multipliers $\bal = (\alpha_1, \alpha_2)$ of $E$ we consider its $(n+1)$-st 
power $\bal^{n+1} := (\alpha_1^{n+1}, \, \alpha_2^{n+1}) = 
(\delta^{(n+1)/2} \alpha^{n+1}, \, \delta^{(n+1)/2} \alpha^{-(n+1)})$.   
In case $(\rC 1)$ of Theorem \ref{thm:ecA}, if $\hat{\nu}(f, \rA_n^{\rt})$ is known  
in some way or another, say, by using FPFs \eqref{eqn:saito2} and \eqref{eqn:tt2}, then the 
second equation in \eqref{eqn:Ai} may be settled to have 
\begin{equation} \label{eqn:P3}
\alpha^{n+1} + \alpha^{-(n+1)} = \hat{P}(\hat{\tau}), \qquad \hat{\tau} := \delta^{1/2} + \delta^{-1/2}, 
\end{equation}
for some rational function $\hat{P}(w) \in \bQ(w)$ as in Hypothesis \ref{hyp:P}, where  
$m = n+1$ and $\alpha$ is replaced by $\alpha^{n+1}$ in formula \eqref{eqn:alpha}.  
Then Theorem \ref{thm:MDR} may be used to examine whether $\bal^{n+1}$ is 
MI, or MD but NR, or resonant. 
This idea will be implemented in later sections, for example, in the proofs of Theorems 
\ref{thm:RHAt} and \ref{thm:appA}.    
\end{remark}
\begin{remark} \label{rem:ecA} 
From dynamical point of view a remark should be in order.     
From \eqref{eqn:gamma} we have $(\gamma_{n+1}, \, \gamma_0^{-1}) = (\alpha_1^{n+1}, \, \alpha_2^{n+1}) = \bal^{n+1}$.    
The geometric progression \eqref{eqn:gamma} for the $m$-th iterate $f^m$ is the $m$-th power, $\gamma_j^m$, 
$0 \le j \le n+1$, of the original progression for $f$.   
Thus Theorem \ref{thm:ecA} implies that cases $(\rC 2)$ and $(\rC 3)$ are stable 
under iterations of $f$, while case $(\rC 1)$ may switch to $(\rC 2)$ or $(\rC 3)$. 
However, $(\rC 1)$ never switches to $(\rC 3)$, if    
\begin{equation} \label{eqn:rem-ecA}
\mbox{neither of the multipliers $\alpha_1$ and $\alpha_2$ is a root of unity in case $(\rC 1)$}.  
\end{equation}
In particular we have $\mu(f^m, \rA_n^{\rt}) = n+1$ for all $m \in \bN$, 
in case $(\rC 2)$ unconditionally, while in case $(\rC 1)$ under condition \eqref{eqn:rem-ecA}.  
In case $(\rC 3)$ we always have $\mu(f^m, \rA_n^{\rt}) \ge n+2$ for all $m \in \bN$. 
In view of formula \eqref{eqn:gamma}, we also remark that if $\bal$ is MI then 
the iterates $f^m$ stay in case $(\rC 1)$ for all $m \in \bN$.  
\end{remark}
\par
{\it Proof of Theorem $\ref{thm:ecA}$}. 
The theorem follow immediately from Lemmas \ref{lem:ecA1}, \ref{lem:ecA2} and \ref{lem:ecA3} below.   
\hfill $\Box$ \par\medskip
To state and prove these lemmas we make some preparations. 
Suppose that $E$ consists of $(-2)$-curves $E_1, \dots, E_n$,  
the dual graph of which is diagram \eqref{eqn:dA}.       
Suppose for the moment that $n \ge 2$. 
By Lemma \ref{lem:ida1} with $\beta := \delta^{1/2}$ and Proposition \ref{prop:cmid}, 
in terms of the inhomogeneous coordinates 
$u_j$ and $v_j$ on $E_j \cong \bP^1$, $u_j v_j = 1$, in \S \ref{ss:cc1}, 
the M\"{o}bius transformation $f|_{E_j} = (f_J \circ h)|_{E_j} : E_j \to E_j$,   
$u_j \mapsto u_j'$ and $v_j \mapsto v_j'$ is given by     
\begin{subequations} \label{eqn:eudvd}
\begin{alignat}{3} 
u_1' &= \gamma_1 (u_1 + b_1), \qquad & 
v_1' &= \frac{\gamma_1^{-1} \, v_1}{1 + b_1 v_1}, \qquad &  j &= 1, 
\label{eqn:eudvd1} \\[1mm]
u_j' &= \gamma_j \, u_j, \qquad &  
v_j' &= \gamma_j^{-1} \, v_j, \qquad & 2 &\le j \le n-1, 
\label{eqn:eudvd2} \\[1mm] 
u_n' &= \frac{ \gamma_n \, u_n }{1 + b_2 u_n }, \qquad & 
v_n' &= \gamma_n^{-1} (v_n + b_2), \qquad &  j  &= n.  
\label{eqn:eudvd3}  
\end{alignat}
\end{subequations}
\par
The points $\{ p_j \} = E_j \cap E_{j+1}$, $1 \le j \le n-1$, are the same as those 
in formula \eqref{eqn:pi} in \S \ref{ss:cc1}, where $p_j$ is given by $(u_{j+1}, v_j) = (0, 0)$ in local coordinates.  
Due to the constants $b_1$ and $b_2$ appearing in \eqref{eqn:eudvd1} and \eqref{eqn:eudvd3}, however, 
the points $p_0 \in E_1$ and $p_n \in E_n$ in formula \eqref{eqn:p0pn} should be redefined as 
$$
\mbox{$p_0 \in E_1$ \, by \, $u_1 = \dfrac{\gamma_1 b_1}{1-\gamma_1}$ \, if \, $\gamma_1\neq 1$}; 
\qquad  
\mbox{$p_n \in E_n$ \, by \, $v_n = \dfrac{\gamma_n^{-1} b_2}{1-\gamma_n^{-1}}$ \, if \, $\gamma_n \neq 1$}.  
$$
Notice that when $\gamma_1 \neq 1$, the new $p_0$ is the unique fixed point of $f$ on $E_1 \setminus \{ p_1\}$.  
When $\gamma_1 = 1$, we do not define $p_0$, as it fails to make sense. 
Similarly, when $\gamma_n = 1$, we do not define $p_n$. 
Formulas \eqref{eqn:eudvd} tell us the following: 
For each $1 \le j \le n-1$, the multipliers of $f$ at $p_j$ are 
$\gamma_{j+1} = \beta^{2 j+1-n} \alpha^{n+1}$ along $E_{j+1}$ and $\gamma_j^{-1} = 
\beta^{n+1-2j} \alpha^{-n-1}$ along $E_j$. 
If $\gamma_1 \neq 1$, the multipliers of $f$ at $p_0$ are  
$\gamma_1 = \beta^{1-n} \alpha^{n+1}$ along $E_1$ and $\gamma_0^{-1} = \beta^{n+1} \alpha^{-n-1}$ 
in a normal direction to $E_1$. 
If $\gamma_n \neq 1$, the multipliers of $f$ at $p_n$ are   
$\gamma_n^{-1} = \beta^{1-n} \alpha^{-n-1}$ along $E_n$ and $\gamma_{n+1} = \beta^{n+1} \alpha^{n+1}$ 
in a normal direction to $E_n$. 
These observations are used to calculate $\hat{\nu}(f, E)$ in Lemmas \ref{lem:ecA1} and \ref{lem:ecA2}.  
\begin{remark} \label{rem:symm} 
Replacing $\alpha$ by $\alpha^{-1}$ in the multipliers $\beta \alpha^{\pm 1}$ of $E$ 
results in the changes $\gamma_j \leftrightarrow \gamma_{n+1-j}^{-1}$, $b_1 \leftrightarrow b_2$ 
and $p_j \leftrightarrow p_{n-j}$.  
Up to this symmetry, we may assume $i \in \{1, \dots, n-1\}$ in case $(\rC 2)$ and 
$\gamma_0 = 1$ in case $(\rC 3)$.  
\end{remark} 
\par
Let $n = 1$. 
Proposition \ref{prop:cmid} then implies $f|_{E_1} = f_J|_{E_1}$ and formula \eqref{eqn:fJn1} shows 
that $f|_{E_1}$ is given by formula \eqref{eqn:eudvd1}, where $\gamma_1 = \alpha^2$. 
We define $p_1 \in E_1$ by $v_1 = 0$, while $p_0 \in E_1$ by $u_1 = \gamma_1 b_1 (1-\gamma_1)^{-1}$ 
if $\gamma_1 \neq 1$. 
When $\gamma_1 = 1$, we do not define $p_0$. 
Then what we mentioned in the case of $n \ge 2$ carries over in this case.           
\begin{lemma} \label{lem:ecA1}   
In case $(\rC 1)$, the exceptional set $E$ contains exactly $n+1$ isolated fixed points $p_0, \dots, p_n$, 
all simple, and no fixed $(-2)$-curve, and formulas \eqref{eqn:Ai} hold true.     
\end{lemma} 
{\it Proof}. 
A careful inspection of formulas \eqref{eqn:eudvd} leads to the assertion about isolated 
fixed points and fixed $(-2)$-curves.     
Then by definitions \eqref{eqn:mufE} and \eqref{eqn:nuhat1} we obtain $\mu(f, E) = n+1$ and  
\begin{align*}
\hat{\nu}(f, E) 
&= \sum_{j=0}^n \dfrac{\beta}{(1-\beta^{2 j+1-n} \alpha^{n+1})(1-\beta^{n+1-2 j} \alpha^{-n-1}) } 
= \dfrac{ \beta (1 + \beta^2 + \cdots + \beta^{2 n}) }{(\beta^{n+1} - \alpha^{n+1}) 
(\beta^{n+1}-\alpha^{-n-1})} 
\\
&= \dfrac{ \beta^{n+1} - \beta^{-n-1} }{(\beta-\beta^{-1}) 
\{ \beta^{n+1} + \beta^{-n-1} -(\alpha^{n+1} + \alpha^{-n-1}) \}} = 
\dfrac{ K_n(\hat{\tau} )}{ H_{n+1}(\hat{\tau}) - \{ \alpha^{(n+1)} + \alpha^{-(n+1)} \} },  
\end{align*}
where $\beta := \delta^{1/2}$ and $\hat{\tau} := \delta^{1/2} + \delta^{-1/2} = \beta + \beta^{-1}$ 
are used.   
This proves formulas \eqref{eqn:Ai}. 
\hfill $\Box$ \par\medskip
Up to the symmetry mentioned in Remark \ref{rem:symm}, case $(\rC 2)$ can be  
divided into two subcases: 
\begin{center}
$(\rC 2 \ra)$ $\gamma_i = 1$ for a unique index $i \in \{1, \dots, n-1 \} \cup \{1 \}$, $b_1 = 0$ if $i = 1$;  \quad    
$(\rC2 \rb)$ $\gamma_1 = 1$ and $b_1 \neq 0$.  
\end{center}
\begin{lemma} \label{lem:ecA2} 
In subcase $(\rC 2\ra)$ the exceptional set $E$ contains just one 
fixed $(-2)$-curve $E_i$ and exactly $n-1$ 
isolated fixed points $p_j$, $j \in \{0, \dots, n \} \setminus \{ i-1, i \}$, all simple, 
and formulas \eqref{eqn:Aii} hold true.            
In subcase $(\rC 2 \rb)$, $E$ contains exactly $n$ isolated fixed points $p_1, \dots, p_n$  
and no fixed $(-2)$-curve, where $p_1$ is double, while $p_2, \dots, p_n$ are simple, and 
formulas \eqref{eqn:Aii} hold true with $i = 1$. 
\end{lemma} 
{\it Proof}. 
In subcase (C2a), a careful inspection of formulas \eqref{eqn:eudvd} leads to 
the assertion about isolated fixed points and fixed curves.   
Since $\gamma_i = 1$, from definition \eqref{eqn:gamma} we have $\alpha^{n+1} = \beta^{n+1 -2 i}$  
and hence  
$$
\gamma_{j+1} = \beta^{2 j+1-n} \cdot \beta^{n+1-2 i} = \beta^{2(j+1-i)} =\delta^{j+1-i},  \quad  
\gamma_j^{-1} = \beta^{n+1-2 j} \cdot \beta^{2 i-n-1} = \beta^{2(i-j)} = \delta^{i-j}. 
$$
Thus by definitions \eqref{eqn:mufE} and \eqref{eqn:nuhat1} we obtain $\mu(f, E) = n+1$ 
and the second formula in \eqref{eqn:Aii} as follows;      
\begin{align*}
\hat{\nu}(f, E) &= \sum_{j=0}^{i-2} \dfrac{\delta^{1/2}}{(1-\delta^{j+1-i})(1-\delta^{i-j})} 
+ \sum_{j=i+1}^n \dfrac{\delta^{1/2}}{(1-\delta^{i-j})(1-\delta^{j+1-i})} + \dfrac{\delta^{1/2}(\delta+1)}{(\delta-1)^2} 
\\[1mm]
&= - \dfrac{\delta^{3/2} (\delta^{i-1}-1)}{(\delta-1)^2(\delta^i-1)} 
- \dfrac{\delta^{3/2} (\delta^{n-i}-1)}{(1-\delta)^2(\delta^{n+1-i}-1)} + \dfrac{\delta^{1/2}(\delta+1)}{(\delta-1)^2} 
= \dfrac{\delta^{1/2} (\delta^{n+1} -1)}{(\delta-1) (\delta^i-1) (\delta^{n+1-i}-1)} 
\\[1mm]
&= \dfrac{\delta^{(n+1)/2} -\delta^{-(n+1)/2} }{(\delta^{1/2}-\delta^{-1/2}) 
(\delta^{i/2}-\delta^{-i/2}) (\delta^{(n+1-i)/2}-\delta^{-(n+1-i)/2})} 
= \dfrac{K_n(\hat{\tau}) }{(H_2(\hat{\tau}) -2) \, K_{i-1}(\hat{\tau}) \, K_{n-i}(\hat{\tau})}.    
\end{align*}
\par
Also in subcase (C2b), a careful inspection of formulas \eqref{eqn:eudvd} leads to 
the assertion about isolated fixed points and fixed $(-2)$-curves.  
In particular, the M\"{o}bius transformation $f|_{E_1}$ has only one fixed point $p_1$ on $E_1 \cong \bP^1$, 
that is, $E_1$ is a {\sl parabolic} curve.       
Recall from \cite[Theorem 6.3]{IT2} that if an isolated fixed point $p$ lies on a parabolic curve, 
then $\mu(f, p) = 2$ and $\nu(f, p) = (\delta + 1)/(\delta-1)^2$. 
We apply this fact to the points $p_1$.  
Since $\gamma_1 = 1$, from definition \eqref{eqn:gamma} we have $\alpha^{n+1} = \beta^{n-1}$ 
and hence  
$\gamma_{j+1} = \beta^{2j} = \delta^j$ and $\gamma_j^{-1} = \beta^{2(1-j)} = \delta^{1-j}$. 
Therefore, in view of definitions \eqref{eqn:mufE} and \eqref{eqn:nuhat1}, we obtain  
$\mu(f, E) = 2 + (n-1) = n+1$ and   
$$
\hat{\nu}(f, E) = \delta^{1/2} \dfrac{\delta+1}{(\delta-1)^2} 
+ \sum_{j=2}^n \dfrac{\delta^{1/2}}{(1-\delta^{1-j})(1-\delta^j)} 
= \dfrac{K_n(\hat{\tau}) }{(H_2(\hat{\tau}) -2) \, K_{n-1}(\hat{\tau})},   
$$ 
which agrees with the second formula in \eqref{eqn:Aii} with $i = 1$.  
Thus Lemma \ref{lem:ecA2} is established. 
\hfill $\Box$ 
\begin{lemma} \label{lem:ecA3}   
In case $(\rC 3)$, $E$ contains exactly $n+1$ isolated fixed points $p_0, \dots, p_n$ 
and no fixed $(-2)$-curve, where $p_0$ is multiple while $p_1, \dots, p_n$ are simple, 
so we have $\mu(f, \rA_n^{\rt}) \ge n+2$.    
\end{lemma}
{\it Proof}. 
A careful inspection of formulas \eqref{eqn:eudvd} shows the assertions.  
In particular, $p_1, \dots, p_n$ are simple, while $p_0$ is multiple, as $f$ has 
multiplier $\gamma_0^{-1} = 1$ at $p_0$ in a normal direction to $E_1$.  
So we have $\mu(f, \rA_n^{\rt}) \ge 2 + n$. \hfill $\Box$  
\begin{remark} \label{rem:A} 
By Theorem \ref{thm:ecA} we always have $\mu(f, \rA_n^{\rt}) \ge n+1$, 
where the minimum $\mu(f, \rA_n^{\rt}) = n+1$ is attained exactly when $(\rC 3)$ is not the case.  
In case $(\rC 3)$ the point $p_0$ or $p_{n+1}$ is an {\sl exceptional fixed point of type} $\rII$, a concept discussed   
in \cite[\S 6]{IT2}, for which the holomorphic local index $\hat{\nu}(f, p)$ is harder to calculate.           
\end{remark}     
\subsection{Other Exceptional Components} \label{ss:oec}
We proceed to an $f$-invariant exceptional component of type $\rA_n$ with reflectional symmetry. 
Following Notation \ref{not:iec}, we denote it by $\rA_n^{\rr}$ collectively.  
Moreover we denote it by $\rA_n^{\rr, m}$, when we consider the $m$-th iterate $f^m$ in place of $f$.  
Note that $\rA_n^{\rr, m}$ has reflectional symmetry if $m$ is odd, and trivial symmetry 
if $m$ is even. 
\begin{proposition} \label{prop:Ar} 
For any $f$-invariant exceptional component of type $\rA_n$, $n \ge 2$, 
with reflectional symmetry,     
\begin{subequations} \label{eqn:Armn} 
\begin{align}
\mu(f^m, \rA_n^{\rr, m}) 
&= \begin{cases} 
1 & \mbox{if $m$ is odd and $n$ is even}, \\[1mm]
2 & \mbox{if $m$ is odd and $n$ is odd}, \\[1mm]
n+1 & \mbox{if $m$ is even},    
\end{cases}
\label{eqn:mnAr} 
\\[2mm]
\hat{\nu}(f^m, \rA_n^{\rr, m}) 
&= \begin{cases} 
\dfrac{1}{H_m(\hat{\tau})} &  \mbox{if $m$ is odd},  \\[3mm]   
\dfrac{K_{n}(H_m(\hat{\tau}))}{H_{m (n+1)}(\hat{\tau}) - 2(-1)^{(m/2)(n+1)}} & \mbox{if $m$ is even}. 
\end{cases}
\label{eqn:nArm} 
\end{align}
\end{subequations}
\end{proposition}
{\it Proof}.  
We use Corollary \ref{cor:ida2}, 
so we write the component $\rA_n^{\rr}$ also as $E$ in accordance with the notation there.   
\par
First we consider the case $m = 1$. 
From definition \eqref{eqn:hH} we have $f = f_J \circ h$, where $h|_E$ is described 
in Proposition \ref{prop:cmid}.    
The map $f$ acts on $E$ essentially in the same manner as $f_J$. 
Namely, although $f$ and $f_J$ may be slightly different on $E_1$ and $E_n$, this 
does not affect the current issue and Corollary \ref{cor:ida2} is available for $f$ 
with $\beta = \delta^{1/2}$.  
When $n = 2 k$ is even, as in case (1) of Corollary \ref{cor:ida2}, the point $p_k$ 
at the intersection $E_k \cap E_{k+1}$ is the only fixed point of $f$ on $E$, which 
is simple and has multipliers $\pm \ri \delta^{1/2}$ (see Figure \ref{fig:excep1}).  
Thus we have  
$$
\mu(f, \rA_n^{\rr}) = 1, \qquad 
\hat{\nu}(f, \rA_n^{\rr}) = \hat{\nu}(f, p_k) = 
\frac{\delta^{1/2}}{(1 - \ri \, \delta^{1/2})(1 + \ri \, \delta^{1/2})} = 
\dfrac{\delta^{1/2}}{ 1 + \delta} = \dfrac{1}{\hat{\tau}},    
$$ 
which verifies the case where $m=1$ and $n$ is even in formulas \eqref{eqn:mnAr} and \eqref{eqn:nArm}.    
When $n = 2 k-1$ is odd, as in case (2) of Corollary \ref{cor:ida2}, the map $f$ has 
exactly two fixed points $q_{\pm}$ on $E$, which lie on $E_k$ (see Figure \ref{fig:excep2}). 
Moreover each of them is simple and has multipliers $-1$ and $-\delta$. 
Thus we have $\mu(f, \rA_n^{\rr}) = 2$ and 
$$
\hat{\nu}(f, \rA_n^{\rr}) = \hat{\nu}(f, q_+) + \hat{\nu}(f, q_-) = 
\frac{\delta^{1/2}}{(1+1)(1 + \delta)} +  \frac{\delta^{1/2}}{(1+1)(1 + \delta)} 
= \dfrac{\delta^{1/2}}{ 1 + \delta} = \dfrac{1}{\hat{\tau}},   
$$ 
which shows the case where $m = 1$ and $n$ is odd in \eqref{eqn:Armn}. 
In either case $\hat{\nu}(f, \rA_n^{\rr})$ takes the same value $1/ \hat{\tau}$. 
\par
Next we proceed to the case where $m$ is an odd number in general.   
In this case $E$ is an $f^m$-invariant exceptional component of type $\rA_n$ with reflectional summetry. 
So the argument in the previous paragraph can be applied to $f^m$ in place of $f$, for which  
$\hat{\tau}$ should be replaced by $H_m(\hat{\tau})$. 
This proves \eqref{eqn:Armn} in this case. 
\par
Finally we consider the case where $m$ is even. 
In this case $E$ is an $f^m$-invariant exceptional component of type $\rA_n$ with trivial symmetry, 
where $\alpha = \ri^m = (-1)^{m/2}$, $\beta = \delta^{m/2}$ and so  
$\gamma_j = \delta^{m(2j-n-1)/2} (-1)^{m(n+1)/2}$ in formula \eqref{eqn:gamma}. 
So, when $n$ is even, we are in case (C1) of Theorem \ref{thm:ecA}. 
Formulas \eqref{eqn:Ai} with $f$ replaced by $f^m$ yield formulas \eqref{eqn:mnAr} and \eqref{eqn:nArm} in this case. 
When $n$ is odd, we are in case (C2) of Theorem \ref{thm:ecA} with $i = (n+1)/2$.  
Note that $2 \le i \le n-1$ as $n \ge 3$.  
Formulas \eqref{eqn:Aii} with $f$ replaced by $f^m$ yield 
$\mu(f^m, \rA_n^{\rr, m}) = n + 1$ and  
$$
\hat{\nu}(f^m, \rA_n^{\rr, m}) = \dfrac{K_n(H_m(\hat{\tau}))}{ \{ H_2(H_m(\hat{\tau})) -2 \} \{ K_{(n-1)/2}(H_m(\hat{\tau})) \}^2}. 
$$
Using the identity $\{ H_2(w) - 2\} \{ K_n(w)\}^2 = H_{2(n+1)}(w) -2$, we can observe that the last equation 
is just the subcase of formula \eqref{eqn:nArm} where $m$ is even and $n$ is odd.    
 \hfill $\Box$
\begin{lemma} \label{lem:DE} 
For any $f$-invariant exceptional component of Dynkin type $\rD$ we have   
\begin{alignat*}{3}
\mu(f, \rD_n^{\rt}) &= n+1, \qquad & 
\hat{\nu}(f, \rD_n^{\rt}) &= \dfrac{1}{\tau -2} 
\left( \dfrac{\tau}{\hat{\tau}} - \dfrac{K_{n-4}(\hat{\tau})}{K_{n-3}(\hat{\tau})}  \right),  \qquad & n &\ge 4, \\[1mm] 
\mu(f, \rD_n^{\rr}) &= n-1, \qquad &
\hat{\nu}(f, \rD_n^{\rr}) &= \dfrac{1}{2} \left( \dfrac{1}{\hat{\tau}} + \dfrac{K_{n-3}(\hat{\tau})}{H_{n-2}(\hat{\tau})} \right), \qquad & n &\ge 4, \\[1mm]   
\mu(f, \rD_4^{\rc}) &= 2, \qquad & 
\hat{\nu}(f, \rD_4^{\rc}) &= \dfrac{\hat{\tau}}{\tau + 1}. & n &=4,    
\end{alignat*}
where $\tau = H_2(\hat{\tau}) = \delta + \delta^{-1}$. 
For any $f$-invariant exceptional component of Dynkin type $\rE$ we have
\begin{alignat*}{3}
\mu(f, \rE_n^{\rt}) &= n+1, \qquad & 
\hat{\nu}(f, \rE_n^{\rt}) &= \dfrac{1}{\tau -2} 
\left( \dfrac{\tau^2 + \tau -1}{\hat{\tau} (\tau +1)} - \dfrac{K_{n-5}(\hat{\tau})}{ K_{n-4}(\hat{\tau})}  \right),  \qquad 
& n &= 6, 7, 8, \\[1mm]
\mu(f, \rE_6^{\rr}) &= 3. \qquad & \hat{\nu}(f, \rE_6^{\rr}) &= \dfrac{\tau+1}{\hat{\tau} \cdot \tau}, \qquad & n &= 6.   
\end{alignat*}
\end{lemma} 
{\it Proof}. 
Except for the case of type $\rD_4^{\rc} $, this lemma is a reformulation of \cite[Lemma 5.1]{IT2} 
in terms of the polynomials $H_n(w)$ and $K_n(w)$. 
Formulas for $\mu(f, \rD_4^{\rc})$ and $\hat{\nu}(f, \rD_4^{\rc})$ are obtained by direct calculations. 
\hfill $\Box$ 
\begin{remark} \label{rem:DE} 
Under iterations of $f$ the symmetry types of an $f$-invariant exceptional component vary like 
$$
\rt \to \rt \to \rt \to \cdots, \qquad \rr \to \rt \to \rr \to \rt \to \cdots, \qquad 
\rc \to \rc \to \rt \to \rc \to \rc \to \rt \to \cdots.  
$$
Formulas in Lemma \ref{lem:DE} for $f^m$ in place of $f$ can be obtained by just following these patterns 
and replacing $\hat{\tau}$ by $H_m(\hat{\tau})$. 
In the case of type $\rA$ with reflectional symmetry, the formulas for the iterates $f^m$  
are given in Proposition \ref{prop:Ar}. 
The case of type $\rA$ with trivial symmetry is most intricate, but even in this case 
the iteration procedure goes well with Theorem \ref{thm:ecA} and Remark \ref{rem:ecA} 
in a generic situation, once case $(\rC 3)$ is ruled out there.     
\end{remark}
\section{Method of Hypergeometric Groups} \label{sec:mhgg}
Our previous papers \cite{IT1,IT2} develop a method for constructing non-projective K3 
surface automorphisms of positive entropy, which we call the {\sl method of hypergometric groups}.  
Leaving full details to the original paper \cite{IT1} and its amplification \cite{IT2}, 
we briefly survey this method.  
It starts with a pair $(\varphi, \psi)$ of monic polynomials of degree $22$ over $\bZ$ 
such that $\varphi(z)$ is anti-palindromic, $z^{11} \varphi(z^{-1}) = - \varphi(z)$, 
$\psi(z)$ is palindromic, $z^{11} \psi(z^{-1}) = \psi(z)$, and the resultant 
of $\varphi(z)$ and $\psi(z)$ is $\pm 1$. 
Let $A$ and $B$ be the companion matrices of $\varphi(z)$ and $\psi(z)$ respectively,  
so that $\varphi(z)$ and $\psi(z)$ are the characteristic polynomials of $A$ and $B$.  
Let $H = \langle A, B \rangle$ be the subgroup of $\GL_{22}(\bZ)$ generated by $A$ and $B$. 
It is a group modeled on the monodromy group of a generalized hypergeometric differential equation of order $22$.  
We find that $C := A^{-1} B \in H$ is a reflection which fixes a hyperplane in $\bQ^{22}$ pointwise 
and sends a nonzero vector $\br \in \bQ^{22}$ to its negative $- \br$. 
We have a free $\bZ$-module of rank $22$, 
$$
L = \langle \br, A \br, \dots, A^{21} \br \rangle_{\bZ} = \langle \br, B \br, \dots, B^{21} \br \rangle_{\bZ},  
$$ 
stable under the action of $H$. 
We can make it into an $H$-invariant even unimodular lattice by providing it with 
the inner product $(A^{i-1} \br, \, A^{j-1} \br) = \xi_{|i-j|}$, where 
$\xi_0 = 2$ and the sequence $\{ \xi_i \}_{i \in \bN}$ is defined by the expansion  
$$
\dfrac{\psi(z)}{\varphi(z)} = 1 + \sum_{i=1}^{\infty} \xi_i \, z^{-i} \qquad 
\mbox{around} \quad z = \infty. 
$$
Turning $\varphi(z)$ and $\psi(z)$ upside down in the expansion gives the Gram matrix 
$(B^{i-1} \br, \, B^{j-1} \br)$ for the $B$-basis. 
\par
Take $F = A$ or $B$ and let $\chi(z)$ be the characterisctic polynomial of $F$, that is, 
$\chi(z) := \varphi(z)$ for $F = A$ and $\chi(z) := \psi(z)$ for $F = B$.  
A suitable choice of the pair $(\varphi, \psi)$ makes the hypergeometric lattice 
$L$ (or its negative) into an even unimodular lattice of signature $(3, 19)$, 
that is, a K3 lattice, equipped with the Hodge structure 
\begin{equation} \label{eqn:HS}
L_{\bC} = H^{1,0} \oplus H^{1,1} \oplus H^{0,2} = \ell \oplus (\ell \oplus \overline{\ell})^{\perp} \oplus \overline{\ell} 
\quad \mbox{of signature $(1, 0) \oplus (1, 19) \oplus (1, 0)$},  
\end{equation}
such that $F$ is a positive Hodge isometry, where $\ell$ is the eigenline of $F$ for the eigenvalue 
$\delta \in S^1 \setminus\{\pm 1\}$ and $\overline{\ell}$ is the complex conjugate to $\ell$; 
see \cite[Theorems 1.2 and 1.3]{IT1}.    
The number $\delta$, which is a root of $\chi(z)$, and its trace $\tau := \delta + \delta^{-1}$ are referred to 
as the {\sl special eigenvalue} and {\sl special trace} (ST for short) for the Hodge structure \eqref{eqn:HS}.   
We remark that $\delta$ is either a root of unity or a non-real conjugate to a Salem number $\lambda$.  
Let us consider the Salem case.  
In this case we have a factorization $\chi(z) = S(z) \cdot \chi_1(z)$, where $S(z)$ is the minimal polynomial 
of $\lambda$ and $\chi_1(z)$ is a product of cyclotomic polynomials or $1$.   
The Picard lattice $\Pic := L \cap H^{1,1}$ is given by  
$$
\Pic = \{ \bv \in L \mid \chi_1(F) \bv = 0 \} = \langle \bs, F \bs, \dots, F^{\rho-1} \bs \rangle_{\bZ},  
\qquad \bs := S(F) \br,  
$$
with Picard number $\rho := \deg \chi_1(z) = 22 - \deg S(z)$. 
The lexicographic order on $\Pic$ relative to the linear basis $\bs, \dots, F^{\rho-1} \bs$ 
gives rise to a root basis $\vD_{\rb}$ of the root system $\vD := \{ \bu \in \Pic \mid (\bu, \bu) = -2 \}$. 
This together with the choice of a positive cone $\cC^+ \subset H^{1,1}_{\bR}$ determines a Weyl chamber 
$\cK := \{ \bv \in \cC^+ \mid \mbox{$(\bv, \bu) > 0$ for any $\bu \in \vD_{\rb}$}\}$,  
which we designate as the ``K\"{a}hler cone".  
In general, the Hodge isometry $F$ may not preserve $\cK$, but there exists a unique Weyl group element 
$w_F$ such that the modified Hodge isometry $\widetilde{F} := w_F \circ F$ preserves $\cK$. 
Let $\widetilde{\chi}(z)$ be the characteristic polynomial of $\widetilde{F}$. 
It has the same Salem factor $S(z)$ as $\chi(z)$, decomposing into  
\begin{equation} \label{eqn:chi}
\widetilde{\chi}(z) = S(z) \cdot \widetilde{\chi}_1(z), 
\end{equation} 
where $\widetilde{\chi}_1(z)$ is a product of cyclotomic polynomials or $1$, 
which may differ from $\chi_1(z)$. 
Notice that $\widetilde{\chi}_1(z)$ is the characteristic polynomial of $\widetilde{F}|_{\Pic}$. 
What is important here is that, once the initial pair $(\varphi, \psi)$ is given concretely,  
all these calculations can be carried out explicitly by computer.  
In particular one can determine the Dynkin type of the root system $\vD$ and detect how 
$\widetilde{F}$ permutes the simple roots in $\vD_{\rb}$.  
\par
Through the period mapping the modified Hodge isometry $\widetilde{F} : L \to L$ lifts to a non-projective 
K3 surface automorphism $f : X \to X$ with entropy $h(f) = \log \lambda$, twist number $\delta(f) = \delta$,  
Picard lattice $\Pic(X) \cong \Pic$, Picard number $\rho(X) = \rho$ and K\"{a}hler cone 
$\cK(X) \cong \cK$; see \cite[Theorem 1.6]{IT1}.  
The constructed automorphism $f : X \to X$ satisfies conditions (A1) and (A2) in \S \ref{sec:intro}. 
There exists a one-to-one correspondence between the $(-2)$-curves in $\cE(X)$ and the simple roots 
in $\vD_{\rb}$. 
The Dynkin type of the exceptional set $\cE(X)$ is that of the root system $\vD$. 
How the automorphism $f$ permutes the $(-2)$-curves in $\cE(X)$ is faithfully 
represented by the action of $\widetilde{F}$ on the simple roots $\vD_{\rb}$. 
So the Dynkin automorphism $\sigma(f, E)$ defined in \eqref{eqn:sigma} can be determined 
explicitly for each $f$-invariant exceptional component $E \subset \cE(X)$. 
\par
Among other possibilities, the article \cite{IT2} employs two setups for the pair $(\varphi, \psi)$. 
In \cite[Setup 3.2]{IT2} we put      
\begin{equation} \label{eqn:setup3.2}
\varphi(z) = S(z) \cdot \varphi_1(z) \quad \mbox{with} \quad \varphi_1(z) = (z-1)(z+1) \cdot C(z); 
\qquad \psi(z) = s(z) \cdot c(z),     
\end{equation}  
where $S(z)$ is a Salem polynomial and $C(z)$ is a product of cyclotomic polynomials or $1$, 
while $s(z)$ is an unramified Salem polynomial and $c(z)$ is a product of unramified 
cyclotomic polynomials or $1$. 
This setup works efficiently for Picard numbers $\rho = 2,4,6, \dots, 16$, that is, for 
$d := \deg S(z) = 20, 18, 16, \dots, 6$, but not for $\rho = 18$ i.e. $d = 4$. 
To cope with the last case, we introduce \cite[Setup 3.4]{IT2}, where $\varphi(z)$ is 
the same as in \eqref{eqn:setup3.2}, but $S(z)$ is specified as $\rS^{(4)}_1(z) = 
z^4-z^3-z^2-z+1$, the minimal polynomial of the smallest Salem number 
$\lambda^{(4)}_1 \approx 1.7220838$ of degree $4$ (see Notation \ref{not:ST}), 
while $\psi(z)$ is an unramified palindromic polynomial of the form 
\begin{equation} \label{eqn:psi}
\psi(z) = z^{22} + c_1 z^{21} + \cdots + c_{10} z^{12} + c_{11} z^{11} + c_{10} z^{10} + 
\cdots + c_1 z + 1 \in \bZ[z],  
\end{equation}
such that $c_j \in \{ 0, \pm 1, \pm 2 \}$ for $j = 1, \dots, 9$, 
the trace polynomial of $\psi(z)$ has ten or eight roots on $(-2, \, 2)$, and the 
resultant of $\psi(z)$ with $\rS^{(4)}_1(z)$ is $\pm 1$; there are $1019$ polynomials of them.  
In both setups the Hodge isometry $F$ is the matrix $A$, so the polynomials $\widetilde{\chi}(z)$ 
and $\widetilde{\chi}_1(z)$ in \eqref{eqn:chi} are expressed as $\widetilde{\varphi}(z)$ and 
$\widetilde{\varphi}_1(z)$.  
These setups give reise to a lot of K3 surface automorphisms of positive entropy, 
which are listed in our database \cite{IT3}. 
\par
We shall examine which of the automorphisms in \cite{IT3} have rotation domains of rank $1$ or of rank $2$. 
\begin{notation} \label{not:ST}
Let $\lambda^{(d)}_i$ be the $i$-th smallest Salem number of degree $d \ge 4$, 
where $d$ is even.  
We denote by $\rS^{(d)}_i(z)$ the minimal polynomial of $\lambda^{(d)}_i$ and by 
$\rST^{(d)}_i(w)$ the trace polynomial of $\rS^{(d)}_i(z)$. 
The roots $\tau_j$ of $\rST^{(d)}_i(w)$ are arranged in decreasing order as 
$\tau_0 > 2 > \tau_1 > \tau_2 > \dots > \tau_{d/2-1} > -2$.   
To be more precise, $\tau_j$ should be written $\tau_{i,j}^{(d)}$, 
but the dependence upon $(d, i)$ is omitted for simplicity of notation.    
The $i$-th cyclotomic polynomial is denoted by $\rC_i(z)$. 
The unusual convention \eqref{eqn:conv} for $\rC_1(z)$ and $\rC_2(z)$ is not employed here, so  
$\rC_1(z) = z-1$, $\rC_2(z) = z+1$, $\rC_3(z) = z^2 + z+ 1$, $\rC_4(z) = z^2+1$, and so forth.    
\end{notation}
\begin{remark} \label{rem:division}
Recall that $\widetilde{\varphi}_1(z)$ is the characteristic polynomial of $\widetilde{A}|_{\Pic}$. 
Let $\widetilde{\varphi}_2(z)$ denote the characteristic polynomial of $\widetilde{A}|_{\Span \vD_{\rb}}$. 
Notice that $\widetilde{\varphi}_1(z)$ is divisible by $\widetilde{\varphi}_2(z)$, since $\Span \vD_{\rb} \subset \Pic$. 
If $\vD_{\rb}$ is partitioned into $\widetilde{A}$-invariant subsets $\vD_{\rb}^{(1)}, \vD_{\rb}^{(2)}, \dots$, 
then $\widetilde{\varphi}_2(z) = \widetilde{\varphi}_2^{(1)}(z) \, \widetilde{\varphi}_2^{(2)}(z) \cdots$, where 
$\widetilde{\varphi}_2^{(i)}(z)$ is the characteristic polynomial of $\widetilde{A}|_{\Span \vD_{\rb}^{(i)}}$. 
Moreover, if $\widetilde{A}$ acts on $\vD_{\rb}^{(i)}$ as a permutation of cycle type $(m, n, \dots)$, then 
$$
\widetilde{\varphi}_2^{(i)}(z) = (z^m-1)(z^n-1) \cdots = \prod_{i|m} \rC_i(z) \prod_{j|n} \rC_j(z) \cdots. 
$$
In some typical cases, the factor $\widetilde{\varphi}_2^{(i)}(z)$ is calculated in the following manner.     
\begin{itemize}
\setlength{\itemsep}{-1pt} 
\item If all elements of $\vD_{\rb}^{(i)}$ are fixed, then $\widetilde{\varphi}_2^{(i)} = \rC_1^k$ where $k$ is the cardinality of the set $\vD_{\rb}^{(i)}$.  
\item If $\vD_{\rb}^{(i)}$ is a root basis of type $\rA_n$, $n \ge 2$, with reflectional symmetry, then 
$\widetilde{\varphi}_2^{(i)} = \rC_1^{\lceil n/2 \rceil} \rC_2^{\lfloor n/2 \rfloor}$.   
\item If $\vD_{\rb}^{(i)}$ is a root basis of type $\rD_n$, $n \ge 4$, with reflectional symmetry, then $\widetilde{\varphi}_2^{(i)} = \rC_1^{n-1} \rC_2$. 
\item If $\vD_{\rb}^{(i)}$ is a root basis of type $\rD_4$ with tricyclic symmetry, then $\widetilde{\varphi}_2^{(i)} = \rC_1^2 \rC_3$. 
\item If $\vD_{\rb}^{(i)}$ is a root basis of type $\rE_6$ with reflectional symmetry, then $\widetilde{\varphi}_2^{(i)} = \rC_1^4 \rC_2^2$.  
\end{itemize} 
For example, if $\vD_{\rb}^{(i)}$ is a root basis of type $\rE_7$ or $\rE_8$, then we have $\widetilde{\varphi}_2^{(i)} = \rC_1^7$ or $\rC_1^8$ 
as one of the top cases. 
\end{remark}     
\subsection{Components of Type A with Trivial Symmetry} \label{ss:Ass} 
A lot of K3 surface automorphisms of positive entropy created by the method of hypergeometric groups 
are presented in our database \cite{IT3}. 
It contains many automorphisms $f : X \to X$ having an $f$-invariant 
exceptional component of type $\rA$ with trivial symmetry.   
Among them, Table \ref{tab:Ass} focuses on those automorphisms 
which have {\sl no fixed points off} the exceptional set $\cE(X)$. 
Our previous article \cite{IT2} employs two setups in applying the method of hypergeometric groups.     
Table \ref{tab:Ass} arises from one of them, that is, from \cite[Setup 3.2]{IT2}. 
For the format of Table \ref{tab:Ass} we refer to  its legend, except for one thing 
in the last column; by RH we mean ``Rotation domain or Hyperbolic set", but 
we postpone the details of its meaning until \S \ref{sec:rd} (see Theorem \ref{thm:RHAt}).      
\begin{longtable}{llllllllllr}
\caption{\setlength{\leftskip}{-16mm}\setlength{\rightskip}{-16mm} 
Having a unique component of type $\rA$ with trivial symmetry; \cite[Setup 3.2]{IT2}. } \label{tab:Ass} 
\\
\hline
      &           &           &           &          &          &      &           &                               &                       &     \\[-3mm]
No. & $\rho$ & $S(z)$ & $C(z)$ & $s(z)$ & $c(z)$ & ST & Dynkin & $\widetilde{\varphi}_1(z)$ & $L(f)$ & RH \\ 
\hline
\endfirsthead
\multicolumn{11}{l}{continued} \\
\hline
      &           &           &           &          &          &      &           &                               &                       &     \\[-3mm]
No. & $\rho$ & $S(z)$ & $C(z)$ & $s(z)$ & $c(z)$ & ST & Dynkin & $\widetilde{\varphi}_1(z)$ & $L(f)$ & RH \\ 
\hline
\endhead
\hline
\endfoot
\hline
\\[-4pt]
\caption*{\setlength{\leftskip}{-25mm}\setlength{\rightskip}{-25mm}
As in formula \eqref{eqn:setup3.2} the pair $(\varphi, \psi)$ is given by $\varphi(z) = (z-1)(z+1) \, S(z) \, C(z)$ 
and $\psi(z) = s(z) \, c(z)$, where $S(z)$ and $s(z)$ are Salem polynomials, 
while $C(z)$ and $c(z)$ are products of cyclotomic polynomials; 
$\rho$ is the Picard number of $X$; ST is the special trace for the Hodge structure \eqref{eqn:HS};  
Dynkin column exhibits the Dynkin type of the root system $\vD$ or of the exceptional set $\cE(X)$;  
$\widetilde{\varphi}_1(z)$ is the cyclotomic component of the characteristic polynomial 
$\widetilde{\varphi}(z) = S(z) \, \widetilde{\varphi}_1(z)$ of $\widetilde{A} := w_A\circ A$ or of $f^*|_{H^2(X, \bC)}$;  
$L(f)$ is the Lefschetz number of $f$, which is equal to $2 + \Tr \widetilde{A}$; see \S \ref{sec:rd} for RH-column. }
\endlastfoot
$1$ & $6$ & $\rS^{(16)}_{2}$ & $\rC_5$ & $\rS^{(22)}_{18}$ & $1$ & $\tau_5$ & $\rA_6$ & $\rC_1^6$ & $7$ & $\rR 2$ \\
$2$ & $6$ & $\rS^{(16)}_{2}$ & $\rC_3\rC_4$ & $\rS^{(22)}_{18}$ & $1$ & $\tau_3$ & $\rA_6$ & $\rC_1^6$ & $7$ & $\rR 2$ \\
$3$ & $6$ & $\rS^{(16)}_{2}$ & $\rC_3\rC_6$ & $\rS^{(22)}_{1}$ & $1$ & $\tau_1$ & $\rA_6$ & $\rC_1^6$ & $7$ & $\rHy$ \\
$4$ & $6$ & $\rS^{(16)}_{2}$ & $\rC_3\rC_6$ & $\rS^{(22)}_{18}$ & $1$ & $\tau_2$ & $\rA_6$ & $\rC_1^6$ & $7$ & $\rR 2$ \\
$5$ & $6$ & $\rS^{(16)}_{2}$ & $\rC_3\rC_6$ & $\rS^{(22)}_{72}$ & $1$ & $\tau_4$ & $\rA_6$ & $\rC_1^6$ & $7$ & $\rR 2$ \\ 
\hline
$6$ & $10$ & $\rS^{(12)}_{1}$ & $\rC_{15}$ & $\rS^{(14)}_{18}$ & $\rC_{24}$ & $\tau_3$ & $\rA_1\oplus \rE_8$ & $\rC_1^9\rC_2$ & $11$ & $\rHy$ \\
$7$ & $10$ & $\rS^{(12)}_{1}$ & $\rC_{20}$ & $\rS^{(22)}_{2}$ & $1$ & $\tau_5$ & $\rA_1\oplus \rE_8$ & $\rC_1^9\rC_2$ & $11$ & $\rR 2$ \\
$8$ & $10$ & $\rS^{(12)}_{1}$ & $\rC_{20}$ & $\rS^{(22)}_{39}$ & $1$ & $\tau_5$ & $\rA_1\oplus \rE_8$ & $\rC_1^9\rC_2$ & $11$ & $\rR 2$ \\
$9$ & $10$ & $\rS^{(12)}_{1}$ & $\rC_{20}$ & $\rS^{(22)}_{71}$ & $1$ & $\tau_5$ & $\rA_1\oplus \rE_8$ & $\rC_1^9\rC_2$ & $11$ & $\rR 2$ \\
$10$ & $10$ & $\rS^{(12)}_{1}$ & $\rC_{24}$ & $\rS^{(6)}_{3}$ & $\rC_{60}$ & $\tau_2$ & $\rA_1\oplus \rE_8$ & $\rC_1^9\rC_2$ & $11$ & $\rR 2$ \\
$11$ & $10$ & $\rS^{(12)}_{1}$ & $\rC_{24}$ & $\rS^{(22)}_{64}$ & $1$ & $\tau_4$ & $\rA_1\oplus \rE_8$ & $\rC_1^9\rC_2$ & $11$ & $\rR 2$ \\
$12$ & $10$ & $\rS^{(12)}_{1}$ & $\rC_{30}$ & $\rS^{(6)}_{1}$ & $\rC_{40}$ & $\tau_3$ & $\rA_1$ & $\rC_1\rC_2\rC_{30}$ & $2$ & $\rR 2$ \\
$13$ & $10$ & $\rS^{(12)}_{1}$ & $\rC_{30}$ & $\rS^{(10)}_{2}$ & $\rC_{42}$ & $\tau_5$ & $\rA_1$ & $\rC_1\rC_2\rC_{30}$ & $2$ & $\rHy$ \\
$14$ & $10$ & $\rS^{(12)}_{1}$ & $\rC_{30}$ & $\rS^{(14)}_{114}$ & $\rC_{24}$ & $\tau_5$ & $\rA_1$ & $\rC_1\rC_2\rC_{30}$ & $2$ & $\rHy$ \\
$15$ & $10$ & $\rS^{(12)}_{1}$ & $\rC_{30}$ & $\rS^{(22)}_{14}$ & $1$ & $\tau_5$ & $\rA_1$ & $\rC_1\rC_2\rC_{30}$ & $2$ & $\rHy$ \\
$16$ & $10$ & $\rS^{(12)}_{1}$ & $\rC_{30}$ & $\rS^{(22)}_{39}$ & $1$ & $\tau_2$ & $\rA_1$ & $\rC_1\rC_2\rC_{30}$ & $2$ & $\rR 2$ \\
$17$ & $10$ & $\rS^{(12)}_{1}$ & $\rC_{30}$ & $\rS^{(22)}_{52}$ & $1$ & $\tau_3$ & $\rA_1$ & $\rC_1\rC_2\rC_{30}$ & $2$ & $\rR 2$ \\
$18$ & $10$ & $\rS^{(12)}_{1}$ & $\rC_{30}$ & $\rS^{(22)}_{89}$ & $1$ & $\tau_4$ & $\rA_1$ & $\rC_1\rC_2\rC_{30}$ & $2$ & $\rR 2$ \\
\hline
$19$ & $12$ & $\rS^{(10)}_{1}$ & $\rC_3\rC_{15}$ & $\rS^{(14)}_{1}$ & $\rC_{20}$ & $\tau_3$ & $\rA_2$ & $\rC_1^2\rC_3\rC_{15}$ & $3$ & $\rR 2$ \\
$20$ & $12$ & $\rS^{(10)}_{1}$ & $\rC_3\rC_{15}$ & $\rS^{(14)}_{3}$ & $\rC_{20}$ & $\tau_1$ & $\rA_2$ & $\rC_1^2\rC_3\rC_{15}$ & $3$ & $\rR 2$ \\
$21$ & $12$ & $\rS^{(10)}_{1}$ & $\rC_3\rC_{15}$ & $\rS^{(22)}_{43}$ & $1$ & $\tau_4$ & $\rA_2$ & $\rC_1^2\rC_3\rC_{15}$ & $3$ & $\rR 2$ \\
$22$ & $12$ & $\rS^{(10)}_{1}$ & $\rC_3\rC_{15}$ & $\rS^{(22)}_{71}$ & $1$ & $\tau_4$ & $\rA_2$ & $\rC_1^2\rC_3\rC_{15}$ & $3$ & $\rR 2$ \\
$23$ & $12$ & $\rS^{(10)}_{1}$ & $\rC_3\rC_{15}$ & $\rS^{(22)}_{72}$ & $1$ & $\tau_4$ & $\rA_2$ & $\rC_1^2\rC_3\rC_{15}$ & $3$ & $\rR 2$ \\
\hline
$24$ & $14$ & $\rS^{(8)}_{1}$ & $\rC_4\rC_7\rC_8$ & $\rS^{(22)}_{72}$ & $1$ & $\tau_2$ & $\rA_7$ & $\rC_1^7\rC_2\rC_4\rC_8$ & $8$ & $\rR 2$ \\ 
\hline
$25$ & $14$ & $\rS^{(8)}_{2}$ & $\rC_{42}$ & $\rS^{(6)}_{1}$ & $\rC_{60}$ & $\tau_2$ & $\rA_1$ & $\rC_1\rC_2\rC_{42}$ & $2$ & $\rR 2$ \\
$26$ & $14$ & $\rS^{(8)}_{2}$ & $\rC_{42}$ & $\rS^{(14)}_{105}$ & $\rC_{30}$ & $\tau_1$ & $\rA_1$ & $\rC_1\rC_2\rC_{42}$ & $2$ & $\rR 2$ \\
\hline
$27$ & $14$ & $\rS^{(8)}_{6}$ & $\rC_{14}\rC_{18}$ & $\rS^{(14)}_{59}$ & $\rC_{15}$ & $\tau_3$ & $\rA_1\oplus \rA_6\oplus \rE_6$ & $\rC_1^8\rC_2^6$ & $6$ & $\rR 2$ \\
\hline
$28$ & $14$ & $\rS^{(8)}_{8}$ & $\rC_{42}$ & $\rS^{(14)}_{306}$ & $\rC_{20}$ & $\tau_2$ & $\rA_2$ & $\rC_1^2\rC_{42}$ & $3$ & $\rR 1$ \\
$29$ & $14$ & $\rS^{(8)}_{8}$ & $\rC_4\rC_5\rC_7$ & $\rS^{(22)}_{18}$ & $1$ & $\tau_2$ & $\rA_{11}$ & $\rC_1^{11}\rC_2\rC_4$ & $12$ & $\rHy$ \\ 
\hline
$30$ & $14$ & $\rS^{(8)}_{10}$ & $\rC_3\rC_5\rC_7$ & $\rS^{(10)}_{1}$ & $\rC_{42}$ & $\tau_1$ & $\rA_6$ & $\rC_1^7\rC_2\rC_3\rC_5$ & $7$ & $\rHy$ \\
\hline
$31$ & $16$ & $\rS^{(6)}_{1}$ & $\rC_3\rC_5\rC_6\rC_7$ & $\rS^{(22)}_{18}$ & $1$ & $\tau_1$ & $\rA_6$ & $\rC_1^7\rC_2\rC_3\rC_5\rC_6$ & $7$ & $\rHy$ \\
\hline
\end{longtable}
Every exceptional component appearing in Table $\ref{tab:Ass}$ is $f$-invariant, because for each entry 
of the table the exceptional set $\cE(X)$ contains no more than one exceptional components of the same Dynkin type.  
We establish the main properties of the automorphisms in Table \ref{tab:Ass} in the following proposition. 
\begin{proposition}[Table $\mathbf{\ref{tab:Ass}}$] \label{prop:Ass} 
In No.~$27$ the $\rA_6$-component and the $\rE_6$-component have reflectional symmetry.   
Any other exceptional component of type $\rA$ has trivial symmetry, where case $(\rC 3)$ 
never occurs in Theorem $\ref{thm:ecA}$. 
Every automorphism $f$ in the table has no fixed point off the exceptional set $\cE(X)$, 
that is, $\mu_{\roff}(f) = 0$.    
\end{proposition}
{\it Proof}. 
We use Notation \ref{not:iec}.   
It is evident that any $\rA_1$-component has trivial symmetry, that is, $\rA_1 = \rA_1^{\rt}$.   
First we consider No.~$27$, which has exceptional set of type $\rA_1 \oplus \rA_6 \oplus \rE_6$. 
By Remark \ref{rem:A}, Proposition \ref{prop:Ar} and Lemma \ref{lem:DE} we have $\mu(f, \rA_1) \ge 2$, 
$\mu(f, \rA_6) \ge 1$ and $\mu(f, \rE_6) \ge 3$, so by definition \eqref{eqn:muof} we have  
$\mu_{\ron}(f) = \mu(f, \rA_1) + \mu(f, \rA_6) + \mu(f, \rE_6) \ge 2 + 1 + 3 = 6$.  
FPF \eqref{eqn:saito2} then implies $0 \le \mu_{\roff}(f) = L(f) - \mu_{\ron}(f) \le 6 -6 = 0$. 
This forces $\mu_{\roff}(f) = 0$ and hence $\mu(f, \rA_1) = 2$, $\mu(f, \rA_6) = 1$ and $\mu(f, \rE_6) = 3$. 
Remark \ref{rem:A}, Proposition \ref{prop:Ar} and Lemma \ref{lem:DE} again show that case $(\rC 3)$ in 
Theorem \ref{thm:ecA} does not occur for $\rA_1$ and we have $\rA_6 = \rA_6^{\rr}$ and $\rE_6 = \rE_6^{\rr}$.    
Next we pick out No.~$29$, which has exceptional set of type $\rA_{11}$.   
By Remark \ref{rem:division}, if $\rA_{11} = \rA_{11}^{\rt}$ then $\tilde{\varphi}_2 = \rC_1^{11}$, 
while if $\rA_{11} = \rA_{11}^{\rr}$ then $\tilde{\varphi}_2 = \rC_1^6 \rC_2^5$.  
Since $\tilde{\varphi}_2$ divides $\tilde{\varphi}_1 = \rC_1^{11} \rC_2 \rC_4$, we must have $\rA_{11} = \rA_{11}^{\rt}$. 
So, by definition \eqref{eqn:muof} and Remark \ref{rem:A}, we have $\mu_{\ron}(f) = \mu(f, \rA_{11}^{\rt}) \ge 12$. 
FPF \eqref{eqn:saito2} then implies $0 \le \mu_{\roff}(f) = L(f) - \mu_{\ron}(f) \le 12 - 12 = 0$, 
which forces $\mu_{\roff}(f) = 0$ and $\mu(f, \rA_{11}^{\rt}) = 12$. 
By Remark \ref{rem:A} case $(\rC 3)$ in Theorem \ref{thm:ecA} does not occur for $\rA_{11}^{\rt}$. 
The other entries can be treated in similar or simpler manners.        
\hfill $\Box$ \par\medskip
\subsection{Components of Type A with Reflectional Symmetry} \label{ss:Ars}
We pick out from database \cite{IT3} and present in Table \ref{tab:Ars1} examples of  
K3 surface automorphisms $f : X \to X$ of positive entropy having an $f$-invariant exceptional 
component of type $\rA$ with reflectional symmetry. 
They also arise from \cite[Setup 3.2]{IT2} and their main properties are established in 
Proposition \ref{prop:Ars1}. 
\begin{longtable}{lllllllllcc}
\caption{\setlength{\leftskip}{-16mm}\setlength{\rightskip}{-16mm} 
Having a component of type $\rA$ with reflectional symmetry; \cite[Setup 3.2]{IT2}.} \label{tab:Ars1}
\\
\hline
     &           &         &           &          &           &      &            &                                &                        &     \\[-3mm]
No. & $\rho$ & $S(z)$ & $C(z)$ & $s(z)$ & $c(z)$ & ST & Dynkin & $\widetilde{\varphi}_1(z)$ & $L(f)$ & $\mu_{\roff}(f)$ \\ 
\hline
\endfirsthead
\multicolumn{11}{l}{continued} \\
\hline
     &           &           &           &          &           &      &        &                                &                        &     \\[-3mm]
No. & $\rho$ & $S(z)$ & $C(z)$ & $s(z)$ & $c(z)$ & ST & Dynkin & $\widetilde{\varphi}_1(z)$ & $L(f)$ & $\mu_{\roff}(f)$ \\ 
\hline
\endhead
\hline
\endfoot
\hline
\\[-4pt] 
\caption*{\setlength{\leftskip}{-20mm}\setlength{\rightskip}{-20mm}
The format of this table is the same as that of Table \ref{tab:Ass} except for the last column,  
which gives the cardinality $\mu_{\roff}(f)$ counted with multiplicities of the fixed points of $f$ off the exceptional set $\cE(X)$; 
No.\! $28$ of this table is the same as No.\! $27$ of Table \ref{tab:Ass}. }
\endlastfoot
$1$ & $4$ & $\rS^{(18)}_{6}$ & $\rC_6$ & $\rS^{(14)}_{4}$ & $\rC_{15}$ & $\tau_8$ & $\rA_1 \oplus \rA_2$ & $\rC_1^2\rC_2^2$ & $4$ & $1$ \\
$2$ & $4$ & $\rS^{(18)}_{6}$ & $\rC_6$ & $\rS^{(22)}_{25}$ & $1$ & $\tau_6$ & $\rA_1 \oplus \rA_2$ & $\rC_1^2\rC_2^2$ & $4$ & $1$ \\ 
\hline
$3$ & $6$ & $\rS^{(16)}_{3}$ & $\rC_{10}$ & $\rS^{(18)}_{1}$ & $\rC_{12}$ & $\tau_1$ & $\rA_4$ & $\rC_1^3\rC_2^3$ & $2$ & $1$ \\
$4$ & $6$ & $\rS^{(16)}_{3}$ & $\rC_{10}$ & $\rS^{(22)}_{2}$ & $1$ & $\tau_1$ & $\rA_4$ & $\rC_1^3\rC_2^3$ & $2$ & $1$ \\
\hline
$5$ & $12$ & $\rS^{(10)}_{1}$ & $\rC_{10}\rC_{14}$ & $\rS^{(22)}_{10}$ & $1$ & $\tau_3$ & $\rA_6$ & $\rC_1^3\rC_2^3\rC_{14}$ & $2$ & $1$ \\
$6$ & $12$ & $\rS^{(10)}_{1}$ & $\rC_{10}\rC_{14}$ & $\rS^{(22)}_{22}$ & $1$ & $\tau_3$ & $\rA_6$ & $\rC_1^3\rC_2^3\rC_{14}$ & $2$ & $1$ \\
$7$ & $12$ & $\rS^{(10)}_{1}$ & $\rC_{10}\rC_{14}$ & $\rS^{(22)}_{57}$ & $1$ & $\tau_2$ & $\rA_6$ & $\rC_1^3\rC_2^3\rC_{14}$ & $2$ & $1$ \\
$8$ & $12$ & $\rS^{(10)}_{1}$ & $\rC_3\rC_6\rC_{14}$ & $\rS^{(22)}_{1}$ & $1$ & $\tau_3$ & $\rA_6$ & $\rC_1^3\rC_2^3\rC_{14}$ & $2$ & $1$ \\
$9$ & $12$ & $\rS^{(10)}_{1}$ & $\rC_3\rC_6\rC_{14}$ & $\rS^{(22)}_{6}$ & $1$ & $\tau_2$ & $\rA_6$ & $\rC_1^3\rC_2^3\rC_{14}$ & $2$ & $1$ \\
$10$ & $12$ & $\rS^{(10)}_{1}$ & $\rC_3\rC_6\rC_{14}$ & $\rS^{(22)}_{10}$ & $1$ & $\tau_3$ & $\rA_6$ & $\rC_1^3\rC_2^3\rC_{14}$ & $2$ & $1$ \\
\hline
$11$ & $14$ & $\rS^{(8)}_{1}$ & $\rC_{28}$ & $\rS^{(10)}_{24}$ & $\rC_{36}$ & $\tau_2$ & $\rA_2$ & $\rC_1\rC_2\rC_{28}$ & $2$ & $1$ \\
$12$ & $14$ & $\rS^{(8)}_{1}$ & $\rC_{28}$ & $\rS^{(14)}_{3}$ & $\rC_{30}$ & $\tau_1$ & $\rA_2$ & $\rC_1\rC_2\rC_{28}$ & $2$ & $1$ \\
$13$ & $14$ & $\rS^{(8)}_{1}$ & $\rC_{28}$ & $\rS^{(14)}_{234}$ & $\rC_{30}$ & $\tau_3$ & $\rA_2$ & $\rC_1\rC_2\rC_{28}$ & $2$ & $1$ \\
$14$ & $14$ & $\rS^{(8)}_{1}$ & $\rC_{28}$ & $\rS^{(18)}_{65}$ & $\rC_{12}$ & $\tau_2$ & $\rA_2$ & $\rC_1\rC_2\rC_{28}$ & $2$ & $1$ \\
$15$ & $14$ & $\rS^{(8)}_{1}$ & $\rC_{28}$ & $\rS^{(18)}_{109}$ & $\rC_{12}$ & $\tau_3$ & $\rA_2$ & $\rC_1\rC_2\rC_{28}$ & $2$ & $1$ \\
\hline
$16$ & $14$ & $\rS^{(8)}_{2}$ & $\rC_7\rC_{14}$ & $\rS^{(22)}_{6}$ & $1$ & $\tau_2$ & $\rA_{13}$ & $\rC_1^7\rC_2^7$ & $3$ & $1$ \\
$17$ & $14$ & $\rS^{(8)}_{2}$ & $\rC_9\rC_{14}$ & $\rS^{(22)}_{5}$ & $1$ & $\tau_3$ & $\rA_6$ & $\rC_1^4\rC_2^4\rC_9$ & $3$ & $2$ \\
$18$ & $14$ & $\rS^{(8)}_{2}$ & $\rC_5\rC_{10}\rC_{12}$ & $\rS^{(22)}_{2}$ & $1$ & $\tau_3$ & $\rA_{10}$ & $\rC_1^5\rC_2^5\rC_{12}$ & $3$ & $2$ \\
\hline
$19$ & $14$ & $\rS^{(8)}_{3}$ & $\rC_{10}\rC_{15}$ & $\rS^{(10)}_{71}$ & $\rC_{28}$ & $\tau_2$ & $\rA_6$ & $\rC_1^3\rC_2^3\rC_{15}$ & $4$ & $3$ \\
$20$ & $14$ & $\rS^{(8)}_{3}$ & $\rC_{10}\rC_{15}$ & $\rS^{(10)}_{71}$ & $\rC_{42}$ & $\tau_3$ & $\rA_6$ & $\rC_1^3\rC_2^3\rC_{15}$ & $4$ & $3$ \\
$21$ & $14$ & $\rS^{(8)}_{3}$ & $\rC_6\rC_7\rC_8$ & $\rS^{(22)}_{72}$ & $1$ & $\tau_1$ & $\rA_2\oplus \rD_{11}$ & $\rC_1^{11}\rC_2^3$ & $11$ & $0$ \\ 
\hline
$22$ & $14$ & $\rS^{(8)}_{5}$ & $\rC_{21}$ & $\rS^{(10)}_{19}$ & $\rC_{12}\rC_{30}$ & $\tau_2$ & $\rA_2$ & $\rC_1\rC_2\rC_{21}$ & $4$ & $3$ \\
$23$ & $14$ & $\rS^{(8)}_{5}$ & $\rC_{21}$ & $\rS^{(18)}_{16}$ & $\rC_{12}$ & $\tau_3$ & $\rA_2$ & $\rC_1\rC_2\rC_{21}$ & $4$ & $3$ \\
$24$ & $14$ & $\rS^{(8)}_{5}$ & $\rC_{21}$ & $\rS^{(18)}_{89}$ & $\rC_{12}$ & $\tau_3$ & $\rA_2$ & $\rC_1\rC_2\rC_{21}$ & $4$ & $3$ \\
$25$ & $14$ & $\rS^{(8)}_{5}$ & $\rC_{21}$ & $\rS^{(22)}_{85}$ & $1$ & $\tau_1$ & $\rA_2$ & $\rC_1\rC_2\rC_{21}$ & $4$ & $3$ \\
$26$ & $14$ & $\rS^{(8)}_{5}$ & $\rC_4\rC_{10}\rC_{18}$ & $\rS^{(22)}_{13}$ & $1$ & $\tau_1$ & $\rA_4$ & $\rC_1^3\rC_2^3\rC_4\rC_{18}$ & $3$ & $2$ \\
$27$ & $14$ & $\rS^{(8)}_{5}$ & $\rC_4\rC_{10}\rC_{18}$ & $\rS^{(22)}_{80}$ & $1$ & $\tau_1$ & $\rA_4$ & $\rC_1^3\rC_2^3\rC_4\rC_{18}$ & $3$ & $2$ \\ 
\hline 
$28$ & $14$ & $\rS^{(8)}_{6}$ & $\rC_{14}\rC_{18}$ & $\rS^{(14)}_{59}$ & $\rC_{15}$ & $\tau_3$ & $\rA_1\oplus \rA_6\oplus \rE_6$ & $\rC_1^8\rC_2^6$ & $6$ & $0$ \\
\hline
$29$ & $14$ & $\rS^{(8)}_{10}$ & $\rC_6\rC_{10}\rC_{18}$ & $\rS^{(22)}_{39}$ & $1$ & $\tau_2$ & $\rA_4$ & $\rC_1^3\rC_2^3\rC_6\rC_{18}$ & $4$ & $3$ \\ 
\hline
$30$ & $14$ & $\rS^{(8)}_{11}$ & $\rC_{26}$ & $\rS^{(6)}_{3}$ & $\rC_{40}$ & $\tau_3$ & $\rA_{12}$ & $\rC_1^7\rC_2^7$ & $3$ & $2$ \\
$31$ & $14$ & $\rS^{(8)}_{11}$ & $\rC_{26}$ & $\rS^{(18)}_{43}$ & $\rC_{12}$ & $\tau_3$ & $\rA_{12}$ & $\rC_1^7\rC_2^7$ & $3$ & $2$ \\
\hline
\end{longtable}
\begin{proposition}[Table $\mbox{\boldmath $\ref{tab:Ars1}$}$] \label{prop:Ars1} 
Every $\rA_n$-component with $n \ge 2$ appearing in Table $\ref{tab:Ars1}$ has reflectional symmetry. 
The $\rD_{11}$-component in No.\! $21$ and the $\rE_6$-component in No.\! $28$ also have reflectional symmetries. 
The cardinality $\mu_{\roff}(f)$ counted with multiplicities of the fixed points off $\cE(X)$ is given in the last column of the table.      
\end{proposition}
{\it Proof}. 
We continue to use Notation \ref{not:iec}. 
First, we consider No.~$21$, which has exceptional set of type $\rA_2 \oplus \rD_{11}$.  
By Remark \ref{rem:A}, Proposition \ref{prop:Ar} and Lemma \ref{lem:DE} we have $\mu(f, \rA_2) \ge 1$ and $\mu(f, \rD_{11}) \ge 10$, 
so by definition \eqref{eqn:muof} we have $\mu_{\ron}(f) = \mu(f, \rA_2) + \mu(f, \rD_{11}) \ge 1 + 10 = 11$.  
FPF \eqref{eqn:saito2} then implies $0 \le \mu_{\roff}(f) = L(f) - \mu_{\ron}(f) \le 11 -11 = 0$. 
This forces $\mu_{\roff}(f) = 0$ and hence $\mu(f, \rA_2) = 1$ and $\mu(f, \rD_{11}) = 10$. 
Again by Remark \ref{rem:A}, Proposition \ref{prop:Ar} and Lemma \ref{lem:DE} we have $\rA_2 = \rA_2^{\rr}$ and $\rD_{11} = \rD_{11}^{\rr}$.    
Secondly, No.~$28$ of Table \ref{tab:Ars1} is the same as No.~$27$ of Table \ref{tab:Ass}, which has already been discussed in Proposition \ref{prop:Ass}. 
Thirdly, we pick out Nos.~$22$--$25$, which has exceptional set of type $\rA_2$. 
By Remark \ref{rem:division}, if $\rA_2 = \rA_2^{\rt}$ then $\widetilde{\varphi}_2 = \rC_1^2$, 
while if $\rA_2 = \rA_2^{\rr}$ then $\widetilde{\varphi}_2 = \rC_1 \rC_2$. 
Since $\widetilde{\varphi}_2$ divides $\widetilde{\varphi}_1 = \rC_1 \rC_2 \rC_{21}$, we must have $\rA_2 = \rA_2^{\rr}$, 
and hence $\mu_{\ron}(f) = \mu(f, \rA_2^{\rr}) = 1$ by Proposition \ref{prop:Ar}. 
FPF \eqref{eqn:saito2} then yields $\mu_{\roff}(f) = L(f) - \mu_{\ron}(f) = 4 - 1 = 3$.       
The remaining entries can be treated in similar manners, except the determination of $\mu_{\roff}(f)$ in Nos.~$1$ and $2$.  
In these entries FPF \eqref{eqn:saito2} and formula \eqref{eqn:mnAr} imply  
$4 = L(f) = \mu_{\ron}(f) + \mu_{\roff}(f) \ge \mu_{\ron}(f) = \mu(f, \rA_1) + \mu(f, \rA_2^{\rr}) = \mu(f, \rA_1) + 1$, 
which means $\mu_{\ron}(f, \rA_1) \le 3$. 
Since $\mu(f, \rA_1) \ge 2$ by Remark \ref{rem:A}, we have $\mu(f, \rA_1) = 2$ or $3$, so  
$\mu_{\ron}(f) = 3$ or $4$, and hence $\mu_{\roff}(f) = 1$ or $0$, respectively. 
In fact, the case of $\mu(f, \rA_1) = 2$ and $\mu_{\roff}(f) = 1$ actually occurs, but the proof of 
this fact requires much space, so it is left to Lemma \ref{lem:A2r1} in Appendix \ref{app:exp}.               
\hfill $\Box$ \par\medskip
To include the case of Picard number $\rho = 18$ in our list, we add Table \ref{tab:Ars2},   
each entry of which has an $f$-invariant exceptional component of type $\rA$ with reflectional symmetry.   
These automorphisms stem from \cite[Setup 3.4]{IT2} and are contained in database \cite{IT3}.  
A similar proof to that of Proposition \ref{prop:Ars1} establishes Table \ref{tab:Ars2}.      
\begin{longtable}{cccccccccc}
\caption{\setlength{\leftskip}{-14mm}\setlength{\rightskip}{-14mm} 
Having a component of type $\rA$ with reflectional symmetry; \cite[Setup 3.4]{IT2}.} \label{tab:Ars2}
\\
\hline
        &            &           &          &               &      &           &                              &                       &       \\[-3mm]          
No. & $\rho$ & $S(z)$ & $C(z)$ & $\psi(z)$ & ST & Dynkin & $\tilde{\varphi}_1(z)$ & $L(f)$ & $\mu_{\roff}(f)$ \\
\hline
\endfirsthead
\multicolumn{10}{l}{continued} \\
\hline
        &           &           &           &              &      &           &                              &                       &      \\[-3mm]          
No. & $\rho$ & $S(z)$ & $C(z)$ & $\psi(z)$ & ST & Dynkin & $\tilde{\varphi}_1(z)$ & $L(f)$ & $\mu_{\roff}(f)$ \\
\hline
\endhead 
\hline
\endfoot
\hline
\\[-4pt] 
\caption*{\setlength{\leftskip}{-13mm}\setlength{\rightskip}{-13mm}
The format of this table is the same as that of Table \ref{tab:Ars1} except for the $\psi(z)$ column, 
where for example the polynomial $\psi(z)$ in No.\! $2$ is the $58$-th one among all the $1091$ polynomials 
with respect to the lexicographic order on the coefficients $(c_1, \dots, c_{11})$ in the expression \eqref{eqn:psi}. 
For each entry the $\rA_2$-component has reflectional symmetry. }
\endlastfoot
$1$ & $18$ & $\rS_1^{(4)}$ & $ \rC_{40}$ & $40$ & $\tau_1$ & $\rA_2$ & $\rC_1 \rC_2 \rC_{40}$ & $3$ & $2$ \\
$2$ & $18$ & $\rS_1^{(4)}$ & $ \rC_{40}$ & $58$ & $\tau_1$ & $\rA_2$ & $\rC_1 \rC_2 \rC_{40}$ & $3$ & $2$ \\
$3$ & $18$ & $\rS_1^{(4)}$ & $ \rC_{40}$ & $515$ & $\tau_1$ & $\rA_2$ & $\rC_1 \rC_2 \rC_{40}$ & $3$ & $2$ \\
$4$ & $18$ & $\rS_1^{(4)}$ & $ \rC_{40}$ & $579$ & $\tau_1$ & $\rA_2$ & $\rC_1 \rC_2 \rC_{40}$ & $3$ & $2$ \\
$5$ & $18$ & $\rS_1^{(4)}$ & $ \rC_{40}$ & $873$ & $\tau_1$ & $\rA_2$ & $\rC_1 \rC_2 \rC_{40}$ & $3$ & $2$ \\
$6$ & $18$ & $\rS_1^{(4)}$ & $ \rC_{48}$ & $692$ & $\tau_1$ & $\rA_2$ & $\rC_1 \rC_2 \rC_{48}$ & $3$ & $2$ \\
$7$ & $18$ & $\rS_1^{(4)}$ & $ \rC_{48}$ & $699$ & $\tau_1$ & $\rA_2$ & $\rC_1 \rC_2 \rC_{48}$ & $3$ & $2$ \\
$8$ & $18$ & $\rS_1^{(4)}$ & $ \rC_{60}$ & $457$ & $\tau_1$ & $\rA_2$ & $\rC_1 \rC_2 \rC_{60}$ & $3$ & $2$ \\
$9$ & $18$ & $\rS_1^{(4)}$ & $ \rC_{60}$ & $699$ & $\tau_1$ & $\rA_2$ & $\rC_1 \rC_2 \rC_{60}$ & $3$ &  $2$ \\
$10$ & $18$ & $\rS_1^{(4)}$ & $ \rC_{60}$ & $744$ & $\tau_1$ & $\rA_2$ & $\rC_1 \rC_2 \rC_{60}$ & $3$ & $2$ \\
$11$ & $18$ & $\rS_1^{(4)}$ & $ \rC_{60}$ & $961$ & $\tau_1$ & $\rA_2$ & $\rC_1 \rC_2 \rC_{60}$ & $3$ & $2$ \\
\end{longtable}
\subsection{Components of Type D or E} \label{ss:DEs}
Our database \cite{IT3} also contains a lot of K3 surface automorphisms $f : X \to X$ having an $f$-invariant 
exceptional component of type $\rD$ or $\rE$, too many to be presented here.  
We restrict our attention to those which have {\sl only one (simple) fixed point off} the exceptional set $\cE(X)$, 
i.e. $\mu_{\roff}(f) = 1$.  
Tables \ref{tab:DEs1} and \ref{tab:DEs2} give such examples.   
\begin{longtable}{lllllllllcr}
\caption{Having components of type $\rD$ or $\rE$; \cite[Setup 3.2]{IT2}.}  \label{tab:DEs1}
\\
\hline
       &          &           &            &          &           &     &           &                                &                   &  \\[-3mm]
No. & $\rho$ & $S(z)$ & $C(z)$ & $s(z)$ & $c(z)$ & ST & Dynkin & $\Tilde{\varphi}_1(z)$ & $L(f)$ & RH \\ 
\hline
\endfirsthead
\multicolumn{11}{l}{continued} \\
\hline
     &         &          &            &          &           &      &            &                               &                    &  \\[-3mm]
No. & $\rho$ & $S(z)$ & $C(z)$ & $s(z)$ & $c(z)$ & ST & Dynkin & $\Tilde{\varphi}_1(z)$ & $L(f)$ & RH \\ 
\hline
\endhead
\hline
\endfoot
\hline
\\[-4pt] 
\caption*{\setlength{\leftskip}{-9mm}\setlength{\rightskip}{-9mm}
The format of this table is the same as that of Table \ref{tab:Ars1} except for the last column.  
We refer to Theorem \ref{thm:rk1-2} for the explanation of the RH column.}
\endlastfoot
$1$ & $6$ & $\rS^{(16)}_{5}$ & $\rC_5$ & $\rS^{(14)}_{3}$ & $\rC_{30}$ & $\tau_3$ & $\rD_6$ & $\rC_1^5\rC_2$ & $6$ & $\rR 2$ \\
$2$ & $6$ & $\rS^{(16)}_{5}$ & $\rC_5$ & $\rS^{(14)}_{12}$ & $\rC_{30}$ & $\tau_6$ & $\rD_6$ & $\rC_1^5\rC_2$ & $6$ & $\rR 2$ \\
$3$ & $6$ & $\rS^{(16)}_{5}$ & $\rC_5$ & $\rS^{(22)}_{10}$ & $1$ & $\tau_5$ & $\rD_6$ & $\rC_1^5\rC_2$ & $6$ & $\rR 2$ \\
$4$ & $6$ & $\rS^{(16)}_{5}$ & $\rC_{10}$ & $\rS^{(10)}_{12}$ & $\rC_{36}$ & $\tau_4$ & $\rD_6$ & $\rC_1^5\rC_2$ & $6$ & $\rR 2$ \\
$5$ & $6$ & $\rS^{(16)}_{5}$ & $\rC_{10}$ & $\rS^{(22)}_{2}$ & $1$ & $\tau_4$ & $\rD_6$ & $\rC_1^5\rC_2$ & $6$ & $\rR 2$ \\
$6$ & $6$ & $\rS^{(16)}_{5}$ & $\rC_{10}$ & $\rS^{(22)}_{3}$ & $1$ & $\tau_4$ & $\rD_6$ & $\rC_1^5\rC_2$ & $6$ & $\rR 2$ \\
$7$ & $6$ & $\rS^{(16)}_{5}$ & $\rC_{10}$ & $\rS^{(22)}_{10}$ & $1$ & $\tau_2$ & $\rD_6$ & $\rC_1^5\rC_2$ & $6$ & $\rR 2$ \\
$8$ & $6$ & $\rS^{(16)}_{5}$ & $\rC_{12}$ & $\rS^{(10)}_{24}$ & $\rC_{42}$ & $\tau_7$ & $\rA_1^{\oplus 2} \oplus \rD_4$ & $\rC_1^3\rC_2\rC_3$ & $3$ & $\rR 2$ \\
$9$ & $6$ & $\rS^{(16)}_{5}$ & $\rC_{12}$ & $\rS^{(10)}_{107}$ & $\rC_{42}$ & $\tau_1$ & $\rA_1^{\oplus 2} \oplus \rD_4$ & $\rC_1^3\rC_2\rC_3$ & $3$ & $\rR 2$ \\
$10$ & $6$ & $\rS^{(16)}_{5}$ & $\rC_{12}$ & $\rS^{(14)}_{12}$ & $\rC_{30}$ & $\tau_2$ & $\rA_1^{\oplus 2} \oplus \rD_4$ & $\rC_1^3\rC_2\rC_3$ & $3$ & $\rR 2$ \\
$11$ & $6$ & $\rS^{(16)}_{5}$ & $\rC_{12}$ & $\rS^{(14)}_{17}$ & $\rC_{15}$ & $\tau_4$ & $\rA_1^{\oplus 2} \oplus \rD_4$ & $\rC_1^3\rC_2\rC_3$ & $3$ & $\rR 2$ \\
$12$ & $6$ & $\rS^{(16)}_{5}$ & $\rC_{12}$ & $\rS^{(14)}_{382}$ & $\rC_{20}$ & $\tau_1$ & $\rA_1^{\oplus 2} \oplus \rD_4$ & $\rC_1^3\rC_2\rC_3$ & $3$ & $\rR 2$ \\
\hline
$13$ & $8$ & $\rS^{(14)}_{1}$ & $\rC_{14}$ & $\rS^{(14)}_{4}$ & $\rC_{24}$ & $\tau_1$ & $\rE_8$ & $\rC_1^8$ & $10$ & $\rR 2$ \\
$14$ & $8$ & $\rS^{(14)}_{1}$ & $\rC_{14}$ & $\rS^{(14)}_{15}$ & $\rC_{24}$ & $\tau_5$ & $\rE_8$ & $\rC_1^8$ & $10$ & $\rR 2$ \\
$15$ & $8$ & $\rS^{(14)}_{1}$ & $\rC_{14}$ & $\rS^{(22)}_{10}$ & $1$ & $\tau_5$ & $\rE_8$ & $\rC_1^8$ & $10$ & $\rR 2$ \\
$16$ & $8$ & $\rS^{(14)}_{1}$ & $\rC_{14}$ & $\rS^{(22)}_{13}$ & $1$ & $\tau_1$ & $\rE_8$ & $\rC_1^8$ & $10$ & $\rR 2$ \\
$17$ & $8$ & $\rS^{(14)}_{1}$ & $\rC_{14}$ & $\rS^{(22)}_{89}$ & $1$ & $\tau_5$ & $\rE_8$ & $\rC_1^8$ & $10$ & $\rR 2$ \\
$18$ & $8$ & $\rS^{(14)}_{1}$ & $\rC_{18}$ & $\rS^{(6)}_{1}$ & $\rC_{60}$ & $\tau_1$ & $\rE_8$ & $\rC_1^8$ & $10$ & $\rR 2$ \\
$19$ & $8$ & $\rS^{(14)}_{1}$ & $\rC_{18}$ & $\rS^{(14)}_{15}$ & $\rC_{24}$ & $\tau_3$ & $\rE_8$ & $\rC_1^8$ & $10$ & $\rR 2$ \\
$20$ & $8$ & $\rS^{(14)}_{1}$ & $\rC_{18}$ & $\rS^{(14)}_{85}$ & $\rC_{24}$ & $\tau_3$ & $\rE_8$ & $\rC_1^8$ & $10$ & $\rR 2$ \\
$21$ & $8$ & $\rS^{(14)}_{1}$ & $\rC_{18}$ & $\rS^{(22)}_{43}$ & $1$ & $\tau_3$ & $\rE_8$ & $\rC_1^8$ & $10$ & $\rR 2$ \\
$22$ & $8$ & $\rS^{(14)}_{1}$ & $\rC_{18}$ & $\rS^{(22)}_{72}$ & $1$ & $\tau_1$ & $\rE_8$ & $\rC_1^8$ & $10$ & $\rR 2$ \\
$23$ & $8$ & $\rS^{(14)}_{1}$ & $\rC_{18}$ & $\rS^{(22)}_{89}$ & $1$ & $\tau_3$ & $\rE_8$ & $\rC_1^8$ & $10$ & $\rR 2$ \\
$24$ & $8$ & $\rS^{(14)}_{1}$ & $\rC_3\rC_{12}$ & $\rS^{(22)}_{61}$ & $1$ & $\tau_1$ & $\rE_8$ & $\rC_1^8$ & $10$ & $\rR 2$ \\
$25$ & $8$ & $\rS^{(14)}_{1}$ & $\rC_4\rC_8$ & $\rS^{(22)}_{13}$ & $1$ & $\tau_5$ & $\rE_8$ & $\rC_1^8$ & $10$ & $\rR 2$ \\
$26$ & $8$ & $\rS^{(14)}_{1}$ & $\rC_4\rC_{10}$ & $\rS^{(10)}_{15}$ & $\rC_{42}$ & $\tau_5$ & $\rE_8$ & $\rC_1^8$ & $10$ & $\rR 2$ \\
$27$ & $8$ & $\rS^{(14)}_{1}$ & $\rC_4\rC_{10}$ & $\rS^{(14)}_{11}$ & $\rC_{24}$ & $\tau_1$ & $\rE_8$ & $\rC_1^8$ & $10$ & $\rR 2$ \\
$28$ & $8$ & $\rS^{(14)}_{1}$ & $\rC_4\rC_{10}$ & $\rS^{(22)}_{13}$ & $1$ & $\tau_4$ & $\rE_8$ & $\rC_1^8$ & $10$ & $\rR 2$ \\
$29$ & $8$ & $\rS^{(14)}_{1}$ & $\rC_4\rC_{12}$ & $\rS^{(22)}_{10}$ & $1$ & $\tau_1$ & $\rE_8$ & $\rC_1^8$ & $10$ & $\rR 2$ \\
$30$ & $8$ & $\rS^{(14)}_{1}$ & $\rC_4\rC_{12}$ & $\rS^{(22)}_{43}$ & $1$ & $\tau_1$ & $\rE_8$ & $\rC_1^8$ & $10$ & $\rR 2$ \\ 
\hline
$31$ & $10$ & $\rS^{(12)}_{1}$ & $\rC_{16}$ & $\rS^{(6)}_{1}$ & $\rC_{60}$ & $\tau_2$ & $\rD_9$ & $\rC_1^8\rC_2^2$ & $9$ & $\rHy$ \\
$32$ & $10$ & $\rS^{(12)}_{1}$ & $\rC_{16}$ & $\rS^{(6)}_{3}$ & $\rC_{60}$ & $\tau_5$ & $\rD_9$ & $\rC_1^8\rC_2^2$ & $9$ & $\rR 2$ \\
$33$ & $10$ & $\rS^{(12)}_{1}$ & $\rC_{16}$ & $\rS^{(10)}_{7}$ & $\rC_{36}$ & $\tau_5$ & $\rD_9$ & $\rC_1^8\rC_2^2$ & $9$ & $\rR 2$ \\
$34$ & $10$ & $\rS^{(12)}_{1}$ & $\rC_{16}$ & $\rS^{(14)}_{18}$ & $\rC_{24}$ & $\tau_5$ & $\rD_9$ & $\rC_1^8\rC_2^2$ & $9$ & $\rR 2$ \\
$35$ & $10$ & $\rS^{(12)}_{1}$ & $\rC_{16}$ & $\rS^{(14)}_{241}$ & $\rC_{15}$ & $\tau_5$ & $\rD_9$ & $\rC_1^8\rC_2^2$ & $9$ & $\rR 2$ \\
$36$ & $10$ & $\rS^{(12)}_{1}$ & $\rC_{16}$ & $\rS^{(22)}_{5}$ & $1$ & $\tau_4$ & $\rD_9$ & $\rC_1^8\rC_2^2$ & $9$ & $\rHy$ \\
$37$ & $10$ & $\rS^{(12)}_{1}$ & $\rC_{16}$ & $\rS^{(22)}_{89}$ & $1$ & $\tau_2$ & $\rD_9$ & $\rC_1^8\rC_2^2$ & $9$ & $\rHy$ \\
$38$ & $10$ & $\rS^{(12)}_{1}$ & $\rC_4\rC_7$ & $\rS^{(6)}_{1}$ & $\rC_{60}$ & $\tau_1$ & $\rD_9$ & $\rC_1^8\rC_2^2$ & $9$ & $\rR 2$ \\
\hline
$39$ & $12$ & $\rS^{(10)}_{1}$ & $\rC_4\rC_{16}$ & $\rS^{(6)}_{1}$ & $\rC_{40}$ & $\tau_2$ & $\rE_6^{\oplus 2}$ & $\rC_1^4\rC_2^4\rC_4^2$ & $1$ & $\rR 2$ \\
$40$ & $12$ & $\rS^{(10)}_{1}$ & $\rC_4\rC_{16}$ & $\rS^{(6)}_{1}$ & $\rC_{60}$ & $\tau_2$ & $\rE_6^{\oplus 2}$ & $\rC_1^4\rC_2^4\rC_4^2$ & $1$ & $\rR 2$ \\
$41$ & $12$ & $\rS^{(10)}_{1}$ & $\rC_4\rC_{16}$ & $\rS^{(10)}_{2}$ & $\rC_{42}$ & $\tau_4$ & $\rE_6^{\oplus 2}$ & $\rC_1^4\rC_2^4\rC_4^2$ & $1$ & $\rHy$ \\
$42$ & $12$ & $\rS^{(10)}_{1}$ & $\rC_4\rC_{16}$ & $\rS^{(22)}_{10}$ & $1$ & $\tau_1$ & $\rE_6^{\oplus 2}$ & $\rC_1^4\rC_2^4\rC_4^2$ & $1$ & $\rR 2$ \\
$43$ & $12$ & $\rS^{(10)}_{1}$ & $\rC_4\rC_{20}$ & $\rS^{(14)}_{2}$ & $\rC_{30}$ & $\tau_4$ & $\rE_6^{\oplus 2}$ & $\rC_1^4\rC_2^4\rC_4^2$ & $1$ & $\rHy$ \\
$44$ & $12$ & $\rS^{(10)}_{1}$ & $\rC_4\rC_{20}$ & $\rS^{(22)}_{1}$ & $1$ & $\tau_4$ & $\rE_6^{\oplus 2}$ & $\rC_1^4\rC_2^4\rC_4^2$ & $1$ & $\rHy$ \\
$45$ & $12$ & $\rS^{(10)}_{1}$ & $\rC_4\rC_{20}$ & $\rS^{(22)}_{14}$ & $1$ & $\tau_2$ & $\rE_6^{\oplus 2}$ & $\rC_1^4\rC_2^4\rC_4^2$ & $1$ & $\rR 2$ \\
$46$ & $12$ & $\rS^{(10)}_{1}$ & $\rC_4\rC_{20}$ & $\rS^{(22)}_{18}$ & $1$ & $\tau_4$ & $\rE_6^{\oplus 2}$ & $\rC_1^4\rC_2^4\rC_4^2$ & $1$ & $\rHy$ \\
$47$ & $12$ & $\rS^{(10)}_{1}$ & $\rC_3\rC_4\rC_6\rC_8$ & $\rS^{(22)}_{72}$ & $1$ & $\tau_3$ & $\rE_6^{\oplus 2}$ & $\rC_1^4\rC_2^4\rC_4^2$ & $1$ & $\rR 2$ \\ 
\hline
$48$ & $14$ & $\rS^{(8)}_{2}$ & $\rC_4\rC_{11}$ & $\rS^{(6)}_{1}$ & $\rC_{60}$ & $\tau_3$ & $\rD_{13}$ & $\rC_1^{12}\rC_2^2$ & $13$ & $\rHy$ \\
$49$ & $14$ & $\rS^{(8)}_{2}$ & $\rC_5\rC_{16}$ & $\rS^{(14)}_{1}$ & $\rC_{30}$ & $\tau_1$ & $\rD_{13}$ & $\rC_1^{12}\rC_2^2$ & $13$ & $\rR 2$ \\
$50$ & $14$ & $\rS^{(8)}_{2}$ & $\rC_8\rC_{24}$ & $\rS^{(10)}_{12}$ & $\rC_{36}$ & $\tau_2$ & $\rD_{13}$ & $\rC_1^{12}\rC_2^2$ & $13$ & $\rR 2$ \\
$51$ & $14$ & $\rS^{(8)}_{2}$ & $\rC_8\rC_{24}$ & $\rS^{(14)}_{26}$ & $\rC_{30}$ & $\tau_2$ & $\rD_{13}$ & $\rC_1^{12}\rC_2^2$ & $13$ & $\rR 2$ \\
$52$ & $14$ & $\rS^{(8)}_{2}$ & $\rC_8\rC_{24}$ & $\rS^{(22)}_{5}$ & $1$ & $\tau_2$ & $\rD_{13}$ & $\rC_1^{12}\rC_2^2$ & $13$ & $\rR 2$ \\
$53$ & $14$ & $\rS^{(8)}_{2}$ & $\rC_8\rC_{24}$ & $\rS^{(22)}_{39}$ & $1$ & $\tau_1$ & $\rD_{13}$ & $\rC_1^{12}\rC_2^2$ & $13$ & $\rR 2$ \\
\hline
$54$ & $14$ & $\rS^{(8)}_{4}$ & $\rC_{26}$ & $\rS^{(6)}_{3}$ & $\rC_{48}$ & $\tau_3$ & $\rD_{14}$ & $\rC_1^{13}\rC_2$ & $14$ & $\rR 2$ \\
$55$ & $14$ & $\rS^{(8)}_{4}$ & $\rC_{26}$ & $\rS^{(10)}_{2}$ & $\rC_{28}$ & $\tau_3$ & $\rD_{14}$ & $\rC_1^{13}\rC_2$ & $14$ & $\rR 2$ \\
$56$ & $14$ & $\rS^{(8)}_{4}$ & $\rC_{26}$ & $\rS^{(14)}_{18}$ & $\rC_{30}$ & $\tau_3$ & $\rD_{14}$ & $\rC_1^{13}\rC_2$ & $14$ & $\rR 2$ \\
\hline
$57$ & $14$ & $\rS^{(8)}_{5}$ & $\rC_7\rC_{18}$ & $\rS^{(10)}_{5}$ & $\rC_{42}$ & $\tau_3$ & $\rE_7 $ & $\rC_1^7\rC_2\rC_{18}$ & $9$ & $\rHy$ \\
$58$ & $14$ & $\rS^{(8)}_{5}$ & $\rC_7\rC_{18}$ & $\rS^{(18)}_{123}$ & $\rC_{12}$ & $\tau_1$ & $\rE_7 $ & $\rC_1^7\rC_2\rC_{18}$ & $9$ & $\rR 2$ \\
$59$ & $14$ & $\rS^{(8)}_{5}$ & $\rC_7\rC_{18}$ & $\rS^{(22)}_{61}$ & $1$ & $\tau_3$ & $\rE_7 $ & $\rC_1^7\rC_2\rC_{18}$ & $9$ & $\rHy$ \\
$60$ & $14$ & $\rS^{(8)}_{5}$ & $\rC_3\rC_{12}\rC_{18}$ & $\rS^{(10)}_{12}$ & $\rC_{28}$ & $\tau_3$ & $\rE_7 $ & $\rC_1^7\rC_2\rC_{18}$ & $9$ & $\rHy$ \\
\hline 
$61$ & $14$ & $\rS^{(8)}_{6}$ & $\rC_8\rC_{15}$ & $\rS^{(22)}_{25}$ & $1$ & $\tau_1$ & $\rD_5\oplus \rE_8$ & $\rC_1^{12}\rC_2^2$ & $14$ & $\rHy$ \\
\hline
$62$ & $14$ & $\rS^{(8)}_{8}$ & $\rC_3\rC_4\rC_{16}$ & $\rS^{(14)}_{2}$ & $\rC_{30}$ & $\tau_3$ & $\rD_4^{\oplus 2}$ & $\rC_1^4\rC_2^4\rC_3\rC_4^2$ & $1$ & $\rHy$ \\
$63$ & $14$ & $\rS^{(8)}_{8}$ & $\rC_3\rC_4\rC_{16}$ & $\rS^{(22)}_{18}$ & $1$ & $\tau_3$ & $\rD_4^{\oplus 2}$ & $\rC_1^4\rC_2^4\rC_3\rC_4^2$ & $1$ & $\rHy$ \\
\hline
$64$ & $16$ & $\rS^{(6)}_{1}$ & $\rC_4\rC_{26}$ & $\rS^{(22)}_{43}$ & $1$ & $\tau_1$ & $\rD_{16}$ & $\rC_1^{16}$ & $18$ & $\rR 2$ \\
$65$ & $16$ & $\rS^{(6)}_{1}$ & $\rC_4\rC_{26}$ & $\rS^{(22)}_{65}$ & $1$ & $\tau_1$ & $\rD_{16}$ & $\rC_1^{16}$ & $18$ & $\rR 2$ \\
$66$ & $16$ & $\rS^{(6)}_{1}$ & $\rC_3\rC_{14}\rC_{18}$ & $\rS^{(14)}_{11}$ & $\rC_{30}$ & $\tau_1$ & $\rA_1^{\oplus 3}\oplus \rE_6$ & $\rC_1^5\rC_2^3\rC_3\rC_{14}$ & $4$ & 
$\rR 2$ \\
$67$ & $16$ & $\rS^{(6)}_{1}$ & $\rC_3\rC_4\rC_9\rC_{12}$ & $\rS^{(14)}_{2}$ & $\rC_{30}$ & $\tau_2$ & $\rD_5^{\oplus 3}$ & $\rC_1^4\rC_2^2\rC_3^4\rC_6$ & $1$ & $\rHy$ \\ 
\hline
\end{longtable}
For some entries of Table \ref{tab:DEs1} the exceptional set $\cE(X)$ contains more than 
one exceptional components of the same Dynkin type; for example, two $\rA_1$-components in Nos.~$8$--$12$ 
and three $\rD_5$-components in No.~$67$.  
Such components are $f$-invariant as a whole, but they may be permuted by $f$. 
It is obvious that every component of type $\rE_7$ or $\rE_8$ in Table $\ref{tab:DEs1}$ has trivial symmetry.   
With these remarks we have the following results.     
\begin{proposition}[Table $\mbox{\boldmath $\ref{tab:DEs1}$}$] \label{prop:DEs}  
In Nos.~$1$--$7$ the $\rD_6$-component has reflectional symmetry.   
In Nos.~$8$--$12$ the two $\rA_1$-components are swapped, while the $\rD_4$-component has tricyclic symmetry. 
In Nos.~$31$--$38$ the $\rD_9$-component has reflectional symmetry. 
In Nos.~$39$--$47$ the two $\rE_6$-components are swapped. 
In Nos.~$48$--$53$ the $\rD_{13}$-component has reflectional symmetry. 
In Nos.~$54$--$56$ the $\rD_{14}$-component has reflectional symmetry. 
In No.~$61$ the $\rD_5$-component has reflectioanl symmetry. 
In Nos.~$62$ and $63$ the two $\rD_4$-components are swapped. 
In Nos.~$64$ and $65$ the $\rD_{16}$-component has trivial symmetry. 
In No.~$66$ the three $\rA_1$-components are cyclically permuted, while the $\rE_6$-component has reflectional symmetry. 
In No.~$67$ the three $\rD_5$-components are cyclically permuted.   
Every entry of the table has exactly one fixed point, which is simple, off the exceptional set $\cE(X)$, that is, $\mu_{\roff}(f) = 1$.     
\end{proposition}
{\it Proof}. 
First let us consider Nos.~$8$--$12$, which have exceptional set of type $\rA_1^{\oplus 2} \oplus \rD_4$.     
Suppose that the two $\rA_1$-components are $f$-invariant. 
By Remark \ref{rem:A} and Lemma \ref{lem:DE} we have   
$\mu(f, \rA_1) \ge 2$ and $\mu(f, \rD_4) \ge 2$, so by definition \eqref{eqn:muof} we have        
$\mu_{\ron}(f) = \mu(f, \rA_1) + \mu(f, \rA_1) + \mu(f, \rD_4) \ge 2 + 2 + 2 = 6$. 
FPF \eqref{eqn:saito2} then implies 
$0 \le \mu_{\roff}(f) = L(f) - \mu_{\roff}(f) \le 2 + 1 -6 = -3$, a contradiction.  
Thus the two $\rA_1$-components must be swapped and the $\rA_1^{\oplus 2}$-part contributes to 
$\widetilde{\varphi}_2$ by $\rC_1 \rC_2$, where $\widetilde{\varphi}_2$ is defined in Remark \ref{rem:division}.   
If $\rD_4 = \rD_4^{\rt}$ then $\widetilde{\varphi}_2 = \rC_1 \rC_2 \cdot \rC_1^4 = \rC_1^5 \rC_2$;  
if $\rD_4 = \rD_4^{\rr}$ then $\widetilde{\varphi}_2 = \rC_1 \rC_2 \cdot \rC_1^3 \rC_2 = \rC_1^4 \rC_2^2$; 
if $\rD_4 = \rD_4^{\rc}$ then $\widetilde{\varphi}_2 = \rC_1 \rC_2 \cdot \rC_1^2 \rC_3 = \rC_1^3 \rC_2 \rC_3$. 
Since $\widetilde{\varphi}_2$ divides $\widetilde{\varphi}_1 = \rC_1^3 \rC_2 \rC_3$, we must have $\rD_4 = \rD_4^{\rc}$ 
and hence $\mu_{\ron}(f) = \mu(f, \rD_4^{\rc}) = 2$ by Lemma \ref{lem:DE}. 
FPF \eqref{eqn:saito2} now gives $\mu_{\roff}(f) = L(f) - \mu(f, \rD_4^{\rc}) = 2 + 1 - 2 = 1$.  
Secondly we look at Nos.~$39$--$47$, which have exceptional set of type $\rE_6^{\oplus 2}$. 
If the two $\rE_6$-components are $f$-invariant, then $\mu_{\ron}(f) = \mu(f, \rE_6) + \mu(f, \rE_6) \ge 3 + 3 = 6$ 
and hence $0 \le \mu_{\roff}(f) = L(f) - \mu_{\ron}(f) \le 2 -1 -6 = -5$, a contradiction. 
So the two $\rE_6$-components must be swapped, hence $\mu_{\ron}(f) = 0$ and 
$\mu_{\roff}(f) = L(f) - \mu_{\ron}(f) = 2-1-0 = 1$. 
Thirdly we consider No.~$66$, which has exceptional set of type $\rA_1^{\oplus 3} \oplus \rE_6$.  
If at least one of the three $\rA_1$-components is $f$-invariant, then 
$\mu_{\ron}(f) \ge \mu(f, \rA_1) + \mu(f, \rE_6) \ge 2 + 3 = 5$ and hence 
$0 \le \mu_{\roff}(f) = L(f) - \mu_{\ron}(f) \le 2 + 2 - 5 = -1$, a contradiction.  
Thus the three $\rA_1$-components must be cyclically permuted and $\mu_{\ron}(f) = \mu(f, \rE_6)$. 
If $\rE_6 = \rE_6^{\rt}$ then $\mu_{\ron}(f) = 7$ and $0 \le \mu_{\roff}(f) = L(f) - \mu_{\ron}(f) \le 2 + 2 - 7 = -3$, 
which is absurd.  
Hence $\rE_6 = \rE_6^{\rr}$ and $\mu_{\ron}(f) = 3$, so $\mu_{\roff}(f) = 2 + 2 -3 = 1$. 
The remaining entries can be treated in similar or simpler manners.                
\hfill $\Box$ 
\begin{longtable}{ccclccllcr}
\caption{Having components of type $\rE$; \cite[Setup 3.4]{IT2}.} \label{tab:DEs2} 
\\
\hline
      &           &           &          &               &      &           &                              &                      &      \\[-3mm]          
No. & $\rho$ & $S(z)$ & $C(z)$ & $\psi(z)$ & ST & Dynkin & $\tilde{\varphi}_1(z)$ & $L(f)$ & RH \\
\hline
\endfirsthead
\multicolumn{10}{l}{continued} \\
\hline
      &           &           &          &               &      &           &                              &                      &      \\[-3mm]          
No. & $\rho$ & $S(z)$ & $C(z)$ & $\psi(z)$ & ST & Dynkin & $\tilde{\varphi}_1(z)$ & $L(f)$ & RH \\
\hline
\endhead 
\hline
\endfoot
\hline
\\[-4pt] 
\caption*{\setlength{\leftskip}{-13mm}\setlength{\rightskip}{-13mm}
The format of this table is the same as that of Table \ref{tab:Ars2} except for RH-column, 
which indicates that the unique fixed point off $\cE(X)$ is a center of a Siegel disk. }
\endlastfoot
$1$ & $18$ & $\rS_1^{(4)}$ & $ \rC_8 \rC_{12} \rC_{30}$ & $523$ & $\tau_1$ & $\rA_2^{\oplus 2}\oplus \rE_6\oplus \rE_8$ & $\rC_1^{13} \rC_2^3 \rC_4$ & $13$ & $\rR 2$ \\
\end{longtable}
\begin{proposition}[Table $\mbox{\boldmath $\ref{tab:DEs2}$}$] \label{prop:DEs2} 
If the two $\rA_2$-components in $\cE(X)$ are denoted by the diagram  
$\xymatrix@C=10pt@R=10pt{ 1 \ar@{-}[r] & 2,} \xymatrix@C=10pt@R=10pt{1' \ar@{-}[r] & 2',}$
then the automorphism $f : X \to X$ induces the cyclic permutation $(1, 1', 2, 2')$ and, in particular,     
swaps the two $\rA_2$-components. 
The $\rE_6$-component has reflectional symmetry and, obviously, the $\rE_8$-component has trivial symmetry. 
The map $f$ has exactly one fixed point off $\cE(X)$, which is simple, i.e. $\mu_{\roff}(f) = 1$.      
\end{proposition}
{\it Proof}. From Lemma \ref{lem:DE} we have $\mu(f, \rE_6) \ge 3$ and $\mu(f, \rE_8) = 9$. 
If each of the two $\rA_2$-components is $f$-invariant, then 
$\mu(f, \rA_2^{\oplus 2}) = \mu(f, \rA_2) + \mu(f, \rA_2) \ge 2 + 2 = 4$ by Remark \ref{rem:A}, 
so $\mu_{\ron}(f) = \mu(f, \rA_2^{\oplus 2}) + \mu(f, \rE_6) + \mu(f, \rE_8) \ge 4 + 3 + 9 = 16$, 
which leads to a contradiction $0 \le \mu_{\roff}(f) = L(f) - \mu_{\ron}(f) \le 13-16 = -3$ by FPF \eqref{eqn:saito2}. 
Thus $f$ swaps the two $\rA_2$-components. 
If $\rE_6 = \rE_6^{\rt}$, then $\mu_{\ron}(f) = \mu(f, \rE_6^{\rt}) + \mu(f, \rE_8) = 7 + 9 = 16$ 
by Lemma \ref{lem:DE}, and hence $0 \le \mu_{\roff}(f) = L(f) - \mu_{\ron}(f) = 13 - 16 = -3$, 
again a contradiction. 
Thus $\rE_6 = \rE_6^{\rr}$ and $\mu_{\ron}(f) = \mu(f, \rE_6^{\rr}) + \mu(f, \rE_8) = 3 + 9 = 12$, 
so $\mu_{\roff}(f) = L(f) - \mu_{\ron}(f) = 13 - 12 = 1$. 
By Remark \ref{rem:division} the $\rE_6^{\rr}$-component and $\rE_8$-component contribute to 
the polynomial $\widetilde{\varphi}_2$ by $\rC_1^4 \rC_2^2$ and $\rC_1^8$ respectively. 
If $f$ acts on the $\rA_2^{\oplus 2}$-part by cycle type $(1, 1') (2, 2')$, then the contribution of 
$\rA_2^{\oplus 2}$ to $\widetilde{\varphi}_2$ is $(\rC_1 \rC_2)^2$, so we have  
$\widetilde{\varphi}_2 = (\rC_1 \rC_2)^2 \cdot \rC_1^4 \rC_2^2 \cdot \rC_1^8 = \rC_1^{14} \rC_2^4$, 
but this does not divide $\widetilde{\varphi}_1 = \rC_1^{13} \rC_2^3 \rC_4$, yet another contradicition. 
Hence the cycle type on the $\rA_2^{\oplus 2}$-part must be $(1, 1', 2, 2')$.         
\hfill $\Box$ 
\section{Rotation Domains} \label{sec:rd}
In this and next sections, putting together all what we have developed so far,   
we shall construct a lot of K3 surface automorphisms of positive entropy having 
rotation domains of rank $1$ and/or of rank $2$. 
Various particular examples from \S \ref{sec:mhgg}, manipulated by general theories  
developed in \S\S \ref{sec:equilin}--\ref{sec:fpf}, are most important.   
\subsection{Reductions to Linear Models}
Detecting rotation domains or hyperbolic sets (or points) often relies on 
reduction to linear models. 
We begin with the simplest linear models around a fixed point (Lemma \ref{lem:rh-splest}) 
and proceed to linear models around an exceptional component (Lemma \ref{lem:rh-lm}). 
A whole space that is a rotation domain is called a {\sl rotation space}.  
\begin{lemma} \label{lem:rh-splest}   
Let $\bc = (c_1, c_2) \in (\bC^{\times})^2$ be a multiplier pair such that $| c_1 c_2| = 1$.  
For the linear map 
\begin{equation} \label{eqn:rh-splest}
\bC^2 \to \bC^2, \quad (z_1, \, z_2) \mapsto (c_1 z_1, \, c_2 z_2),  
\end{equation}
\begin{enumerate}
\setlength{\itemsep}{-1pt}
\item[$(\rR 1)$] if $|c_1| = |c_2| = 1$, $\bc$ is MD and $c_1 c_2$ is not a root of unity, then $\bC^2$ is a rotation space of rank $1$,  
\item[$(\rR 2)$] if $|c_1| = |c_2| = 1$ and $\bc$ is MI, then $\bC^2$ is a rotation space of rank $2$,    
\item[$(\rHy)$] if $|c_1| = |c_2|^{-1} \neq 1$, then the origin $(0, 0)$ is a hyperbolic fixed point and the coordinate axes are 
the stable and unstable curves flowing into and emanating from the origin respectively.   
\end{enumerate}
\end{lemma}
{\it Proof}. 
Assertions $(\rR 1)$ and $(\rR 2)$ follow from the definition of a rotation domain and an 
elementary result in Diophantine approximations, while $(\rHy)$ is trivial from the definition of a hyperbolic fixed point.  
\hfill $\Box$ \par\medskip
For a nontrivial finite subgroup $G$ of $\SU(2)$, let $(Y, E)$ denote the minimal 
resolution of the quotient singularity $\bC^2/G$ constructed in \S\ref{sec:rs}, where $E$ is the exceptional set in $Y$. 
It was denoted by $(X, E)$ in \S\S \ref{sec:rs}--\ref{sec:lmnec}, but in this section $(X, E)$   
is reserved for the pair of a K3 surface $X$ and an exceptional component $E$ on $X$.   
The two sets $E \subset X$ and $E \subset Y$ are often identified, 
so using the same symbol $E$ would not be confusing.     
Given a normalizer $J \in N_{\GL}$ of $G$ whose determinant has modulus $1$, 
we consider its linear model map $f_J : (Y, E) \to (Y, E)$ discussed in \S\S\ref{sec:nmr}--\ref{sec:lmnec}.  
We assume that $J$ is diagonalizable when $G$ is a group of order $2$, so that 
$J$ is diagonalizable in any case by Theorem \ref{thm:normalizer}.   
We denote the eigenvalues of $J$ by $\bal = (\alpha_1, \alpha_2) \in (\bC^{\times})^2$.  
\begin{lemma} \label{lem:rh-lm} 
For the map $f_J : (Y, E) \to (Y, E)$ mentioned above the following assertions hold.  
\begin{enumerate}
\setlength{\itemsep}{-1pt}
\item[$(\ri)$] When $G$ is a group of Dynkin type $\rA_n$, $n \in \bN$, and $\sigma_J \in \Aut(\rA_n)$ is trivial,  
\begin{enumerate}
\setlength{\itemsep}{-1pt}
\item[$(\rR 1)$] if $|\alpha_1| = |\alpha_2| = 1$, $\bal$ is MD and $\det J$ is not a root of unity, then $Y$ is a rotation space of rank $1$,  
\item[$(\rR 2)$] if $|\alpha_1| = |\alpha_2| = 1$ and $\bal$ is MI, then $Y$ is a rotation space of rank $2$,    
\item[$(\rHy)$] if $|\alpha_1| = |\alpha_2|^{-1} \neq 1$, then $E$ is a hyperbolic set to the effect that a stable curve flows into $E$ and 
an unstable curve emanates from $E$.    
\end{enumerate}
\item[$(\rii)$] 
Suppose that $\det J$ is not a root of unity. 
Then $Y$ is a rotation space of rank $1$ in the following cases:   
\begin{enumerate} 
\setlength{\itemsep}{-1pt}
\item $G$ is a group of Dynkin type $\rA_n$, $n \ge 2$, and $\sigma_J \in \Aut(\rA_n) \cong \bZ_2$ is nontrivial; 
\item $G$ is a group of Dynkin type $\rD$ or $\rE$, in which case $\sigma_J$ may be arbitrary. 
\end{enumerate}  
\end{enumerate}
\end{lemma} 
{\it Proof}. 
Assertion $(\ri)$. 
We express the eigenvalues of $J$ as $(\alpha_1, \alpha_2) = (\beta \alpha, \, \beta \alpha^{-1})$.  
Since $J$ is diagonalizable, we may assume that $J$ is of the form \eqref{eqn:fa} in Lemma \ref{lem:ida1}, 
so $f_J$ is given by formula \eqref{eqn:f3}. 
This means that on the chart $W_i \cong \bC^2$ the map $f_J$ takes the normal form \eqref{eqn:rh-splest} 
with $c_1 := \beta^{2 i-n-1} \alpha^{n+1}$ and $c_2 := \beta^{n+3-2 i} \alpha^{-n-1}$. 
Notice that $(c_1, c_2)$ is MI, if and only if $(\alpha, \beta)$ is MI, if and only if  $(\alpha_1, \alpha_2)$ is MI.  
Notice also that $c_1 c_2 = \beta^2 = \alpha_1 \alpha_2 = \det J$. 
Thus assertions $(\rR 1)$ and $(\rR 2)$ of Lemma \ref{lem:rh-lm} follow from those of Lemma \ref{lem:rh-splest}. 
In case $(\rHy)$, since $|\det J | = |\beta|^2 = 1$, we have $|\alpha| \neq 1$. 
By Lemma \ref{lem:ida1}, if $|\alpha| > 1$ then $E_0$ is the stable curve flowing into $E$ at $p_0$ with multiplier $(\beta/\alpha)^{n+1}$, 
and $E_{n+1}$ is the unstable curve emanating from $E$ at $p_{n}$ with multiplier $(\beta \alpha)^{n+1}$.  
If $|\alpha| < 1$, the stable and unstable curves alternate. 
This gives assertion $(\rHy)$. 
\par
Assertion $(\rii)$. 
To address the present issue we may replace $f_J$ by a power of $f_J$.   
In case (a) the matrix $J$ is of the form \eqref{eqn:fb} in Lemma \ref{lem:ida2}, so $J^2$ takes the form 
\eqref{eqn:fa} with $\alpha = -1$ and $\beta^2$ in replace of $\beta$ in Lemma \ref{lem:ida1}. 
Now case $(\rR1)$ of assertion $(\ri)$ can be applied to $(f_J)^2 = f_{J^2}$ to have the conclusion. 
In case (b) we may assume $\sigma_J = 1$ by replacing $J$ with $J^2$ or $J^3$ according to  
whether $\sigma_J$ has order $2$ or $3$.  
Recall from Figure \ref{fig:dDE} that the exceptional set $E$ consists of the core $E_0$ and three arms   
$\rE^{\nu}$, $\nu = \re, \rf, \rv$. 
By Lemma \ref{lem:order1}, case $(\rR 1)$ of Lemma \ref{lem:rh-splest} applies 
with $c_1 := (\det J)^{-i}$ and $c_2 := (\det J)^{i+1}$ in a local chart around the fixed point 
$p_i^{\nu}$, $0 \le i \le n^{\nu}$. 
\hfill $\Box$ \par\medskip
We turn our attention to a K3 surface automorphism $f : X \to X$ satisfying 
conditions (A1) and (A2) in \S \ref{sec:intro}. 
Given an $f$-invariant exceptional component $E \subset \cE(X)$, it is interesting 
to ask whether and when $E$ falls into one of the cases $(\rR 1), (\rR 2), (\rHy)$  
in Problem \ref{prob:rd}.  
Our answer is the following result.  
\begin{theorem} \label{thm:rd} 
Let $E$ be an $f$-invariant exceptional component. 
If $E$ is of type $\rA$ with reflectional symmetry, or if $E$ is of type $\rD$ or $\rE$ 
with arbitrary symmetry, then $E$ always falls into case $(\rR 1)$. 
On the other hand, when $E$ is of type $\rA$ with trivial symmetry, 
all of the three cases $(\rR 1), (\rR 2), (\rHy)$ can actually occur. 
\end{theorem}
{\it Proof}. 
If $E$ is of type $\rA$ with reflectional symmetry, or if $E$ is of type $\rD$ or $\rE$, 
then Theorem \ref{thm:lm} says that the map $f$ is equivariantly linearizable to 
its linear model $f_J$. 
Now the first half of the theorem follows from part $(\rii)$ of Lemma \ref{lem:rh-lm}.  
To prove the second half of the theorem, we put the idea in Remark \ref{rem:ecA0} into action. 
\par
Let $\bal = (\alpha_1, \alpha_2)$ be the multipliers of $E$.   
Using FPFs \eqref{eqn:saito2} and \eqref{eqn:tt2}, we find a rational function 
$\hat{P}(w) \in \bQ(w)$ that satisfies equation \eqref{eqn:P3}, from which 
we proceed to the function $Q(w)$ in formula \eqref{eqn:Q2}.   
We use Theorem \ref{thm:MDR} to confirm that the $(n+1)$-st power 
$\bal^{n+1} = (\alpha_1^{n+1}, \, \alpha_2^{n+1})$ is MI, or MD but NR, hence 
so is $\bal$ itself.     
It follow from Theorem \ref{thm:elinz} that $F$ is $G$-equivariantly linearizable to $J$, 
where $F$ is the function in diagram \eqref{cd:rtlm}. 
Namely, $f$ is equivariantly linearizable to its linear model $f_J$.  
We apply part $(\ri)$ of Lemma \ref{lem:rh-lm} to $f_J$ and conclude that 
$E$ falls in case $(\rR 2)$ if $\bal$ is MI and $Q(\tau) \in (-2, \, 2)$; 
in case $(\rHy)$ if $\bal$ is MI and $Q(\tau) > 2$; and in case $(\rR 1)$ if 
$\bal$ is MD but NR, respectively, where $\tau = \delta + \delta^{-1}$ is the special trace of $f$.    
This scenario works well with the K3 surface automorphisms in Table \ref{tab:Ass}, 
realizing all the three cases $(\rR 1), (\rR 2), (\rHy)$. 
Details of these contents are presented in Theorem \ref{thm:RHAt} and its proof. 
The second half of the theorem is established in this manner. \hfill $\Box$ 
\begin{example} \label{ex:rArDE}
Tables \ref{tab:Ars1} and \ref{tab:Ars2} present K3 surface automorphisms having an exceptional component of type $\rA$ 
with reflectional symmetry, while Tables \ref{tab:DEs1} and \ref{tab:DEs2} exhibit ones having exceptional components of 
type $\rD$ or $\rE$. 
By Theorem \ref{thm:rd} each of these components is contained in a (periodic) rotation domain of rank $1$.   
\end{example}
\subsection{Detecting Rotation Domains} \label{ss:drd} 
We are now in a position to detect rotation domains around an exceptional 
component of type $\rA$ with trivial symmetry, by using examples from Table \ref{tab:Ass}. 
The main result in this direction is Theorem \ref{thm:RHAt}.  
In order to deal with most entries in Table \ref{tab:Ass} in a unified manner, we begin 
by the following lemma.   
\begin{lemma} \label{lem:rd-At}
Suppose that $\cE(X)$ consists of only one exceptional component of type $\rA_n$ 
with trivial symmetry and that $f$ has no fixed point off $\cE(X)$.  
Then in case $(\rC 1)$ of Theorem $\ref{thm:ecA}$ we have 
\begin{equation} \label{eqn:C1}
\alpha^{n+1} + \alpha^{-(n+1)} = \hat{P}_n(\hat{\tau}), \qquad 
\hat{P}_n(w) := H_{n+1}(w) - w^{-1} K_n(w),      
\end{equation} 
where $\hat{\tau} := \delta^{1/2} + \delta^{-1/2}$. 
In case $(\rC 2)$ the number $\hat{\tau}$ must be a root of the algebraic equation   
\begin{equation} \label{eqn:C2C3}
w  \{ H_2(w) - 2 \} K_{i-1}(w) K_{n-i}(w) - K_n(w) = 0,  
\end{equation}
where $i \in \{ 1, \dots, n\}$ is the index appearing in the statement of case $(\rC 2)$ 
in $\S \ref{ss:linec}$.   
\end{lemma} 
{\it Proof}. 
Following Notation \ref{not:iec} we denote the unique exceptional component by $\rA_n^{\rt}$. 
Since $f$ has no fixed point off $\cE(X)$ and only one exceptional component $\rA_n^{\rt}$, 
we have $\hat{\nu}_{\roff}(f) = 0$ and FPF \eqref{eqn:tt2} reads $\hat{\tau} =\hat{\nu}_{\ron}(f) = \hat{\nu}(f, \rA_n^{\rt})$. 
In case (C1) formula \eqref{eqn:Ai} in Theorem \ref{thm:ecA} allows us to rewrite the equation $\hat{\tau} = \hat{\nu}(f, \rA_n^{\rt})$ 
in the form \eqref{eqn:C1}.  
Similarly, in case $(\rC 2)$ formula \eqref{eqn:Aii} recasts $\hat{\tau} = \hat{\nu}(f, \rA_n^{\rt})$ 
into equation \eqref{eqn:C2C3}. \hfill $\Box$ \par\medskip
Remark \ref{rem:HK} allows us to express the function $\hat{P}_n(w)$ in \eqref{eqn:C1} in the 
form \eqref{eqn:P2} with $m = n+1$, that is,  
$$
\hat{P}_n(w) = P_n(H_2(w)) \quad \mbox{when $n$ is odd}; \qquad 
\hat{P}_n(w) = w^{-1} P_n(H_2(w)) \quad \mbox{when $n$ is even}, 
$$
for some polynomial $P_n(w)$.  
Explicitly, according to the parity of $n$, we have 
$$
P_{2 j-1}(w) = H_j(w) - K_{j-1}(w), \qquad P_{2 j}(w) = (w+1) K_j(w) - (w+3) K_{j-1}(w), \qquad j \in \bN. 
$$
Similarly, the algebraic equation \eqref{eqn:C2C3} for $\hat{\tau}$ can be converted to an algebraic equation 
$R^{(n)}_i(w) = 0$ for $\tau := H_2(\hat{\tau}) = \delta + \delta^{-1}$. 
Explicit formulas for $R^{(n)}_i(w)$ exhibit four patterns depending on the parities of $n$ and $i$, 
which are somewhat too lengthy to be included here. 
Instead, some examples of $R^{(n)}_i(w)$ are given in Example \ref{ex:rd-At}. 
Due to the symmetry $R^{(n)}_i(w) = R^{(n)}_{n+1-i}(w)$, the index $i$ may vary within $1 \le i \le (n+1)/2$.   
\begin{example} \label{ex:rd-At} 
Entries in Table \ref{tab:Ass} admit exceptional components of type $\rA_n^{\rt}$ for $n = 1, 2, 6, 7, 11$. 
For those values of $n$ we illustrate the polynomials mentioned above.   
Regarding the polynomials $P_n(w)$ we find         
\begin{alignat*}{2}
P_1(w) &= w-1, \qquad P_2(w) = w^2-3, \qquad & P_6(w) &= w^4 - 5 w^2 - w + 3, \\
P_7(w) &= w^4 - w^3 -4 w^2 + 2 w + 2, \qquad & P_{11}(w) &= w^6 - w^5 - 6 w^4 + 4 w^3 + 9 w^2 -3 w - 2.  
\end{alignat*}
For $n = 1, 2, 6, 7$, the polynomials $R^{(n)}_i = R^{(n)}_i(w)$, $1 \le i \le (n+1)/2$, are calculated as  
\begin{alignat*}{4}
R^{(1)}_1 &= w-3, \qquad & R^{(2)}_1 &= \rST^{(4)}_6, \qquad & R^{(6)}_1 &= \rST^{(8)}_{47}, \qquad & R^{(6)}_2 &= \rST^{(8)}_{21}, \\
R^{(6)}_3 &= \rST^{(8)}_{16}, \qquad & R^{(7)}_1 &= \rST^{(8)}_{46}, \qquad & R^{(7)}_2 &= \rST^{(8)}_{19}, \qquad & R^{(7)}_3 &= \rST^{(8)}_{14}.   
\end{alignat*}
As for $n = 11$, $R^{(11)}_1 = w^6-2 w^5 -6 w^4 + 9 w^3+ 9 w^2 - 8 w - 2$, which is a Salem trace polynomial, and   
$$
R^{(11)}_2 = \rST^{(12)}_{95}, \quad R^{(11)}_3 = \rST^{(6)}_6 \cdot \rCT_3 \cdot \rCT_{10}, \quad  R^{(11)}_4 = \rST^{(4)}_1 \cdot \rCT_4 \cdot \rCT_9, \quad 
R^{(11)}_5 = \rST^{(8)}_{11} \cdot \rCT_8.   
$$ 
\end{example}
\begin{theorem} \label{thm:RHAt} 
For any entry of Table $\ref{tab:Ass}$ the exceptional component of type $\rA$ 
with trivial symmetry falls into one of the cases $(\rR 1), (\rR 2), (\rHy)$.  
$\mathrm{RH}$-column in the table indicates which of the three occurs for 
each individual entry.     
It tells us that all of the three cases actually occur; No.~$28$ is the only entry 
lying in case $(\rR 1)$, see also Remark $\ref{rem:RHAt}$. 
In Nos.~$7$--$11$ a rotation domain of rank $2$ and a rotation domain of rank $1$ 
coexist, where the former domain contains an exceptional component of type $\rA_1$,  
while the latter contains one of type $\rE_8$.         
\end{theorem}
{\it Proof}. 
For every automorphism $f$ in Table \ref{tab:Ass}, since $\hat{\nu}_{\roff}(f) = 0$, FPF \eqref{eqn:tt2} becomes   
$\hat{\tau} = \hat{\nu}_{\ron}(f)$.   
We apply Theorem \ref{thm:ecA} to the current situation, 
where case $(\rC 3)$ in the theorem is already ruled out by Proposition \ref{prop:Ass}. 
Case $(\rC 2)$ can also be ruled out in the following manner.    
First we consider Nos.~$6$--$11$ of Table \ref{tab:Ass}, which have dynamical degree $\lambda^{(12)}_1$ and 
exceptional set of type $\rA_1 \oplus \rE_8$. 
If $(\rC2)$ is the case, then it follows from formula \eqref{eqn:Aii} with $n = 1$ and a formula for $\rE_8$ 
in Lemma \ref{lem:DE} that the FPF $\hat{\tau} = \hat{\nu}_{\ron}(f)$ reads   
$$
\hat{\tau} = \hat{\nu}_{\ron}(f) = \hat{\nu}(f, \rA_1) + \hat{\nu}(f, \rE_8) = \dfrac{\hat{\tau}}{\tau-2} 
+ \dfrac{1}{\tau-2} \left( \dfrac{\tau^2 + \tau -1}{\hat{\tau} (\tau+1)} - \dfrac{K_3(\hat{\tau})}{K_4(\hat{\tau})} \right),  
$$ 
which leads to an algebraic equation $R(\tau) := \tau^5 - 9 \tau^3 - 9 \tau^2 + 5 \tau + 5 = 0$, 
where $R(w)$ is a Salem polynomial of degree $5$. 
Since $R(w)$ is different from $\rST^{(12)}_1(w)$, 
we have $R(\tau) \neq 0$ and hence case $(\rC2)$ cannot occur. 
Secondly we consider No.~$27$, which has dynamical degree $\lambda^{(8)}_6$ and 
exceptional set of type $\rA_1 \oplus \rA_6^{\rr} \oplus \rE_6^{\rr}$. 
If $(\rC2)$ is the case, then by some formulas in Proposition \ref{prop:Ar} and Lemma \ref{lem:DE}  
the FPF $\hat{\tau} = \hat{\nu}_{\ron}(f)$ reads   
$$
\hat{\tau} = \hat{\nu}_{\ron}(f) = \hat{\nu}(f, \rA_1) + \hat{\nu}(f, \rA_6^{\rr}) + \hat{\nu}(f, \rE_6^{\rr}) 
= \dfrac{\hat{\tau}}{\tau-2} + \dfrac{1}{\hat{\tau}} + \dfrac{\tau+1}{\hat{\tau} \cdot \tau},    
$$ 
which leads to an algebraic equation $R(\tau) := \tau^3 - 3 \tau^2 - 3 \tau + 2 = 0$, 
where $R(w)$ is a Salem polynomial of degree $3$. 
Since $R(w)$ is different from $\rST^{(8)}_6(w)$, we have $R(\tau) \neq 0$ and hence case $(\rC2)$ cannot occur. 
Thirdly, for all entries of Table \ref{tab:Ass} other than Nos.~$6$--$11$ and No.~$27$, 
with explicit formulas for $R^{(n)}_i(w)$ in Example \ref{ex:rd-At}, similar reasonings to those 
mentioned above show that $(\rC 2)$ is not the case at all in Theorem \ref{thm:ecA}.    
\par
Now assertion $(\rC 1)$ of Theorem \ref{thm:ecA} applies to every $\rA_n^{\rt}$-component in Table \ref{tab:Ass}. 
Let $\delta^{1/2} \alpha^{\pm 1}$ be its multipliers. 
In Nos.~$6$--$11$, formula \eqref{eqn:Ai} and a formula for $\rE_8$ in 
Lemma \ref{lem:DE} imply that FPF $\hat{\tau} = \hat{\nu}_{\ron}(f)$ becomes 
$$
\hat{\tau} = \hat{\nu}_{\ron}(f) = \hat{\nu}(f, \rA_1) + \hat{\nu}(f, \rE_8) =
\dfrac{\hat{\tau}}{\tau - (\alpha^2 + \alpha^{-2} )} + 
\dfrac{1}{\tau-2} \left( \dfrac{\tau^2 + \tau -1}{\hat{\tau} (\tau+1)} - \dfrac{K_3(\hat{\tau})}{K_4(\hat{\tau})} \right).   
$$
This equation is settled as $\alpha^2 + \alpha^{-2} = \hat{P}(\hat{\tau})$, 
where $\hat{P}(w) = P(H_2(w))$ with $P(w)$ given by    
$$
P(w) = \dfrac{w^6 - 7 w^4 - w^3 + 13 w^2 + 3 w -4}{w^5 + w^4 - 5 w^3 - 5 w^2 + 4 w + 3}.  
$$
In No.~$27$, formula \eqref{eqn:Ai} and formulas for $\rA_6^{\rr}$ and $\rE_6^{\rr}$ in \S \ref{ss:oec} imply that 
FPF $\hat{\tau} = \hat{\nu}_{\ron}(f)$ becomes 
$$
\hat{\tau} = \hat{\nu}_{\ron}(f) = \hat{\nu}(f, \rA_1) + \hat{\nu}(f, \rA_6^{\rr}) + \hat{\nu}(f, \rE_6^{\rr}) 
= \dfrac{\hat{\tau}}{\tau-(\alpha^2 + \alpha^{-2})} + \dfrac{1}{\hat{\tau}} + \dfrac{\tau+1}{\hat{\tau} \cdot \tau}.    
$$
This equation is settled as $\alpha^2 + \alpha^{-2} = \hat{P}(\hat{\tau})$, 
where $\hat{P}(w) = P(H_2(w))$ with $P(w)$ given by    
$$
P(w) = \dfrac{w(w^2-w-3)}{(w-1)(w+1)}. 
$$
For each entry other than Nos.~$6$--$11$ and No.~$27$, formula 
\eqref{eqn:C1} is available with the associated polynomial $P(w) = P_n(w)$ in Example \ref{ex:rd-At}, 
depending on the value $n = 1, 2, 6,7, 11$ of the $\rA_n^{\rt}$-component.   
\par
\begin{longtable}{llccccccccccc}
\caption{\setlength{\leftskip}{-16mm}\setlength{\rightskip}{-16mm} 
Examinations of the entries in Table \ref{tab:Ass}. } \label{tab:rAss} 
\\
\hline
      &           &                &     &              &              &             &              &              &              &              &            & \\[-3mm]
No. & Dynkin & $\lambda$ & AI & $\tau_0$ & $\tau_1$ & $\tau_2$ & $\tau_3$ & $\tau_4$ & $\tau_5$ & $\tau_6$ & $\tau_7$ & MID \\ 
\hline
\endfirsthead
\multicolumn{13}{l}{continued} \\
\hline
      &           &                &     &              &              &              &              &              &             &              &            &   \\[-3mm]
No. & Dynkin & $\lambda$ & AI & $\tau_0$ & $\tau_1$ & $\tau_2$ & $\tau_3$ & $\tau_4$ & $\tau_5$ & $\tau_6$ & $\tau_7$ & MID \\ 
\hline
\endhead
\hline
\endfoot
\hline
\\[-4pt]
\caption*{\setlength{\leftskip}{-20mm}\setlength{\rightskip}{-20mm}
The first column refers to the entry numbers in Table \ref{tab:Ass}; the second column exhibits 
the Dynkin type of the exceptional set $\cE(X)$; the third one gives the dynamical degree $\lambda$ of $f$; 
for the remaining we refer to the second to the last paragraph in the proof of Theorem \ref{thm:RHAt}. }
\endlastfoot
$1$--$5$ & $\rA_6^{\rt}$ & $\lambda^{(16)}_{2}$ & Y & Y & N & Y & Y & Y & Y & N & N & MI \\
$6$--$11$ & $\rA_1^{\rt} \oplus \rE_8^{\rt}$ & $\lambda^{(12)}_{1}$ & Y & Y & Y & Y & N & Y & Y &  & & MI \\
$12$--$18$ & $\rA_1^{\rt}$ & $\lambda^{(12)}_{1}$ & Y & Y & Y & Y & Y & Y & N &  & & MI \\
$19$--$23$ & $\rA_2^{\rt}$ & $\lambda^{(10)}_{1}$ & Y & Y & Y & N & Y & Y & & & & MI \\
$24$ & $\rA_7^{\rt}$ & $\lambda^{(8)}_{1}$ & Y & Y & N & Y & Y &  &  &  & & MI \\ 
$25$, $26$ & $\rA_1^{\rt}$ & $\lambda^{(8)}_{2}$ & Y & Y & Y & Y & N &  &  &  & & MI \\
$27$ & $\rA_1^{\rt} \oplus \rA_6^{\rr} \oplus \rE_6^{\rr}$ & $\lambda^{(8)}_{6}$ & Y & Y & N & Y & Y &  &  &  & & MI \\
$28$ & $\rA_2^{\rt}$ & $\lambda^{(8)}_{8}$ & Y & Y & Y & Y & Y &  &  &  & & MD \\
$29$ & $\rA_{11}^{\rt}$ & $\lambda^{(8)}_{8}$ & Y & Y & N & N & Y &  &  &  & & MI \\
$30$ & $\rA_6^{\rt}$ & $\lambda^{(8)}_{10}$ & Y & Y & N & Y & Y &  &  &  & & MI \\
$31$ & $\rA_6^{\rt}$ & $\lambda^{(6)}_{1}$ & Y & Y & N & Y &  &  &  &  & & MI \\
\hline
\end{longtable}
As in definition \eqref{eqn:Q2} of \S \ref{sec:pp} we set $Q(w) := P^2(w)-2$ if $n+1$ is even, 
and $Q(w) := (w+2)^{-1} P^2(w)-2$ if $n+1$ is odd. 
Putting necessary information from Table \ref{tab:Ass} into the function $Q(w)$, 
we have the results in Table \ref{tab:rAss}. 
Here AI-column concerns the question of whether $Q(\tau)$, with $\tau$ being the special trace of $f$, 
is an algebraic integer (AI for short) or not; the answer is yes (Y) for every entry. 
The $\lambda$-column gives the dynamical degree of $f$, which is a Salem number $\lambda^{(d)}_i$.   
Let $\tau_0, \dots, \tau_{d/2-1}$ be the roots of the Salem trace polynomial $\rST^{(d)}_i(w)$ arranged as 
$\tau_0 > 2 > \tau_1 > \cdots > \tau_{d/2-1} > -2$.      
The $\tau_j$-column gives yes (Y) or no (N) answer as to whether $Q(\tau_j) \in [-2, \, 2]$ or not, i.e. $Q(\tau_j) > 2$.     
Then the results in MID-column follow from Theorem \ref{thm:MDR} and Corollary \ref{cor:MDR},  
where MI (resp. MD) means that $\bal^{n+1} = (\alpha_1^{n+1}, \, \alpha_2^{n+1})$ and hence 
$\bal = (\alpha_1, \alpha_2) = (\delta^{1/2} \alpha, \, \delta^{1/2} \alpha^{-1})$ itself is MI (resp. MD). 
In particular, No.~$28$ is the only entry that falls in the MD case, for which 
$$
Q(w) = (w+2)^{-1} (w^2-3)^2 -2, \qquad \alpha^6 + \alpha^{-6} = Q(\tau) = -1 \quad 
\mbox{for every root $\tau$ of $\rST^{(8)}_8(w)$}, 
$$ 
and hence $\alpha$ is an $18$-th root of unity.    
Now the results in RH-column of Table \ref{tab:Ass} follow from Theorem \ref{thm:elinz}, Proposition \ref{prop:eltlm}, 
Lemma \ref{lem:rh-lm}, Table \ref{tab:rAss} and ST-column in Table \ref{tab:Ass}. 
\par
An inspection of Table \ref{tab:Ass} shows that in Nos.~$7$--$11$ the $\rA_1$-component is containd in 
a rotation domain of rank $2$. 
For these entries the $\rE_8$-component is contained in a rotation domain of rank $1$ by Theorem \ref{thm:rd}. 
Hence the coexistence assertion in the theorem follows.  \hfill $\Box$ 
\begin{remark} \label{rem:RHAt} 
The automorphism $f$ in No.~$28$ of Table \ref{tab:Ass} has dynamical degree 
$\lambda^{(8)}_8 = (\lambda^{(8)}_1)^2$ and exceptional set $E$ of type $\rA_2$ with trivial symmetry. 
We know from Theorem \ref{thm:RHAt} that $(f, E)$ lies in case $(\rR 1)$.   
We show that it is {\sl genuinely} of case $(\rR 1)$ in the sense of Remark \ref{rem:genuine}. 
Otherwise, $f$ is the square of an automorphism $g$ of dynamical degree $\lambda^{(8)}_1$  
having exceptional set $E$ with reflectional symmetry. 
We find from Table \ref{tab:Ass} that $f^*|_{H^2(X)}$ has characteristic polynomial 
$\tilde{\varphi}(z) = \rS^{(8)}_8(z) \cdot \rC_1^2(z) \cdot \rC_{42}(z)$, for which   
$\tilde{\varphi}(z^2) = \rS^{(8)}_1(z) \cdot \rS^{(8)}_1(-z) \cdot \rC_1^2(z) \cdot \rC_2^2(z) \cdot \rC_{84}(z)$.   
If $f = g^2$ then $\tilde{\varphi}(z^2)$ is divisible by the characteristic polynomial of $g^*|_{H^2(X)}$,  
which takes the form $S^{(8)}_1(z) \cdot C(z)$, where $C(z)$ is a product of cyclotomic polynomials with $\deg C(z) = 14$.  
Then $C(z)$ must divide $\rC_1^2(z) \cdot \rC_2^2(z) \cdot \rC_{84}(z)$, but this is impossible because 
$\rC_{84}(z)$ is an irreducible polynomial of degree $24$.     
\end{remark}
\par
Next we illustrate the coexistence of rotation domains of ranks $1$ and $2$  
by using Tables \ref{tab:DEs1} and \ref{tab:DEs2}.  
\begin{theorem} \label{thm:rk1-2}
For each entry of Tables $\ref{tab:DEs1}$ and $\ref{tab:DEs2}$, every exceptional 
component of type $\rD$ or $\rE$ is contained in a rotation domain of rank $1$.  
Moreover, the unique fixed point off the exceptional set $\cE(X)$ is either a center 
of a rotation domain of rank $2$ or a hyperbolic fixed point, as indicated by $\rR 2$ 
or $\rHy$ in $\mathrm{RH}$-column of the tables. 
In the former case a rotation domain of rank $1$ and a rotation domain of rank $2$ 
coexist in $X$.   
\end{theorem}
{\it Proof}. 
The first part of the theorem is a consequence of Example \ref{ex:rArDE}. 
For the second part we recall from Proposition \ref{prop:DEs} that for each entry of Table \ref{tab:DEs1} 
the map $f : X \to X$ has a unique fixed point $p$ off $\cE(X)$, which is simple. 
If its multipliers are $\delta^{1/2} \alpha^{\pm 1}$, then 
$\hat{\nu}_{\roff}(f) = \hat{\nu}(f, p) = \{ \hat{\tau} - (\alpha + \alpha^{-1})\}^{-1}$ 
and FPF \eqref{eqn:tt2} gives  
\begin{equation} \label{eqn:rk1-2}
\hat{\tau} = \hat{\nu}_{\ron}(f) + \hat{\nu}_{\roff}(f) = \hat{\nu}_{\ron}(f) + \dfrac{1}{\hat{\tau} - (\alpha + \alpha^{-1})}. 
\end{equation}
The Dynkin types of the $f$-invariant exceptional components (not of the entire exceptional set $\cE(X)$) 
can be read off from Proposition \ref{prop:DEs} and are shown in Dynkin column of Table \ref{tab:rk1-2}. 
This information together with Lemma \ref{lem:DE} enables us to calculate $\hat{\nu}_{\ron}(f)$ explicitly. 
Then equation \eqref{eqn:rk1-2} is settled as $\alpha + \alpha^{-1} = \hat{P}(\hat{\tau}) = P(\tau)/\hat{\tau}$ 
for some rational function $P(w) \in \bQ(w)$; cf. formulas \eqref{eqn:P1} and \eqref{eqn:P2} in \S \ref{sec:pp}. 
Here the explicit formulas for $P(w)$ are given in Table \ref{tab:rk1-2}. 
We put $Q(w) := (w+2)^{-1} P^2(w) - 2$ as in definition \eqref{eqn:Q2}. 
With this function we can make a table similar to Table \ref{tab:rAss}, which together with Corollary \ref{cor:MDR} 
tells us that the multiplier pair  $\bal = (\alpha_1, \alpha_2)$ is algebraic and MI for every entry in Table \ref{tab:DEs1}. 
Thus, by Theorem \ref{thm:Poschel} and \ref{thm:tnt}, in a neighborhood of $p$ the map $f$ is biholomorphically conjugate to 
the linear model \eqref{eqn:rh-splest} with $(c_1, c_2) = (\alpha_1, \alpha_2)$. 
Now the results in RH-column of Table \ref{tab:DEs1} follow from  
Lemma \ref{lem:rh-splest}; we have $\rR 2$ if $Q(\tau_j) \in (-2, \, 2)$  
and $\rHy$ if $Q(\tau_j) > 2$, when the special trace $\tau$ of $f$ is $\tau_j$.  
The automorphism in Table \ref{tab:DEs2} can be treated in a similar manner.     
\hfill $\Box$ 
\begin{longtable}{lll}
\caption{Explicit formulas for $P(w)$ in the proof of Theorem \ref{thm:rk1-2}. } \label{tab:rk1-2} 
\\
\hline
      &     &                 \\[-3mm]
No.  & Dynkin  & $P(w)$      \\ 
\hline
\endfirsthead
\multicolumn{3}{l}{continued} \\
\hline
      &    &                  \\[-3mm]
No. & Dynkin &  $P(w)$      \\ 
\hline
\endhead
\hline
\endfoot
\hline
\\[-4pt]
\caption*{\setlength{\leftskip}{-24mm}\setlength{\rightskip}{-24mm}
The first column refers to the entry numbers in Table \ref{tab:DEs1}; the second column exhibits 
the Dynkin type of the $f$-invariant exceptional components in $X$; the third column gives the formula for $P(w)$. }
\endlastfoot
$1$--$7$ & $\rD_6^{\rr}$ & $\dfrac{(w+2)(w^3-3 w-1)}{(w+1)(w^2-3)}$ \\[2mm]
$8$--$12$ & $\rD_4^{\rc}$ & $\dfrac{w^2+w-1}{w}$ \\[2mm]
$13$--$30$ & $\rE_8^{\rt}$ & $\dfrac{w^6+2w^5-6w^4-17w^3-4w^2+9w+2}{w^5+w^4-6w^3-9w^2+4w+3}$ \\[2mm]
$31$--$38$ & $\rD_9^{\rr}$ & $\dfrac{(w+2)(w^4-w^3-3w^2+w+1)}{(w-2)(w+1)(w^2+w-1)}$ \\[2mm]
$39$--$47$, $62$, $63$, $67$ & $\emptyset$ & $w+1$ \\[2mm] 
$48$--$53$ & $\rD_{13}^{\rr}$ & $\dfrac{(w+1)(w+2)(w^5-2w^4-3w^3+6w^2-1)}{(w^2-2)(w^4-4w^2-w+1)}$ \\[2mm]
$54$--$56$ & $\rD_{14}^{\rr}$ & $\dfrac{(w+2)(w^7-7w^5-w^4+13w^3+3w^2-5w-1)}{(w^2+w-1)(w^5-6w^3-w^2+8w+3)}$ \\[2mm]
$57$--$60$ & $\rE_7^{\rt}$ & $\dfrac{(w+2)(w^4-w^3-3w^2+w+1)}{w^4-4w^2-w+1}$ \\[2mm]
$61$ & $\rD_5^{\rr} \oplus \rE_8^{\rt}$ & $\dfrac{(w+2)(w^6-2w^5-5w^4+8w^3+7w^2-6w-1)}{(w^2-3)(w^4-w^3-3w^2+w+1)}$ \\[2mm]
$64$, $65$ & $\rD_{16}^{\rt}$ & $\dfrac{(w+2)(w^8-2w^7-6w^6+11w^5+11w^4-16w^3-7w^2+5w+1)}{(w^3-3w-1)(w^5-w^4-5w^3+4w^2+5w-3)}$ \\[2mm]
$66$ & $\rE_6^{\rr}$ & $\dfrac{(w-1)(w+1)(w+2)}{w^2+w-1}$ \\[1mm]
\hline
\end{longtable}
\section{Periodic Cycles} \label{sec:pc} 
Under conditions  (A1) and (A2) in \S \ref{sec:intro}, the K3 surface $X$ contains  
only a finite number of exceptional components. 
Moreover, any automorphism of a Dynkin diagram of type $\rA$, $\rD$ or $\rE$ 
has order $1$, $2$ or $3$ (not respectively).   
So there exists the smallest number $m_0 \in \bN$ such that for any exceptional 
component $E \subset \cE(X)$ the map $f^{m_0}$ preserves $E$ and induces 
the trivial Dynkin automorphism $\sigma(f^{m_0}, E) = 1$.     
Once we know how $f$ acts on $\cE(X)$, it is relatively easy to detect the way in which 
all iterates $f^m$, $m \in \bN$, act on $\cE(X)$.   
The pattern of these actions is periodic in $m$ with period $m_0$. 
In particular, in a generic situation, it is not hard to see  how the numbers $\mu_{\ron}(f^m)$ 
and $\hat{\nu}_{\ron}(f^m)$ vary as $m$ moves (Remark \ref{rem:DE}), 
where $\mu_{\ron}(f)$ and $\hat{\nu}_{\ron}(f)$ are defined in \eqref{eqn:muof} and \eqref{eqn:nuhat2}.    
\par
A {\sl periodic cycle} of {\sl primitive} period $m$ is an $f$-orbit $C = \{p, f(p), \dots, f^{m-1}(p) \}$ 
consisting of distinct $m$ points such that $f^m(p) = p$.    
Since $\cE(X)$ is $f$-invariant, we have either $C \subset \cE(X)$ or $C \subset X \setminus \cE(X)$, 
that is to say, $C$ is either {\sl on} or {\sl off} $\cE(X)$. 
By definition, the {\sl multiplicity} of $C$ is the local index $\mu(f^m, p)$ at a representative point $p \in C$     
and the {\sl multipliers} of $C$ are those of $f^m$ at the point $p$.  
They are well defined, because they do not depend on the choice of $p$.   
We say that $C$ is {\sl simple} if its multiplicity is $1$, or equivalently if neither of its multipliers is $1$. 
We say that $C$ is non-resonant (NR) if its multipliers are NR. 
\par
Given any $k \in \bN$, we consider the following two situations and introduce a condition $(\mathrm{PC})_k$ 
as follows.  
\begin{enumerate}
\setlength{\itemsep}{-1pt} 
\item[$(k)_0$] There is {\sl no} periodic cycle of primitive period $k$ off $\cE(X)$. 
\item[$(k)_1$] There is {\sl exactly one} periodic cycle $C_k$ of primitive period $k$ off $\cE(X)$.  
In addition $C_k$ is {\sl simple}.   
\end{enumerate} 
\par\medskip\noindent 
{\bf Condition} $\mbox{\boldmath $(\mathrm{PC})_k$}$ 
Either $(k)_0$ or $(k)_1$ should occur and in the latter case $C_k$ should be NR.         
\par\medskip
Let $\{ \bar{\mu}(k) \}_{k \in \bN}$ be the unique sequence of integers defined recursively by  
\begin{equation} \label{eqn:lambda}
\mu_{\roff}(f^k) = \sum_{j \div k} \bar{\mu}(j), \qquad \mbox{or directly by} \qquad 
\bar{\mu}(k) = \sum_{j \div k} \Mob(k/j) \cdot \mu_{\roff}(f^j),     
\end{equation}
where the sums extend over all (positive) divisors $j$ of $k$ and $\Mob$ is the M\"{o}bius 
function in elementary arithmetic.  
We say that a divisor $j$ of $k$ is {\sl proper} if $j < k$, in which case 
we write $j \pd k$.  
In particular $1$ has no proper divisor. 
The multipliers of $C_k$ can and shall be expressed as $\delta^{k/2} \alpha_k^{\pm 1}$ 
for some number $\alpha_k \in \bC^{\times}$.      
\begin{lemma} \label{lem:pc} 
Let $m \in \bN$ and suppose that conditions $(\mathrm{PC})_k$ hold for all proper divisors $k$ of $m$. 
Then,      
\begin{enumerate} 
\setlength{\itemsep}{-1pt} 
\item for each proper divisor $k$ of $m$ we have either $\bar{\mu}(k) = 0$ or $\bar{\mu}(k) = k$,        
\item for $\ve = 0, 1$, the case $(m)_{\ve}$ occurs if and only if $\bar{\mu}(m) = \ve m$.   
\end{enumerate}
If either case $(m)_0$ or $(m)_1$ occurs, then FPF \eqref{eqn:tt2} for the map $f^m$ takes the form   
\begin{equation} \label{eqn:tt-pc}
H_m( \hat{\tau} ) = \hat{\nu}_{\ron} (f^m) + 
\sum_{k \pd m} \dfrac{ \bar{\mu}(k)}{H_m(\hat{\tau}) - H_{m/k} (\alpha_k + \alpha_k^{-1}) } 
+ \dfrac{\bar{\mu}(m)}{H_m(\hat{\tau}) - (\alpha_m + \alpha_m^{-1}) }.    
\end{equation}
\end{lemma}
{\it Proof}. 
For each $j \in \bN$, let $\cP_j$ be the set of all $f$-periodic points of primitive period $j$ off $\cE(X)$.  
For each $(k, j) \in \bN^2$ with $j \div k$, let $\bar{\mu}(k, j)$ denote the sum of local indices 
$\mu(f^k, p)$ over all points $p \in \cP_j$. 
Then we have    
\begin{equation} \label{eqn:mubar}
\mu_{\roff}(f^k) = \sum_{j \pd k} \bar{\mu}(k, j) + \bar{\mu}(k, k) \qquad \mbox{for every} \quad k \in \bN.     
\end{equation}
\par 
{\bf Claim 1}. For any $j \in \bN$ and $\ve = 0, 1$, the case $(j)_{\ve}$ occurs if and only if $\bar{\mu}(j, j) = \ve j$. 
\par\medskip
Notice that $\bar{\mu}(j, j)$ is $j$ times the cardinality, counted with multiplicity, of 
all periodic cycles of primitive period $j$ off $\cE(X)$. 
Claim 1 is clear from this characterization of the number $\bar{\mu}(j, j)$ and the definition of $(j)_{\ve}$.  
\par\medskip  
{\bf Claim 2}. For every divisor $k$ of $m$ we have $\bar{\mu}(k) = \bar{\mu}(k, k)$.  
\par\medskip
Let $k \div m$ and $j \pd k$ so that $j \pd m$.    
From $(\mathrm{PC})_j$ and Claim 1 we have either $\bar{\mu}(j, j) = 0$ or $\bar{\mu}(j, j) = j$.   
In the former case $\cP_j$ is empty, so $\bar{\mu}(k, j) = 0 = \bar{\mu}(j, j)$.  
In the latter case $\cP_j = C_j$ and, as $C_j$ is NR by assumption $(\mathrm{PC})_j$, neither of the 
multipliers of $C_j$ is a root of unity.  
Thus $\mu(f^k, p) = \mu((f^j)^{k/j}, p) = 1$ for all $p \in C_j$, and so $\bar{\mu}(k, j) = j = \bar{\mu}(j, j)$.    
In either case we have $\bar{\mu}(k, j) = \bar{\mu}(j, j)$ for every $j \pd k$. 
Thus equation \eqref{eqn:mubar} becomes  
$$
\mu_{\roff}(f^k) = \sum_{j \pd k} \bar{\mu}(j, j) + \bar{\mu}(k, k) = \sum_{j \div k} \bar{\mu}(j, j) \qquad \mbox{for every $k \div m$}. 
$$
This together with definition \eqref{eqn:lambda} shows that the two sequences $\{ \bar{\mu}(k) \}_{k | m}$ and 
$\{ \bar{\mu}(k, k) \}_{k | m}$ satisfy the same recurrence relation, so they must be the same sequence.  
This proves Claim 2. 
\par
Assertions (1) and (2) readily follow from Claims $1$ and $2$.  
Suppose that $(m)_{\ve}$ is the case. 
Then FPF \eqref{eqn:tt2} for the map $f^m$ reads $H_m(\hat{\tau}) = \hat{\nu}_{\ron}(f^m) + \hat{\nu}_{\roff}(f^m)$, 
where $\hat{\mu}_{\roff}(f^m)$ is given by the middle and last terms in the right-hand side of \eqref{eqn:tt-pc}, 
which are the contributions of the periodic cycles $C_k$ for $k \div m$ with $\bar{\mu}(k) = k$.  
Here the denominators $H_m(\hat{\tau}) -H_{m/k}(\alpha_k + \alpha_k^{-1})$ are nonzero, 
since $C_k$ is NR for $k \pd m$ and simple for $k = m$.   
\hfill $\Box$ 
\begin{induction} \label{induction} 
Let $m \in \bN$. 
Suppose that for all proper divisors $k$ of $m$ the conditions $(\mathrm{PC})_k$ are satisfied 
and the numbers $\alpha_k + \alpha_k^{-1}$ with $\bar{\mu}(k) = k$ are known. 
With this hypothesis, we try to deduce condition $(\mathrm{PC})_m$ and find the value 
of $\alpha_m + \alpha_m^{-1}$ in case $\bar{\mu}(m) = m$.  
In view of Lemma \ref{lem:pc} the first step is to confirm that    
\begin{equation} \label{eqn:test} 
\mbox{$\bar{\mu}(k) = 0$ or $\bar{\mu}(k) = k$ for every divisor $k$ of $m$}. 
\end{equation}
\begin{itemize}
\setlength{\itemsep}{-1pt} 
\item If $\bar{\mu}(m) = 0$, we are done as Lemma \ref{lem:pc}.(2) shows that we are in case $(m)_0$.  
Confirmation of equation \eqref{eqn:tt-pc} supports the correctness of the preceding 
calculations of $\alpha_k + \alpha_k^{-1}$ for proper divisors $k$ of $m$.   
\item If $\bar{\mu}(m) = m$, we are in case $(m)_1$ and led to 
the next step of solving equation \eqref{eqn:tt-pc} in the form 
$$
\alpha_m + \alpha_m^{-1} = \hat{P}_m(\hat{\tau}), 
$$
where $\hat{P}_m(w) = P_m(H_2(w))$ for even $m$ and $\hat{P}_m(w) = w^{-1} P_m(H_2(w))$ for odd $m$  
with some function $P_m(w)$ as in formulas \eqref{eqn:P1} and \eqref{eqn:P2} in Hypothesis \ref{hyp:P} 
and Remark \ref{rem:Q}.  
Theorem \ref{thm:MDR} is then used to examine whether the multipliers    
$\delta^{m/2} \alpha_m^{\pm 1}$ of the periodic cycle $C_m$ are NR. 
If so, our trial is successful. 
\end{itemize}     
\end{induction}
\par
There is a recipe to carry out the test \eqref{eqn:test} in the context of 
the method of hypergeometric groups.    
\begin{recipe} \label{recipe} 
Let $f : X \to X$ be a K3 surface automorphism arising from the modified 
Hodge isometry $\widetilde{A}$ for a hypergeometric group 
$H = \langle A, B \rangle$ as described in \S \ref{sec:mhgg}.     
Then our recipe proceeds as follows.  
\begin{enumerate}
\setlength{\itemsep}{-1pt} 
\item The Lefschetz numbers $L(f^m) = 2 + \Tr((f^m)^*|{H^2(X)}) = 2 + \Tr(\widetilde{A}^m)$ 
for $m \in \bN$ are given by  
\begin{equation} \label{eqn:lnm}
L(f^m) = 2+ \mbox{the coefficient of $z^m$ in the Maclaurin expansion of 
$- z \dfrac{d}{d z} \log \widetilde{\varphi}(z)$},  
\end{equation} 
where $\widetilde{\varphi}(z) = S(z) \cdot \widetilde{\varphi}_1(z)$ is the characteristic 
polynomial of $\widetilde{A}$; see \cite[formula (41)]{IT2}.  
\item Calculate the numbers $\mu_{\ron}(f^m)$, $m \in \bN$, following the recipe 
in Remark \ref{rem:DE}, where $\mu_{\ron}(f)$ is defined by formula \eqref{eqn:muof} in \S \ref{sec:fpf}.   
When an exceptional component of type $\rA$ with trivial symmetry is involved, 
this step is rather intricate but manageable, once case $(\rC 3)$ is ruled out in Theorem \ref{thm:ecA}.  
\item FPF \eqref{eqn:saito2} for the $m$-th iterate $f^m$ yields   
$\mu_{\roff}(f^m) = L(f^m) - \mu_{\ron}(f^m)$. 
\item Use formula \eqref{eqn:lambda} to produce the sequence $\{ \bar{\mu}(k) \}_{k|m}$ from  
$\{ \mu_{\roff}(f^k) \}_{k | m}$ and check the condition \eqref{eqn:test}.   
\end{enumerate}
\end{recipe}
\begin{example} \label{ex:npc} 
To illustrate how Recipre \ref{recipe} works, we apply it to the automorphisms $f$ in 
Nos.~$11$--$15$ of Table \ref{tab:Ars1}, which have exceptional set of type $\rA_2^{\rr}$.   
Step (1) of Recipe \ref{recipe} with polynomial 
$\widetilde{\varphi} = \rS^{(8)}_1\cdot \rC_1 \rC_2 \rC_{28}$ 
yields the second row of Table \ref{tab:npc}. 
In step (2) we have $\mu_{\ron}(f^m) = 2 + (-1)^m$ for $m \in \bN$,  
which follow from formula \eqref{eqn:mnAr} in Proposition \ref{prop:Ar}.  
Step (3) amounts to subtracting the third row from the second one in Table \ref{tab:npc}. 
Step (4) is just to pass from the fourth row to the fifth one via formula \eqref{eqn:lambda}.      
The last row shows that the test \eqref{eqn:test} is successful for all values of $m$ up to 
$14$ except for $m = 13$.   
Thus it is worth trying Induction \ref{induction} for those values of $m$.  
It turns out that this inductive procedure is successful, leading to Theorem \ref{thm:rdpc} below. 
\begin{longtable}{lccccccccccccccc}
\caption{ 
Recipe \ref{recipe} applied to Nos.~$11$--$15$ in Table \ref{tab:Ars1}. } \label{tab:npc} 
\\
\hline
      &       &       &       &       &       &       &       &       &       &        &         &         &        &        &        \\[-3mm]
$m$ & $1$ & $2$ & $3$ & $4$ & $5$ & $6$ & $7$ & $8$ & $9$ & $10$ & $11$ & $12$ & $13$ & $14$ & $15$ \\ 
\hline
\endfirsthead
\multicolumn{16}{l}{continued} \\
\hline
      &       &       &      &       &       &       &       &       &       &         &        &         &        &        &         \\[-3mm]
$m$ & $1$ & $2$ & $3$ & $4$ & $5$ & $6$ & $7$ & $8$ & $9$ & $10$ & $11$ & $12$ & $13$ & $14$ & $15$ \\ 
\hline
\endhead
\hline
\endfoot
\hline
\\[-4pt]
\caption*{\setlength{\leftskip}{-5mm} \setlength{\rightskip}{-5mm} The bottom row shows whether the test \eqref{eqn:test} is succesful (S) or failed (F).  }
\endlastfoot
$L(f^m)$ & $2$ & $6$ & $5$ & $6$ & $7$ & $9$ & $9$ & $6$ & $14$ & $21$ & $13$ & $21$ & $28$ & $27$ & $40$ \\
$\mu_{\ron}(f^m)$ & $1$ & $3$ & $1$ & $3$ & $1$ & $3$ & $1$ & $3$ & $1$ & $3$ & $1$ & $3$ & $1$ & $3$ & $1$ \\
$\mu_{\roff}(f^m)$ & $1$ & $3$ & $4$ & $3$ & $6$ & $6$ & $8$ & $3$ & $13$ & $18$ & $12$ & $18$ & $27$ & $24$ & $39$ \\
$\bar{\mu}(m)$ & $1$ & $2$ & $3$ & $0$ & $5$ & $0$ & $7$ & $0$ & $9$ & $10$ & $11$ & $12$ & $26$ & $14$ & $30$ \\
S/F & S & S & S & S & S & S & S & S & S & S & S & S & F & S & F \\
\hline
\end{longtable}
\end{example}
\par
Maps in Example \ref{ex:npc} are simpler to the effect that they contain 
no exceptional component of type $\rA$ with trivial symmetry. 
In Appendix \ref{app:exp} we illustrate Recipe \ref{recipe} by using examples containing such a component.  
\begin{example} \label{ex:recipe} 
For a lot of automorphisms, Recipe \ref{recipe} can be used to confirm that condition  
\eqref{eqn:test} is satisfied by a good many primitive periods $m$.  
We present some examples of them from Tables \ref{tab:Ass}, \ref{tab:Ars1} and \ref{tab:DEs1}.  
\begin{itemize}
\setlength{\itemsep}{-1pt}
\item Table \ref{tab:Ass}: Nos.~$1$--$5$, all $m \le 10$; Nos.~$6$--$11$, all $m \le 14$; Nos.~$12$--$18$, all $m \le 15$; Nos.~$19$--$23$, all $m \le 22$; 
No.~$24$, all $m \le 12$. 
\item Table \ref{tab:Ars1}: Nos.~$3$, $4$, all $m \le 10$; Nos.~$5$--$10$, all $m \le 22$; Nos.~$11$--$15$, all $m \le 14$ but $m = 13$; 
No.~$16$, all $m \le 9$. 
\item Table \ref{tab:DEs1}: Nos.~$1$--$7$, all $m \le 10$; Nos.~$8$--$12$, all $m \le 10$; Nos.~$13$--$30$, all $m \le 20$ but $m = 19$; 
Nos.~$31$--$38$, all $m \le 14$; Nos.~$39$--$47$, all $m \le 22$; Nos.~$48$--$53$, all $m \le 9$.    
\end{itemize}  
\end{example}
\par
If Induction \ref{induction} goes well until reaching $(\mathrm{PC})_m$ with case $(m)_1$, 
then each point $p$ of the periodic cycle $C_m$ is a NR fixed point of $f^m$ with algebraic 
multipliers, hence $f^m$ is linearizable near $p$ to its Jacobian matrix by Theorems \ref{thm:Poschel} and \ref{thm:tnt}.  
Then Lemma \ref{lem:rh-splest} is used to examine whether $p$ is a center of a (periodic) rotation domain, 
or a hyperbolic periodic point. 
We illustrate this process by continuing with Example \ref{ex:npc}. 
Note that the special trace of $f$ is among $\tau_1, \tau_2, \tau_3$, where 
$\tau_0 > \tau_1 > \tau_2 > \tau_3$ are the roots of $\rST^{(8)}_1(w)$ with $\tau_0 > 2$.  
\begin{longtable}{llllllllllll}
\caption{ \setlength{\leftskip}{-10mm}\setlength{\rightskip}{-10mm} 
Periodic cycles and rotation domains for Nos.~$11$--$15$ in Table \ref{tab:Ars1}. } \label{tab:rdpc} 
\\
\hline
                      &      &      &       &       &       &       &       &         &        &         &          \\[-3mm]
prim. period & $m$ & $1$ & $2$ & $3$ & $5$ & $7$ & $9$ & $10$ & $11$ & $12$ & $14$   \\ 
\hline
\endfirsthead
\multicolumn{12}{l}{continued} \\
\hline
                      &       &      &       &       &       &       &       &        &         &         &         \\[-3mm]
prim. period & $m$ & $1$ & $2$ & $3$ & $5$ & $7$ & $9$ & $10$ & $11$ & $12$ & $14$  \\ 
\hline
\endhead
\hline
\endfoot
\hline
\\[-4pt]
\caption*{\setlength{\leftskip}{-15mm} \setlength{\rightskip}{-15mm} $\rR 1$ and $\rR 2$ indicate that  
points in the periodic cycle are centers of rotation domains of ranks $1$ and $2$ respectively, while $\rHy$ says that 
they are hyperbolic periodic points. }
\endlastfoot
special trace & $\tau_1$ & $\rR 2$ & $\rHy$ & $\rR 1$ & $\rHy$ & $\rR 2$ & $\rHy$ & $\rHy$ & $\rHy$ & $\rR 2$ & $\rR 2$  \\
special trace & $\tau_2$ & $\rHy$ & $\rR 2$ & $\rR 1$ & $\rR 2$ & $\rHy$ & $\rHy$ & $\rR 2$ & $\rHy$ & $\rR 2$ & $\rHy$ \\
special trace & $\tau_3$ & $\rR 2$ & $\rR 2$ & $\rR 1$ & $\rHy$ & $\rHy$ & $\rHy$ & $\rHy$ & $\rHy$ & $\rHy$ & $\rHy$  \\
\hline
\end{longtable}
\begin{theorem} \label{thm:rdpc} 
Each automorphism $f  : X \to X$ in Nos.~$11$--$15$ of Table $\ref{tab:Ars1}$ has a 
unique periodic cycle $C_m$ of primitive period $m = 1, 2, 3, 5, 7, 9, 10, 11, 12, 14$ off $\cE(X)$, 
which is NR. 
Each points of $C_m$ is a center of a (periodic) rotation domain of rank $2$ $(\rR 2)$ or of rank $1$ $(\rR 1)$,   
or a hyperbolic periodic point $(\rHy)$, as indicated in Table $\ref{tab:rdpc}$. 
Moreover we know from Example $\ref{ex:rArDE}$ that $\cE(X)$ is contained in a rotation domain of rank $1$.     
\end{theorem}
{\it Proof}. 
Since $f$ has exceptional set of type $\rA_2^{\rr}$, we have $\hat{\nu}_{\ron}(f^m) = 
\hat{\nu}(f^m, \rA_2^{\rr, m})$ in the notation of \S \ref{ss:oec}, where   
\begin{equation} \label{eqn:A2r}
\hat{\nu}(f^m, \rA_2^{\rr, m}) = \dfrac{1}{H_m(\hat{\tau})} \quad \mbox{if $m$ is odd}; \qquad 
\hat{\nu}(f^m, \rA_2^{\rr, m}) = \dfrac{K_2(H_m(\hat{\tau}))}{H_{3 m}(\hat{\tau}) - 2(-1)^{m/2}} \quad \mbox{if $m$ is even},   
\end{equation}
by formula \eqref{eqn:nArm} with $n = 2$ in Proposition \ref{prop:Ar}. 
Solving equation \eqref{eqn:tt-pc}, with \eqref{eqn:A2r} in mind, for the unknown $\alpha_m + \alpha_m^{-1}$ 
recursively along the divisors of $m$, we find that the functions $P_m(w)$ 
in Induction \ref{induction} are given by  
\begin{align*}
P_1(w) &= \dfrac{w(w+2)}{w+1}, \hspace{30mm} P_2(w) = \dfrac{w^4-w^3-3w^2+w+1}{w(w^2-w-1)}, \\[2mm]
P_3(w) &= \dfrac{(w-1)(w+1)^2 \, \rST^{(8)}_1(w)}{w(w^2-w-1)(w^2+w-1)} \quad \mbox{with} \quad \rST^{(8)}_1(w) = w^4 -4 w^2 - w + 1, \\[2mm]
P_5(w) &= \dfrac{(w^2-w-1) (w^9 + 4 w^8 -18 w^6 - 20 w^5 + 5 w^4 + 6 w^3 - 12 w^2 - 8 w -1)}{w (w+1) (w^3-3w-1) (w^3+w^2-2w-1)}, 
\end{align*}  
and so on. 
Formulas for $P_m(w)$, $m = 7, 9, 10, 11, 12, 14$, are too messy to be given here, so left to Appendix \ref{app:rfp}. 
\par
Define the function $Q_m(w)$ from $P_m(w)$ as in definition \eqref{eqn:Q2} in Remark \ref{rem:Q} and apply Corollary \ref{cor:MDR} to it. 
For each $m = 1, 2, 5, 7, 9, 10, 11, 12, 14$, we have $Q_m(\tau_j) > 2$ for at least one index $j = 1, 2,3$, 
so the multipliers $\delta^{m/2} \alpha_m^{\pm 1}$ of $C_m$ are MI, in particular NR, where the case 
$Q_m(\tau_j) > 2$ leads to $\rHy$, while the case $Q_m(\tau_j) \in (-2, \, 2)$ leads to $\rR 2$ in Table \ref{tab:rdpc}. 
Things are different for $m = 3$, because $P_3(\tau_j) = 0$ and $Q_3(\tau_j) = -2$ for $j = 0, 1, 2, 3$, and hence  
$\alpha_3 + \alpha_3^{-1} = 0$, that is, $\alpha_3 = \pm \ri$. 
We are in case (I) of Theorem \ref{thm:MDR} with multipliers $\pm \ri \delta^{3/2}$ being MD but NR. 
Thus each point of $C_3$ is a center of a periodic rotation domain of rank $1$.  
\hfill $\Box$ \par\medskip  
The automorphisms in Theorem \ref{thm:rdpc} are remarkable to the effect that they admit  
a rotation domain of rank $1$ with {\sl center off the exceptional set}.  
In this respect the following exampls are also interesting.   
\begin{longtable}{ccclclcclcr}
\caption{
Examples \`{a} la \cite[Setup 3.2]{IT2} due to Y.~Takada. } \label{tab:takada} 
\\
\hline
      &           &           &           &          &          &      &           &                               &                       &     \\[-3mm]
No. & $\rho$ & $S(z)$ & $C(z)$ & $s(z)$ & $c(z)$ & ST & Dynkin & $\widetilde{\varphi}_1(z)$ & $L(f)$ & RH \\ 
\hline
\endfirsthead
\multicolumn{11}{l}{continued} \\
\hline
      &           &           &           &          &          &      &           &                               &                       &     \\[-3mm]
No. & $\rho$ & $S(z)$ & $C(z)$ & $s(z)$ & $c(z)$ & ST & Dynkin & $\widetilde{\varphi}_1(z)$ & $L(f)$ & RH \\ 
\hline
\endhead
\hline
\endfoot
\hline
\\[-4pt]
\caption*{\setlength{\leftskip}{-25mm}\setlength{\rightskip}{-25mm} 
The format of this table is the same as that of Table \ref{tab:DEs1}. }
\endlastfoot
$1$ & $14$ & $\rS^{(8)}_{17}$ & $\rC_8 \rC_{16}$ & $\rS^{(10)}_{1}$ & $\rC_{12} \rC_{30}$ & $\tau_2$ & $\emptyset$ & $\rC_1 \rC_2 \rC_8 \rC_{16}$ & $1$ & $\rR 1$ \\
$2$ & $14$ & $\rS^{(8)}_{17}$ & $\rC_8 \rC_{24}$ & $\rS^{(14)}_{17}$ & $\rC_{30}$ & $\tau_1$ & $\emptyset$ & $\rC_1 \rC_2 \rC_8 \rC_{24}$ & $1$ & $\rR 1$ \\
\hline
\end{longtable}
\begin{example} \label{ex:takada} 
The method of hypergeometric groups \`{a} la \cite[Setup 3.2]{IT2} produces two automorphisms $f : X \to X$ 
in Table \ref{tab:takada}, which have dynamical degree $\lambda^{(8)}_{17}$ and empty exceptional set. 
Since $\mu_{\roff}(f) = L(f) - \mu_{\ron}(f) = 1 - 0 = 1$, the map $f$ has a unique fixed point $p \in X$, 
which is simple, and the function $Q(w)$ in \eqref{eqn:Q2} is $Q(w) = (w^2-3)/(w+2)$. 
We observe that $Q(\tau_j) = 2 \cos(2 \pi/5)$ for $j = 0, 3$, and $Q(\tau_j) = 2 \cos(6\pi/5)$ for $j = 1, 2$, 
where $\tau_0 > \tau_1 > \tau_2 > \tau_3$ are the roots of $\rST^{(8)}_{17}(w) = w^4 + w^3 - 5 w^2 - 7 w - 1$.  
Thus we are in case (I) of Theorem \ref{thm:MDR} and the multipliers of $f$ at $p$ are algebraic, MD but NR.  
In conclusion $f$ admits a rotation domain of rank $1$ (with center at $p$), while $X$ contains no exceptional set.   
These automorphisms were found by a computer search due to Y.~Takada, after a disscussion 
with the author about the possibility of such examples. 
\end{example} 
\begin{appendix} 
\renewcommand{\thetable}{\Alph{section}.\arabic{table}}
\renewcommand{\thefigure}{\Alph{section}.\arabic{figure}}
\makeatletter
\@addtoreset{table}{section}
\@addtoreset{figure}{section}
\makeatother
\section{Study of Further Examples} \label{app:exp}
For a further illustration of Induction \ref{induction} as well as for a supplement to Proposition \ref{prop:Ars1},   
we study in detail the automorphisms $f : X \to X$ in Nos.~$1$ and $2$ of Table \ref{tab:Ars1}, 
which have dynamical degree $\lambda^{(18)}_6$ and exceptional set of type $\rA_1 \oplus \rA_2^{\rr}$. 
They are interesting because they contain an $\rA_1$-component (necessarily having trivial 
symmetry), which is intricate to handle under iterations of $f$ in the sense of Remarks \ref{rem:ecA} and \ref{rem:DE}. 
\begin{theorem} \label{thm:appA} 
In both of Nos.~$1$ and $2$ in Table $\ref{tab:Ars1}$, the $\rA_1$-component 
is contained in a rotation domain of rank $2$, while the $\rA_2$-component is 
contained in a rotation domain of rank $1$. 
For each $m = 1, 7, 8, 9, 10, 11$, there is exactly one periodic cycle $C_m$ 
of primitive period $m$ off $\cE(X)$. 
In No.~$1$ the periodic cycles $C_1$ and $C_{10}$ serve as centers of rotation domains of rank $2$, 
while for each $m = 7, 8, 9, 11$, $C_m$ is a hyperbolic periodic cycle.  
In No.~$2$ the unique point in $C_1$ is a center of a rotation domain of rank $2$, 
while for each $m = 7, 8, 9, 10, 11$, $C_m$ is a hyperbolic periodic cycle.    
\end{theorem}
\par
The proof of Theorem \ref{thm:appA} is completed after Lemma \ref{lem:A2r3} is proved. 
\begin{lemma} \label{lem:A2r1}  
We have $\mu(f, \rA_1) = 2$, hence for the $\rA_1$-component the case $(\rC 3)$ 
does not occur in Theorem $\ref{thm:ecA}$.  
\end{lemma}
{\it Proof}. 
It was shown in the proof of Proposition \ref{prop:Ars1} that either 
$\mu(f, \rA_1) = 2$ or $\mu(f, \rA_1) = 3$ held.     
We rule out the latter possibility.  
Suppose the contrary that $\mu(f, \rA_1) = 3$ is the case, so that the $\rA_1$-component 
falls into case (C3) of Theorem \ref{thm:ecA} with $n = 1$.  
Then this component contains two fixed points $p$ and $q$ with local indices $\mu(f, p) = 2$ and $\mu(f, q) = 1$. 
The multipliers of $f$ at $p$ are $\delta$ along $\rA_1$ and $1$ in a normal direction to $\rA_1$, while     
those at $q$ are $\delta^{-1}$ along $\rA_1$ and $\delta^2$ in a normal direction to $\rA_1$.  
Notice that $p$ is an exceptional fixed point of type $\rII$ in the terminology of \cite[\S 6]{IT2}.   
From \cite[Theorem 6.8]{IT2} we have for every $m \in \bN$,  
$$
\mu(f^m, p) = 2, \qquad 
\nu(f^m, p) = \dfrac{m-1 + (m+1) \delta^m + \theta (1 + \delta + \cdots + \delta^{m-1})}{m(1- \delta^m)^2}
$$
for some constant $\theta \in \bC$ independent of $m$. 
A little manipulation of $\nu(f^m, p)$ yields  
$$
\hat{\nu}(f^m, p) := \delta^{m/2} \nu(f^m, p) 
= \dfrac{m H_m(\hat{\tau}) + \hat{\theta} K_{m-1}(\hat{\tau})}{m \{ H_{2 m}(\hat{\tau}) -2 \}}, 
\qquad \hat{\theta} := \delta^{1/2} + \delta^{-1/2} (\theta - 1). 
$$
Since $f^m$ has multipliers $\delta^{(m/2) \pm (3m/2)}$ at $q$, we have $\mu(f^m, q) = 1$ and 
$$
\hat{\nu}(f^m, q) = \dfrac{1}{H_m(\hat{\tau}) - H_{3 m}(\hat{\tau})} = - \frac{1}{H_m(\hat{\tau}) \, \{H_{2m} (\hat{\tau}) -2\}}, 
$$
where $H_3(w) -w = w \{ H_2(w) -2 \}$ is used. 
Therefore $\mu(f^m, \rA_1) = \mu(f^m, p) + \mu(f^m, q) = 2+1 = 3$ 
and   
$$
\hat{\nu}(f^m, \rA_1) = \hat{\nu}(f^m, p) + \hat{\nu}(f^m, q) = \dfrac{1}{H_{2 m}(\hat{\tau})-2} 
\left\{ H_m(\hat{\tau}) - \dfrac{1}{H_m(\hat{\tau})} + \dfrac{\hat{\theta}}{m} K_{m-1}(\hat{\tau}) \right\}, 
\quad m \in \bN.   
$$
\par
By formula \eqref{eqn:mnAr} in Proposition \ref{prop:Ar} we have $\mu(f^m, \rA_2^{\rr, m}) = 2 + (-1)^m$ 
and so $\mu_{\ron}(f^m) = \mu(f^m, \rA_1) + \mu(f^m, \rA_2^{\rr, m}) = 3 + \{ 2 + (-1)^m \} = 5 + (-1)^m$ 
for every $m \in \bN$.   
Moreover, $\hat{\nu}(f^m, \rA_2^{\rr, m})$ is given by formulas \eqref{eqn:A2r}. 
All these formulas enable us to calculate $\hat{\nu}_{\ron}(f^m) 
= \hat{\nu}(f^m, \rA_1) + \hat{\nu}(f^m, \rA_2^{\rr, m})$ explicitly.  
For $m = 1$ we have $L(f) = 4 = \mu_{\ron}(f)$, hence $\mu_{\roff}(f) = L(f) - \mu_{\ron}(f) = 0$ and 
so $\hat{\nu}_{\roff}(f) = 0$.  
FPF \eqref{eqn:tt2} now reads $\hat{\tau} = \hat{\nu}_{\ron}(f)$. 
Solving it for the unknown $\hat{\theta}$ gives $\hat{\theta} = (\tau + 1)(\tau-3)/\hat{\tau}$.   
For $m = 2$ we have $L(f^2) = 6$ from formula \eqref{eqn:lnm} with 
$\widetilde{\varphi} = \rS^{(18)}_6 \rC_1^2 \rC_2^2$. 
Since $\mu_{\ron}(f^2) = 5 + (-1)^2 = 6$, we have $\mu_{\roff}(f^2) = L(f^2) - \mu_{\ron}(f^2) = 6 - 6 = 0$ 
and so $\hat{\nu}_{\roff}(f^2) = 0$. 
FPF \eqref{eqn:tt2} applied to $f^2$ then becomes $H_2(\hat{\tau}) = \hat{\nu}_{\ron}(f^2)$. 
Substituting $\hat{\theta} = (\tau + 1)(\tau-3)/\hat{\tau}$ into this equation, we obtain an 
algebraic equation $2 \tau^4 + \tau^3 - 7 \tau^2 - \tau + 1 = 0$.  
But this is absurd, as $\tau$ is a root of $\rST^{(16)}_6(w)$. 
Hence $\mu(f, \rA_1) = 3$ is impossible and we must have $\mu(f, \rA_1) = 2$. 
\hfill $\Box$ 
\begin{lemma} \label{lem:A2r2} 
For the $\rA_1$-component the case $(\rC 2)$ also does not occur in Theorem $\ref{thm:ecA}$.  
\end{lemma}
{\it Proof}. 
Formula \eqref{eqn:lnm} with $\widetilde{\varphi} = \rS^{(18)}_6 \rC_1^2 \rC_2^2$ evaluates 
$L(f^m)$, $m \in \bN$, as in the second row of Table \ref{tab:A2r}. 
Contrary to the assertion, suppose that $(\rC 2)$ is the case.  
Then by Remark \ref{rem:ecA} we have $\mu(f^m, \rA_1) = 2$ for every $m \in \bN$. 
Thus $\mu_{\ron}(f^m) = \mu(f^m, \rA_1) + \mu(f^m, \rA_2^{\rr, m}) = 2 + \{ 2 + (-1)^m \}= 4 + (-1)^m$  
as in the third row of Table \ref{tab:A2r}, 
so the numbers $\mu_{\roff}(f^m)$, $m \in \bN$, are given as in the fourth row. 
In particular, for each $m = 1, \dots, 6$, the unique $f^m$-fixed point $p_1$ off $\cE(X)$ is simple. 
Let $\delta^{1/2} \beta^{\pm 1}$ be the multipliers of $f$ at $p_1$ and put $B := \beta + \beta^{-1}$. 
Then for $m = 1, \dots, 6$, we have $\hat{\nu}_{\roff}(f^m) = \{ H_m(\hat{\tau}) - H_m(B) \}^{-1}$ 
with denominator being nonzero. 
Since $\hat{\nu}(f^m, \rA_1) = H_m(\hat{\tau}) \, \{ H_{2 m}(\hat{\tau}) - 2\}^{-1}$ by 
formula \eqref{eqn:Aii} in Theorem \ref{thm:ecA}, FPF \eqref{eqn:tt2} becomes
$$
H_m(\hat{\tau}) = \dfrac{H_m(\hat{\tau})}{H_{2 m}(\hat{\tau}) -2} + \hat{\nu}(f^m, \rA_2^{\rr, m}) 
+ \dfrac{1}{H_m(\hat{\tau}) - H_m(B)}, \qquad m = 1, \dots, 6,  
$$
where $\hat{\nu}(f^m, \rA_2^{\rr, m})$ is given by formula \eqref{eqn:A2r}. 
For $m=1$, solving this equation for the unknown $B$, we have 
$B = \hat{\tau} (\tau^2- 3 \tau -2)/(\tau^2- 2 \tau -4)$.  
Substituting this into the equation for $m = 2$, we find that $\tau$ satisfies the algebraic equation 
$\tau^5 - 2 \tau^4 - 8 \tau^3 + 8 \tau^2 + 11 \tau - 1 = 0$. 
This contradicts the fact that $\tau$ is a root of $\rST^{(18)}_6(w)$. 
Therefore for the $\rA_1$-component case (C2) does not occur in Theorem \ref{thm:ecA}.      
\hfill $\Box$ \par\medskip  
It follows from Lemmas \ref{lem:A2r1} and \ref{lem:A2r2} that the $\rA_1$-component falls   
in case (C1) of Theorem \ref{thm:ecA}. 
Let $\delta^{1/2} \alpha^{\pm 1}$ be the multipliers of the $\rA_1$-component.  
Formula \eqref{eqn:Ai} with $n = 1$ yields   
$\hat{\nu}(f, \rA_1) = \hat{\tau} \, \{ H_2(\hat{\tau}) - (\alpha^2 + \alpha^{-2}) \}^{-1}$.  
Since $\mu_{\roff}(f) = 1$, the map $f$ has a unique fixed point $p_1$ off $\cE(X)$, 
which is simple and whose multipliers can be expressed as $\delta^{1/2} \beta_1^{\pm 1}$ 
for some $\beta_1 \in \bC^{\times}$. 
We shall evaluate $A := \alpha^2 + \alpha^{-2}$ and $B_1 := \beta_1 + \beta_1^{-1}$ in terms of $\tau$. 
\begin{longtable}{lccccccccccccccc}
\caption{ 
Recipe \ref{recipe} applied to Nos.~$1$ and $2$ in Table \ref{tab:Ars1}. } \label{tab:A2r} 
\\
\hline
      &       &       &       &       &       &       &       &       &       &        &         &         &        &        &        \\[-3mm]
$m$ & $1$ & $2$ & $3$ & $4$ & $5$ & $6$ & $7$ & $8$ & $9$ & $10$ & $11$ & $12$ & $13$ & $14$ & $15$ \\ 
\hline
\endfirsthead
\multicolumn{16}{l}{continued} \\
\hline
      &       &       &      &       &       &       &       &       &       &         &        &         &        &        &         \\[-3mm]
$m$ & $1$ & $2$ & $3$ & $4$ & $5$ & $6$ & $7$ & $8$ & $9$ & $10$ & $11$ & $12$ & $13$ & $14$ & $15$ \\ 
\hline
\endhead
\hline
\endfoot
\hline
\\[-4pt]
\caption*{\setlength{\leftskip}{-5mm} \setlength{\rightskip}{-5mm} The bottom row shows whether the test \eqref{eqn:test} is succesful (S) or failed (F).  }
\endlastfoot
$L(f^m)$ & $4$ & $6$ & $4$ & $6$ & $4$ & $6$ & $11$ & $14$ & $13$ & $16$ & $15$ & $30$ & $30$ & $41$ & $49$ \\
$\mu_{\ron}(f^m)$ & $3$ & $5$ & $3$ & $5$ & $3$ & $5$ & $3$ & $5$ & $3$ & $5$ & $3$ & $5$ & $3$ & $5$ & $3$ \\
$\mu_{\roff}(f^m)$ & $1$ & $1$ & $1$ & $1$ & $1$ & $1$ & $8$ & $9$ & $10$ & $11$ & $12$ & $25$ & $27$ & $36$ & $46$ \\
$\bar{\mu}(m)$ & $1$ & $0$ & $0$ & $0$ & $0$ & $0$ & $7$ & $8$ & $9$ & $10$ & $11$ & $24$ & $26$ & $28$ & $45$ 
\\
S/F & S & S & S & S & S & S & S & S & S & S & S & F & F & F & F \\
\hline
\end{longtable}
\begin{lemma} \label{lem:A2r3} 
In the situation and notation mentioned above we have 
\begin{equation} \label{eqn:ABa}
A := \alpha^2 + \alpha^{-2} = -(\tau^2-3)(\tau^2-\tau-1), \qquad 
B_1 := \beta_1 + \beta_1^{-1} = \dfrac{\hat{\tau}(\tau-1)(\tau^4-4 \tau^2 + 2)}{\tau^5 - 5 \tau^3 + 6 \tau + 1}. 
\end{equation}
Moreover, the pairs $\delta^{1/2} \alpha^{\pm 1}$ and $\delta^{1/2} \beta_1^{\pm 1}$ are MI,   
all data in Table $\mathrm{\ref{tab:A2r}}$ is correct, $\mu(f^m, p_1) = 1$ and 
\begin{equation} \label{eqn:nuonpa}
\hat{\nu}_{\ron}(f^m) = \dfrac{H_m(\hat{\tau})}{H_{2 m}(\hat{\tau}) - H_m(A)} + \hat{\nu}(f^m, \rA_2^{\rr, m}), 
\qquad 
\hat{\nu}(f^m, p_1) = \dfrac{1}{H_m(\hat{\tau}) - H_m(B_1)},
\end{equation}
for every $m \in \bN$, where $\hat{\nu}(f^m, \rA_2^{\rr, m})$ is explicitly given   
by formula \eqref{eqn:A2r}. 
\end{lemma}
{\it Proof}. 
Although Table \ref{tab:A2r} is not yet established, 
we already know that $L(f^m)$-row and the column of $m = 1$ are correct  in Table \ref{tab:A2r}.  
Since $\mu_{\roff}(f) = 1$, we have $\mu_{\roff}(f^2) \ge 1$. 
Remark \ref{rem:A} and formula \eqref{eqn:mnAr} in Proposition \ref{prop:Ar} imply 
$\mu(f^2, \rA_1) \ge 2$ and $\mu(f^2, \rA_2^{\rr, 2}) = 3$, thus FPF \eqref{eqn:saito2} yields 
$$
6 = L(f^2) = \mu_{\ron}(f^2) + \mu_{\roff}(f^2) = \mu(f^2, \rA_1) + \mu(f^2, \rA_2^{\rr, 2}) 
+ \mu_{\roff}(f^2) \ge 2 + 3 + 1 = 6, 
$$
which forces $\mu(f^2, \rA_1) = 2$, $\mu_{\ron}(f^2) = 5$ and $\mu_{\roff}(f^2) = 1$. 
This implies that the column of $m = 2$ in Table \ref{tab:A2r} is correct. 
A similar reasoning shows that the columns of $m = 3, \dots, 6$ are also correct. 
In particular, for $m = 1, \dots, 6$, we find that $p_1$ is the only $f^m$-fixed point off $\cE(X)$, 
which is simple, and $\hat{\nu}(f^m, p_1)$ is given by the second formula in \eqref{eqn:nuonpa}.      
Moreover, as in the proof of Lemmas \ref{lem:A2r1} and \ref{lem:A2r2}, for each $m = 1, \dots, 6$, we can show 
that the $\rA_1$-component falls in case (C1) of Theorem \ref{thm:ecA} relative to the map $f^m$.  
Hence by formula \eqref{eqn:Ai} in Theorem \ref{thm:ecA}, $\hat{\nu}_{\ron}(f^m)$ is 
given by the first formula in \eqref{eqn:nuonpa}.  
FPF \eqref{eqn:tt2} applied to $f^m$ now reads  
\begin{equation} \label{eqn:tt1a}
H_m(\hat{\tau}) = \hat{\nu}_{\ron}(f^m) + \hat{\nu}(f^m, p_1), \qquad m = 1, \dots, 6.  
\end{equation}
\par
Solving equation \eqref{eqn:tt1a} with $m = 1$ for the unknown $B_1$, we have 
\begin{equation} \label{eqn:Ba}
B_1 = \hat{\tau} \dfrac{(\tau+1)(\tau-2)- A \tau}{\tau^2-2- A(\tau + 1)}. 
\end{equation}
Substituting formula \eqref{eqn:Ba} into equations \eqref{eqn:tt1a} with $m = 2, 3$, we have  
two algebraic equations for $A$,  
$$
\{ A + (\tau^2-3)(\tau^2-\tau-1) \} \, R_m(\tau; A) = 0, \qquad m = 2, 3, 
$$
where $R_2(\tau; A)$ and $R_3(\tau; A)$ are certain polynomials in $A$ of degrees $1$ and $3$,  
respectively, with coefficients in the field $\bQ(\tau)$. 
We observe that these polynomials have no roots in common, 
so we have $A = -(\tau^2-3)(\tau^2-\tau-1)$, the first formula in \eqref{eqn:ABa}. 
Substituting this into \eqref{eqn:Ba} yields the second formula in \eqref{eqn:ABa}. 
\par
Let $\tau_0 > \tau_1> \cdots > \tau_8$ be the roots of $\rST^{(18)}_6(w)$. 
Substituting $\tau = \tau_j$ into formulas \eqref{eqn:ABa} gives the table  
\begin{equation} \label{tab:YNa}
\mbox{\begin{tabular}{lccccccccc}
\hline
$\tau_j$ & $\tau_0$ & $\tau_1$ & $\tau_2$ & $\tau_3$ & $\tau_4$ & $\tau_5$ & $\tau_6$ & $\tau_7$ & $\tau_8$ \\ 
\hline
$A$ & Y & Y & Y & Y & N & Y & Y & Y & Y \\
$B_1$ & Y & Y & N & Y & Y & N & Y & Y & Y \\
\hline
\end{tabular}}
\end{equation}
where the signs Y and N for $A$ indicate $|A| < 2$ and $|A| > 2$ respectively  
and the same rule applies to $B_1$. 
By Corollary \ref{cor:MDR} the multipliers $\delta^{1/2} \alpha^{\pm 1}$ of the $\rA_1$-component  
and those at $p_1$, i.e. $\delta^{1/2} \beta_1^{\pm 1}$ are MI. 
By Remark \ref{rem:ecA} the $\rA_1$-component stays in case (C1) of Theorem \ref{thm:ecA} 
relative to $f^m$ for all $m \in \bN$. 
In particular we have $\mu_{\ron}(f^m) = 
\mu(f^m, \rA_1) + \mu(f^m, \rA_2^{\rr, m}) = 2 +\{ 2 + (-1)^m \} = 4 + (-1)^m$ for all $m \in \bN$, 
hence all data in Table \ref{tab:A2r} is correct. 
This establishes Table \ref{tab:A2r}.  
Now $\mu(f^m, p_1) = 1$ and formulas \eqref{eqn:nuonpa} hold for all $m \in \bN$. 
\hfill $\Box$ \par\medskip
{\it Proof of Theorem} \ref{thm:appA}. 
From Table \ref{tab:Ars1} the special trace $\tau$ of $f$ is either $\tau_8$ or $\tau_6$, 
for which $A$ and $B_1$ have Y in table \eqref{tab:YNa}. 
So Lemmas \ref{lem:rh-lm} and \ref{lem:rh-splest} combined with Theorems \ref{thm:elinz} and \ref{thm:Poschel} 
imply that the $\rA_1$-component is contained in a rotation domain of rank $2$ and the point $p_1$ is a center of 
a rotation domain of rank $2$.  
On the other hand, the $\rA_2$-component is contained in a rotation domain of rank $1$ 
by Theorem \ref{thm:rd}.   
\par
It is clear from Table \ref{tab:A2r} that for each $m = 7, 8, 9, 10, 11$, the map $f$ 
has a unique periodic cycle $C_m$ of prmitive period $m$ off $\cE(X)$, which is simple.  
Let $\delta^{m/2} \beta_m^{\pm 1}$ be the multipliers of $C_m$. 
Since $C_1 = \{ p_1\}$ is the only periodic cycle off $\cE(X)$ 
whose primitive period is a proper divisor of $m$, FPF \eqref{eqn:tt2} takes the form   
\begin{equation*} \label{eqn:tt2a}
H_m(\hat{\tau}) = \hat{\nu}_{\ron}(f^m) + \hat{\nu}(f^m, p_1) + 
\dfrac{m}{H_m(\hat{\tau}) - (\beta_m + \beta_m^{-1})}, \qquad m = 7, 8, 9, 10, 11,  
\end{equation*}
where the first and second terms in the right-hand side are given by formulas \eqref{eqn:nuonpa}. 
Solving this equation, we have an expression $\beta_m + \beta_m^{-1} = \hat{P}_m(\hat{\tau})$ as in formula 
\eqref{eqn:P1} with \eqref{eqn:P2}, to which we can associate a rational function $Q_m(w)$ 
as in definition \eqref{eqn:Q2}. 
Explicit formula for $Q_m(w)$ is too long to be included here, but for each $m = 7, 8, 9, 10, 11$,  
it tells us that $Q_m(\tau_6) > 2$, so the multipliers $\delta^{m/2} \beta_m^{\pm 1}$ are MI 
by Corollary \ref{cor:MDR}. 
For the automorphism $f$ of No.~$2$, which has special trace $\tau_6$, the periodic cycle $C_m$ is hyperbolic 
for every $m = 7, 8, 9, 10, 11$. 
On the other hand, we have $Q_{10}(\tau_8) \in (-2, \, 2)$ and $Q_m(\tau_8) > 2$ for $m = 7, 8, 9, 11$.   
Thus for the automorphism $f$ of No.~$1$, which has special trace $\tau_8$, we have the result 
stated in the theorem.     
\hfill $\Box$
\section{Rational Functions in the Proof of Theorem \ref{thm:rdpc}} \label{app:rfp} 
We present the explicit formulas for the rational functions $P_m(w)$, $m = 7, 9, 10, 11, 12, 14$, 
which are used in the proof of Theorem \ref{thm:rdpc} but not given there.     
In the following table we have $P_m(w) = P_m^{\rn}(w)/P_m^{\rd}(w)$.   
\begingroup
\begin{align*} 
P_7^{\rn}(w) 
&= (1 - 2 w - w^2 + w^3)  \\
&\phantom{=} \cdot (-1 - 12 w - 54 w^2 - 99 w^3 - 3 w^4 + 271 w^5 + 365 w^6 - 19 w^7- 403 w^8  \\
&\phantom{==} - 203 w^9 + 170 w^{10} + 156 w^{11} - 15 w^{12} - 43 w^{13} - 6 w^{14} + 4 w^{15} + w^{16}), 
\\[2mm]
P_7^{\rd}(w) &= w^2 (1 + w) (-1 - w + w^2) (-1 + w + w^2) (-1 - 3 w + w^3) (3 + 4 w - 5 w^2 - 5 w^3 + w^4 + w^5), 
\\[2mm] 
P_9^{\rn}(w) 
&= (-1 - 3 w + w^3) \\
&\phantom{=} \cdot (1 - w - 18 w^2 + 26 w^3 + 177 w^4 - 155 w^5 - 877 w^6 + 578 w^7 + 2643 w^8 - 1906 w^9 \\ 
&\phantom{==} - 6494 w^{10} + 3032 w^{11} + 11592 w^{12} - 911 w^{13} - 12700 w^{14} - 2671 w^{15} + 8178 w^{16} + 3452 w^{17} \\
&\phantom{==} - 3005 w^{18} - 1888 w^{19} + 555 w^{20} + 547 w^{21} - 19 w^{22} - 82 w^{23} - 9 w^{24} + 5 w^{25} + w^{26}), 
\\[2mm] 
P_9^{\rd}(w) 
&= w^2 (1 + w)^2 (-1 - w + w^2)^2 (-1 + w + w^2)^2 (-1 - 2 w + w^2 + w^3) \\ 
&\phantom{=} \cdot (1 - 4 w^2 + w^4) (3 + 4 w - 5 w^2 - 5 w^3 + w^4 + w^5),  
\\[2mm] 
P_{10}^{\rn}(w) 
&= 1 + 8 w - 27 w^2 - 370 w^3 - 276 w^4 + 5394 w^5 + 11335 w^6 - 36005 w^7 - 120516 w^8 + 112088 w^9 \\ 
&\phantom{=} + 681870 w^{10} - 42243 w^{11} - 2449139 w^{12} - 996059 w^{13} + 6113807 w^{14} + 4422993 w^{15} \\
&\phantom{=} - 11165936 w^{16} - 10979308 w^{17} + 15220163 w^{18} + 18819973 w^{19} - 15459172 w^{20} - 23624673 w^{21} \\
&\phantom{=}  + 11543667 w^{22} + 22236295 w^{23} - 6168625 w^{24} - 15899433 w^{25} + 2209896 w^{26} + 8707773 w^{27} \\
&\phantom{=} - 415589 w^{28} - 3667536 w^{29} - 40943 w^{30} + 1186599 w^{31} + 55414 w^{32} - 292645 w^{33} - 18927 w^{34} \\
&\phantom{=}  + 54092 w^{35} + 3734 w^{36} - 7264 w^{37} - 455 w^{38} + 670 w^{39} + 32 w^{40} - 38 w^{41} - w^{42} + w^{43},  
\\[2mm]
P_{10}^{\rd}(w) 
&= (-1 + w) w (1 + w)^2 (-1 - 3 w + w^3) (1 - 4 w^2 + w^4) \\ 
&\phantom{=} \cdot  (-1 - 5 w + 21 w^2 + 123 w^3 - 143 w^4 - 1196 w^5 + 215 w^6 + 5979 w^7 \\
&\phantom{==} + 1521 w^8 - 17237 w^9 - 7888 w^{10} + 30870 w^{11} + 16806 w^{12} - 35909 w^{13}  - 20105 w^{14}  \\
&\phantom{==} + 28075 w^{15} + 14821 w^{16} - 15091 w^{17} - 6994 w^{18} + 5614 w^{19} + 2122 w^{20}  \\
&\phantom{==} - 1425 w^{21} - 401 w^{22} + 236 w^{23} + 43 w^{24} - 23 w^{25} - 2 w^{26} + w^{27}), 
\\[2mm]
P_{11}^{\rn}(w) 
&= (-1 + 3 w + 3 w^2 - 4 w^3 - w^4 + w^5) \\
&\phantom{=} \cdot (1 - w - 84 w^2 - 475 w^3 - 1034 w^4 - 181 w^5 + 3468 w^6 + 5329 w^7 - 2924 w^8 - 14179 w^9 \\
&\phantom{==} - 5306 w^{10} + 18178 w^{11} + 15958 w^{12} - 12145 w^{13} - 18265 w^{14} + 3097 w^{15} + 11958 w^{16}  + 1442 w^{17} \\
&\phantom{==} - 4691 w^{18} - 1499 w^{19} + 1030 w^{20} + 533 w^{21} - 94 w^{22} - 90 w^{23} - 4 w^{24} + 6 w^{25} + w^{26}), 
\\[2mm]
P_{11}^{\rd}(w) 
&= w (1 + w) (-1 - w + w^2)^2 (-1 + w + w^2)^2 (-1 - 3 w + w^3) (-1 - 2 w + w^2 + w^3) \\ 
&\phantom{=} \cdot (1 - 3 w - 18 w^2 - 10 w^3 + 21 w^4 + 11 w^5 - 12 w^6 - 6 w^7 + 2 w^8 + w^9), 
\\[2mm] 
P_{12}^{\rn}(w)
&= -20 - 160 w + 476 w^2 + 6950 w^3 + 3764 w^4 - 117154 w^5 - 239193 w^6 + 893485 w^7  \\ 
&\phantom{=} + 3249683 w^8 - 1360366 w^9 - 19468339 w^{10} - 29662466 w^{11} + 22026364 w^{12}  \\ 
&\phantom{=} + 251855097 w^{13} + 481582562 w^{14} - 844985968 w^{15} - 3888644749 w^{16} + 302366509 w^{17} \\ 
&\phantom{=} + 16183296599 w^{18} + 9173623610 w^{19} - 44382464862 w^{20} - 44516363187 w^{21}  \\ 
&\phantom{=} + 86650822625 w^{22} + 121777600435 w^{23} - 124615231567 w^{24} - 233699234030 w^{25}  \\
&\phantom{=} + 133401139929 w^{26} + 339610170404 w^{27} - 104500466296 w^{28}- 388842857445 w^{29}  \\
&\phantom{=} + 55054266142 w^{30} + 359460832770 w^{31} - 11827612288 w^{32} - 272638478144 w^{33}  \\
&\phantom{=} - 10460123031 w^{34} + 171430947927 w^{35} + 14357358820 w^{36} - 89892019367 w^{37} \\ 
&\phantom{=} - 9943938296 w^{38} + 39387929930 w^{39} + 4887824270 w^{40} - 14399955526 w^{41}  \\ 
&\phantom{=} - 1841628256 w^{42} + 4371118862 w^{43} + 545303524 w^{44} - 1092492530 w^{45} - 127460120 w^{46}  \\
&\phantom{=} + 222049100 w^{47} + 23328765 w^{48} - 36058143 w^{49} - 3278148 w^{50} + 4560665 w^{51} + 341598 w^{52} \\
&\phantom{=} - 432442 w^{53} - 24864 w^{54} + 28892 w^{55} + 1128 w^{56} - 1212 w^{57} - 24 w^{58} + 24 w^{59},  
\\[2mm] 
P_{12}^{\rd}(w)
&= (-1 - w + w^2) \\
&\phantom{=} \cdot (2 + 14 w - 110 w^2 - 958 w^3 + 1749 w^4 + 26028 w^5 + 3082 w^6 - 380478 w^7 - 436880 w^8  \\ 
&\phantom{==} + 3358025 w^9 + 6467144 w^{10} - 18645884 w^{11} - 51827938 w^{12} + 63449893 w^{13}   \\ 
&\phantom{==} + 268432193 w^{14} - 105324660 w^{15} - 961844232 w^{16} - 104600789 w^{17} + 2474871557 w^{18} \\
&\phantom{==}  + 1117734403 w^{19} - 4691077974 w^{20} - 3351543459 w^{21} + 6687134676 w^{22} + 6284364526 w^{23}  \\
&\phantom{==} - 7302285375 w^{24} - 8427089759 w^{25} + 6213029709 w^{26} + 8542446972 w^{27} - 4182301586 w^{28} \\
&\phantom{==} - 6750878162 w^{29} + 2255525881 w^{30} + 4237514457 w^{31} - 982464866 w^{32} - 2134619098 w^{33}  \\
&\phantom{==} + 346205827 w^{34} + 865936298 w^{35}  - 98057541 w^{36} - 282103298 w^{37} + 21974692 w^{38}  \\ 
&\phantom{==} + 73159465 w^{39} - 3795259 w^{40} - 14873041 w^{41}+ 485199 w^{42} + 2314363 w^{43} \\
&\phantom{==} - 43060 w^{44} - 265650 w^{45} + 2360 w^{46} + 21164 w^{47} - 60 w^{48} - 1044 w^{49} + 24 w^{51}), 
\\[2mm]
P_{14}^{\rn}(w) 
&= 1 + 18 w + 48 w^2 - 1150 w^3 - 10584 w^4 - 16856 w^5 + 201582 w^6 + 1051039 w^7 - 245359 w^8 \\
&\phantom{=} - 15254219 w^9 - 30725714 w^{10} + 96191858 w^{11} + 436744053 w^{12} - 91245445 w^{13} \\
&\phantom{=} - 3150357080 w^{14} - 3254256234 w^{15}  + 13985975722 w^{16} + 30302036439 w^{17} - 36475611459 w^{18} \\
&\phantom{=}  - 156936591344 w^{19}  + 18041518157 w^{20} + 561751929661 w^{21} + 317438627252 w^{22} \\
&\phantom{=}  - 1457659086048 w^{23} - 1710734068488 w^{24} + 2699238148398 w^{25} + 5300481871754 w^{26} \\ 
&\phantom{=} - 3180939933930 w^{27} - 11730513662541 w^{28} + 882982500346 w^{29} + 19786697628112 w^{30} \\ 
&\phantom{=} + 5635816908727 w^{31} - 26191426046071 w^{32} - 15245985143865 w^{33} + 27603595017537 w^{34} \\
&\phantom{=}  + 24123674061494 w^{35} - 23267158000794 w^{36} - 28266220500808 w^{37}  + 15573787714152 w^{38} \\
&\phantom{=} + 26347932663471 w^{39} - 8040754054189 w^{40} - 20219254808844 w^{41} + 2913757422045 w^{42} \\ 
&\phantom{=} + 13039306145362 w^{43} - 436814784895 w^{44} - 7166401637458 w^{45} - 309896054469 w^{46} \\
&\phantom{=}  + 3390734778271 w^{47} + 322618975980 w^{48} - 1391029125965 w^{49} - 177724816066 w^{50} \\ 
&\phantom{=}  + 496967831823 w^{51} + 72107362132 w^{52} - 154846433912 w^{53} - 23305413552 w^{54} \\ 
&\phantom{=} + 42007015414 w^{55} + 6178786202 w^{56} - 9872582209 w^{57} - 1358882276 w^{58} \\ 
&\phantom{=} + 1993500799 w^{59}  + 248224654 w^{60} - 341791956 w^{61} - 37420814 w^{62}  \\
&\phantom{=} + 48983724 w^{63}  + 4590102 w^{64} - 5747538 w^{65} - 447302 w^{66} + 536800 w^{67}  \\ 
&\phantom{=} + 33335 w^{68}- 38319 w^{69} - 1784 w^{70} + 1961 w^{71} + 61 w^{72} - 64 w^{73} - w^{74} + w^{75},  
\\[2mm]
P_{14}^{\rd}(w)
&= w (1 + w)^2 (-1 - w + w^2) (-1 - 3 w + w^3) \\ 
&\phantom{=} \cdot (1 + 12 w - 2 w^2 - 668 w^3 - 2814 w^4 + 8166 w^5 + 80212 w^6 + 50896 w^7 -  975110 w^8  \\
&\phantom{==} - 2205443 w^9 + 5973625 w^{10} + 25470134 w^{11} - 12450873 w^{12} - 166424288 w^{13} \\ 
&\phantom{==} - 86034985 w^{14} + 700635446 w^{15} + 867601239 w^{16} - 1961974794 w^{17} - 3961632148 w^{18} \\ 
&\phantom{==} + 3507867687 w^{19} + 11849187936 w^{20} - 2971372911 w^{21} - 25625376296 w^{22} \\
&\phantom{==}  - 3227909676 w^{23} + 42039951115 w^{24} + 16153466649 w^{25} - 53956431186 w^{26} \\
&\phantom{==}   - 31238294251 w^{27} + 55464309358 w^{28} + 40715106798 w^{29} - 46599607721 w^{30} \\
&\phantom{==} - 39989365139 w^{31} + 32616742126 w^{32} + 30974868105 w^{33} - 19372520142 w^{34} \\
&\phantom{==}  - 19381444121 w^{35} + 9926067055 w^{36} + 9935616227 w^{37} - 4438547736 w^{38} \\
&\phantom{==}  - 4206912815 w^{39} + 1737989039 w^{40} + 1476621872 w^{41} - 592318987 w^{42} - 429453035 w^{43} \\ 
&\phantom{==} + 173279152 w^{44} + 103003098 w^{45} - 42710554 w^{46} - 20177661 w^{47} + 8682712 w^{48} \\ 
&\phantom{==} + 3177299 w^{49} - 1420971 w^{50} - 392392 w^{51} + 181736 w^{52} + 36584 w^{53} \\
&\phantom{==} - 17441 w^{54} - 2419 w^{55} + 1179 w^{56} + 101 w^{57} - 50 w^{58} - 2 w^{59} + w^{60}). 
\end{align*}
\endgroup
\end{appendix}
\par\vspace{3mm} 
{\bf Acknowledgments}. 
This work was supported by JSPS KAKENHI Grant Number JP22K03365.    
The author is indebted to Yuta Takada for providing him with Example \ref{ex:takada}.    

\end{document}